\documentclass[11pt,letterpaper]{amsart}
\usepackage[marginparwidth=75pt]
{geometry}

\usepackage{mathtools,amssymb,amsthm,mathrsfs,color,lineno,paralist,graphicx,float}
\usepackage[colorlinks,
linkcolor=blue,
anchorcolor=green,
citecolor=blue,
]{hyperref}

\usepackage{amssymb,amsthm,amsmath}
\usepackage[T1]{fontenc}
\usepackage[utf8]{inputenc}
\usepackage{color}
\usepackage[svgnames]{xcolor}
\usepackage{ifpdf}
\usepackage[all]{xy}

\usepackage{tikz-cd} 
\usepackage{bbm}

\usepackage{calc}
\newcommand{\Aa}{{\overset{\circ}{A}}}
\renewcommand{\leq}{\leqslant}
\renewcommand{\geq}{\geqslant}
\newcommand{\hide}[1]{}

\newcommand{\rmm}[1]{\mathrm{#1}} 
\newcommand\pa[1]{\left(#1\right)}

\newcommand\av[1]{\left|#1\right|}
\newcommand\op[1]{\operatorname{#1}}
\newcommand\diag[1]{\operatorname{diag}\left(#1\right)}

\usepackage[marginparwidth=75pt]
{geometry}
\usepackage{pgfplots}
\usepackage{mathtools,amssymb,amsthm,mathrsfs,color,lineno,paralist,graphicx,float}
\usepackage[colorlinks,
linkcolor=blue,
anchorcolor=green,
citecolor=blue,
]{hyperref}

\usepackage{tikz}
\usepackage{pgfplots}
\pgfplotsset{compat=1.18}
\usepackage{graphicx}

\usepackage{enumitem}
\usepackage[T1]{fontenc}
\usepackage[utf8]{inputenc}
\usepackage{subfig} 

\usepackage[justification = centering, labelsep =period]{caption} 

\usepackage{calc}
\definecolor{bleu1}{RGB}{0,57,128}
\def\bleu1{\color{bleu1}}

\usepackage{etoolbox}
\patchcmd{\section}{\normalfont}{\normalfont \bleu1}{}{}
\patchcmd{\subsection}{\normalfont}{\normalfont \bleu1}{}{}
\patchcmd{\subsubsection}{\normalfont}{\normalfont \bleu1}{}{}

\author{Xianzhe Li}
\address{Department of Mathematics, University of California, Berkeley, CA 94720, USA} 
 \email{xianzhe@berkeley.edu}

\author{Disheng Xu}
\address{School of Science, Great Bay University and Great bay institute for advanced study, 
Songshan Lake International Innovation Entrepreneurship Community A5, Dongguan 523000, China}
\email{xudisheng@gbu.edu.cn}

\author{Qi Zhou}
\address{
Chern Institute of Mathematics and LPMC, Nankai University, Tianjin 300071, China
}
\email{qizhou@nankai.edu.cn}

\newtheorem{theorem}{Theorem}[section]
\newtheorem{lemma}[theorem]{Lemma}
\newtheorem{proposition}[theorem]{Proposition}
\newtheorem{corollary}[theorem]{Corollary}
\newtheorem{conjecture}[theorem]{Conjecture}

\theoremstyle{definition}

\newtheorem{definition}[theorem]{Definition}

\newtheorem{remark}[theorem]{Remark}

\theoremstyle{plain}

\numberwithin{equation}{section}

\newcommand{\I}{{\mathbf I}}
\newcommand{\OO}{{\mathcal O}}
\newcommand{\mc}[1]{\mathcal{#1}}
\newcommand{\PP}{\mathbb P}
\newcommand{\A}{{\mathbb A}}

\newcommand{\N}{{\mathbb N}}
\newcommand{\Q}{{\mathbb Q}}
\newcommand{\R}{{\mathbb R}}
\newcommand{\C}{\mathbb{C}}
\newcommand{\T}{{\mathbb T}}
\newcommand{\V}{{\mathbb V}}

\newcommand{\Z}{{\mathbb Z}}

\newcommand{\e}{\varepsilon} 


\makeatletter
\def\saveenum{\xdef\@savedenum{\the\c@enumi\relax}}
\def\resetenum{\global\c@enumi\@savedenum}
\makeatother

\title[Stability of DTMP]{Monotonicity, 
Global symplectification
and the stability of Dry Ten Martini Problem}

\def \V{\mathcal{V}}

\newcommand{\CC}{\mathbb{C}}

\DeclareMathOperator{\GL}{GL}
\DeclareMathOperator{\Her}{Her}

\newcommand{\U}{\mathbf{U}}

\newcommand{\Lag}{\operatorname{Lag}}

\usepackage{tikz}
\usetikzlibrary{decorations.pathreplacing, decorations.pathmorphing, calc, arrows.meta, patterns}
\definecolor{lightblue}{RGB}{173,216,230}
\definecolor{darkblue}{RGB}{0,0,139}
\definecolor{lightred}{RGB}{255,182,193}
\definecolor{darkred}{RGB}{139,0,0}

\definecolor{basecolor}{RGB}{230,240,255}
\definecolor{fibercolor}{RGB}{255,230,230}
\definecolor{connection}{RGB}{180,0,0}
\definecolor{transport}{RGB}{0,100,0}
\definecolor{holonomy}{RGB}{0,0,180}

\newcommand{\tnorm}[1]{\left\vert\mkern-1.5mu\left\vert\mkern-1.5mu\left\vert #1 
                            \right\vert\mkern-1.5mu\right\vert\mkern-1.5mu\right\vert}

\begin{document}

\begin{abstract}

We prove that, for every irrational frequency and every analytic Type I
potential, each supercritical spectral energy satisfying the gap-labelling
condition is an endpoint of an open spectral gap. This establishes the
conjecture of Ge--Jitomirskaya--You \cite{GJY,You} in the supercritical
regime. Consequently, the
``all gaps open'' property of the supercritical almost Mathieu operator
persists under sufficiently small analytic perturbations.

The main ingredient is a global symplectification of the center bundle
that preserves quantitative monotonicity. This allows us to study gap
opening through the center dynamics of the dual long-range operator,
which has no natural Schr\"odinger form. We first establish the result
for trigonometric polynomial potentials and then pass to general analytic
potentials by controlling the dependence on the truncation dimension.
The proof combines a discrete Hellmann--Feynman identity,
dimension-free Aubry duality in weighted analytic norms, and a
quantitative cone argument based on pre-monotonicity. These estimates
ensure that the gaps survive in the analytic limit.

Our results establish analytic stability of the Dry Ten Martini Problem
in the supercritical regime and give a partial answer to a question of
M.~Shamis on the persistence of periodic spectral gaps.
\end{abstract}
		\maketitle{\Large }

\section{Introduction}

\subsection{TKNN Theory and the Robustness of the Dry Ten Martini Problem}

Following von Klitzing's discovery of the integer quantum Hall effect \cite{Klitzing1980}, Thouless and coauthors provided a profound theoretical framework \cite{TKNN1982}, demonstrating that the quantized Hall conductance is inherently determined by a topological invariant---the Chern number---a breakthrough that earned Thouless a share of the 2016 Nobel Prize in Physics. Within this framework, the Gap Labelling Theorem (GLT) for quasiperiodic operators assigns precise labels to spectral gaps, corresponding directly to these Chern numbers. Consequently, the assertion that ``all spectral gaps are open'' carries the direct physical implication that all predicted topological quantum Hall phases are strictly realized in the corresponding system. Notably,  some arguments  in the seminal TKNN work \cite{TKNN1982} implicitly rely on the assumption that all spectral gaps of the almost Mathieu operator (AMO) are open, particularly as the coupling parameter $\lambda$ is varied. The AMO is defined as
\begin{equation}\label{amo}
(H_{\lambda,\alpha,x}u)_n = u_{n+1} + u_{n-1} + 2\lambda \cos(2\pi(n\alpha +x))u_n,
\end{equation}
and was originally introduced by Peierls \cite{Pe} to model an electron on a two-dimensional lattice in a perpendicular magnetic field \cite{Ha,R}.

A natural question is whether all gaps remain open under small analytic perturbations of the almost Mathieu potential. This is relevant to  physical models, where the potential need not be exactly a cosine.
Gap stability is also related to the robustness of topological phases
and quantized conductance under perturbations, as discussed by
Schneider \textit{et al.} \cite{GBSS}:
``A key feature of topological phases is the robustness
to local perturbations that leads to quantized physical
observables characterized by a bulk topological invariant,
a famous example being the integer quantum Hall effect ”.
In this work, we prove that, the ``all spectral gaps are open'' property holds for an open set of quasiperiodic operators,  for any irrational frequency.

To state our results precisely, we first recall the framework of the Gap Labelling Theorem (GLT), which provides the proper context for defining what it means for ``all spectral gaps to be open.'' Consider one-frequency analytic quasiperiodic Schr\"odinger operators on $\ell^2(\mathbb{Z})$, defined by
\begin{equation}\label{eq:finalH}
    (H_{v,\alpha,x}u)_n = u_{n+1} + u_{n-1} + v(x+n\alpha)u_n,
\end{equation}
where $\alpha \in \mathbb{R}\setminus\mathbb{Q}$ is the frequency, $x \in \mathbb{T}$ is the phase, and $v \in C^\omega(\mathbb{T},\mathbb{R})$ is the analytic potential. The spectrum $\Sigma_{v,\alpha}$ of $H_{v,\alpha,x}$ is a compact subset of $\mathbb{R}$ with no isolated points and is independent of $x$ for irrational $\alpha$. Any connected component of $\mathbb{R}\setminus\Sigma_{v,\alpha}$ is called a \emph{spectral gap}.

The integrated density of states (IDS) $N_{v,\alpha}(E)$ is defined as
\[
N_{v,\alpha}(E) = \int_{\mathbb{T}} \mu_{v,\alpha,x}(-\infty,E]\,dx,
\]
where $\mu_{v,\alpha,x}$ is the spectral measure of $H_{v,\alpha,x}$ associated with $\delta_0$. The IDS is continuous, strictly increasing on the spectrum, and locally constant on its complement. The GLT \cite{BLT,gaplabel} asserts that in each bounded spectral gap, there exists a unique nonzero integer $k \in \mathbb{Z}\setminus\{0\}$ such that
$N_{v,\alpha}(E) = k\alpha \mod \mathbb{Z}.$
This integer $k$ is called the \emph{label} of the spectral gap.

At the 1981 AMS Annual Meeting, Mark Kac famously posed the question of whether AMO ``has all its gaps there'', jokingly offering ten martinis for its solution. Barry Simon later reformulated this into two distinct problems: the \emph{Ten Martini Problem}, which asserts that the spectrum $\Sigma_{\lambda,\alpha}$ of the AMO is a Cantor set for all irrational $\alpha$ and $\lambda \neq 0$; and the more difficult \emph{Dry Ten Martini Problem} (DTMP), which asks whether for every nonzero integer $k \in \mathbb{Z}\setminus\{0\}$, there exists a spectral gap with label $k$.

The Ten Martini Problem was ultimately resolved by Avila and Jitomirskaya \cite{AJ05}, with partial results having been obtained in \cite{AK06,BS,HS,L,P}. However, the Dry Ten Martini Problem remained open, recent result of Avila-You-Zhou \cite{avila2016dry}, showed that all spectral gaps of $H_{\lambda,\alpha,x}$ are open, except possibly in the critical case $(\lambda,\beta(\alpha))=(\pm1,0)$. Here,
\begin{equation*}
  \beta(\alpha) = \limsup_{n\to\infty} \frac{\log q_{n+1}}{q_n},
\end{equation*}
where $p_n/q_n$ are the continued fraction convergents of $\alpha$. For related developments and partial results, see \cite{AJ05,AJ08,CEY,P,LY2015}.

Recently, Argentieri and Avila \cite{AA} discovered a striking phenomenon: in the subcritical regime ($0 < \lambda < 1$), the TMP and DTMP are \emph{not} robust. There exists arbitrarily small analytic perturbations which collapse the gaps and produce interval spectrum. In  contrast, we establish that the DTMP exhibits strict topological rigidity in the supercritical regime:

\begin{theorem}\label{thm:example-thm}
Let $\alpha \in \mathbb{R}\setminus\mathbb{Q}$, and \( f \in C^\omega(\mathbb{T}, \mathbb{R}) \). For every $\lambda > 1$, there exists $\delta_0 = \delta_0(\lambda,\alpha, f) > 0$ such that for all $0 < \delta < \delta_0$, every spectral gap of the perturbed operator $H_{2\lambda \cos + \delta f,\alpha,x}$ is open.
\end{theorem}

From a physical perspective, the resolution of the DTMP provides a complete set of topological invariants, while gap stability ensures that each invariant corresponds to a measurable, quantized observable. From a mathematical perspective, stability is particularly subtle in the Liouvillean regime ($\beta(\alpha) > 0$), where spectral gaps are tiny and non-uniformly distributed \cite{ALSZ2024abominable},  exponentially and sub-exponentially small gaps  (in $q_n$) coexist \cite{avila2016dry}. This intricate structure makes the persistence of every labelled gap under a fixed perturbation a separate issue from the persistence of Cantor spectrum. The latter was established in the Type I setting by Ge--Jitomirskaya--You \cite{GJY,Martini}; the precise scope of their results is recalled below. For DTMP with rough potential, one can consult \cite{ABD,DGY,BBL,GeWXu} and  the reference therein. 

\subsection{The Evolution of DTMP}

To illuminate the mechanism behind Theorem \ref{thm:example-thm}, it is necessary to trace the methodological evolution of the DTMP:

\smallskip

\noindent\textbf{Periodic approximation.}
For Liouvillean frequencies, the method of periodic approximation provides a natural approach. 
Combining this with $C^*$-algebraic techniques, Choi-Elliott-Yui \cite{CEY} established 
the DTMP for $\lambda=1$ when $\beta(\alpha)$ is sufficiently large. 
Avila and Jitomirskaya \cite{AJ05} subsequently refined these estimates to the methodological 
limit, extending the result to $\beta(\alpha) > |\ln \lambda|$. 
However, the periodic approximation method appears to encounter fundamental limitations 
for weakly Liouvillean frequencies. As demonstrated in \cite{CEY,AJ05}, this approach can only 
guarantee gaps that are exponentially small in $q_n$—an issue we will revisit later 
in our discussion.

\medskip

\noindent\textbf{Aubry duality and Moser-Poschel argument.}
For Diophantine frequencies $\alpha$, the dynamical systems approach to the DTMP  relies on the \emph{reducibility} of associated quasiperiodic cocycles. 
Recall that two cocycles $(\alpha, A^i)$, $i=1,2$, are said to be \emph{conjugated} if there exists $B \in C^\omega(\mathbb{T}, \mathrm{PSL}(2, \mathbb{R}))$ such that 
$A^1(x) = B(x + \alpha)A^2(x)B(x)^{-1}.$
A cocycle $(\alpha, A)$ is called \emph{reducible} if it is conjugated to a constant cocycle. 

The eigenvalue equation $H_{v,\alpha,x}u = Eu$ gives rise to the  cocycle $(\alpha, S_E^v)$, where
\begin{equation} \label{eq:transfer_matrix}
    S_E^{v}(x) = \begin{pmatrix} E - v (x) & -1 \\ 1 & 0 \end{pmatrix}.
\end{equation}
In the Diophantine case, if $E$ is a potential gap boundary according to the GLT, then Eliasson's reducibility result \cite{E92} implies that $(\alpha, S_E^{2\lambda \cos})$ is analytically reducible to a parabolic cocycle of the form $\begin{pmatrix} 1 & d \\ 0 & 1 \end{pmatrix}$. 
Aubry duality \cite{GJLS1997Duality} subsequently ensures that $d \neq 0$, while the Moser--P\"oschel argument \cite{MP84} verifies the openness of the gap.
Using this framework, Puig \cite{P} resolved the DTMP in the perturbative regime for sufficiently small $\lambda$; he later extended the result \cite{P06} to generic small analytic potentials. Nevertheless, this approach has two  limitations: it requires the potential to be sufficiently small, and it fundamentally depends on the frequency being Diophantine.

\medskip

\noindent\textbf{Avila's global theory.}
To extend the result to the global regime, Avila and Jitomirskaya \cite{AJ08} established the \textit{almost reducibility} of the cocycle $(\alpha, S_E^{2\lambda\cos})$ by exploiting the almost localization of the dual operator. (Recall that a cocycle $(\alpha, A)$ is almost reducible if its analytic conjugacy class contains a constant cocycle.) This reduction of global dynamics to local behavior allowed them to resolve the DTMP for Diophantine frequencies and all $\lambda \neq 1$ \cite{AJ08}. The result was later extended by Avila \cite{avila2008absolutely} to the case $\beta(\alpha) = 0$ with $\lambda \neq 1$.

For energies in the spectrum, almost reducibility implies the strong vanishing of the Lyapunov exponent, and moreover this approach removes the requirement that $\lambda$ be small. This line of investigation culminated in Avila's Almost Reducibility Conjecture (ARC) and his subsequent development of the global theory \cite{avila0,AJ08}. 
Let $L(E)$ denote the Lyapunov exponent of the cocycle $(\alpha, S_E^v)$. A cornerstone of Avila's theory is the analysis of the complexified Lyapunov exponent
$
L_\varepsilon(E) = L\bigl(\alpha, S_E^v(\cdot+i\varepsilon)\bigr),$
which he showed is an even, convex, piecewise‑linear function of $\varepsilon$ with integer slopes. The acceleration is defined as the integer
\[
\omega(E)=\lim_{\varepsilon\to0^+} \frac{L_\varepsilon(E)-L_0(E)}{2\pi\varepsilon}\in \mathbb{Z}.
\]
This structure naturally partitions the spectrum into three dynamical regimes:
\begin{itemize}
    \item \textbf{Subcritical:} $L(E)=0$ and $\omega(E)=0$.
    \item \textbf{Critical:} $L(E)=0$ and $\omega(E)\geq 1$.
    \item \textbf{Supercritical:} $L(E)>0$.
\end{itemize}
The ARC asserts that subcriticality (i.e., strong vanishing of the Lyapunov exponent in the spectrum) implies almost reducibility. The conjecture was fully proved by Avila \cite{AviAlmost,Avi2023KAM}. See also \cite{GeARC,HS1} for trigonometric polynomial potentials with Diophantine frequencies. The ARC played a central role in the subsequent study of the spectral theory of quasiperiodic Schr\"odinger operators \cite{AYZ,avila2016dry,LYZZ2017Asymptotics,GY2020Arithmetic,GJ,GJY,GYZ2024Structured,LYZ2024Exact}.

\medskip

\noindent\textbf{Quantitaive Aubry duality and Generalized  Moser-P\"oschel argument.}
However, this left the case of Liouvillean frequencies unresolved. The useful insight came from a fundamental shift in perspective. Previous spectral analysis \cite{avila2008absolutely,AJ08,BJ02}, all proceeded from the supercritical to the subcritical regime, moreover it was long believed within the community that Aubry duality was ineffective for Liouvillean frequencies. Avila--You--Zhou \cite{avila2016dry} reversed this direction: they began dynamically on the subcritical side, leveraging the ARC to reach conclusions about the supercritical regime. To handle Liouvillean frequencies effectively, their method required the development of a \textit{quantitative Aubry duality} combined with a \textit{quantitative almost reducibility} analysis. This analysis in fact provided the first all-frequency approach to problems concerning Cantor spectrum.

This approach is novel even for Diophantine frequencies, as it yields good control over reducible parabolic cocycles, which in turn leads to exponential lower bounds on the decay of spectral gaps. Recently, Ge-You-Zhou \cite{GYZ2024Structured} developed a sharp version of quantitative Aubry duality, which provides the precise exponential decay rates for spectral gaps of the almost Mathieu operator.
With the quantitative Aubry duality established, to open the gap for any frequency, \cite{avila2016dry} develops the \textit{generalized Moser--P\"oschel argument}, which  assumes only quantitative almost reducibility rather than full reducibility. It proves the existence of some \(\tau' \in \mathbb{R} \setminus \{0\}\) such that the cocycle \((\alpha, S_{E+\tau'}^{2\lambda \cos})\) is uniformly hyperbolic, while satisfying \(N_{\lambda,\alpha}(E+\tau') = N_{\lambda,\alpha}(E)\). This directly implies the openness of the corresponding spectral gap.

\medskip

\noindent\textbf{Quantitative Avila's Global Theory.}
While these techniques provided the first all-frequency resolution for the AMO, they encounter a severe structural limitation when generalized: they fundamentally rely on the natural Schr\"odinger structure of the equations. This renders them highly effective in the subcritical regime of general analytic quasiperiodic operators.
For instance, combining Avila's solution of ARC
with the quantitative Moser–P\"oschel argument developed in \cite{avila2016dry}, one can show
that the DTMP holds generically in the subcritical regime for all irrational
frequencies \cite{LZ}.  However, these techniques are structurally inadequate for the supercritical regime.

For the spectral analysis in the  supercritical regime, Bourgain, Goldstein, and Schlag
\cite{BouGreens,BG,GS1,GS11} made significant progress. To handle resonances effectively,
however, their methods necessarily exclude certain Diophantine frequencies.
For example, Goldstein and Schlag \cite{GS11} proved that Cantor spectrum holds
for almost every Diophantine frequency in the supercritical regime.

To obtain results at a fixed frequency—or better, for all frequencies—a new
approach is required. The key development is the \textit{Quantitative Avila
Global Theory} introduced by Ge-Jitomirskaya-You-Zhou \cite{GJYZ}.
This exciting direction of progress is to use duality techniques, which gives
new information about the supercritical region \cite{Avi2023KAM}. Let us explain this in detail.

Assume \( v_d(\theta) = \sum_{k=-d}^{d} \hat{v}_k e^{2\pi i k\theta} \) is a
trigonometric polynomial of degree \(d\). The dual operator of \eqref{eq:finalH}
takes the form
\begin{equation}\label{eq:finalL}
    (L_{v_d,\alpha,x}u)_n = \sum_{k=-d}^{d} \hat{v}_k u_{n+k} + 2\cos\bigl(2\pi (x+n\alpha)\bigr) u_n.
\end{equation}
Define the matrices
\[
C_d = \begin{pmatrix}
\hat{v}_d & \cdots & \hat{v}_1 \\
0 & \ddots & \vdots \\
0 & 0 & \hat{v}_d
\end{pmatrix},
\qquad
V_d(x) = \begin{pmatrix}
2\cos (2\pi x_0) & \hat{v}_{-1} & \cdots & \hat{v}_{-d+1} \\
\hat{v}_1 & 2\cos (2\pi x_1) & \ddots & \vdots \\
\vdots & \ddots & \ddots & \hat{v}_{-1} \\
\hat{v}_{d-1} & \cdots & \hat{v}_1 & 2\cos (2\pi x_{d-1})
\end{pmatrix},
\]
where \( x_j = x + j\alpha \). Then \( L_{v_d,\alpha,\theta} \) coincides with the
generalized matrix-valued Schr\"odinger operator
\( \widehat H_{V,\alpha,\theta} \) acting on \( \ell^2(\mathbb{Z}, \mathbb{C}^d) \):
\begin{equation} \label{eq:hatH}
    (\widehat H_{V_d,\alpha,x}\mathbf{u})_n = C_d \mathbf{u}_{n+1} + V_d(x+dn\alpha) \mathbf{u}_n + C_d^* \mathbf{u}_{n-1}.
\end{equation}
This induces a cocycle \( (d\alpha, A^{E,d}) \) with
\[
A^{E,d}(x)=
\begin{pmatrix}
C_d^{-1}\bigl(E - V_d(x)\bigr) & -C_d^{-1}C_d^* \\[4pt]
I_d & 0
\end{pmatrix}.
\]

A fundamental observation of \cite{GJYZ} is that \( \omega(E) = k > 0 \) if and only
if \( (d\alpha, A^{E,d}) \) is \emph{partially hyperbolic} with a \(2k\)-dimensional
center. To extend this characterization to all energies—including subcritical and
uniformly hyperbolic cocycles—and to obtain stability results (note that the
acceleration is only upper semicontinuous), Ge, Jitomirskaya, and You \cite{GJY} introduced
the concept of \emph{$T$-acceleration}. Let \( \varepsilon_1(E) \) be the first
turning point of the complexified Lyapunov exponent \( L_\varepsilon(E) \) as a
function of the imaginary part \( \varepsilon \). The $T$-acceleration at energy
\( E \) is defined by
\[
\overline{\omega}(E)=\lim_{\varepsilon\to\varepsilon_1^+}
\frac{L_\varepsilon(E)-L_{\varepsilon_1}(E)}{2\pi(\varepsilon-\varepsilon_1)},
\]
with \( \overline{\omega}(E)=0 \) if no such turning point exists. We say that
\( E \) is a \emph{Type I} energy for the operator \eqref{eq:finalH} if
\( \overline{\omega}(E)=1 \). The property of being Type I is stable in the
analytic topology\footnote{One must take care with the analytic radius.}, and
\( \overline{\omega}(E)=1 \) holds if and only if \( (d\alpha, A^{E,d}) \) is
partially hyperbolic with a \(2\)-dimensional center.

An analytic potential is of Type I if every energy in its spectrum is
of Type I. The corresponding operator is called a Type I operator.
Ge--Jitomirskaya--You \cite{GJY} proved Cantor spectrum for even
trigonometric-polynomial Type I potentials at every irrational frequency,
without arithmetic restrictions. Their proof developed higher-dimensional
Kotani theory and simplicity of point spectrum. The extension to general analytic Type I potentials,
announced in \cite{GJY} in 2023, is proved in their manuscript
\emph{The Robust Ten Martini Problem} \cite{Martini}, which was largely
completed in 2025.

Beyond Cantor spectrum, one may ask whether every gap allowed by gap
labelling is open. This stronger conclusion was conjectured for Type I
operators in \cite{GJY}; see also Conjecture~4.10 in \cite{You}.

\begin{conjecture}[\cite{GJY,You}]
All gaps of Type I operators are open.
\end{conjecture}

In the present paper, we establish this conjecture in the supercritical
regime:

\begin{theorem}\label{main-thm-analytic}
 Let $\alpha\in\mathbb{R}\setminus\mathbb{Q}$ and let $v$ be an analytic
Type I potential. Define
\[
\Sigma_{v,\alpha}^{+}
:=
\{E\in\Sigma_{v,\alpha}:L(E)>0\}.
\]
Then every energy $E_k\in\Sigma_{v,\alpha}^{+}$ satisfying
$N_{v,\alpha}(E_k)\equiv k\alpha\pmod{\mathbb{Z}}$
for some integer $k$ is an endpoint of an open spectral gap.
\end{theorem}

In the subcritical regime, the conjecture was proved in \cite{GJY} for trigonometric polynomial potentials under the assumption $\beta(\alpha)=0$. As noted in \cite[Remarks 1.9 and 1.11]{GJY}, the result is expected to extend to arbitrary irrational frequencies and general analytic Type I potentials; this extension will appear in \cite{GJY2}. See also \cite{L} for a recent alternative proof.
However, small analytic perturbations of the subcritical AMO need not remain of Type I. This property is preserved when
the perturbations are analytic in a strip of width greater than
$-\ln\lambda/(2\pi)$. The loss of the Type I property in narrower
strips allows the counterexample of Argentieri and Avila \cite{AA}.

\subsection{Periodic Approximation and Spectral Rigidity: A question of M. Shamis}
\label{q:shamis}

In the trigonometric polynomial case, a natural initial strategy to tackle Theorem \ref{main-thm-analytic} would be to adapt the generalized Moser-P\"oschel argument developed in \cite{avila2016dry}. However, such an approach is  obstructed by the lack of an underlying Schr\"odinger structure in the dual cocycle. While one might attempt a symplectic block diagonalization followed by brute-force perturbative estimates, this would not only obscure the intrinsic dynamics but also entail degenerate technicalities (one can consult Remark \ref{mp-co} for related discussions). This conceptual barrier necessitates a paradigm shift: instead of forcing a perturbative Schr\"odinger framework onto the dual system, we ask whether a more robust, geometric mechanism governs the survival of spectral gaps. This perspective motivates our return to the method of periodic approximation, fundamentally reshaping the approach of Avila-You-Zhou \cite{avila2016dry}. 
Two key distinctions separate \cite{avila2016dry} from the earlier work of Choi–Elliott–Yui \cite{CEY}. First, while there indeed exist gaps that are exponentially small in $q_n$ \cite{CEY}, the method in \cite{avila2016dry} establishes sub-exponential estimates for gaps. Second, whereas the $C^*$-algebraic technique of \cite{CEY} is operator-theoretic, the approach in \cite{avila2016dry} is inherently dynamical: it exploits the projective action of $\mathrm{SL}(2,\mathbb{R})$ cocycles.
In this paper, we adapt and significantly extend the method of \cite{avila2016dry} to prove the DTMP in a robust, geometric manner.
We also explain why we adopt the method of periodic approximation, which is partly motivated by a question of M.~Shamis. 

For a periodic Schr\"odinger operator with frequency $p/q$, it is well known that
$G_k(p/q) = \bigcap_{x\in\T} G_k(p/q,x)$, where $G_k(p/q,x)$ denotes the spectral gap of the operator $H_{v,p/q,x}$ with labeling $k$.  
In general, if the band functions oscillate wildly with respect to $x$, it may happen that $S_+= \bigcup_{x\in\T}\sigma(p/q,x)$ fills an entire interval, while $G_k(p/q)$ may vanish entirely.
For the AMO, such a scenario is impossible, indeed, by Chambers' formula, the $k$-th gap satisfies 
$G_k(p/q)= G_k(p/q,x_0)$
for some  $x_0\in\T$, so the gap persists uniformly in $x$.
During the preparation of this manuscript, M. Shamis raised the question of whether periodic gaps inevitably collapse as the phase $x$ varies.  The present work offers a partial answer to her inquiry, relying on Lemma \ref{lem:gap-esti} and Lemma \ref{lem:local-stability} to characterize two complementary regimes of spectral stability: Lemma \ref{lem:gap-esti} establishes a quantitative lower bound on the gap width $G_k(p/q,x)$, estimating the separation between the endpoints $E_+^{k}(x)$ and $E_-^{k}(x)$. In contrast, Lemma \ref{lem:local-stability} establishes the spatial clustering of the spectrum, controlling the variation of the endpoint $E_\pm^{k}(x)$ across different base points $x, x'$. This approach is 
applicable to any Type I trigonometric polynomial potentials (covering both the subcritical and supercritical cases).

\subsection{Methodology: Global Symplectification and Infinite-Dimensional Limit}

Our proof builds on the duality approach to global theory developed in \cite{GJYZ}, the Type I framework introduced in \cite{GJY}, and the theory of the dual center developed in \cite{GJ}. Building on \cite{GJYZ}, Ge-Jitomirskaya-You \cite{GJY} introduced $T$-acceleration and the Type I class, isolating the natural global regime corresponding to a two-dimensional dual center and formulating the Type I Ten Martini program. The identification of the two-dimensional center and its hidden subcriticality in \cite{GJ} are essential inputs. For general analytic potentials, we also use the intrinsic center construction and convergence established there.

Our contribution is to develop the quantitative and geometric tools needed to turn this structure into gap-opening estimates. These include a global symplectic trivialization preserving monotonicity, quantitative Aubry duality with bounds uniform in the truncation dimension, and control of projective dynamics under premonotonicity. Together, these tools yield quantitative nondegeneracy and gap estimates that remain effective as the trigonometric-polynomial truncations converge to a general analytic potential.

\subsubsection{The Geometric Framework: Global Symplectification for monotonic cocycles and Projective Dynamics}

According to the recently developed quantitative global theory \cite{GJYZ}, for a trigonometric polynomial potential of degree $d$, the dual cocycle $(d\alpha, A^{E,d})$ is partially hyperbolic with a two-dimensional center bundle. A (Hermitian) symplectic quasi-periodic cocycle $A$ is said to admit a \emph{partially hyperbolic splitting} if there exists an $A$-invariant dominated splitting
$$E^u(x)\oplus E^c(x)\oplus E^s(x)$$
where $A$ uniformly expands $E^u$ and uniformly contracts $E^s$, and the hyperbolic bundle $E^u\oplus E^s$ is symplectically orthogonal to the center bundle $E^c$. 

To extract spectral information, we must reduce the dynamics to this center bundle. However, the local coordinates of the reduced system do not naturally preserve a strict symplectic and Schr\"odinger structure. The natural substitute for the Schr\"odinger structure is monotonicity. This concept was originally introduced for $\mathrm{SL}(2,\R)$ cocycles by Avila and Krikorian \cite{AK}, and was later extended to higher-dimensional symplectic settings in \cite{Xu2019Density} and to the Hermitian symplectic setting in \cite{WAXZ}. Concretely, a $C^1$ one-parameter family of cocycles $A_t(x)\in C^1\bigl(\I\times\T,\mathrm{HSp}(2d,\psi)\bigr)$ is called \emph{monotonic} if, for every isotropic vector $v\in\C^{2d}$ and every base point $x$, the quadratic form
$$\Psi_{A_t}(v,x):=\psi\bigl(A_t(x)v, \partial_t A_t(x)v\bigr)$$
is strictly positive.

To overcome this, we introduce a
\textbf{global symplectification for monotonic cocycles} (Proposition \ref{prop:globalmonoframe}), whose proof we refer to as the \textbf{Holonomy‐Driven Parallel Transport}. This proof adapts techniques from geometric ideas to the setting of dynamical systems via the following three steps: 
  \begin{itemize} \item \textbf{Graph‐Transform Connection.} We first construct a near‐parameter “graph‐transform’’ connection $\widehat\nabla_{t,s,\alpha} \colon E_{s,\alpha}^c \;\longrightarrow\; E_{t,\alpha}^c, $ defined as the graph of a bundle map over a fixed reference splitting (Lemma \ref{lem:graphy}).
  \item \textbf{Symplectic Correction.} Next, we apply a quantitative implicit‐function theorem to perturb $\widehat\nabla_{t,s,\alpha}$ into an exactly symplectic connection, thereby producing the desired symplectic parallel transport (Lemma \ref{lem:N_s(t)}). 
  
  \item \textbf{Path‐Ordered Holonomy.} Finally, we trivialize globally by subdividing the spectral‐parameter interval $\I$ into finitely many subintervals and concatenating the local transports (Lemma \ref{center-frame}). The resulting holonomy matrix $R_{(t,s)}(x)$ records the path‐dependence in direct analogy with curvature in classical differential geometry. \end{itemize}

  To define the global symplectification, our approach parallels the familiar extension of ordinary‐differential‐equation solutions via parallel transport: one seeks a symplectic trivialization for the non‐autonomous linear system. The connection $\widetilde\nabla_{t,s,\alpha}$ generalizes this ODE parallel transport, and the key estimate (Lemma \ref{lem:localR_t}) \[ \bigl\|R_{t,s,\alpha}(x)\bigr\|_{C^0} \;\le\; e^{\eta\,\lvert t-s\rvert}, \quad \eta=\eta(\I), \] ensures that the holonomy matrices remain uniformly bounded, which in turn yields the desired global trivialization.

For Type I operators, Ge-Jitomirskaya \cite{GJ}
previously introduced a pair of fibered rotation numbers
for the two-dimensional dual center and studied their
relation to the IDS.
Li-Wu \cite{LW} developed a rotation-number
framework for Hermitian-symplectic cocycles in arbitrary
dimension, establishing a direct correspondence with
the IDS of generalized Schr\"odinger operators. The present paper uses a discrete-time formulation of
this framework to study gap opening from a
Hermitian-symplectic perspective, thereby providing
what can be viewed as a natural extension of the
projective-dynamical approach of \cite{avila2016dry}.

Once the reduction to a two-dimensional center is in place, monotonicity is then used again to prove that these labelled gaps are in fact open (Lemma \ref{lem:gap-esti}).
To establish the local stability of gap endpoints (Lemma~\ref{lem:local-stability}), we use direct projective-action estimates. Subcriticality controls the dependence of the projective action on the base point $x$, while monotonicity provides a uniform lower bound on its energy derivative. Together, these estimates control the variation of periodic gap endpoints with the phase.

\subsubsection{From Trigonometric Polynomials to Analytic Potentials}

Extending the above geometric framework from trigonometric polynomial potentials to general analytic potentials requires approximating $v$ by its truncations $v_d$ and passing to the limit $d\to\infty$. The main obstruction is that several quantitative estimates obtained for the finite-dimensional dual cocycle depend implicitly on the truncation dimension. This manifests itself in two distinct ways:
\begin{enumerate}
    \item \textbf{Analytic Degeneration.}
    Quantitative Aubry duality is one of the principal ingredients in the proof of the DTMP \cite{avila2016dry}. However, a direct extension to the $2d$-dimensional dual cocycle leads to constants that deteriorate with $d$. This couples the periodic approximation requirements to the truncation scale, causing the error bounds to degenerate and preventing a rigorous uniform limit.

      \item \textbf{Loss of Strict Monotonicity:} 
     When approximating an analytic potential by trigonometric polynomials of degree $d$, establishing the DTMP requires showing that the spectral gaps of the approximants are substantially larger than $e^{-O(d)}$. 
    For high-degree truncations, the conjugated central cocycle becomes merely premonotonic.
    In this context, tracking the monotonicity of the projective phase dynamics, specifically, iterating the dual cocycle to establish uniform hyperbolicity makes the estimate exceptionally delicate.
    Its one-step projective action is no longer strictly increasing with respect to the energy parameter $E$. A naive iteration fails because any initial negative phase drift is actively accelerated by the unperturbed parabolic fixed point, which acts as a repeller and destroys the invariant cone field necessary for uniform hyperbolicity.
\end{enumerate}

\subsubsection{Taming the Infinity: Quantitative Stability of the  Center Cocycle}

To resolve the dimension dependence, the first step is to identify an object that admits a stable limit as $d\to\infty$. A fundamental observation of \cite{GJ} is that, although the  dual cocycles have increasing dimension, their two-dimensional center cocycles converge. This limiting object, referred to as the \emph{intrinsic center cocycle}, provides a canonical finite-dimensional reduction of the infinite-range dual system. Our quantitative analysis is built upon this convergence and consists of the following four ingredients.

\begin{itemize}[leftmargin=1cm]
  \item \textbf{Dimension-Free Aubry Duality via Weighted Norms:} A principal ingredient is a quantitative Aubry-duality estimate for approximate scalar reductions, valid at irrational frequencies and suitable rational approximants, with bounds uniform in the truncation dimension (Proposition~\ref{dim-duality}). We use the weighted analytic norm $\tnorm{f}_{\epsilon} = \sum_{k=-d}^{d-1} e^{-\xi|k|}\|f(k,\cdot)\|_{\epsilon}$. These weights combine the analytic decay of the potential with the subexponential spatial growth of the center frame to control both accumulated component errors and normalization estimates without losses depending on the dimension. This allows the duality estimate to remain effective in the limit of general analytic potentials.

 \item \textbf{Center Estimates Independent of Dimension.} The key observation is that the bounds for the center coordinate changes
depend only on the two center columns, without involving the stable or
unstable frames
(Lemmas~\ref{lem:detQ} and~\ref{lem:localR_t}).
On energy intervals shrinking to the target gap endpoint, the exponential
decay of the hopping coefficients gives uniform bounds for the symplectic
pairings of these columns (Lemma~\ref{lem:psi-bound}); Cauchy estimates
give the corresponding energy derivative bounds.
These bounds control the symplectic correction and transport independently
of the truncation dimension
(Lemma~\ref{lem:N_s(t)} and Proposition~\ref{prop:10.3}).
We construct the transport in one reference dimension, set it equal to
the identity at the target energy, and apply the same two-dimensional
change of coordinates to all larger truncations. The reduced cocycles
then converge, with uniform bounds for the scalar phases and energy
dependence (Proposition~\ref{prop:10.3}). At the target energy, the limiting
cocycle is independent of the reference dimension and the transport mesh.
Thus the almost-reducibility data can be fixed first
(Lemma~\ref{lem:unified-ar}), and the mesh can then be refined to control
the monotonicity error without changing these data.
    \item \textbf{The Discrete Hellmann-Feynman Identity:} To bypass the lack of one-step monotonicity, we apply a discrete analogue of the Hellmann-Feynman theorem (Lemma \ref{lem:pre-mono}). By projecting the $2d$-step variations onto the Gram matrix of the center basis, we extract a strict $e^{-o(d)}$ lower bound for the projective action, despite the intermediate premonotonicity.
      \item \textbf{Unified Cone Field Argument via Pre-monotonicity:} 
    We develop a cone argument adapted to pre-monotonic cocycles. Instead of requiring one-step monotonicity, we control the intermediate phase evolution quantitatively and prove that the orbit remains inside the invariant cone until the monotonicity of the $2d$-th iterate becomes effective (Lemma~\ref{lem:orbit}). Since the argument is independent of the arithmetic properties of $\alpha$, it applies uniformly to both the Diophantine ($\beta(\alpha)=0$) and Liouvillean ($\beta(\alpha)>0$) regimes.
\end{itemize}

\subsection{Organization of the Paper}
The paper is organized as follows. Section~\ref{sec:pre} introduces the necessary preliminaries and notation. Section~\ref{sec:global-symp} constructs a global symplectic trivialization, which is refined in Section~\ref{sec:global-symp-mono} to preserve monotonicity, a key property for the subsequent analysis. Section~\ref{sec:proj-rotationnumber} develops the projective-action tools used in the spectral analysis, following \cite{LW,avila2016dry}.  Section~\ref{sec:periodicapprox} establishes gap labeling and derives quantitative estimates for spectral gaps, while Section~\ref{sec:aubry} develops a quantitative version of Aubry duality. These ingredients are combined in Section~\ref{sec:endpf} to prove Theorem~\ref{main-thm-analytic} for trigonometric polynomial potentials, and the proof is completed for general analytic potentials in Section~\ref{analytic-proof}. Finally, the Appendix provides proofs for several standard results that are well-known to the community but lack explicit references in the literature.

\section{Preliminaries}\label{sec:pre}

Unless otherwise noted, we use \(\|\cdot\|\) to denote the spectral norm.

\subsection{Linear cocycles and Lyapunov exponents}

Let $ G $ be a Lie group. Given $ \alpha \in \mathbb{R} \setminus \mathbb{Q} $ and $ A \in C^\omega(\mathbb{T}, G) $, we define the $G$-valued one-frequency cocycle $ (\alpha, A) $ as follows:
$$
(\alpha, A) \colon \left\{
\begin{array}{rcl}
\mathbb{T} \times \mathbb{C}^{m} & \to & \mathbb{T} \times \mathbb{C}^{m} \\[1mm]
(x, v) & \mapsto & (x + \alpha, A(x) \cdot v)
\end{array}
\right. .
$$
The iterates of $ (\alpha, A) $ take the form $ (\alpha, A)^n = (n\alpha, A_n) $, where
$$
A_n(x) :=
\begin{cases}
A(x + (n-1)\alpha) \cdots A(x + \alpha) A(x), & n \geq 0 \\[1mm]
A^{-1}(x + n\alpha) A^{-1}(x + (n+1)\alpha) \cdots A^{-1}(x - \alpha), & n < 0
\end{cases} .
$$
We denote the Lyapunov exponents of $ (\alpha, A) $ by $ L_1(A) \geq L_2(A) \geq \ldots \geq L_m(A) $, ordered according to their multiplicities. These exponents are defined as follows:
$$
L_k(A) = \lim_{n \to \infty} \frac{1}{n} \int_{\mathbb{T}} \ln \sigma_k(A_n(x)) \, dx,
$$
where $ \sigma_1(A_n) \geq \ldots \geq \sigma_m(A_n) $ denote the singular values (the eigenvalues of $ \sqrt{A_n^* A_n} $).

Let $(\alpha, A)$ be a cocycle. An $A$-invariant continuous splitting
\[
\mathbb{C}^{d} = E^1(x) \oplus \cdots \oplus E^k(x)
\]
is called \textit{dominated} if there exists an integer $n \geq 1$ such that for any $x$, any $1 \leq j < k$, and any unit vectors $w_j \in E^j(x)$ and $w_{j+1} \in E^{j+1}(x)$, we have:
\[
\| A_n(x) w_j \| > \| A_n(x) w_{j+1} \|.
\]
A cocycle $(\alpha, A)$ is said to be \textit{partially hyperbolic} if there exists an $A$-invariant dominated continuous splitting
\[
\mathbb{C}^{d} = E^u(x) \oplus E^c(x) \oplus E^s(x),
\]
where $E^u$, $E^c$, and $E^s$ denote the unstable, center, and stable bundles, respectively. Moreover, there exist constants $C > 0$ and $c > 0$ such that for every $n \geq 0$, the following contraction properties hold:
\begin{align*}
    \|A_n(x)v\| &\leq C e^{-cn}\|v\|, \quad \text{for } v \in E^s(x), \\
    \|A_n(x)^{-1}v\| &\leq C e^{-cn}\|v\|, \quad  \text{for } v \in E^u(x+n\alpha).
\end{align*}
In the special case where the center bundle vanishes (i.e., the phase space decomposes as $\mathbb{C}^{d} = E^u(x) \oplus E^s(x)$), the cocycle $(\alpha, A)$ is said to be \textit{uniformly hyperbolic}.

\subsection{(Hermitian) symplectic group actions}\label{sec:symp}

In this paper, we focus primarily on the cases where $ G = \mathrm{Sp}(2d, \mathbb{F}) $(where  $\mathbb{F} = \mathbb{R}$  or  $\mathbb{C}$), $ \mathrm{HSp}(2d) $, or $ \mathrm{SL}(2, \mathbb{R}) $.
The symplectic group, denoted $ \mathrm{Sp}(2d,\mathbb{F}) $, consists of matrices $ M $ that preserve the standard symplectic matrix $ \mathcal{J}_{2d} $, which is defined as:
$$
M^{T} \mathcal{J}_{2d} M = \mathcal{J}_{2d}, \quad \text{where} \quad \mathcal{J}_{2d} = \begin{pmatrix} 0 & -I_d \\ I_d & 0 \end{pmatrix}.
$$
The Hermitian-symplectic group is the subgroup of complex general linear matrices that preserve the standard symplectic matrix in the Hermitian sense:
$$
\mathrm{HSp}(2d) := \left\{ M \in \mathrm{GL}(2d, \mathbb{C}) : M^{*} \mathcal{J}_{2d} M = \mathcal{J}_{2d} \right\},
$$
where $ M^{*} = \overline{M}^{\top} $ denotes the conjugate transpose. The associated sesquilinear, skew-Hermitian symplectic form $ \omega_d $ is given by:
$$
\omega_d(u, v) = u^{*} \mathcal{J}_{2d} v \quad \text{for } u, v \in \mathbb{F}^{2d}.
$$
When the underlying form is a (Hermitian) skew-symmetric form 
$\psi$ with associated structure matrix $\mc{S}$ congruent to $\mc{J}_{2d}$ (the only case considered in this paper), 
we indicate this dependence by writing $\mathrm{HSp}(2d,\psi)$. In the notable special case where $ \psi $ is the Hermitian form induced by the diagonal matrix
$
I_{d,d} = \operatorname{diag}(I_d, -I_d),
$
the corresponding preserving group is the pseudo-unitary group:
$$
\mathbf{U}(d, d) = \{ M \in \mathrm{GL}(2d, \mathbb{C}) : M^{*} I_{d,d} M = I_{d,d} \}.
$$

Let \(\operatorname{Lag}(\mathbb{F}^{2d},\omega_d)\) denote the Lagrangian Grassmannian of \((\mathbb{F}^{2d},\omega_d)\).  A \emph{Lagrangian frame} of a given \(\Lambda\in\operatorname{Lag}(\mathbb{F}^{2d},\omega_d)\) is an injective linear map
\(\mathcal L:\mathbb{F}^{d}\to\Lambda\) of the form 
$\mathcal L=\begin{psmallmatrix} X \\ Y \end{psmallmatrix},$
where \(X\) and \(Y\) are \(d\times d\) matrices over \(\mathbb F\) satisfying \(X^*Y=Y^*X\) and \(\operatorname{rank}\mathcal L=d\).  Two frames \(\mathcal L_1,\mathcal L_2\) are declared equivalent, \(\mathcal L_1\sim\mathcal L_2\), if \(\mathcal L_1=\mathcal L_2 R\) for some \(R\in\mathrm{GL}(d,\mathbb F)\).  In this way one may represent a Lagrangian subspace \(\Lambda\) by the equivalence class
\(\Lambda=\begin{psmallmatrix} X\\ Y\end{psmallmatrix}_{\sim}\), and for brevity we will often write simply
\(\Lambda=\begin{bsmallmatrix} X\\ Y\end{bsmallmatrix}\).
The natural actions of \(\mathrm{Sp}(2d,\R)\) and \(\mathrm{HSp}(2d)\) leave \(\operatorname{Lag}(\R^{2d},\omega_d)\) and \(\operatorname{Lag}(\C^{2d},\omega_d)\) invariant, respectively.  

Consider the Cayley element
\begin{equation}\label{eq:Cayley}
\mathcal{C}:=\frac{1}{\sqrt{2i}}\begin{pmatrix} I_d & -iI_d\\[4pt] I_d & iI_d \end{pmatrix}.
\end{equation}
For any \(2d\times 2d\) complex matrix \(A\) set \(\Aa:=\mathcal{C}A\mathcal{C}^{-1}\).  Conjugation by \(\mathcal C\) identifies the real groups \(\mathrm{Sp}(2d,\R)\), \(\mathrm{HSp}(2d)\) with the subgroups \(\mathbf{U}(d,d)\cap\mathrm{Sp}(2d,\C)\) and \(\mathbf{U}(d,d)\), respectively (more precisely, \(A\mapsto\Aa\) is a Lie-group isomorphism onto its image).

Introduce the following `chart' of the Lagrangian Grassmannian:
\begin{align*}
\overset{\circ}{\operatorname{Lag}}(\R^{2d},\omega_d)
&=\Big\{\begin{pmatrix} M\\[2pt] I_d\end{pmatrix}_{\sim}: M\in \mathbf{U}(d)\cap\operatorname{Sym}_d(\C)\Big\},\\[4pt]
\overset{\circ}{\operatorname{Lag}}(\C^{2d},\omega_d)
&=\Big\{\begin{pmatrix} M\\[2pt] I_d\end{pmatrix}_{\sim}: M\in \mathbf{U}(d)\Big\}.
\end{align*}
The Cayley element \(\mathcal C\) induces a diffeomorphism (respecting the relevant group actions) from \(\operatorname{Lag}(\mathbb F^{2d},\omega_d)\) onto the corresponding \(\overset{\circ}{\operatorname{Lag}}(\mathbb F^{2d},\omega_d)\): if
\(\Lambda=\begin{bsmallmatrix} X\\ Y\end{bsmallmatrix}\) then
\[
\overset{\circ}{\Lambda}=\mathcal C\Lambda=\begin{bsmallmatrix} I_d\\[2pt] W_\Lambda\end{bsmallmatrix},
\qquad\text{where}\qquad
W_\Lambda=(X+iY)(X-iY)^{-1}.
\]
Hence the actions of \(\mathbf{U}(d,d)\cap\mathrm{Sp}(2d,\C)\) and \(\mathbf{U}(d,d)\) preserve \(\overset{\circ}{\operatorname{Lag}}(\R^{2d},\omega_d)\) and \(\overset{\circ}{\operatorname{Lag}}(\C^{2d},\omega_d)\), respectively.  

 \begin{lemma}[{\cite[Lemma 4.7]{LW}}]\label{lemWlambda}
	Let $ \Lambda=\begin{bsmallmatrix}
		X\\Y
	\end{bsmallmatrix} $ be a Lagrangian subspace. Then we have
	\begin{enumerate}
		\item 
		$ \det W_\Lambda=\det(X+i Y)\det(X^*+i Y^*) \det (X^*X+Y^*Y)^{-1}$. 
		\item $ W_\Lambda z=z $ if and only if $ Y^*z=0 $. In particular, the eigenspace of $ W_\Lambda $ associated to $ 1 $ agrees with the kernel of $ Y^* $. \label{lemWlambda:2}
        \item $ W_\Lambda z=-z $ if and only if $ X^*z=0 $. In particular, the eigenspace of $ W_\Lambda $ associated to $ -1 $ agrees with the kernel of $ X^* $. \label{lemWlambda:3}
	\end{enumerate}
\end{lemma}  

\subsection*{Direct product operation ($\diamond$)}

The direct product operation $\diamond$ is defined for both Lagrangian subspaces and Hermitian symplectic matrices as follows.

\vspace{0.5\baselineskip}

\noindent
Let $\Lambda_d=\begin{bmatrix} X_d\\ Y_d \end{bmatrix} \in \mathrm{Lag}(\mathbb{C}^{2d}, \omega_d)$ and $\Lambda_m=\begin{bmatrix} X_m\\ Y_m \end{bmatrix} \in \mathrm{Lag}(\mathbb{C}^{2m}, \omega_m)$. Their direct product $\Lambda_d\diamond \Lambda_m \in \mathrm{Lag}(\mathbb{C}^{2d+2m}, \omega_{d+m})$ is defined as:
\[
\Lambda_d\diamond \Lambda_m=\begin{bmatrix}
    X_d & 0\\
    0 & X_m\\
    Y_d & 0\\
    0 & Y_m
\end{bmatrix}.
\]
Let $A=\begin{pmatrix} X_d & Z_d \\ Y_d & W_d \end{pmatrix} \in \mathrm{HSp}(2d)$ and $B=\begin{pmatrix} X_m & Z_m \\ Y_m & W_m \end{pmatrix} \in \mathrm{HSp}(2m)$. Their direct product $A\diamond B \in \mathrm{HSp}(2d+2m)$ is defined as:
\[
A\diamond B=\begin{pmatrix}
    X_d & 0 & Z_d & 0\\
    0 & X_m & 0 & Z_m\\
    Y_d & 0 & W_d & 0\\
    0 & Y_m & 0 & W_m
\end{pmatrix}.
\]

\subsection{Aubry duality}

Let $\alpha\in\mathbb R$, $v,w:\T\to \R$ be real analytic functions with Fourier coefficients $\{\hat{v}_k\}_{k\in\Z}$ and $\{\hat{w}_k\}_{k\in\Z}$. Consider the Hilbert space $\mathcal H = L^2(\mathbb T \times \mathbb Z)$, equipped with the norm
\begin{equation}\label{eq:norm_def}
\lVert\Psi\rVert_{\mathcal H}^2 = \sum_{n\in\mathbb Z} \int_{\mathbb T} |\Psi(x,n)|^2\, dx < \infty.
\end{equation}
For each $x\in\mathbb T$, let $L_x$ be the operator on $\ell^2(\mathbb Z)$ defined by
\begin{equation}\label{eq:Lx_def}
(L_x u)_n = \sum_{k\in\mathbb Z} \hat v_k\, u_{n+k} + w(x+n\alpha)\, u_n.
\end{equation}
Let $\mathcal L = \int_{\mathbb T}^{\oplus} L_x\, dx$ be the associated direct integral operator on $\mathcal H$. A standard fact (see also \cite{JM2012} for details) is that its spectrum satisfies
\begin{equation}\label{eq:spec_union}
\sigma(\mathcal{L}) = \bigcup_{x\in\mathbb T}\sigma(L_x).
\end{equation}
The \textit{dual operator} $\mathcal L^\ast = \int_{\mathbb T}^{\oplus} L_x^\ast\, dx$ is similarly defined via the family
\[
(L_x^\ast u)_n = \sum_{k\in\mathbb Z} \hat w_k \, u_{n+k} + v(-x-n\alpha)\, u_n .
\]
Aubry duality is implemented by the map $\mathcal U:\mathcal H\to\mathcal H$
\begin{equation*}
(\mathcal U\Psi)(x,n) = \sum_{m\in\mathbb Z} e^{-2\pi i m (x+n\alpha)} \int_{\mathbb T} e^{-2\pi i n\theta}\, \Psi(\theta,m)\, d\theta .
\end{equation*}
A direct computation verifies that $\mathcal U$ is a unitary operator and induces the unitary equivalence $\mathcal U^{*}\, \mathcal L\, \mathcal U = \mathcal L^\ast$, which implies the spectral identity
$\sigma(\mathcal{L}) = \sigma(\mathcal{L}^*).$

\subsection{The integrated density of states}\label{sec:ids-revisited}

Consider $L_x=L_{v,w,\alpha,x}$ acting on $\ell^{2}(\mathbb Z)$,
\[
(L_{v,w,\alpha,x} u)_n
=
\sum_{k\in\mathbb Z} \hat v_k\, u_{n+k}
+
w(x+n\alpha)\, u_n .
\]
For any $\Lambda\subseteq\mathbb Z$, let
$P_\Lambda : \ell^2(\mathbb Z) \to \ell^2(\mathbb Z)$
denote the canonical orthogonal projection. On the standard basis
$\delta_n$,
\[
P_\Lambda \delta_n =
\begin{cases}
\delta_n, & n\in\Lambda,\\[4pt]
\mathbf 0, & n\notin\Lambda.
\end{cases}
\]
Define the finite-volume (Dirichlet) restriction by
\[
L_x^\Lambda
:=
\bigl(P_\Lambda L_x P_\Lambda^*\bigr)\big|_{\Lambda\times\Lambda},
\qquad
L_x^N := L_x^{[1,N]}.
\]
Let $E_1^N\le \cdots \le E_N^N$ be the eigenvalues of $L_x^N$ (counted with multiplicity).  
The integrated density of states (IDS) of $L_x$ is defined by
\begin{equation}\label{eq:idsdef}
N_x^L(E)
:=
\lim_{N\to\infty}
\frac{1}{N}
\#\{\, j : E_j^N \le E \,\}.
\end{equation}
When $\alpha$ is irrational, this limit converges uniformly in $x$ and is independent of $x$.

Specially, assume $v=v_d$ is a real trigonometric polynomial of degree $d$, so that  
$L_{v,w,\alpha,x}$ can be rewritten as a strip operator on $\ell^2(\mathbb Z,\mathbb C^d)$
of the form \eqref{eq:hatH}, which we denote by $\widehat H_x$.
For $\Lambda\subseteq\mathbb Z$, let
$P_\Lambda : \ell^2(\mathbb Z,\mathbb C^d)
\to
\ell^2(\mathbb Z,\mathbb C^d)$
be the canonical projection. On the basis vectors $\delta_{n,i}$ (the $i$-th standard basis at site $n$),
\[
P_\Lambda \delta_{n,i} =
\begin{cases}
\delta_{n,i}, & n\in\Lambda,\\[4pt]
\mathbf 0, & n\notin\Lambda.
\end{cases}
\]
Define the finite-volume restriction of $\widehat H_x$ by
\[
\widehat H_x^\Lambda
:=
\bigl(P_\Lambda \widehat H_x P_\Lambda^*\bigr)\big|_{\Lambda\times\Lambda},
\qquad
\widehat H_x^N := \widehat H_x^{[1,N]}.
\]
Note that $\widehat H_x^N = L_x^{dN}$.
The IDS of $\widehat H_x$ is given by
\begin{equation}\label{eq:idsdef-strip}
N_x^{\widehat H}(E)
:=
\lim_{N\to\infty}
\frac{1}{dN}
\#\{\, j : E_j^{dN} \le E \,\},
\end{equation}
which necessarily coincides with $N_x^L(E)$.

\section{Global Symplectification  for  Hermitian Symplectic cocycles}\label{sec:global-symp}

As proved in \cite[Proposition~3.2]{WAXZ}, any analytic one-frequency quasiperiodic
Hermitian-symplectic cocycle that is partially hyperbolic admits an analytic
block diagonalization.  
Here, we strengthen this result by showing that parameterized quasiperiodic
Hermitian-symplectic cocycles can be globally analytically block diagonalized.

To make this precise, for any $c, h > 0$, we denote the complex neighborhoods $\mathbb{T}_h := \{x \in \mathbb{C}/\mathbb{Z} : |\text{Im}\, x| < h\}$ and $\I(c) := \{t \in \mathbb{C} : \operatorname{dist}(t, \I) < c\}$. Let $C^\omega_{c,h}(\I \times \mathbb{T}, *)$ be the Banach space of $*$-valued functions on $\I\times\T$ that extend holomorphically to $\I(c) \times \mathbb{T}_h$, equipped with the norm
\[
\|A\|_{c,h} = \sup_{(t,x) \in \I(c) \times \mathbb{T}_h} \|A(t,x)\|.
\]
Similarly, for a given interval $\mathcal{O} \subset \mathbb{R}$, we define $C^{\omega}_{c,0,h}(\I \times \mathcal{O} \times \mathbb{T}, *)$ as the Banach space of $*$-valued functions on $\I \times \mathcal{O} \times \mathbb{T}$ that are jointly continuous in $(t, \alpha, x)$ and extend holomorphically to $ \I(c)\times\T_h$ in $(t,x)$, under the norm
\[
\|A\|_{c,0,h} = \sup_{(t,\alpha,x) \in \I(c) \times \mathcal{O} \times \mathbb{T}_h} \|A(t,\alpha,x)\|.
\]
With these conventions in place, we obtain the following result:
\begin{proposition}\label{prop:globalsympframe}
Let \( \I \subset \mathbb{R} \) be a compact interval, and let  
\( A_t(x)\in C_{c,h_0}^\omega(\I \times \mathbb{T}, \mathrm{HSp}(2d)) \). Assume the cocycle \( (\alpha_0, A_t(\cdot+i\epsilon)) \) is partially hyperbolic with a $2m$-dimensional center bundle for all $t\in\I$ and $|\epsilon|<h_0$. 
Then for any $0<h<h_0$, there exist an open neighborhood $\mathcal{O}'=\OO'(\I,\alpha_0,h)$ of $\alpha_0$ and $\delta'(\I,\alpha_0,h)>0$, $h'(\I,\alpha_0,h)\in (0,h] $ such that for all $\alpha\in \OO'$, the following hold:
\begin{enumerate}
    \item The cocycle \( (\alpha, A_t(\cdot+i\epsilon)) \) is partially hyperbolic with a $2m$-dimensional center bundle for all $t\in\I$ and $|\epsilon|<h$.
    \item The center bundle \( E^c_{t,\alpha} \) admits a canonical symplectic basis of vectors
        $
        u^{t,\alpha}_{\pm i} \in C^\omega_{\delta',0,h'}(\I \times \OO' \times \mathbb{T}, \mathbb{C}^{2d})$, for $ i = d-m+1, \dots, d.$
    \item The hyperbolic bundle \( E^u_{t,\alpha} \oplus E^s_{t,\alpha'} \) admits a canonical symplectic basis of vectors $
        u^{t,\alpha}_{\pm i} \in C^\omega_{\delta',0,h}(\I \times  \OO'\times \mathbb{T},\mathbb{C}^{2d})$, for $ i = 1, \dots, d-m.
        $
\end{enumerate}
\end{proposition}

\subsection{Main ingredients} 

To block diagonalize a partially hyperbolic Hermitian--symplectic cocycle,
the proof proceeds in two steps.  
First, one shows that the center bundle and the stable/unstable bundles admit
holomorphic trivializations with a fixed analytic radius.  
Second, one constructs a symplectic basis adapted to this splitting.
Here we present the main ingredients of the argument; they are of independent
interest and will be used repeatedly throughout the paper.

\subsubsection{Holomorphic Structure of the Grassmannian}
To obtain holomorphic trivializations with a fixed analytic radius,
the key ingredient concerns the holomorphic structure of the Grassmannian. Let $\mathcal{M}_k(d)$ denote the set of $d \times k$ complex matrices of rank $k$. Let $G(k,d)$ denote the set of $k$-dimensional subspaces of $\mathbb{C}^d$, viewed as a compact complex manifold (the Grassmannian) with its unique holomorphic structure \cite{AJS2014}. We are interested in the trivialization of the Grassmannian, and the proof follows the spirit of the Oka-Grauert principle (topologically trivial holomorphic vector bundles over Stein manifolds are holomorphically trivial) \cite{Forstneric2017}.

\begin{proposition}\label{lem:ajslem}
Let $\Omega\subset\C^n$ be a contractible Stein domain and $\delta>0$.
Suppose $
u:\Omega\times\T\longrightarrow G(k,d)$ 
is a map which extends holomorphically to $\Omega\times\T_\delta$.  Then $u$ admits a holomorphic lift
$\widetilde u:\Omega\times\T_\delta\longrightarrow \mathcal{M}_k(d).$ Equivalently, for every $(x,z)\in\Omega\times\T_\delta$ the columns of $\widetilde u(x,z)$ form a holomorphically varying basis of the subspace $u(x,z)\subset\C^d$.
\end{proposition}

\begin{proof}
The map \(u:\Omega\times\T_\delta\to G(k,d)\) determines a holomorphic rank-\(k\) subbundle 
$E \;=\; u^*(\mathcal T)\;\subset\;(\Omega\times\T_\delta)\times\C^d,$
where \(\mathcal T\to G(k,d)\) is the tautological holomorphic vector bundle and \(u^*(\mathcal T)\) is its pullback. Thus \(E\) is a holomorphic vector bundle of rank \(k\) over the complex manifold \(\Omega\times\T_\delta\).

By assumption, $\Omega$ is Stein. Since $\T_\delta$ is an open Riemann surface, it is also Stein. The product of Stein manifolds is Stein, so $\Omega\times\T_\delta$ is Stein. See \cite[Ch.~5]{Forstneric2017}.
On the other hand,
the strip \(\T_\delta\) is homotopy equivalent to the circle \(\T\). Since \(\Omega\) is contractible by hypothesis, the product \(\Omega\times\T_\delta\) is homotopy equivalent to \(\T\). Complex vector bundles over a space \(X\) are classified (up to isomorphism) by homotopy classes of maps \(X\to BU(k)\) (the classifying space for rank-\(k\) complex bundles) \cite{HatcherVB}. Because \(BU(k)\) is simply connected (indeed \(\pi_1(BU(k))\cong\pi_0(U(k))=0\)), any continuous map \(S^1\to BU(k)\) is null-homotopic. Hence every complex vector bundle over \(\T\)  is topologically trivial \cite{HatcherVB,MS1974}. Thus \(E\) is topologically trivial.

By the Oka-Grauert principle the holomorphic bundle \(E\) is holomorphically trivial. Therefore there exists a global holomorphic bundle isomorphism
$
\Phi:E\longrightarrow (\Omega\times\T_\delta)\times\C^k.$
Let \(\{e_1,\dots,e_k\}\) denote the standard basis of \(\C^k\). For each \(i=1,\dots,k\) define the holomorphic section $
s_i(x,z):=\Phi^{-1}\bigl((x,z),e_i\bigr)\in E_{(x,z)}\subset\C^d.$
The sections \(s_1,\dots,s_k\) form a global holomorphic frame of \(E\). Assembling these sections as the columns of a \(d\times k\) matrix gives a holomorphic map
$
\widetilde u(x,z):=\bigl(s_1(x,z)\;\; s_2(x,z)\;\;\cdots\;\; s_k(x,z)\bigr)\in\mathcal M_k(d),
$
and by construction the column span of \(\widetilde u(x,z)\) equals the fiber \(E_{(x,z)}=u(x,z)\). Because the sections are holomorphic and pointwise linearly independent, \(\widetilde u\) is a holomorphic lift of \(u\) with full-rank matrices at every point, as required.
\end{proof}

\subsubsection{Non-degeneracy of the Krein Matrix}

To find symplectic frames in each bundle, the key tool is the so-called \textbf{Krein matrix}. Let $V \subset \mathbb{C}^{2d}$ be a Hermitian symplectic subspace, and let $\{v_1, \dots, v_{2d}\}$ be a basis of $V$. Following M. Krein, for any $u, v \in V$, we define the \textbf{Krein form} by $i u^* \mathcal{J}_{2d}\, v$ \cite{Long}. Consequently, we call the matrix
\begin{equation}\label{GRAM MATRIX}
    G(v_1, \dots, v_{2n}) = \frac{1}{i}
    \begin{bmatrix}
        v_1^* \mathcal{J}_{2d}\, v_1 & v_1^* \mathcal{J}_{2d}\, v_2 & \cdots & v_1^* \mathcal{J}_{2d}\, v_{2d} \\
        v_2^* \mathcal{J}_{2d}\, v_1 & v_2^* \mathcal{J}_{2d}\, v_2 & \cdots & v_2^* \mathcal{J}_{2d}\, v_{2d} \\
        \vdots & \vdots & \ddots & \vdots \\
        v_{2d}^* \mathcal{J}_{2d}\, v_1 & v_{2d}^* \mathcal{J}_{2d}\, v_2 & \cdots & v_{2d}^* \mathcal{J}_{2d}\, v_{2d}
    \end{bmatrix}
\end{equation}
the \textbf{Krein matrix} of this basis. The following observation establishes the non-degeneracy of the holomorphic extension of the Krein matrix.
\begin{lemma}\label{lem:Key-holom-extend}
    Let $\Omega \subset \mathbb{C}^n$ be a contractible Stein domain with $\Omega = \Omega^*$, where $\Omega^* = \{\bar z : z \in \Omega\}$.
    Let $E: \Omega\times \mathbb{T}_\delta \to G(2m, 2d)$ and $F: \Omega\times \mathbb{T}_\delta \to G(2d-2m,2d)$ be holomorphic distributions satisfying
    \[
        \mathbb{C}^{2d} = E(z) \oplus F(z) \quad \text{for all } z \in \Omega\times \mathbb{T}_\delta.
    \]
    Suppose that for real parameters $z \in (\Omega \cap \mathbb{R}^n )\times \mathbb{T}$, $F(z)$ is the symplectic orthogonal complement of $E(z)$. Then:
    \begin{enumerate}
        \item For every $z\in\Omega\times\mathbb T_\delta$ we have
        $F(z)=(E(\bar z))^{\perp_\omega}$.
  \item If $M(z)=[v_1(z)\ \cdots\ v_{2m}(z)]$ is any global holomorphic frame of $E$, then the associated Krein matrix
        \[
        G(z)=\tfrac{1}{i}\,M(\bar z)^* \mathcal{J}_{2d} M(z)
        \]
        is a holomorphic map $\Omega\times\mathbb T_\delta\to\mathrm{GL}(2m,\mathbb C)$.
    \end{enumerate}
\end{lemma}

\begin{proof}
    Since \(E\) and \(F\) are topologically trivial, Proposition \ref{lem:ajslem} implies they are holomorphically trivial. Let \(f\) and \(g\) be arbitrary holomorphic sections of \(F\) and \(E\) respectively. Define the pairing  holomorphic function
    \[
        P(z):= \omega\bigl( g(\bar z), f(z) \bigr)= g(\bar z)^{*} \mathcal{J}_{2d} f(z) , \quad z \in \Omega \times \mathbb{T}_\delta.
    \]
By hypothesis $P\equiv0$ on the real slice $(\Omega\cap\mathbb R^n)\times\mathbb T$, which is a maximal totally real submanifold; hence the identity principle in several complex variables \cite{Shabat92} implies $P\equiv0$ on the full domain.
    This implies $
        F(z) \subset \bigl(E(\bar z)\bigr)^{\perp_\omega}.$
    Since \(\dim F(z) = 2d-2m = \dim E(\bar z)^{\perp_\omega}\), we verify (1).

    For (2),  $G(z)$ is holomorphic extension of  $\tfrac{1}{i}\, M(x)^* J_{2d} M(x),$ which is non-degenerate on the real slice. 
    Suppose instead that \(\det G(z_0)=0\) for some \(z_0 \in \Omega \times \mathbb{T}_\delta\). Then there exists a nonzero vector \(v_0 \in \mathbb{C}^{2m}\) such that
    \[
        w^* G(z_0)\, v_0 = 0 \quad \text{for all } w\in \mathbb{C}^{2m}.
    \]
    Define \(u_0 := M(z_0)v_0 \in E(z_0)\). Since $M(z_0)$ has full rank, $u_0 \neq 0$. The degeneracy condition is equivalent to
    \[
        \omega\bigl(u_0, M(\bar z_0)w \bigr) = 0 \quad \text{for all } w \in \mathbb{C}^{2m},
    \]
    which implies \( u_0 \in E(\bar z_0)^{\perp_\omega} \).
    Using the result from (1), we have \( E(\bar z_0)^{\perp_\omega} = F(z_0) \).
    Thus, \( u_0 \in E(z_0) \cap F(z_0) \).
    Since $\mathbb{C}^{2d} = E(z_0) \oplus F(z_0)$, the intersection must be $\{0\}$, so $u_0 = 0$, which is a contradiction.
    The lemma follows.
\end{proof}

\subsection{Holomorphic trivializations}

We first establish the holomorphic trivializations of the center, stable(unstable) bundle with a fixed analytic radius:

\begin{proposition}\label{prop:globalframe}
Let \( \I \subset \mathbb{R} \) be a compact interval, and let  
\( A_t(x)\in C^\omega_{c,h_0}(\I \times \mathbb{T}, \mathrm{HSp}(2d)) \). Assume the cocycle \( (\alpha_0, A_t(\cdot+i\epsilon)) \) is partially hyperbolic with a $2m$-dimensional center bundle for all $t\in\I$ and $|\epsilon|<h_0$. 
Then for any $0<h<h_0$, there exist $\delta=\delta(\I,\alpha_0,h)>0$, an open neighborhood $\mathcal{O}=\OO(\I,\alpha_0,h) \subset \mathbb{T}$ of $\alpha_0$ such that for $\alpha \in\OO$, the following hold:
\begin{enumerate}
    \item The cocycle \( (\alpha, A_t(\cdot+i\epsilon)) \) is partially hyperbolic with a $2m$-dimensional center bundle for all $t\in\I$ and $|\epsilon|<h$.
    \item The center bundle \( E^c_{t,\alpha} \) admits a canonical basis of vectors
        $
        u^{t,\alpha}_{\pm i} \in C^\omega_{\delta,0,h}(\I  \times\OO\times \mathbb{T}, \mathbb{C}^{2d})$, for $ i = d-m+1, \dots, d.$
    \item The hyperbolic bundle \( E^u_{t,\alpha} \oplus E^s_{t,\alpha} \) admits a canonical basis of vectors $
        u^{t,\alpha}_{\pm i} \in C^\omega_{\delta,0,h}(\I  \times\OO\times \mathbb{T}, \mathbb{C}^{2d})$, for $ i = 1, \dots, d-m.
        $
\end{enumerate}
\end{proposition}

\begin{proof}
To apply Proposition~\ref{lem:ajslem} and deduce Proposition~\ref{prop:globalframe} we only need a uniform regularity statement for the invariant bundles near \(\I\times\{\alpha_0\}\times\T\).  The following lemma gives us desired regularity:

\begin{lemma}\label{prop:regu-bundle}
Let \( \I \subset \mathbb{R} \) be a compact interval, and let  
\( A_t(x)\in C_{c,h}^\omega(\I \times \mathbb{T}, \mathrm{HSp}(2d)) \).  
Assume that the cocycle \( (\alpha_0, A_t(\cdot + i\epsilon)) \) is partially hyperbolic with a $2m$-dimensional center bundle for all \( t \in \I \) and all \( |\epsilon| < h_0 \). Then for any $0<h<h_0$, there exists $\delta=\delta(\I,\alpha_0,h)>0$ and an open neighborhood \( \mathcal{O}=\OO(\I,\alpha_0,h) \subset \mathbb{T} \) of \( \alpha_0 \) such that the partially hyperbolic splitting  
\(
E^s_{t,\alpha} \oplus E^c_{t,\alpha} \oplus E^u_{t,\alpha}
\)
satisfies  
\[
E^s_{t,\alpha},\, E^u_{t,\alpha} \in C^\omega_{\delta,0,h}(\I \times \OO\times \mathbb{T}, G(d-m,2d)), 
\qquad  
E^c_{t,\alpha} \in C^\omega_{\delta,0,h}(\I \times \OO\times\mathbb{T}, G(2m,2d)),
\]

\end{lemma}
\begin{remark}
   The proof is actually the classical cone argument, where the holomorphicity part is essentially the same as \cite[Lemma 6.2]{AJS2014}, for completeness we sketch the proof in Appendix \ref{regu-bund}.  
\end{remark}

Once we have this, let's finish the proof, while we provide the proof only for $E^c_{t,\alpha}(x)$; the remaining cases follow in the same way. For any $ 0<h<h_0$, let $\tilde h=\frac{h+h_0}{2}$, by Lemma \ref{prop:regu-bundle} and Proposition \ref{lem:ajslem}, there exists $\tilde \delta=\tilde \delta(\I,\alpha_0,\tilde h)>0$ and a  neighborhood $\tilde{ \mathcal{O}}=\tilde{\OO}(\I,\alpha_0,\tilde h) \subset \mathbb{T}$ of $\alpha_0$ such that: (i) For any $\alpha\in\tilde{\OO}$, we have continuous splitting $\mathbb{C}^{2d} = E_{t,\alpha}^c(x) \oplus E_{t,\alpha}^s(x) \oplus E_{t,\alpha}^u(x)$; (ii). For $\alpha_0$, there exist frames holomorphically on $\I({\tilde\delta})\times\T_{\tilde h}$:
\[
\{\bar{u}_{-i}^{t,\alpha_0}(x)\}_{i=1}^{d-m} \subset E_{t,\alpha_0}^s(x), \quad 
\{\bar{u}_{i}^{t,\alpha_0}(x)\}_{i=1}^{d-m} \subset E_{t,\alpha_0}^u(x),
\]
and 
\[
\{\bar{u}_{\pm i}^{t,\alpha_0}(x)\}_{i=d-m+1}^{d} \subset E_{t,\alpha_0}^c(x).
\]

Take $\delta=\tilde \delta/2$,  By compactness of $\{\alpha_0\}\times \overline{\I(\delta)}\times \overline{\T_{h}}$ and continuity,
there exists a small neighborhood $\OO=\OO(\I,\alpha_0,h) \subset \tilde{\OO}$ 
such that for all $\alpha \in \OO$, the bundle $E_{t,\alpha}^c(x) $ is close to $E_{t,\alpha_0}^c(x)$, uniformly in $(t,x)$, ensuring transverse intersection:
\[
E_{t,\alpha_0}^c(x)  \cap (E_{t,\alpha}^s(x) \oplus E_{t,\alpha}^u(x)) = \{0\}.
\]

Consequently, let $\mathbb{P}_{E_{t,\alpha}^c(x)}$ be the projection onto $E_{t,\alpha}^c(x)$ along $E_{t,\alpha}^s(x) \oplus E_{t,\alpha}^u(x)$. Since the splitting depends continuously on $\alpha$ and holomorphically on $(t,x)$, $\mathbb{P}_{E_{t,\alpha}^c(x)}$ depends continuously on $\alpha$ and holomorphically on $(t,x)$. Moreover, since the bundle $E_{t,\alpha}^c(x)$ is defined over $\mathbb{T}$, the projection automatically satisfies the periodicity condition $\mathbb{P}_{E_{t,\alpha}^c(x+1)} = \mathbb{P}_{E_{t,\alpha}^c(x)}$.

We now construct the frame for $E_{t,\alpha}^c(x)$ by projecting the reference frame at $\alpha_0$. Define the new vector fields:
\begin{equation}\label{eq:proj_frame}
    u_{\pm i}^{t,\alpha}(x) := \mathbb{P}_{E_{t,\alpha}^c(x)} \left( \bar{u}_{\pm i}^{t,\alpha_0}(x) \right), \quad i=d-m+1,\dots,d.
\end{equation}
Then $u_{\pm i}^{t,\alpha}(x)$ are continuous in $\alpha \in \mathcal{O}$ and holomorphic in $(t,x) \in \I(\delta)\times\T_{h}$. The periodicity in $x$ is preserved as both the projection operator $\mathbb{P}_{E_{t,\alpha}^c(x)}$ and the reference frame $\bar{u}_{\pm i}^{t,\alpha_0}(x)$ are periodic.

It remains to verify that $\{u_{\pm i}^{t,\alpha}(x)\}_i$ forms a frame for $E_{t,\alpha}^c(x)$. For all $\alpha \in \mathcal{O}$, the transversality condition ensures that $E_{t,\alpha_0}^c(x)$ intersects the kernel of the projection $\mathbb{P}_{t,\alpha}^c(x)$ trivially. Hence, the projection restricted to $E_{t,\alpha_0}^c(x)$ is a linear isomorphism onto $E_{t,\alpha}^c(x)$, mapping the reference basis $\{\bar{u}_{\pm i}^{t,\alpha_0}(x)\}_i$ to a basis $\{u_{\pm i}^{t,\alpha}(x)\}_i$. 
\end{proof}

\subsection{Proof of Proposition \ref{prop:globalsympframe}}

We construct the required frames in two steps:

\vspace{0.2cm}
\noindent \textbf{Step I: Symplectic frame for $E_{t,\alpha}^s \oplus E_{t,\alpha}^u$.} \\
For any $0<h<h_0$, by Proposition \ref{prop:globalframe}, there exists $\delta=\delta(\I,\alpha_0,h)>0$ and a  neighborhood $\mathcal{O}=\OO(\I,\alpha_0,h) \subset \mathbb{T}$ of $\alpha_0$ such that one can select frames in $C^\omega_{\delta,0,h}(\I\times\OO\times\T,\C^{2d})$:
\[
\{\bar{u}_{-i}^{t,\alpha}(x)\}_{i=1}^{d-m} \subset E_{t,\alpha}^s(x), \quad 
\{\bar{u}_{i}^{t,\alpha}(x)\}_{i=1}^{d-m} \subset E_{t,\alpha}^u(\alpha,x),
\]
forming the $2d \times 2(d-m)$ matrix:
\[
M_{\alpha}(t,x) = \left[ \bar{u}_1^{t,\alpha}, \dots, \bar{u}_{d-m}^{t,\alpha} \mid \bar{u}_{-1}^{t,\alpha}, \dots, \bar{u}_{-(d-m)}^{t,\alpha} \right](x).
\]

Consider the associated Krein matrix
\[
G_{\alpha}(t,x) = \tfrac{1}{i}\, M_{\alpha}(t,x)^* \mathcal{J}_{2d} M_{\alpha}(t,x).
\]
Note that
\[
G_{\alpha}(t,z)
    = \tfrac{1}{i}\, M_{\alpha}(\bar t,\bar z)^* \mathcal{J}_{2d} M_{\alpha}(t,z)
    = \tfrac{1}{i}
    \begin{pmatrix}
        A_{\alpha}(t,z) & -K_{\alpha}(t,z) \\
        K_{\alpha}(\bar t,\bar z)^* & D_{\alpha}(t,z)
    \end{pmatrix}
\]
gives the holomorphic extension of \(G_{\alpha}(t,x)\) to the domain
\(\I(\delta)\times\T_h\). 
\begin{lemma}\label{lem:nondeg-krein}
    We have
    \[
        G_{\alpha}(t,z)=\frac{1}{i}\begin{pmatrix}
            0 & -K_{\alpha}(t,z)\\[4pt]
            K_{\alpha}(\bar t,\bar z)^* & 0
        \end{pmatrix}
    \]
where    $K_{\alpha}(t,z)\in C^\omega\bigl(\I(\delta)\times\T_h,\GL(d-m,\C)\bigr)$.
\end{lemma}

\begin{proof}
On the real slice $M=\I\times\T$, the subspaces $E^u_{t,\alpha}(z)$ and
$E^s_{t,\alpha}(z)$ are isotropic and their direct sum is Hermitian symplectic
\cite{WXZ}. Hence the diagonal blocks of the Krein matrix vanish:
\[
A_{\alpha}(t,z)=D_{\alpha}(t,z)=0 \qquad \forall\,(t,z)\in M.
\]
Since $A_{\alpha}$ and $D_{\alpha}$ are holomorphic on the connected domain
$\I(\delta)\times\T_h$, the identity principle yields
\[
A_{\alpha}\equiv D_{\alpha}\equiv0 \quad\text{on }\I(\delta)\times\T_h.
\]
Thus $G_{\alpha}$ has the stated off-diagonal form. By
Lemma~\ref{lem:Key-holom-extend}, $G_{\alpha}(t,z)$ is invertible throughout
$\I(\delta)\times\T_h$. Because the diagonal blocks vanish identically, this
invertibility is equivalent to the nonsingularity of $K_{\alpha}(t,z)$, which
completes the proof.
\end{proof}

Consequently, $K_{\alpha}^{-1}(t,x)\in C^\omega_{\delta,0,h}\bigl(\I\times\OO\times\T,\GL(d-m,\C)\bigr) $, we have:
\[
\begin{pmatrix}
    I_{d-m} & 0 \\
    0 & (K_{\alpha}^{-1}(t,x))^*
\end{pmatrix}
G_{\alpha}(t,x)
\begin{pmatrix}
    I_{d-m} & 0 \\
    0 & K_{\alpha}^{-1}(t,x)
\end{pmatrix}
= \tfrac{1}{i} \begin{pmatrix}
    0 & -I_{d-m} \\
    I_{d-m} & 0
\end{pmatrix}.
\]
Define the adjusted frame:
\[
\widetilde{M}_{\alpha}(t,x) := M_{\alpha}(t,x) \begin{pmatrix}
    I_{d-m} & 0 \\
    0 & K_{\alpha}^{-1}(t,x)
\end{pmatrix} = \left[ u_1^{t,\alpha}, \dots, u_{d-m}^{t,\alpha} \mid u_{-1}^{t,\alpha}, \dots, u_{-(d-m)}^{t,\alpha} \right](x),
\]
where $\{u_{\pm i}^{t,\alpha}\}$ forms a symplectic basis for $E_{t,\alpha}^u \oplus E_{t,\alpha}^s$ with $\omega_d(u_{-i}^{t,\alpha}, u_j^{t,\alpha}) = \delta_{ij}$.

\vspace{0.2cm}
\noindent \textbf{Step II: Symplectic frame for $E_{t,\alpha}^c$.} \\

Again by Proposition \ref{prop:globalframe}, there exists $\delta=\delta(\I,\alpha_0,h)>0$ and a  neighborhood $\mathcal{O}=\OO(\I,\alpha_0,h) \subset \mathbb{T}$ of $\alpha_0$ such that one can select frames for the center bundle in $C^\omega_{\delta,0,h}(\I\times\OO\times\T,\C^{2d})$:
$\{\bar{u}_{\pm i}^{t,\alpha}(x)\}_{i=d-m+1}^d \subset E_{t,\alpha}^c(x),$
forming the $2d \times 2m$ matrix:
\[
N_{\alpha}(t,x) = \left[ \bar{u}_{d-m+1}^{t,\alpha}, \dots, \bar{u}_d^{t,\alpha} \mid \bar{u}_{-(d-m+1)}^{t,\alpha}, \dots, \bar{u}_{-d}^{t,\alpha} \right](x).
\]
By  Lemma \ref{lem:Key-holom-extend}, one concludes that the associated Krein matrix
\[
H_{\alpha}(t,x)
    := \tfrac{1}{i}\, N_{\alpha}(t,x)^* \mathcal{J}_{2d} N_{\alpha}(t,x)
    \in C^\omega_{\delta,0,h}(\I \times\OO \times \mathbb{T}, \mathrm{GL}(2m,\mathbb{C})\cap \rmm{Her}(2m,\C))
\]
is invertible on the whole domain \(\I(\delta)\times\OO\times\T_h\).
To proceed, We will use the following Sylvester inertia theorem, a generalization of \cite[Theorem 1.3]{WXZ}:

\begin{lemma}\label{hetong}
Let $\I\subset\mathbb R$ be a compact interval, 
\(
H_{\alpha}(t,x)\in C^\omega_{\delta,0,h}\bigl(\I\times\OO\times\T,\GL(2m,\CC)\cap \Her(2m,\C)\bigr)
\)
and 
is invertible on \(\mathcal{D}:=\I(\delta)\times\OO\times\T_h\). 
Then, there exist constants \(\delta'\in(0,\delta]\) and \(h'\in(0,h]\), an open neighborhood \(\OO'\subset\OO\) of \(\alpha_0\), a holomorphic map
\(
P_{\alpha}(t,x)\in C^\omega_{\delta',0,h'}\bigl(\I\times\OO'\times\T,\GL(2m,\CC)\bigr),
\)
and an integer \(p\in\{0,\dots,2m\}\) such that for all \((t,\alpha,z)\in \mathcal{D}':=\I(\delta')\times\OO'\times\T_{h'}\), we have:
\[
P_{\alpha}(\bar t,\bar z)^{*}\,H_{\alpha}(t,z)\,P_{\alpha}(t,z)
=\operatorname{diag}\bigl(I_p,\,-I_{2m-p}\bigr).
\]
\end{lemma}

\begin{proof}
As \(H_{\alpha}(t,x)\) is Hermitian and invertible for every real pair \((t,\alpha,x)\), the positive inertia index (numbers of positive  eigenvalues, denoted by $p$) is constant on the compact set \(\I\times\{\alpha_0\}\times\T\), and there exists \(\eta_0>0\) such that for all \((t,x)\in\I\times\T\)
\[
\sigma\bigl(H_{\alpha_0}(t,x)\bigr)\subset (-\infty,-\eta_0]\cup[\eta_0,\infty).
\]
By continuity of \(H_\alpha\) in \((t,\alpha,z)\) and compactness of \(\I\times\{\alpha_0\}\times\overline{\T_{h}}\), we may choose
\(\delta'\in(0,\delta]\), \(h'\in(0,h]\) and a neighborhood \(\OO'\) of \(\alpha_0\) so that for all
\((t,\alpha,z)\in\mc D'\) the spectrum of \(H_\alpha(t,z)\) stays uniformly away from the imaginary axis:
there is \(\eta\in(0,\eta_0/2)\) with
\begin{equation}\label{eigen-dis}   
\sigma\bigl(H_{\alpha}(t,z)\bigr)\subset
\{w\in\CC:\Re w\le -\eta\}\ \cup\ \{w\in\CC:\Re w\ge \eta\}.
\end{equation}

Let \(\Gamma_+\) (resp.\ \(\Gamma_-\)) be a fixed positively oriented contour contained in \(\{ \Re w>\eta/2\}\) (resp.\ \(\{ \Re w<-\eta/2\}\)) and enclosing the right (resp.\ left) spectral component of \(H_\alpha(t,x)\) for real \(x\). By \eqref{eigen-dis} these contours lie in the resolvent set of \(H_\alpha(t,z)\) for all \((t,\alpha,z)\in\mc D'\). Define the Riesz projections
\[
\Pi^\pm_\alpha(t,z)
:=\frac{1}{2\pi i}\int_{\Gamma_\pm} (\zeta I - H_\alpha(t,z))^{-1}\,d\zeta.
\]
It is clear  the images \(E^\pm_\alpha(t,z):=\mathrm{Ran}\,\Pi^\pm_\alpha(t,z)\) define holomorphic subbundles of the trivial bundle over \(\I_{\delta'}\times\T_{h'}\) which depend continuously on \(\alpha\). By Proposition~\ref{lem:ajslem} and Proposition~\ref{prop:globalframe}, there exist matrices
\[
V^+_\alpha(t,z)\in C^\omega_{\delta',0,h'}(\I\times\OO'\times\T,\CC^{2m\times p}),\qquad
V^-_\alpha(t,z)\in C^\omega_{\delta',0,h'}(\I\times\OO'\times\T,\CC^{2m\times (2m-p)}),
\]
whose columns give bases of \(E^+_\alpha\) and \(E^-_\alpha\), respectively, for every \((t,\alpha,z)\in\mc D'\).

 Put
\(V_\alpha(t,z):=\big[V^+_\alpha(t,z)\ \ V^-_\alpha(t,z)\big]\in\GL(2m,\CC)\).
Then the matrix
\[
F_\alpha(t,z):=V_\alpha(\bar t,\bar z)^*\,H_{\alpha}(t,z)\,V_\alpha(t,z)
\]
has holomorphic entries.
On the real slice, \(\Pi^\pm_\alpha(t,x)\) are exactly the positive / negative spectral subspaces of the Hermitian matrix \(H_\alpha(t,x)\), then the subspaces $E^\pm_\alpha(t,x)$ are
$H_{\alpha}(t,x)$--orthogonal, hence $F_\alpha(t,x)$ is block diagonal.
Again by the identity principle, 
\[
F_\alpha(t,z)=V_\alpha(\bar t,\bar z)^*\,H_{\alpha}(t,z)\,V_\alpha(t,z)=\diag{M^+_\alpha(t,z), -M^-_\alpha(t,z) }
\]
 is  analytic and block diagonal on the whole strip. Moreover, \(M^+_\alpha(x)\) and \(M^-_\alpha(x)\) are positive definite for real \(x\).

By continuity and the spectra of \(M^\pm_\alpha(t,z)\) are contained in a fixed simply connected open subset of the right half-plane for all \((t,\alpha,z)\in\mc D'\). In particular none of the eigenvalues of \(M^\pm_\alpha\) vanishes on \(\mc D'\), and the spectral sets admit a single holomorphic branch of the square root.
Choose a fixed contour \(\Gamma_{s}\) contained in the right half-plane and enclosing \(\sigma(M^\pm_\alpha(t,z))\) for all \((t,\alpha,z)\in\mc D'\). Define, for example,
\[
S^\pm_\alpha(t,z) \;:=\; \frac{1}{2\pi i}\int_{\Gamma_{s}} \zeta^{-1/2}\,(\zeta I - M^\pm_\alpha(t,z))^{-1}\,d\zeta.
\]
 By the holomorphic functional calculus, \(S^\pm_\alpha\in C^\omega_{\delta',0,h'}(\I\times\OO'\times\T,\GL(\cdot,\CC))\) satisfies \(S^\pm_\alpha (\bar t,\bar z)^* S^\pm_\alpha(t,z)=(M^\pm_\alpha(t,z))^{-1}\) on \(\mc D'\) (again using  the identity principle), and on the real slice \(S^\pm_\alpha(t,x)\) is the usual positive definite inverse square root.

Set $
P_\alpha(t,z):=\bigl[\,V^+_\alpha(t,z)S^+_\alpha(t,z)\ \ ;\ \ V^-_\alpha(t,z)S^-_\alpha(t,z)\,\bigr]\in C^\omega_{\delta',0,h'}(\I\times\OO'\times\T,\GL(2m,\CC)).$
A direct computation shows that $P_\alpha(\bar t,\bar z)^* H_\alpha(t,z) P_\alpha(t,z) = \diag{I_p,-I_{2m-p}}$.
\end{proof}

\begin{remark}\label{rem:uniformityofh}
By construction, the parameters $\delta'$ and $h'$ depend  on $m$, and they are
chosen so that the eigenvalues $\sigma\bigl(H_{\alpha}(t,z)\bigr)$ satisfy
\eqref{eigen-dis} uniformly for all $(t,\alpha,z)\in\mathcal D'$.
\end{remark}

 As showed in \cite[Proposition 5.2]{WXZ}, in the center bundle,  the signature (the difference of positive inertia index and nagetive positive inertia index) is zero, i.e. $p=m$. 
Hence, by Lemma \ref{hetong},
there exists $P_{\alpha}(t,x) \in C^\omega_{\delta',0,h'}(\I \times\OO' \times \mathbb{T}, \mathrm{GL}(2m, \mathbb{C}))$ such that:
\[
P_{\alpha}(t,x)^* H_{\alpha}(t,x) P_{\alpha}(t,x) = \begin{pmatrix}
    I_m & 0 \\
    0 & -I_m
\end{pmatrix}.
\]
Let $Q \in \mathrm{GL}(2m, \mathbb{C})$ be a fixed transformation satisfying:
\[
Q^* \begin{pmatrix} 0 & iI_m \\ -iI_m & 0 \end{pmatrix} Q = \begin{pmatrix}
    I_m & 0 \\
    0 & -I_m
\end{pmatrix}.
\]
Define the symplectic frame:
\[
\widetilde{N}_{\alpha}(t,x) := N_{\alpha}(t,x) P_{\alpha}(t,x) Q^{-1} = \left[ u_{d-m+1}^{t,\alpha}, \dots, u_d^{t,\alpha} \mid u_{-(d-m+1)}^{t,\alpha}, \dots, u_{-d}^{t,\alpha} \right](x),
\]
where  this adjusted frame forms a canonical symplectic basis.

\section{Global Symplectification  for Monotonic Families of Hermitian Symplectic cocycles}\label{sec:global-symp-mono}

As a consequence of Proposition~\ref{prop:globalsympframe}, if  the cocycle \( (\alpha_0, A_t(\cdot+i\epsilon)) \) is partially hyperbolic with a $2m$-dimensional center bundle for all $t\in\I$ and $|\epsilon|<h_0$ with $\I$ compact, then 
there exists an open neighborhood \( \mathcal{O}'\subset \mathbb{T} \) of $\alpha_0$, $\delta'>0 $, $h'>0 $ and
$
    B_{t,\alpha}(x) \in C^\omega_{\delta',0,h'}(\I \times\OO' \times \mathbb{T}, \mathrm{HSp}(2d))$
such that 
\begin{equation}\label{eq:not-modif-conj}
    B_{t,\alpha}(x+\alpha)^{-1} A_t(x) B_{t,\alpha}(x)
    = H_{t,\alpha}(x) \diamond C_{t,\alpha}(x),
\end{equation}
where $
    H_{t,\alpha}(x) = 
    \begin{pmatrix}
        H_{t,\alpha}^+(x) & \\
        & H_{t,\alpha}^-(x)
    \end{pmatrix}
    \in C^\omega_{\delta',0,h'}(\I \times\OO \times \mathbb{T}, \mathrm{HSp}(2d-2m)),$ 
    $ 
    C_{t,\alpha}(x) \in C^\omega_{\delta',0,h'}(\I \times\OO \times \mathbb{T}, \mathrm{HSp}(2m))$.
However, if the family \( A_t \) is monotonic, the reduced cocycle \( C_{t,\alpha}(x) \) does not in general preserve monotonicity.  

Recall that a one-parameter family of cocycles  $A_t(x)\in C^1(\mathbf I\times\T,\rmm{HSp}(2d,\psi))$ is \emph{monotonic} \cite{WAXZ}, if for every isotropic vector \( v \in \mathbb{C}^{2d} \),
\[
    \Psi_{A_t}(v,x) 
    := \psi\!\left( A_t(x) v, \ \partial_t A_t(x) v \right) > 0.
\]
The next result shows that monotonicity can be recovered through a suitable modification.

\begin{proposition}\label{prop:globalmonoframe}
Let \( \I \subset \mathbb{R} \) 
be a compact interval. 
Let  \( A_t(x)\in C^\omega_{c,h_0}(\I \times \mathbb{T}, \mathrm{HSp}(2d)) \) be monotonic, and 
 the cocycle \( (\alpha_0, A_t(\cdot+i\epsilon)) \) is partially hyperbolic with a $2$-dimensional center bundle for all $t\in\I$ and $|\epsilon|<h_0$. 
Then there exists a closed neighborhood \( \mathcal{O}\subset \mathbb{T} \) of $\alpha_0$, $h'>0$ and  $
    \widetilde{B}_{t,\alpha}(x) \in C^\omega_{0,h'}(\OO\times\mathbb{T}, \mathrm{HSp}(2d)),$
piecewise \( C^\omega \) with respect to \( t \in \I \)\footnote{We adopt a similar Banach space notation for $(\alpha, x)$, omitting the $t$-analyticity from the previous definition for triples; piecewise analyticity is assumed to be continuous.}, such that
\[
    \widetilde{B}_{t,\alpha}(x+\alpha)^{-1} A_t(x) \widetilde{B}_{t,\alpha}(x)
    = H_{t,\alpha}(x) \diamond \widetilde{C}_{t,\alpha}(x),
\]
where $
    \widetilde{C}_{t,\alpha}(x) \in C^\omega_{0,h'}(\OO\times\mathbb{T},\mathrm{HSp}(2))
$
is piecewise \( C^\omega \) and piecewise monotonic in \( t \in \I \).  
\end{proposition}
\begin{remark}
    {Quantitative bounds for the monotonicity are established in Lemma \ref{lem:LUbd}.}
\end{remark}

	\begin{proof}

Fix any $0<h<h_{0}$. Once $\alpha_{0}$, $\I$, and $h$ are fixed, the constants $\delta'$, $h'>0$, and the open neighborhood $\OO'$ provided by Proposition~\ref{prop:globalsympframe} are thereby determined. Without loss of generality, one can slightly shrink $\OO'$ to a closed neighborhood $\OO$ of $\alpha_0$.
Then the proof will be divided into four steps.

\vspace{0.2cm}
\noindent \textbf{ Step I: Graph Transform}

Let $\{u_{\pm i}^{s,\alpha}(x)\}_{i=1}^d \subset C^\omega_{\delta',0,h'}(\I\times\OO\times\T, \C^{2d}) $ be the symplectic frame constructed in Proposition~\ref{prop:globalsympframe}, where $\{u_i^{s,\alpha}(x)\}_{i=1}^{d-1}$ spans the stable bundle $E_{s,\alpha}^{s}$,  while  $\{u_{-i}^{s,\alpha}(x)\}_{i=1}^{d-1}$ spans the unstable bundle $E_{s,\alpha}^{u}$.
 
For any fixed $ s \in  \I $, we have the splitting $\mathbb{C}^{2d} = E_{s,\alpha}^c (x)\oplus E_{s,\alpha}^s(x) \oplus E_{s,\alpha}^u(x)$.  By robustness of dominated splitting, 
there exists a uniform neighborhood $\I_s \subset \I$ (with size controlled by the derivative bounds of the frame) 
such that for all $t \in \I_s$, the center bundle $E_{t,\alpha}^c(x)$ is sufficiently close to $E_{s,\alpha}^c(x)$, uniformly in $(\alpha,x)\in \OO\times\T_{h'}$, ensuring transverse intersection:
\[
E_{t,\alpha}^c(x) \cap (E_{s,\alpha}^s(x) \oplus E_{s,\alpha}^u(x)) = \{0\},
\text{ for all }(\alpha,x)\in \OO\times \T_{h'}.
\]
Consequently, $E_{t,\alpha}^c$ is the graph of a linear map $\Phi_{t,s,\alpha} : E_{s,\alpha}^c \to E_{s,\alpha}^s \oplus E_{s,\alpha}^u$. 
Let $ \mathbb{P}_{E_{s,\alpha}^c}:\C^{2d}\to E_{s,\alpha}^c $ and $ \mathbb{P}_{E_{s,\alpha}^{s}\oplus E_{s,\alpha}^{u}}:\C^{2d}\to E_{s,\alpha}^{s}\oplus E_{s,\alpha}^{u} $ be the projections determined by the direct sum decomposition, then 
\begin{equation}\label{eq:graphTC}
    \Phi_{t,s,\alpha} = \left( \PP_{E_{s,\alpha}^s \oplus E_{s,\alpha}^u} \big|_{E_{t,\alpha}^c} \right) \circ \left( \PP_{E_{s,\alpha}^c} \big|_{E_{t,\alpha}^c} \right)^{-1}.
\end{equation}
Define the \textit{graph transform}:
\[
\widehat{\nabla}_{t,s,\alpha} : E_{s,\alpha}^c \to E_{t,\alpha}^c, \quad v \mapsto v + \Phi_{t,s,\alpha}(v).
\]
Its inverse is the projection $\widehat{\nabla}_{t,s,\alpha}^{-1} = \PP_{E_{s,\alpha}^c} |_{E_{t,\alpha}^c}$.  A quantitative version is the following:
\begin{lemma}\label{lem:graphy}
For each \(s\in \I\), there exists a constant \(\delta_1>0\), independent of \(s\), and a uniform neighborhood
  $\I_s= [s-\delta_1,\,s+\delta_1]\,\cap\,\I$, such that 
 the graph transform connection   $\widehat{\nabla}_{t,s,\alpha}$
 is well-defined in  $\I_s$, which is $C^\omega$ in $t$ and $s$. Moreover,  there exists $\bar{C}=\bar{C}(\delta',\|u_{\pm d}^{t,\alpha}(x)\|_{\delta',\OO,h'} )>0$, independent of $s$, such that\begin{equation}\label{eq:derivative-bound}
  \bigl\|\partial_t^k\widehat{\nabla}_{t,s,\alpha} u_{\pm d}^{s,\alpha}\bigr\|_{\I_s,\OO,h'}=\bigl\|\partial_t^k\Phi_{t,s,\alpha} u_{\pm d}^{s,\alpha} \bigr\|_{\I_s,\OO,h'} \;\le\;\bar C,\  k=1,2
  \quad\text{for all }s \in\I.
\end{equation}
where we denote 
$\|A\|_{\I,\OO,h}= \sup_{t\in \I, \alpha\in \OO,x\in \T_h}   \|A(t,\alpha,x)\|  $. 
\end{lemma}
\begin{proof}

To establish this, define the $2d \times 2d$ matrix:
\[
M_s^{t,\alpha}(x) = \left[ u_1^{s,\alpha}(x), \dots, u_{d-1}^{s,\alpha}(x), u_d^{t,\alpha}(x) \mid u_{-1}^{s,\alpha}(x), \dots, u_{-(d-1)}^{s,\alpha}(x), u_{-d}^{t,\alpha}(x) \right].
\]  
By compactness of $\I$ and analyticity of the symplectic frame, we have
\begin{lemma}\label{lem:detQ}
For each \(s\in \I\), there exists $\delta_1 = \delta_1(\delta',\|u_{\pm d}^{t,\alpha}(x)\|_{\delta',\OO,h'} ) > 0$, independent of \(s\), such that for all $t\in [s-\delta_1,s+\delta_1]$ and $(\alpha,x) \in \OO \times \mathbb{T}_{h'}$,
$ M_s^{t,\alpha}(x)$ is invertible.
\end{lemma}
\begin{proof}

Using the fact $M_s^{s,\alpha}(x)^{-1}M_s^{s,\alpha}(x)=I_{2d}$ and analyticity, a direct computation shows that \[
\det M_s^{s,\alpha}(x)^{-1}M_s^{t,\alpha}(x)=\det  Q_s^{t,\alpha}(x),
\] where $Q_s^{t,\alpha}(x)$ represents the restriction of the coordinate transition to the center subspace, that is,
\begin{equation}\label{eq:Q-def}
(\mathbf{0}_{2d-2} \diamond I_2) M_s^{s,\alpha}(x)^{-1} M_s^{t,\alpha}(x) (\mathbf{0}_{2d-2} \diamond I_2) = \mathbf{0}_{2d-2} \diamond Q_s^{t,\alpha}(x).
\end{equation}
More precisely, by symplecticity, $M_s^{s,\alpha}(x)^{-1} = -\mathcal{J}_{2d} M_s^{s,\alpha}(\bar{x})^* \mathcal{J}_{2d}$. Then we have
\begin{equation}\label{eq:representQ}
Q_s^{t,\alpha}(x) = \begin{pmatrix}
\omega_d\bigl(u_{-d}^{s,\alpha}(\bar{x}),\, u_d^{t,\alpha}(x)\bigr) & \omega_d\bigl(u_{-d}^{s,\alpha}(\bar{x}),\, u_{-d}^{t,\alpha}(x)\bigr) \\[4pt]
-\omega_d\bigl(u_d^{s,\alpha}(\bar{x}),\, u_d^{t,\alpha}(x)\bigr) & -\omega_d\bigl(u_d^{s,\alpha}(\bar{x}),\, u_{-d}^{t,\alpha}(x)\bigr)
\end{pmatrix}.
\end{equation}

Since $t \mapsto \det Q_s^{t,\alpha}(x)$ is analytic, for $ t\in \I(\delta')$, write
\[
\det Q_s^{t,\alpha}(x) - 1 = \int_s^t \frac{d}{d\tau} \det Q_s^{\tau,\alpha}(x) \, d\tau,
\]
the integral taken along the segment from $s$ to $t$. By the Cauchy estimate,
\[
\sup_{\xi  \in \I(\delta'/2)} \bigl| \partial_t|_{t=\xi} \det Q_s^{t,\alpha}(x) \bigr|
\le \frac{4}{\delta'} \sup_{\tau \in \I(\delta')} \bigl| \det Q_s^{\tau,\alpha}(x) \bigr|.
\]
Set
\[
A_d(\delta',h') := \frac{4}{\delta'} \sup_{s\in\I} \| \det Q_s^{t,\alpha}(x) \|_{\delta',\OO,h'}.
\]
 Choosing $\delta_1 = \frac{1}{2A_d} $, we obtain
\[
|\det Q_s^{t,\alpha}(x)| \ge 1 - |t-s| A_d \geq \frac{1}{2}, \qquad t\in [s-\delta_1,s+\delta_1],
\]
uniformly for $s \in \I$ and $(\alpha,x) \in \mathcal{O} \times \mathbb{T}_{h'}$.
\qedhere
\end{proof}

Shrinking $\I_s$ to $(s - \delta_1, s + \delta_1)$ ensures uniform transversality.
For any non-zero  $w \in E_{t,\alpha}^c(x)$, $w' \in E^c_{s,\alpha}(x)$, consider its representation in the frame:
\[
w = M_s^{t,\alpha}(x) 
\begin{pmatrix}
    0_{d-1} & \bar{x} & 0_{d-1} & \bar{y}
\end{pmatrix}^{\top}, \quad (\bar{x},\bar{y}) \in \C^2 \setminus \{0\}.
\]
\[
w' = M_s^{s,\alpha}(x) 
\begin{pmatrix}
    0_{d-1} & x' & 0_{d-1} & y'
\end{pmatrix}^{\top}, \quad (x',y') \in \C^2 \setminus \{0\}.
\]
The projections onto $E_{s,\alpha}^c$ and $E_{s,\alpha}^s \oplus E_{s,\alpha}^u$ are given by:
\begin{eqnarray*}
\mathbb{P}_{E_{s,\alpha}^c}(w) &=& M_s^{s,\alpha}(x) (\mathbf{0}_{2d-2} \diamond I_2) M_s^{s,\alpha}(x)^{-1} w, \\
\mathbb{P}_{E_{s,\alpha}^s \oplus E_{s,\alpha}^u}(w) &=& M_s^{s,\alpha}(x) (I_{2d-2} \diamond \mathbf{0}_2) M_s^{s,\alpha}(x)^{-1} w.
\end{eqnarray*}
The inverse of the restricted projection is then:
\[
(\mathbb{P}_{E_{s,\alpha}^c}|_{E_{t,\alpha}^c})^{-1}(w') = M_s^{t,\alpha}(x) (\mathbf{0}_{2d-2} \diamond Q_s^{t,\alpha}(x)^{-1}) M_s^{s,\alpha}(x)^{-1} w'.
\]

The graph map $\Phi_{t,s,\alpha}$ acts as:
\begin{align}
\Phi_{t,s,\alpha}(w') &= \mathbb{P}_{E_{s,\alpha}^s \oplus E_{s,\alpha}^u} \left( (\mathbb{P}_{E_{s,\alpha}^c}|_{E_{t,\alpha}^c})^{-1}(w') \right) \notag\\
&= M_s^{s,\alpha}(x) (I_{2d-2} \diamond \mathbf{0}_2) M_s^{s,\alpha}(x)^{-1} M_s^{t,\alpha}(x) (\mathbf{0}_{2d-2} \diamond Q_s^{t,\alpha}(x)^{-1}) M_s^{s,\alpha}(x)^{-1} w'.\label{eq:represent-of-Phi}
\end{align}
Note by \eqref{eq:Q-def} together with \eqref{eq:represent-of-Phi}, the graph map $\Phi_{t,s,\alpha}$ is also independent of the choice of frame in the stable, unstable, or center bundle.

Again by the symplectic relation $M_s^{s,\alpha}(x)^{-1} = -\mathcal{J}_{2d} M_s^{s,\alpha}(\bar{x})^* \mathcal{J}_{2d}$, we have
\begin{equation}\label{eq:leftinverse}
    (\mathbf{0}_{2d-2} \diamond I_2)M_s^{s,\alpha}(x)^{-1}
= \begin{pmatrix}0\\ (u_{-d}^{s,\alpha}(\bar{x}))^* \mathcal{J}_{2d}\\ 0\\ -(u_d^{s,\alpha}(\bar{x}))^* \mathcal{J}_{2d}\end{pmatrix}.
\end{equation}
Thus
\begin{equation}\label{eq:ProjeE_c}
    \mathbb{P}_{E_{s,\alpha}^c} = \left[u_d^{s,\alpha}(x) (u_{-d}^{s,\alpha}(\bar{x}))^* - u_{-d}^{s,\alpha}(x) (u_d^{s,\alpha}(\bar{x}))^*\right] \mathcal{J}_{2d},
\end{equation}
which involves only the center columns $u_{\pm d}^{s,\alpha}$ evaluated at $x$ and $\bar{x}$. Moreover, we have
\begin{equation}
    \|\mathbb{P}_{E_{s,\alpha}^c}\|
    \leq 2 \|u^{s,\alpha}_d(x)\|
    \|u^{s,\alpha}_{-d}(x)\|
\end{equation}

From \eqref{eq:represent-of-Phi} and $\mathbb{P}_{E_{s,\alpha}^s \oplus E_{s,\alpha}^u} = I_{2d} - \mathbb{P}_{E_{s,\alpha}^c}$ we have
\[
\Phi_{t,s,\alpha} = (I_{2d} - \mathbb{P}_{E_{s,\alpha}^c}) \,
M_s^{t,\alpha}(x) \, (\mathbf{0}_{2d-2} \diamond Q_s^{t,\alpha}(x)^{-1}) \, M_s^{s,\alpha}(x)^{-1}.
\]
Let $w' \in E_{s,\alpha}^c$ with $\|w'\|=1$, \eqref{eq:leftinverse} gives
\[
M_s^{s,\alpha}(x)^{-1} w' = \begin{pmatrix}0\\ \omega_d(u_{-d}^{s,\alpha}(\bar{x}), w') \\ 0\\ -\omega_d(u_d^{s,\alpha}(\bar{x}), w')\end{pmatrix}.
\]

Let $ \mathfrak{U}^{t,\alpha}(x) =\begin{pmatrix}
   u_d^{t,\alpha}(x) &  u_{-d}^{t,\alpha}(x))
\end{pmatrix}  $ and $z=(\omega_d(u_{-d}^{s,\alpha}(\bar{x}), w'), -\omega_d(u_d^{s,\alpha}(\bar{x}), w'))^\top$,
\begin{equation}\label{eq:represent-of-Phi-2}
    \Phi_{t,s,\alpha}(w') = (I_{2d} - \mathbb{P}_{E_{s,\alpha}^c}) \, \mathfrak{U}^{t,\alpha}(x)  \, Q_s^{t,\alpha}(x)^{-1} z .
\end{equation}
Lemma~\ref{lem:detQ} gives $|\det Q_s^{t,\alpha}(x)| \ge \frac12$ for $t\in\I_s$, so $\|Q^{-1}\| \le 2\|Q\|$.  Together with the explicit structure for $\mathbb{P}_{E^c_{s,\alpha}}$, $\mathfrak{U}^{t,\alpha}(x) $ and $Q^{t,\alpha}_s(x)$, we obtain
\[
\|\Phi_{t,s,\alpha}u_{\pm d}^{s,\alpha}\|_{\I_s,\OO,h'}
\le C\bigl( \|u_{\pm d}^{t,\alpha}\|_{\delta',\OO,h'}
\bigr),
\]

For the $t$-derivative, analyticity of $t \mapsto \Phi_{t,s,\alpha}$ (as a composition of analytic objects) allows the Cauchy estimate ($k=1,2$) for $t\in \I_s$
\[
\|\partial_t^k \Phi_{t,s,\alpha}u_{\pm d}^{s,\alpha}\|_{\I_s,\OO,h'}
\leq  \bar C\bigl(\delta', \|u_{\pm d}^{t,\alpha}\|_{\delta',\OO,h'}
\bigr).
\]
This completes the proof.
\end{proof}

This \textit{graph transform} is particularly valuable because, as established in \cite[Theorem 4.2]{WAXZ}, it preserves the monotonicity of the center bundle --- a property that will be further exploited in subsequent steps.

\vspace{0.2cm}
\noindent \textbf{ Step II: Local Symplectification via parallel transport}

However $\widehat{\nabla}_{t,s,\alpha}$ is not necessarily symplectic, we require a \textbf{symplectic correction}  to obtain a symplectic frame in $E_{t,\alpha}^c$. Fix $s\in\I$, for any $t\in \I_s$,
consider the image of the center basis under the connection:
\[
\widehat{u}_{\pm d}^{t,\alpha}(x) = \widehat{\nabla}_{t,s,\alpha} u_{\pm d}^{s,\alpha}(x) \in E_{t,\alpha}^c.
\]
The Krein matrix satisfies:
\[
\left( \widehat{u}_{d}^{t,\alpha}(x) \quad \widehat{u}_{-d}^{t,\alpha}(x) \right)^* \mathcal{J}_{2d} \left( \widehat{u}_{d}^{t,\alpha}(x) \quad \widehat{u}_{-d}^{t,\alpha}(x) \right) =: J_s(t,\alpha,x),
\]
with $J_s(s,\alpha,x) = \mathcal{J}_2$ and 
\begin{equation}\label{van-de}
\partial_t J_s(t,\alpha,x)\big|_{t=s} = 0.
\end{equation}This vanishing derivative follows because $\frac{d}{dt}\big|_{t=s} \widehat{u}_{\pm d}^{t,\alpha}(x) \in E_{s,\alpha}^s \oplus E_{s,\alpha}^u$, which is symplectically orthogonal to $E_{s,\alpha}^c$. Moreover, note that by \eqref{eq:represent-of-Phi-2}, $J_s(t,\alpha, x)
$ can be expressed in entries relative to the frame $\{u_{\pm d}^{t,\alpha}(x)\}$. More precisely, its derivatives is controlled by \begin{equation}\label{eq:partial-J_s}
    \bigl\|\partial_t^k J_s\bigr\|_{\I_s,\OO,h'}\leq \bar C'\bigl(\delta', \|u_{\pm d}^{t,\alpha}\|_{\delta',\OO,h'}
\bigr), \qquad k=1,2.
\end{equation}

\begin{lemma}\label{lem:N_s(t)}
For each \( s \in \I \), there exists a constant \( 0 < \delta_2 < \delta_1 \), independent of \(s\), and a uniform neighborhood
  $\I'_s \;=\; [s-\delta_2,\,s+\delta_2] \,\cap\, \I$
together with a map
  $N_s \;\in\; C^\omega_{\delta_2,0,h'}\!\bigl(\I'_s\times\OO \times \T,\, \mathrm{GL}(2,\C)\bigr)$
such that, for all \((t,\alpha,x) \in \I'_s \times \OO \times \T\),
\[
  N_s(t,\alpha,x)^{*}\,J_s(t,\alpha,x)\,N_s(t,\alpha,x)\;=\;\mathcal{J}_2,
  \qquad
  N_s(s,\alpha,x) = I_2,
  \qquad
  \partial_t N_s(s,\alpha,x) = 0.
\]
Moreover, there exists an absolute constant \( C > 0 \) such that the following estimates hold:
\begin{align}
  \label{eq:derivative-1}
  \bigl\|\partial_t N_s\bigr\|_{\I'_s,\OO,0}
  &\;\le\;
  C\,\bigl\|\partial_t J_s\bigr\|_{\I'_s,\OO,0}, \\[6pt]
  \label{eq:derivative-2}
  \bigl\|\partial_t^2 N_s\bigr\|_{\I'_s,\OO,0}
  &\;\le\;
  C\,\Bigl(\,\bigl\|\partial_t J_s\bigr\|_{\I'_s,\OO,0}
       \;+\;\bigl\|\partial_t^2 J_s\bigr\|_{\I'_s,\OO,0}\Bigr).
\end{align}
\end{lemma}

\begin{proof}
The proof follows from the Quantitative Implicit Function Theorem.

\begin{theorem}[\cite{BB,De}]\label{thm6.1}(Quantitative Implicit Function Theorem)
Let \(X,Y,Z\) be Banach spaces, and let \(U \subset X\), \(V \subset Y\) be neighborhoods of \(x_0\) and \(y_0\), respectively.  
Fix \( \varepsilon,\delta>0 \) and define 
\( X_\varepsilon := \{x \in X : \|x-x_0\|_X < \varepsilon\} \) and 
\( Y_\delta := \{y \in Y : \|y-y_0\|_Y < \delta\} \).  
Let \( \Psi \in C^1(U \times V, Z) \).  
Suppose \( \Psi(x_0,y_0) = 0 \) and that \( D_y \Psi(x_0,y_0) \in \mathcal{L}(Y,Z) \) is invertible.  
If
\[ \sup_{\overline{X_{\varepsilon}}}\|\Psi(x,y_0)\|_Z\leq\frac{\delta}{2\|D_y\Psi(x_0,y_0)^{-1}\|_{\mathcal{L}(Z,Y)}}, \]\[ \sup_{\overline{X_{\varepsilon}}\times\overline{Y_\delta}}\|\mathrm{Id}_Y-D_y\Psi(x_0,y_0)^{-1}D_y\Psi(x,y)\|_{\mathcal{L}(Y,Y)}\leq\frac{1}{2}, \]
then there exists \( y \in C^1(X_{\varepsilon}, \overline{Y_\delta}) \) such that \( \Psi(x,y(x)) = 0 \).
\end{theorem}

Define the Banach spaces $\mathcal{J}$ and $\mathcal{A}$ of analytic matrix-valued functions equipped with the norm
\[
\|A\|_{c,0,h} 
= \sup_{( t,\alpha,z) \in \I(c)\times\OO\times \T_h}
   \|A(t,\alpha,z)\|
\]
as follows:
\[
\begin{aligned}
\mathcal{J} &= \Bigl\{ J \in C_{c,0,h}^\omega(\I \times \OO \times \T, \mathrm{M}(2,\C)) : J( t,\alpha,z) = -J^*(\bar t,\alpha,\bar z) \Bigr\}, \\
\mathcal{A} &= \Bigl\{ A \in C_{c,0,h}^\omega(\I \times \OO \times \T, \mathrm{M}(2,\C)) : \mathcal{J}_2 A(t,\alpha,z)=A(\bar t,\alpha,\bar z)^*\mathcal{J}_2   \Bigr\}.
\end{aligned}
\]

Consider the \(C^1\)-map \(F: \mathcal{J} \times \mathcal{A} \to \mathcal{J}\) defined by
\[
F(J,A)(t,\alpha,z) = A^*(\bar t,\alpha,\bar z) J(t,\alpha,z) A(t,\alpha,z) - \mathcal{J}_2.
\]
At \((J,A) = (\mathcal{J}_2,I_2)\), we have \(F(\mathcal{J}_2,I_2) = 0\). The Fréchet derivative with respect to \(A\) is
\[
D_A F(\mathcal{J}_2,I_2)[H] = H^* \mathcal{J}_2 + \mathcal{J}_2 H=2\mathcal{J}_2 H,\footnote{Here $H^*=H^*(\bar t,\alpha,\bar z)$ for $H=H(t,\alpha,z)$.}
\qquad \forall H \in \mathcal{A},
\]
which is invertible with inverse bounded by $\|(D_A F(\mathcal{J}_2, I_2))^{-1}\|_{\mathcal{L}(\mathcal{J},\mathcal{A})} \leq \frac{1}{2} \|\mathcal{J}_2^{-1}\|$. At general \((J,A)\), one has
\[
D_A F(J,A)[H]
= H^* J A + A^* J H.
\]

Fix \( \varepsilon = \delta < \tfrac{1}{6} \). Define the neighborhoods
\[
K_\varepsilon = \bigl\{ J \in \mathcal{J} : \|J-\mathcal{J}_2\|_{c,0,h} < \varepsilon \bigr\}, 
\qquad
H_\delta = \bigl\{ A \in \mathcal{A} : \|A-I_2\|_{c,0,h} < \delta \bigr\}.
\]
The hypotheses of Theorem~\ref{thm6.1} are satisfied:
\[ 2\|(D_A F(\mathcal{J}_2, I_2))^{-1}\| \cdot \sup_{J \in \overline{K_\varepsilon}} \|F(J, I_2)\| \leq \|\mathcal{J}_2^{-1}\| \cdot \varepsilon \leq \delta, 
\] 
\[ 
\begin{aligned}
    &\sup_{J \in \overline{K_\varepsilon}, A \in \overline{H_\delta}}\left\| \rmm{Id}_\mathcal{A} - (D_A F(\mathcal{J}_2, I_2))^{-1} D_A F(J, A) \right\| \\
    &\qquad\qquad\qquad\leq \|\mathcal{J}_2^{-1}\|( \|J - \mathcal{J}_2\|_{\I,0, h}+ 2\|A - I_2\|_{\I,0, h} ) \leq \frac{1}{2}. 
\end{aligned}
\]
hence there exists a unique solution 
\(A(J) \in C^1(K_\varepsilon,\overline{H}_\delta)\) with \(F(J,A(J))=0\).

By \eqref{eq:partial-J_s}, there exists \(0 < \delta_2(\delta',\|u_{\pm d}^{t,\alpha}(x)\|_{\delta',\OO,h'} )) < \delta_1\) such that
\[
\|J_s(t,\alpha,x)-\mathcal{J}_2\|_{\I'_s,\OO,h'} < \varepsilon,
\qquad \I'_s = [s-\delta_2,\,s+\delta_2]\cap \I.
\]
Reducing $\delta_2<\frac12\min\{\delta_1,\delta'\}$ uniformly in $s$ if necessary, the same bound holds on the complex $\delta_2$-neighborhood of $\I'_s$, with the same dependence on the preceding data.
Thus, setting
$N_s(t,\alpha,x) := A(J_s(t,\alpha,x)) \in C^\omega_{\delta_2,0,h'}(\I'_s \times \OO \times \T, \mathrm{M}(2,\C)),$
we obtain
\begin{equation}\label{eq:sym}
N_s(t,\alpha,x)^*\,J_s(t,\alpha,x)\,N_s(t,\alpha,x) \;=\; \mathcal{J}_2,
\end{equation}
with \(\|N_s-I_2\|_{\I'_s,0,h'} < \tfrac{1}{6}\), so $N_s\in \mathrm{GL}(2,\C)$.

Since \(J_s(s,\alpha,x) = \mathcal{J}_2\), we have \(N_s(s,\alpha,x) = I_2\).  
Differentiating \eqref{eq:sym} at \(t=s\) gives
\[
(\partial_t N_s|_{t=s})^* \mathcal{J}_2 + \mathcal{J}_2 (\partial_t N_s|_{t=s}) + \partial_t J_s|_{t=s} = 0.
\]
Thus, by invertibility of \(D_A F(\mathcal{J}_2,I_2)\), one concludes \(\partial_t N_s|_{t=s} = 0\).

Differentiating \eqref{eq:sym} twice yields
\[
(\partial_t^2 N_s|_{t=s})^* \mathcal{J}_2 + \mathcal{J}_2 (\partial_t^2 N_s|_{t=s}) + \partial_t^2 J_s|_{t=s} + R = 0,
\]
where \(R\) depends only on first derivatives of \(N_s\) and \(J_s\).  
Since \(\partial_t N_s|_{t=s} = \partial_t J_s|_{t=s} = 0\), one has \(R=0\).  
Hence,
\[
(\partial_t^2 N_s|_{t=s})^* \mathcal{J}_2 + \mathcal{J}_2 (\partial_t^2 N_s|_{t=s})
= -\partial_t^2 J_s|_{t=s}.
\]
The operator \(H \mapsto H^* \mathcal{J}_2 + \mathcal{J}_2 H\) is invertible on \(\mathcal{A}\), giving
\[
\|\partial_t^2 N_s|_{t=s}\| 
\;\leq\; \|(D_A F(\mathcal{J}_2,I_2))^{-1}\| \cdot \|\partial_t^2 J_s|_{t=s}\|.
\]
The uniform first and second derivative bounds for the implicit map $A(J)$ and the chain rule give $\|\partial_tN_s\|\leq C\|\partial_tJ_s\|$ and $\|\partial_t^2N_s\|\leq C(\|\partial_t^2J_s\|+\|\partial_tJ_s\|^2)$ on $\I'_s\times\OO\times\T$. Since $\partial_tJ_s(s)=0$, reducing $\delta_2$ uniformly if necessary ensures $\|\partial_tJ_s\|\leq1$, yielding \eqref{eq:derivative-1}-\eqref{eq:derivative-2}.
\end{proof}

Within Lemma~\ref{lem:N_s(t)}, for each \(s \in \I\) we define the \textbf{symplectic parallel transport} on 
\(\I_s'=(s-\delta_2,s+\delta_2)\cap\I\) by
\[
\begin{array}{rcl}
	\widetilde{\nabla}_{t,s,\alpha}: E^c_{s,\alpha}(\cdot) &\longmapsto& E_{t,\alpha}^c(\cdot), \\[1mm]
	\begin{pmatrix} u_{d}^{s,\alpha}(x) \quad u_{-d}^{s,\alpha}(x) \end{pmatrix}
	&\longmapsto& 
	\widehat{\nabla}_{t,s,\alpha}\begin{pmatrix} u_{d}^{s,\alpha}(x) & u_{-d}^{s,\alpha}(x) \end{pmatrix}
	\,N_s(t,\alpha,x),
\end{array}
\]
where \(\widetilde{\nabla}_{t,s,\alpha}\begin{pmatrix} u_d^{s,\alpha}(x) & u_{-d}^{s,\alpha}(x) \end{pmatrix}\) satisfies the symplectic condition
\begin{equation}\label{sym-con}
N_s(t,\alpha,x)^*\!
\Big[\widehat{\nabla}_{t,s,\alpha}\begin{pmatrix} u_d^{s,\alpha}(x) & u_{-d}^{s,\alpha}(x) \end{pmatrix}\Big]^{\!*}
\mathcal{J}_{2d}\,
\widehat{\nabla}_{t,s,\alpha}\begin{pmatrix} u_d^{s,\alpha}(x) & u_{-d}^{s,\alpha}(x) \end{pmatrix}
N_s(t,\alpha,x)
= \mathcal{J}_2.
\end{equation}
In particular, \(\widetilde{\nabla}_{t,s,\alpha}u_{\pm d}^{s,\alpha}(x)=u_{\pm d}^{t,\alpha}(x)\) when \(t=s\).
Recall
\[
\mathfrak{U}^{s,\alpha}(x):=\begin{pmatrix} u_d^{s,\alpha}(x) & u_{-d}^{s,\alpha}(x) \end{pmatrix}.
\]

For \(t\in \I_s'\), the \textbf{local holonomy matrix} 
$ R_{t,s,\alpha}(x)\in C^\omega_{\delta',0,h'}\bigl(\I_s'\times\OO\times\T,\mathrm{HSp}(2)\bigr)$
is defined by
\begin{equation}\label{eq:local-holonomy}
\widetilde{\nabla}_{t,s,\alpha}\mathfrak{U}^{s,\alpha}(x)
=\mathfrak{U}^{t,\alpha}(x)\,R_{t,s,\alpha}(x),
\end{equation}
with \(R_{s,s,\alpha}(x)=I_2\). Thus \(R_{t,s,\alpha}(x)\) quantifies the deviation of the transported frame relative to the reference frame at the same \((\alpha,x)\).

\begin{lemma}\label{lem:localR_t}
Given \(0<\epsilon<\delta_2\), there exists \(\eta=\eta\bigl(\delta', \|u_{\pm d}^{t,\alpha}\|_{\delta',\OO,h'}
\bigr)>0\) such that for all \(|t-s|<\epsilon\), $t,s\in\I$,
\[
\|\partial_tR_{t,s,\alpha}\|_{\OO,h'}\;\le\;\eta,\qquad\|R_{t,s,\alpha}\|_{\OO,h'}\;\le\; e^{\eta\,\epsilon},
\]
where \(\|R_{t,s,\alpha}\|_{\OO,h}:=\sup_{(\alpha,x)\in\OO\times\T_{h}}\|R_{t,s,\alpha}(x)\|\).
\end{lemma}

\begin{proof}
Completely analogous to \eqref{eq:representQ}, the matrix $R_{t,s,\alpha}(x)$ is given by the symplectic inner products
\begin{equation}\label{eq:R_tsalpha}
    R_{t,s,\alpha}(x)=\begin{pmatrix}
    \omega_d(u_{-d}^{t,\alpha}(\bar{x}), \widetilde{\nabla}_{t,s,\alpha}u^{s,\alpha}_d(x)) & \omega_d(u_{-d}^{t,\alpha}(\bar{x}), \widetilde{\nabla}_{t,s,\alpha}u^{s,\alpha}_{-d}(x))\\
    -\omega_d(u_{d}^{t,\alpha}(\bar{x}), \widetilde{\nabla}_{t,s,\alpha}u^{s,\alpha}_d(x))&-\omega_d(u_d^{t,\alpha}(\bar{x}), \widetilde{\nabla}_{t,s,\alpha}u^{s,\alpha}_{-d}(x)) 
\end{pmatrix}
\end{equation}
By \eqref{eq:derivative-bound}, \eqref{eq:partial-J_s} and Lemma~\ref{lem:N_s(t)}, the quantities 
\(\partial_t\!\big(\widehat{\nabla}_{t,s,\alpha}(\alpha)\,u_{\pm d}^{s,\alpha}(x)\big)\) and 
\(\partial_t N_s(t,\alpha,x)\) are uniformly controlled (in \((\alpha,x)\)) by some $\bar C''=\bar C''\bigl(\delta', \|u_{\pm d}^{t,\alpha}\|_{\delta',\OO,h'}
\bigr)>0$.  
Consequently, for all \(|t-s|<\epsilon\), $t,s\in\I$, $\|\partial_t R_{t,s,\alpha}\|_{\OO,h'}\le \eta$. 
Integrating $R_{t,s,\alpha}=I_2+\int_s^t \partial_\tau R_{\tau,s,\alpha}\,d\tau$ yields
\[
\|R_{t,s,\alpha}\|_{\OO,h'} \le 1+\eta|t-s| \le e^{\eta|t-s|} \le e^{\eta\epsilon},
\]
which completes the proof.
\end{proof}

Lemmas~\ref{lem:graphy} and~\ref{lem:N_s(t)} ensure that the connection \(\widetilde{\nabla}_{t,s,\alpha}\) is well defined on \(\I_s'\) relative to the reference frame \(\{u_{\pm i}^{t,\alpha}(x)\}_{i=1}^d\).  
Via the local holonomy \(R_{t,s,\alpha}(x)\), this construction extends globally to produce a global symplectification (see Lemma~\ref{center-frame}).

\vspace{0.2cm}
\noindent \textbf{ Step III: Global Symplectification by holonomy}

Given $ \delta_2>\epsilon > 0$, the interval $\I$ is covered by $N \leq \left\lfloor \frac{|\I|}{\epsilon} \right\rfloor + 2$ closed sub-intervals $\{\widehat{\I}_s = [i_s^-, i_s^+]\}_{s=1}^N$ with $0 < |\widehat{\I}_s| < \epsilon$ and $i_s^+ = i_{s+1}^-$ for $s = 1, \dots, N-1$. 
Within the framework of \textbf{symplectic parallel transport}, we define a global frame for the center bundle as follows, as shown in Figure \ref{fig:placeholder}. First, fix the initial center space at $i_1^-$:
$
E_{i_1^-,\alpha}^c = \operatorname{span}_{\mathbb{C}}\{u_d^{i_1^-,\alpha}(x), u_{-d}^{i_1^-,\alpha}(x)\},
$
using the frame comes from Proposition \ref{prop:globalsympframe}. We then propagate this frame along the partitioned interval:

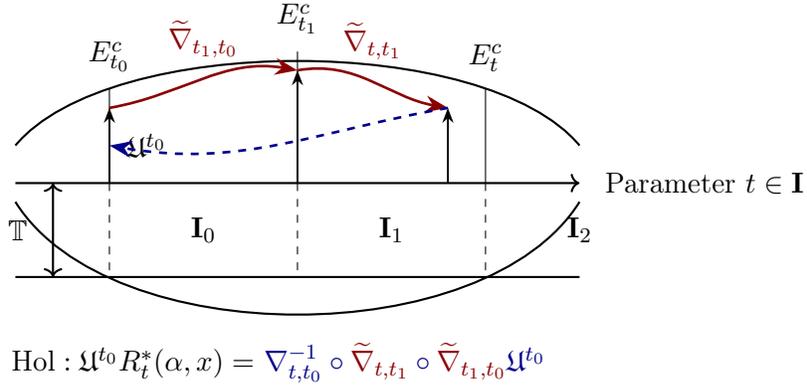
\begin{figure}
    \centering
    \begin{tikzpicture}[
    scale=2.5,
    vector/.style={-{Stealth[length=3mm]}, line width=1pt},
    parallel path/.style={decorate, decoration={snake, amplitude=.4mm, segment length=2mm}},
    basis/.style={-{Stealth[length=2mm]}, line width=0.7pt},
    annotation/.style={font=\footnotesize}
]
\begin{scope}[xshift=4.5cm, yshift=-1cm]
    \draw[thick, ->] (0,0) -- (3,0) node[right, xshift=2mm] {Parameter $t \in \I$};
    \draw[thick] (0,-0.5) -- (3,-0.5);
    \draw[thick, <->] (0.2,-0.5) -- (0.2,0) node[midway, left, xshift=-2mm] {$\mathbb{T}$};
    
    \foreach \x in {0.5,1.5,2.5} {
        \draw[dashed] (\x,0) -- (\x,-0.5);
    }
    \node at (1,-0.25) {$\I_{0}$};
    \node at (2,-0.25) {$\I_{1}$};
    \node at (3,-0.25) {$\I_{2}$};
    
    \draw[thick] (0,0.2) .. controls (0.5,0.8) and (2.5,0.8) .. (3,0.2)
                  (0,-0.1) .. controls (0.5,-0.9) and (2.5,-0.9) .. (3,-0.1);
    
    \draw (0.5,0) -- (0.5,0.5) node[above, yshift=1mm] {$E^c_{t_0}$};
    \draw[basis] (0.5,0) -- (0.5,0.4) node[midway, right, xshift=1mm] {$\mathfrak{U}^{t_0}$};
    
    \draw (1.5,0) -- (1.5,0.7) node[above, yshift=1mm] {$E^c_{t_1}$};
    \draw[basis] (1.5,0) -- (1.5,0.6);
    
    \draw (2.5,0) -- (2.5,0.5) node[above, yshift=1mm] {$E^c_{t}$};
    \draw[basis] (2.3,0) -- (2.3,0.4);
    
    \draw[vector, darkred] (0.5,0.4) to[out=10, in=170] node[midway, above, yshift=1.5mm] {$\widetilde{\nabla}_{t_1,t_0}$} (1.5,0.6);
    \draw[vector, darkred] (1.5,0.6) to[out=10, in=170] node[midway, above, yshift=1.5mm] {$\widetilde{\nabla}_{t,t_1}$} (2.3,0.4);
    
    \draw[vector, darkblue, dashed] (2.3,0.4) to[out=190, in=-10] node[midway, below, yshift=-25mm]
        {${\color{black}\text{Hol} : \mathfrak{U}^{t_0} R_t^*(\alpha, x) =\  } { \color{darkblue}{\nabla}^{-1}_{t,t_0}} \circ {\color{darkred}{\widetilde{\nabla}}_{t,t_1}} \circ {\color{darkred}{\widetilde{\nabla}}_{t_1,t_0}}\mathfrak{U}^{t_0}$} (0.5,0.2);
\end{scope}
\end{tikzpicture}
    \caption{Parallel Transport and Holonomy}
    \label{fig:placeholder}
\end{figure}

\begin{enumerate}
    \item For $t \in \widehat{\I}_1 = [i_1^-, i_1^+]$, define 
    \[
    \widetilde{u}_{\pm d}^{t,\alpha}(x) := \widetilde{\nabla}_{t, i_1^-}\, u_{\pm d}^{\,i_1^-,\alpha}(x).
    \]
    \item For $t \in \widehat{\I}_s = [i_s^-, i_s^+]$ ($s \geq 2$), define 
    \[
    \widetilde{u}_{\pm d}^{t,\alpha}(x) := 
    \widetilde{\nabla}_{t, i_s^-} \circ\ 
    \widetilde{\nabla}_{i_{s-1}^+, i_{s-1}^-}\circ\ \cdots \circ\ 
    \widetilde{\nabla}_{i_1^+, i_1^-}\, u_{\pm d}^{\,i_1^-,\alpha}(x).
    \]
\end{enumerate}

With the help of local holonomy, we obtain:

\begin{lemma}\label{center-frame}
The frame 
\begin{equation}\label{eq:frak-U}
    \widetilde{\mathfrak{U}}^{t,\alpha}(x) := 
\begin{pmatrix} \widetilde{u}_d^{t,\alpha}(x) & \widetilde{u}_{-d}^{t,\alpha}(x) \end{pmatrix}
\in C^{\omega}_{0,h'}\!\bigl(\OO \times \T, \mathrm{M}(2d,2)\bigr)
\end{equation}
forms a \textbf{canonical symplectic basis} for the center bundle $E_{t,\alpha}^c(\cdot)$, which is continuous in $t \in \I$ and piecewise analytic in $t$ on each subinterval $\widehat{\I}_s$.
\end{lemma}

\begin{proof}
By Lemma~\ref{lem:graphy} and Lemma~\ref{lem:N_s(t)}, the vectors $\widetilde{u}_{\pm d}^{t,\alpha}(x)$ are well defined on $\widehat{\I}_1$.  
By iteration, it suffices to check well-definedness on $\widehat{\I}_2$.  
Using the local holonomy relation \eqref{eq:local-holonomy}, we have
\[
    \widetilde\nabla_{i_2^-,i_1^-,\alpha} \mathfrak{U}^{\,i_1^-,\alpha}(x)
 \;=\;  \mathfrak{U}^{\,i_2^-,\alpha}(x)\, R_{i_2^-,i_1^-,\alpha}(x).
\]
Applying \eqref{eq:local-holonomy} again for $t \in \widehat{\I}_2$,
\begin{eqnarray*}
\widetilde{\mathfrak{U}}^{t,\alpha}(x)
   &=& \widetilde\nabla_{t,i_2^-,\alpha}\circ\,
      \widetilde\nabla_{i_2^-,i_1^-,\alpha}\mathfrak{U}^{\,i_1^-,\alpha}(x) \\
   &=& \widetilde\nabla_{t,i_2^-,\alpha}\mathfrak{U}^{\,i_2^-,\alpha}(x)\, R_{i_2^-,i_1^-,\alpha}(x) 
   = \mathfrak{U}^{t,\alpha}(x)\, R_{t,i_2^-,\alpha}(x)\,R_{i_2^-,i_1^-,\alpha}(x).
\end{eqnarray*}
Hence $\widetilde{\mathfrak{U}}^{t,\alpha}(x)$ is analytic in $(t,x)\in\widehat{\I}_2\times\T_{h'} $ and continuous in $\alpha \in \OO$.  
Moreover, the symplectic condition \eqref{sym-con} ensures that at each partition point the transported frame spans $E_{i_2^-,\alpha}^c(\cdot)$, providing a consistent reference for subsequent transports.
\end{proof}

Therefore, the family $t \mapsto \widetilde{u}_{\pm d}^{t,\alpha}(x)$ defines a \textbf{path‐ordered holonomy} along the broken path
\[
  i_1^- \;\longrightarrow\; i_1^+= i_2^- \;\longrightarrow\;\cdots\;\longrightarrow\; i_{s-1}^+= i_s^- \;\longrightarrow\; t.
\]
Composing over all segments yields the \textbf{total holonomy}
\begin{equation}\label{eq:total-holonomy}
 \widetilde{\mathfrak{U}}^{t,\alpha}(x)
   = \widetilde\nabla_{t,i_s^-,\alpha}\circ\cdots\circ\,
     \widetilde\nabla_{i_1^+,i_1^-,\alpha} \mathfrak{U}^{\,i_1^-,\alpha}(x)
   = \mathfrak{U}^{t,\alpha}(x)\,\widetilde{R}_{t,\alpha}(x),
\end{equation}
where
\[
  \widetilde{R}_{t,\alpha}(x) :=
  R_{t,i_s^-,\alpha}(x)\,R_{i_{s}^-,i_{s-1}^-,\alpha}(x)\,\cdots\,R_{i_2^-,i_1^-,\alpha}(x),
\]
which is piecewise $C^\omega$ in $t$, continuous in $\alpha$ and analytic in $x$, and continuous on $\I$.
Indeed,  let  $\nabla_{t,s,\alpha}$ denote the transformation induced by the reference frame, i.e.
\[
\nabla_{t,s,\alpha} u_{\pm d}^{s,\alpha}(x) = u_{\pm d}^{t,\alpha}(x).
\]
Then the total holonomy $\widetilde{R}_{t,\alpha}(x)$ satisfies
\begin{equation}\label{hol}
    \mathfrak{U}^{\,t_0,\alpha}(x)\, \widetilde{R}_{t,\alpha}(x) \;=\;
    \nabla^{-1}_{t,t_0,\alpha} \circ\ \widetilde{\nabla}_{t,t_1,\alpha} \circ\ 
    \widetilde{\nabla}_{t_1,t_0,\alpha} \mathfrak{U}^{\,t_0,\alpha}(x),
\end{equation}
as illustrated schematically in Figure~\ref{fig:placeholder}.  
See also Remark~\ref{rem:holo} for further discussion of holonomy.

{\bf \noindent Substep IV: Verifying monotonicity in center }

Once the symplectic basis for the center bundle $E_{t,\alpha}^{c}$ is constructed, we define 
\[
\widetilde{B}_{t,\alpha}(x) 
= \Big[\, u_1^{t,\alpha}(x), \dots, u_{d-1}^{t,\alpha}(x), \;\widetilde{u}_d^{t,\alpha}(x) 
\;\big|\;
u_{-1}^{t,\alpha}(x), \dots, u_{-(d-1)}^{t,\alpha}(x), \;\widetilde{u}_{-d}^{t,\alpha}(x)\,\Big].
\]
By Proposition~\ref{prop:globalsympframe} and Lemma~\ref{center-frame},  
we obtain that $
\widetilde{B}_{t,\alpha}(\cdot) \in C^\omega_{0,h'}(\OO\times\mathbb{T}, \mathrm{HSp}(2d)),
$
piecewise $C^\omega$ and continuous in $t \in \I$, such that  
\begin{equation}\label{eq:conj}
\widetilde{B}_{t,\alpha}(x+\alpha)^{-1}\,A_{t,\alpha}(x)\,\widetilde{B}_{t,\alpha}(x) 
= H_{t,\alpha}(x) \;\diamond\; \widetilde{C}_{t,\alpha}(x),    
\end{equation}
with $\widetilde{C}_{t,\alpha}(x) \in C^\omega_{0,h'}(\OO\times\mathbb{T}, \mathrm{HSp}(2))$ piecewise $C^\omega$ and continuous in $t\in\I$. 

\medskip 

We also remark that if one uses the reference frame 
\[
B_{t,\alpha}(x) 
= \Big[\, u_1^{t,\alpha}(x), \dots, u_{d-1}^{t,\alpha}(x), \;u_d^{t,\alpha}(x) 
\;\big|\;
u_{-1}^{t,\alpha}(x), \dots, u_{-(d-1)}^{t,\alpha}(x), \;u_{-d}^{t,\alpha}(x)\,\Big],
\] 
as in Proposition~\ref{prop:globalsympframe}, then by the total holonomy relation \eqref{eq:total-holonomy} we have 
\begin{equation}\label{tb-b}
    \widetilde{B}_{t,\alpha}(x)
    \;=\; B_{t,\alpha}(x)\,\big(I_{2d-2}\diamond \widetilde{R}_{t,\alpha}(x)\big).
\end{equation}

For any isotropic vector \( v \in \mathbb{C}^{2d} \), define (identify $\psi=\omega_d$ below)
\[
\Psi_{A_t}(v,x) := \psi \big( A_t(x) v,\ \partial_t A_t(x) v \big).
\]
It is clear that \( \Psi_{A_t}(v,x) \) is continuous in \( (t,v,x) \).  
By compactness, there exist real functions \( m(t,x), M(t,x)  \) such that
\[
m(t,x) := \inf_{\|v\|=1} \Psi_{A_t}(v,x), 
\quad 
M(t,x) := \sup_{\|v\|=1} \Psi_{A_t}(v,x).
\]
\begin{lemma}\label{lem:LUbd}
There exists
$C'''=C'''\bigl(\delta', \|u_{\pm d}^{t,\alpha}\|_{\delta',\OO,h'}
\bigr)>0$ such that for any $0<\epsilon<\delta_2$, after excluding a finite set of points where differentiability fails, for every unit isotropic vector \( v \in \mathbb{C}^{2} \), 
one has
\begin{equation}\label{eq:LUbd}
    \begin{aligned}
        (C''')^{-1}\left(\inf_{y\in\T}\lambda_{\min}W_{t,\alpha}(y)\right) &\,m(t,x) - C''' \|C_{t,\alpha}(\cdot)\|_0^2\,\epsilon\\
    &< \Psi_{\widetilde C_{t,\alpha}}(v,x)<C''' \|\lambda_{\max}W_{t,\alpha}(\cdot)\|_0 \,M (t,x)
      + C''' \|C_{t,\alpha}(\cdot)\|_0^2\,\epsilon.
    \end{aligned}
\end{equation}
where $W_{t,\alpha}(x)={\mathfrak{U}}^{t,\alpha}(x)^* {\mathfrak{U}}^{t,\alpha}(x)$.
\end{lemma}

\begin{proof}

Let \( u = \begin{pmatrix} u_1, u_2 \end{pmatrix}^{\top} \in \mathbb{C}^2 \) be a unit isotropic vector, and set 
\[
\bar{u} := (0, \dots, 0, u_1, 0, \dots, 0, u_2)^\top \in \mathbb{C}^{2d},
\]
where \( u_1 \) and \( u_2 \) occupy the $d$-th and $2d$-th entries, respectively.  
Define
\[
\bar{v} := \widetilde{B}_{t,\alpha}(x) \bar{u}, 
\quad 
\bar{w} := A_t(x)\, \widetilde{B}_{t,\alpha}(x)\, \bar{u}, 
\quad 
\bar{x} := \widetilde{B}_{t,\alpha}(x+\alpha)^{-1} A_{t,\alpha}(x)\, \widetilde{B}_{t,\alpha}(x)\, \bar{u}.
\]
Then $\bar{u},\bar{v},\bar{w},\bar{x}$ are isotropic vectors in $(\C^{2d},\psi)$.

For $t \in \overset{\circ}{\widehat{\I}}_s$, a direct computation yields:
\begin{equation}\label{eq:Psi-C-chain}
    \begin{aligned}
    \Psi_{\widetilde{C}_{t,\alpha}}(u,x)
    :&= \psi\left( \widetilde{C}_{t,\alpha}(x) u,\, \partial_t \widetilde{C}_{t,\alpha}(x) u \right) \\
    &= \psi \left( \widetilde{B}_{t,\alpha}(x) \bar{u},\, \partial_t \widetilde{B}_{t,\alpha}(x) \bar{u} \right)
    - \psi \left( \widetilde{B}_{t,\alpha}(x+\alpha) \bar{x},\, \partial_t \widetilde{B}_{t,\alpha}(x+\alpha) \bar{x} \right) \\
    &\quad + \psi \left( A_t(x)\bar{v},\, (\partial_t A_t)(x) \bar{v} \right).
\end{aligned}
\end{equation}
Furthermore, as indicated in \eqref{eq:conj}, the vector $\bar{x}$ is expressed in the form $$ \bar{x} := (0, \dots, 0, x_1, 0, \dots, 0, x_2)^\top \in \mathbb{C}^{2d},$$ we reduce to estimating terms of the form
\[
\psi\!\left( \widetilde{\mathfrak{U}}^{t,\alpha}(x) v,\ \partial_t \widetilde{\mathfrak{U}}^{t,\alpha}(x) v \right), 
\quad v\in \C^2,
\]
which is addressed in the following quantitative estimate.

\begin{lemma}\label{lem:quanti-tildeB}
   There exists $\eta'=\eta'\bigl(\delta', \|u_{\pm d}^{t,\alpha}\|_{\delta',\OO,h'}
\bigr)>0$ such that for any $t \in (i_s^-, i_s^+) $, we have estimate:
\begin{align}
    \sup_{(x,\alpha, v)\in \T \times \OO \times S}\left| \psi \left( \widetilde{\mathfrak{U}}^{t,\alpha} v, \partial_t \widetilde{\mathfrak{U}}^{t,\alpha} v \right) \right| &\leq e^{4|\I|\eta} \eta' \epsilon.
\end{align}
where $S:=\{ v\in \C^2: \|v\|=1 \}$.
\end{lemma}
\begin{proof}
First we prove for any $s\in \I$, $t\in 
(s-\epsilon,s+\epsilon)\cap \I_s^{'}$, there exists $\eta'\bigl(\delta', \|u_{\pm d}^{t,\alpha}\|_{\delta',\OO,h'}
\bigr)>0$ such that
\begin{equation}\label{eq:delta-epsilon}
    \sup_{(x,\alpha,v)\in\T\times \OO\times S}\left|\psi(\widetilde{\nabla}_{t,s,\alpha} \begin{pmatrix} u_{d}^{s,\alpha}(x) & u_{-d}^{s,\alpha}(x) \end{pmatrix}\,v, \partial_{t} \widetilde{\nabla}_{t,s,\alpha} \begin{pmatrix} u_{d}^{s,\alpha}(x) & u_{-d}^{s,\alpha}(x) \end{pmatrix}\,v )\right|<\eta' \epsilon.
\end{equation}
 To see this, 
define 
\[
f_s(t,\alpha,v,x) := \psi\left(
\widetilde{\nabla}_{t,s,\alpha}\mathfrak{U}^{s,\alpha}(x)v,\quad
\partial_t \widetilde{\nabla}_{t,s,\alpha}\mathfrak{U}^{s,\alpha}(x)v
\right).
\]
Using the decomposition
\[
\widetilde{\nabla}_{t,s,\alpha}\mathfrak{U}^{s,\alpha}(\cdot)
=\mathfrak{U}^{s,\alpha}N_s(t,\alpha,\cdot) + \Phi_{t,s,\alpha}\mathfrak{U}^{s,\alpha}N_{s}(t,\alpha,\cdot),
\]
we expand $f_s(t,\alpha,v,x) $ as (omit $\alpha$ from notation for clarity)
\begin{align*}
& f_s(t,\alpha,v,\cdot) \\ 
&= \underbrace{
\psi\left(
\mathfrak{U}^s(\cdot)N_s(t,\cdot) v,\quad
\mathfrak{U}^s(\cdot)\partial_t N_s(t,\cdot) v
\right)
}_{\text{(I)}} 
\quad + \underbrace{
\psi\left(
\Phi_{t,s}\mathfrak{U}^s(\cdot)N_s(t,\cdot) v,\quad
\mathfrak{U}^s(\cdot)\partial_t  N_s(t,\cdot) v
\right)
}_{\text{(II)}} \\
&\quad + \underbrace{
\psi\left(
\Phi_{t,s}\mathfrak{U}^s(\cdot)N_s(t,\cdot) v,
\Phi_{t,s}\mathfrak{U}^s(\cdot)\partial_t N_s(t,\cdot) v
\right)
}_{\text{(III)}} 
\quad + \underbrace{
\psi\left(
\mathfrak{U}^s(\cdot)N_s(t,\cdot) v,
\Phi_{t,s}\mathfrak{U}^s(\cdot)\partial_t N_s(t,\cdot) v
\right)
}_{\text{(IV)}} \\
&\quad + \underbrace{
\psi\left(
\Phi_{t,s}\mathfrak{U}^s(\cdot)N_s(t,\cdot) v,
\partial_t \left[\Phi_{t,s}\mathfrak{U}^s(\cdot)\right] N_s(t,\cdot) v
\right)
}_{\text{(V)}} 
\quad + \underbrace{
\psi\left(
\mathfrak{U}^s(\cdot)N_s(t,\cdot) v,
\partial_t \left[\Phi_{t,s}\mathfrak{U}^s(\cdot)\right] N_s(t,\cdot) v
\right)
}_{\text{(VI)}}.
\end{align*}
Since \( E_{s,\alpha}^s \oplus E_{s,\alpha}^u \) is symplectic orthogonal to \( E_{s,\alpha}^c \), the terms (II), (IV), and (VI) vanish. Moreover, when \( t = s \), we have \( \Phi_{s,s} = 0 \) and 
by Lemma \ref{lem:N_s(t)}, we have $\partial_t\bigl|_{t=s}N_s(t,\cdot)=0$, and therefore 
\( f_s(s,\alpha,v,\cdot) = 0 \).
Then one can write  
\[
	f_s(t,\alpha,v,\cdot)=f_s'(\xi,\alpha,v,\cdot)(t-s), \  \xi\in (s-\epsilon,s+\epsilon)
\] 
where $ f_s'(\xi,\alpha,v,\cdot) $ is controlled by the first and second derivatives of $ u_{\pm d}^{t,\alpha}( x) $, by \eqref{eq:derivative-bound} \eqref{eq:partial-J_s}  and Lemma \ref{lem:N_s(t)} and hence $ \sup_{( \alpha, x,v)\in \OO\times \T \times S}f_s'(\xi,\alpha,v,x) $ is uniformly bounded by some $\eta'\bigl(\delta', \|u_{\pm d}^{t,\alpha}\|_{\delta',\OO,h'}
\bigr)>0$.

Thus, for $t \in (i_s^-, i_s^+)$, again by the using total holonomy  \eqref{eq:total-holonomy},  we can have estimate: 
\begin{eqnarray*}
\frac{1}{\|v\|^2} \left| \psi \left( \widetilde{\mathfrak{U}}^t v, \partial_t \widetilde{\mathfrak{U}}^t v \right) \right| 
&=& \frac{1}{\|v\|^2}  \left| \psi \left( \widetilde{\nabla}_{t,i_s^-} \mathfrak{U}^{i_s^-}\widetilde R_{i_s^-} v, \partial_t \widetilde{\nabla}_{t,i_s^-} \mathfrak{U}^{i_s^-} \widetilde R_{i_s^-} v\right)\right| \\
&=& \frac{\|\widetilde R_{i_s^-}v\|^2}{\|v\|^2}  \frac{1}{\|\widetilde R_{i_s^-}v\|^2}\left| \psi \left( \widetilde{\nabla}_{t,i_s^-} \mathfrak{U}^{i_s^-} \widetilde R_{i_s^-} v, \partial_t \widetilde{\nabla}_{t,i_s^-} \mathfrak{U}^{i_s^-}  \widetilde R_{i_s^-}(x) v\right)\right|\\
&\leq& \|\widetilde R_{i_s^-}\|_0^2 \sup_{v\in \C^2-\{0\}} \frac{1}{\|v\|^2} \left| \psi \left( \widetilde{\nabla}_{t,i_s^-} \mathfrak{U}^{i_s^-} v, \partial_t \widetilde{\nabla}_{t,i_s^-} \mathfrak{U}^{i_s^-}  v\right)\right|\\
&\leq&  e^{4|\I|\eta} \eta' \epsilon,
\end{eqnarray*}
where the last estimate follows by \eqref{eq:delta-epsilon} and  Lemma \ref{lem:localR_t}. 
\end{proof}

Applying Lemma~\ref{lem:quanti-tildeB}, we estimate:
\begin{align}
    \Psi_{\widetilde{C}_{t,\alpha}}(u,x) 
    &> \psi \left( A_t(x) \bar{v},\, (\partial_t A_t)(x) \bar{v} \right) \notag\\
    &\quad - \left| \psi \left( \widetilde{B}_{t,\alpha}(x) \bar{u},\, \partial_t \widetilde{B}_{t,\alpha}(x) \bar{u} \right) \right| 
    - \left| \psi \left( \widetilde{B}_{t,\alpha}(x+\alpha) \bar{x},\, \partial_t \widetilde{B}_{t,\alpha}(x+\alpha) \bar{x} \right) \right| \label{eq:4.28}\\
    &> \|\bar{v}\|^2m(t,x) 
    - \Big(1+ \|\bar{x}\|^2\Big) e^{4|\I|\eta} \eta' \epsilon \notag
\end{align}
 Then by \eqref{eq:not-modif-conj}, \eqref{eq:total-holonomy} and \eqref{tb-b}, we have
\begin{align*}
    \|\bar{x}\|^2&=\|\widetilde{R}_{t,\alpha}(x+\alpha
)^{-1} C_{t,\alpha}(x)\widetilde{R}_{t,\alpha}(x
)u\|^2\leq \|\widetilde{R}_{t,\alpha}\|_0^4\|C_{t,\alpha}(\cdot)\|_0^2,\\
\|\bar{v}\|^2&=\|(\widetilde{R}_{t,\alpha}u)^*\,\mathfrak{U}^{t,\alpha}(x)^*\mathfrak{U}^{t,\alpha}(x)\,\widetilde{R}_{t,\alpha}u\|\geq \|\widetilde{R}_{t,\alpha}\|_0^{-2}\inf_{y\in\T}\lambda_{\min}W_{t,\alpha}(y),\\
\|\bar{v}\|^2&=\|(\widetilde{R}_{t,\alpha}u)^*\,\mathfrak{U}^{t,\alpha}(x)^*\mathfrak{U}^{t,\alpha}(x)\,\widetilde{R}_{t,\alpha}u\|\leq \|\widetilde{R}_{t,\alpha}\|_0^{2}\|\lambda_{\max}W_{t,\alpha}(\cdot)\|_0.
\end{align*}
Thus by Lemma \ref{lem:localR_t}, \eqref{eq:LUbd} holds.
\end{proof}

Taking in the construction in Step III. To complete the whole proof, it suffices to verify its monotonicity properties. 
By compactness and the positive definiteness of $W_{t,\alpha}(x)$,
$\inf_{x\in\T}\lambda_{\min}W_{t,\alpha}(x)$ is uniformly bounded away from zero,
while $\|\lambda_{\max}W_{t,\alpha}(\cdot)\|_0$ and
$\|C_{t,\alpha}(\cdot)\|_0$ are uniformly bounded for
$(t,\alpha)\in\I\times\OO$. Hence there exists $ C_*\bigl(\I,\OO,\delta', \|u_{\pm d}^{t,\alpha}\|_{\delta',\OO,h'}
\bigr)>0 $ such that \begin{equation}\label{eq:C_*mono}
     C_*^{-1} m(t,x) - C_*\epsilon<\Psi_{\widetilde{C}_{t,\alpha}}(u,x) <C_* M (t,x)
      + C_* \epsilon.
\end{equation}
Monotonicity of $A_t$ means $m(t,x)>c>0$,
the proof is thus completed by choosing $\epsilon$ in Lemma \ref{lem:LUbd} small enough.
\end{proof}

\begin{remark}\label{rem:holo}
The holonomy matrices $R_{t,s,\alpha}(x)$ play a pivotal role in our construction. We clarify their relationship to classical holonomy theory \cite{SS}.  For a given connection $\nabla$ of the bundle, the classical  \textbf{holonomy group} $\rmm{Hol}_b(\nabla)$ at $b \in \mathcal{B}$ is:
    \[
    \rmm{Hol}_b(\nabla) = \{ P_\gamma \in \GL(\mathcal{E}_b) \mid \gamma \ \text{is a closed loop based at}\ b \},
    \]
    which measures path-dependence of parallel transport.
    
    While not classical closed-loop holonomy, the local holonomy  $R_{t,s,\alpha}(x)$ is a \textbf{generalized holonomy operator} that:
    \begin{itemize}
     
  \item \textbf{Path-dependence}: $\widetilde{R}_{t,\alpha}$ measures how parallel transport deviates from a reference due to curvature, the \textit{core feature} of holonomy.
    
    \item \textbf{Curvature bounds}: The estimate Lemma \ref{lem:localR_t} encodes curvature data of $\widetilde\nabla_{t,s,\alpha}$.

  \item \textbf{Monodromy interpretation}: Relative to the fixed initial frame $\mathfrak{U}^{s,\alpha}$, 
    $\widetilde{R}_{t,\alpha}(x)$
     exhibits the usual \textbf{monodromy} interpretation of holonomy as the mismatch between initial and transported frames, as was shown in \eqref{hol} and Figure \ref{fig:placeholder}.
\end{itemize}
\end{remark}

\begin{remark}\label{rem:premono}
Recall that a one-parameter family of cocycles 
\(A_t(x)\in C^1(\mathbf{I}\times\mathbb{T}, \mathrm{HSp}(2d,\psi))\) is called
\textit{premonotonic} if some iterate of it is monotonic.
Note that Proposition~\ref{prop:globalmonoframe} also holds for premonotonic cocycles.
Indeed, assume that there exists a transformation \(\widetilde B_{t,\alpha}\)
which block-diagonalizes an iterate \((\alpha,A_t)^n\) of the cocycle.
Then the same transformation also block-diagonalizes the original cocycle \(A_t\), that is,
\begin{equation}\label{eq:block-diag-A_t}
\widetilde B_{t,\alpha}(x+\alpha)^{-1} A_t(x)\widetilde B_{t,\alpha}(x)
= \widehat{H}_{t,\alpha}(x)\diamond \widehat{C}_{t,\alpha}(x),
\end{equation}
where \(\widehat{C}_{t,\alpha}(x)\in C^\omega_{0,h'}(\OO\times\mathbb{T},\mathrm{HSp}(2))\)
is piecewise analytic in \(t\in\mathbf{I}\).
Moreover, since the invariant splitting associated with \((\alpha,A_t)\)
coincides with that of its iterate \((\alpha,A_t)^n\), we have
\begin{equation}\label{eq:iterate-2}
(\alpha,\widehat{H}_{t,\alpha})^n
=
(n\alpha, H_{t,\alpha}),\qquad (\alpha,\widehat{C}_{t,\alpha})^n
=
(n\alpha,\widetilde{C}_{t,\alpha}).
\end{equation}
Therefore, the piecewise monotonicity of
\((\alpha,\widetilde{C}_{t,\alpha})\) implies the piecewise premonotonicity of
\((\alpha,\widehat{C}_{t,\alpha})\). {Moreover, Lemma \ref{lem:LUbd} remains valid for $(\alpha,A_t)$, since it does not assume monotonicity a priori. }
\end{remark}

\section{Projective action and the fibred rotation number}\label{sec:proj-rotationnumber}

\subsection{Projective action}
Recall that \(\psi\) is a Hermitian-symplectic form on \(\C^{2d}\) with structure matrix
\(\mathcal S\), i.e. $ 
\psi(u,v)=u^* \mathcal S v,$
and \(\mathcal S\) is congruent to the standard symplectic matrix \(\mathcal J_{2d}\),
i.e. \(\mathcal S=\mathcal P^*\mathcal J_{2d}\mathcal P\) for some invertible
\(\mathcal P\in\mathrm{GL}(2d,\C)\).  By Section~\ref{sec:symp}, the Cayley element (replaced when necessary by
\(\mathcal C\mathcal P\)) induces an identification
\[
\operatorname{Lag}(\C^{2d},\psi)\;\cong\;\mathbf U(d),
\]
so that every Lagrangian subspace \(\Lambda\) is represented by a unitary matrix
\(W_\Lambda\in\mathbf U(d)\).  If \(A\in\mathrm{HSp}(2d,\psi)\), then \(A\) preserves
\(\operatorname{Lag}(\C^{2d},\psi)\) and the correspondence
 $ W_\Lambda \;\longmapsto\; W_{A\Lambda}$
realizes the projective action of \(A\) on \(\mathbf{U}(d)\).

In what follows we adopt the common convention that our Hermitian-symplectic cocycle
\(A\) is homotopic to the identity.  Fix once and for all a homotopy
\(\{A^s\}_{s\in[0,1]}\subset\mathrm{HSp}(2d,\psi)\) with \(A^0=\mathrm{id}\) and
\(A^1=A\).  The chosen homotopy determines a unique continuous lift of the
projective action to the universal cover \(\widetilde{\mathbf U(d)}\).  Concretely,
the determinant map \(\det:\mathbf U(d)\to S^1\) has kernel \(\mathbf{SU}(d)\),
and \(\mathbf{SU}(d)\) is simply connected for \(d\ge2\) (the trivial case \(d=1\)
is immediate).  Lifting the \(S^1\)-coordinate to its universal cover \(\R\)
while keeping the \(\mathbf{SU}(d)\)-factor yields the identification
\(\widetilde{\mathbf U(d)}\cong \R\times\mathbf{SU}(d)\).  We therefore obtain
a lift
\[
F_A:\R\times\mathbf{SU}(d)\longrightarrow \R\times\mathbf{SU}(d)
\]
of the projective action which is homotopic to the identity and is determined by the
chosen homotopy \(\{A^s\}\).

That is, if a unitary matrix \( W_\Lambda \in \mathbf{U}(d) \) lifts to a pair \( (x_\Lambda, S_\Lambda) \in \mathbb{R} \times \mathbf{SU}(d) \) with $W_\Lambda = e^{2\pi i x_\Lambda} S_\Lambda$, then \( F_A \) maps it to
\[
  F_A(x_\Lambda, S_\Lambda) 
   := (x_{A\Lambda}, S_{A\Lambda}),\]
which lifts \( W_{A\Lambda} \). A priori $x_{A\Lambda}$ may not be uniquely-defined, but as $A$ is homotopic to identity, we could assume that along the deformation path $x_{A\Lambda}$ will continuously move to $x_\Lambda$, obviously such $x_{A\Lambda}$ is unique,
so \(F_A\) is well defined. Using the lift \(F_A\) define the phase map \(\rho:\operatorname{Lag}(\C^{2d},\psi)\to\R\)
by choosing the real coordinate of any lift:
\[
\rho(\Lambda):=dx_\Lambda\qquad\text{where }W_\Lambda=e^{2\pi i x_\Lambda}S_\Lambda.
\]
By construction \(\rho(\Lambda)\) is a real-valued lift of \(\det W_\Lambda\), i.e.
\(\pi(\rho(\Lambda))=\det W_\Lambda\) where \(\pi:\R\to S^1\) is \(\pi(s)=e^{2\pi i s}\).
The lift $\rho$ is not unique --- it is only defined up to an integer,
but the difference
\[
  \phi \equiv \phi[A](\Lambda) = \rho(A\Lambda) - \rho(\Lambda)
\]
does not depend on the choice of lift and hence descends to a well-defined real-analytic function $\operatorname{Lag}(\mathbb{C}^{2d},\psi)\to \R$. {In the following, we just call $\phi$ the phase of the cocycle $A$.} Different choices of deformation path lead to phases $\phi$ that differ by integers, and the integers are  determined by the homotopy, not by the choice of lift.  We refer to \cite[Section 4]{LW} for further details.

Moreover, the following uniform bound holds:
\begin{lemma}[{\cite[Lemma~3.6]{LW}}]\label{lem:minmax} $
  \max_{\Lambda} \phi[A](\Lambda) \;-\; \min_{\Lambda} \phi[A](\Lambda) \;<\; d.$
\end{lemma}

Given any \(A\in C^0(\R,\rmm{HSp}(2d,\psi))\), it is always homotopic to the identity in the sense that there exists a homotopy \(\{A^s(x)\}_{s\in[0,1]}\subset C^0(\R,\mathrm{HSp}(2d,\psi))\) with \(A^0(x)=\mathrm{id}\) and
\(A^1(x)=A\). Consequently one may make a choice of 
$\phi_x(\Lambda)\equiv\phi[A(x)](\Lambda)
$
depending continuously on \(x\in\R\).

By contrast, if \(A\in C^0(\T,\mathrm{HSp}(2d,\psi))\), a globally continuous choice of \(\phi_x(\Lambda)\) on \(\T\) need not exist.  Indeed, when viewed on \(\R\) any chosen lift satisfies
$\phi_{x+1}(\Lambda)=\phi_x(\Lambda)+m,$
for some integer \(m\in\Z\) which does not depend on the choice of lift but reflects the topological degree (winding number) of the loop \(x\mapsto\det\!\bigl(W_{\Lambda}^{-1}W_{A(x)\Lambda}\bigr)\).  Thus individual phase lifts \(\phi_x(\Lambda)\) are only defined on \(\T\) up to an integer ambiguity; nevertheless the difference of two such lifts is free of that ambiguity, as the next lemma states.

\begin{lemma}\label{non-homo}
Fix \(A\in C^0(\T,\mathrm{HSp}(2d,\psi))\) and two Lagrangian frames
\(\Lambda_1, \Lambda_2 \in\operatorname{Lag}(\mathbb{C}^{2d},\psi)\)
the difference
\[ \eta[A(x)](\Lambda_1,\Lambda_2) := \phi_x(\Lambda_1) - \phi_x(\Lambda_2), \]
is well-defined as a continuous function on \(\T\), and independent of the choices of lifts.
\end{lemma}

If \(A\in C^0(\T,\mathrm{HSp}(2d,\psi))\) is homotopic to a constant, so the integer obstruction discussed above is zero. Equivalently, there exists a consistent choice of 
\(\phi_x(\Lambda)=\phi[A(x)](\Lambda)\)
which depends continuously on \(x\in\T\). In this case one may define the fibred rotation number of the cocycle using such a choice.

\subsection{Fibred rotation number.}

Let $\alpha\in \R$, \(A\in C^0(\T,\rmm{HSp}(2d,\psi))\)  homotopic to a constant.  Let \(\phi_{x,1}=\phi_x\), we may define 
\(\phi_{x,k}\equiv\phi[A_k(x)]\) for all \(k\in\Z\) by imposing the cocycle rule
\begin{equation}\label{eq:cocycle}
\phi_{x,k+l}(\Lambda)
=\phi_{x,k}(\Lambda)+\phi_{x+k\alpha,l}\big(A_k(x)\Lambda\big)\qquad(k,l\in\Z,\ \Lambda\in\op{Lag}(\C^{2d},\psi)).
\end{equation}
In particular, we obtain for \(k\ge1\)
\[
\phi_{x,k}(\Lambda)=\sum_{m=0}^{k-1}\phi_{x+m\alpha}\big(A_m(x)\Lambda\big).
\]
The fibred rotation number is defined by the averaged limit
\[
\rho_x\equiv \rho_x(\alpha,A)
:=\lim_{k\to\infty}\frac{1}{k}\,\phi_{x,k}(\Lambda).
\]
The limit exists and is uniform in \(\Lambda\) (hence determines an element of \(\R/\Z\)) \cite{LW}.   From \eqref{eq:cocycle} one checks immediately that
\(\rho_{x+\alpha}=\rho_x\).  Consequently, when \(\alpha\in\R\setminus\Q\) the function
\(x\mapsto \rho_x\) is constant and we shall then denote its value simply by \(\rho\).

If \(A^t\in C^0(\I\times\T,\rmm{HSp}(2d,\psi))\) is a continuous family (in parameter \(t\in\I\)) and always homotopic to a constant, then one can choose a branch \(\phi[A^t(x)]\) so that
\((t,x)\mapsto\phi[A^t(x)]\) depends continuously on \(\I\times\T\).
With this choice the associated objects \(\phi_{x,k}^t\) and the averaged quantities
\(\rho_x(t)\) are well defined and depend continuously on \(t\) (and on \(x\) and \(A^t\)).

\subsection{Monotonicity}

In this subsection, we aim to explore the monotonicity of the projective action, with a particular focus on its relation to the monotonicity of the cocycle, while the key link is the monotonicity of  Lagrangian paths. 
Without loss of generality, we restrict our discussion to the standard form $\psi = \omega_d$, 
as the argument is coordinate invariant. 
To begin with, we start with the following basic algebraic observation:
\begin{lemma}[\cite{howardmaslov,LW}]\label{lem:Wdot-phi}
    Let \( \Lambda_t = \begin{psmallmatrix} X_t \\ Y_t \end{psmallmatrix} \) be a $C^1$-smooth family of Lagrangian frames in $\op{Lag}(\C^{2d},\omega_d)$. Then the time derivative of $ W_{\Lambda_t}$ satisfies  $$
    \partial_t W_{\Lambda_t} \;=\; i\, W_{\Lambda_t}\,\Omega_{\Lambda_t},$$
    where $\Omega_{\Lambda_t}: \op{Lag}(\C^{2d},\omega_d) \to \op{Her}_d(\C)$ is   given by 
    \begin{eqnarray*}
        \Omega_{\Lambda_t}= -2\,[(X_t-iY_t)^{-1}]^*
      \mathcal{M}(t)
      (X_t-iY_t)^{-1}
    \end{eqnarray*}
    with $\mathcal{M}(t)=\Lambda_t^* \mc{J}_{2d}\, \partial_t \Lambda_t$.
    Moreover, let $\phi(\Lambda_t)$ be any continuous lift  of $\det W_{\Lambda_t}$, i.e.\ $\pi(\phi(\Lambda_t))=\det W_{\Lambda_t}$, then 
    \[
      \partial_t \phi(\Lambda_t) \;=\; \frac{1}{2\pi}\,\op{tr} \Omega_{\Lambda_t}.
    \]
\end{lemma}

We recall the following definition due to Ekeland \cite{Ekeland}: 
\begin{definition}
    The  Lagrangian paths $\Lambda_t$ is said to be monotonic, if $\mathcal{M}(t)=\Lambda_t^* \mc{J}_{2d}\, \partial_t \Lambda_t$ is definite (strictly positive- or strictly negative-definite)
\end{definition}

Monotonicity of  Lagrangian paths is quite important. First, a basic property is the following:
\begin{lemma}[{\cite[Lemma 3.11]{HOWARD2016}}]\label{lem:lag-monotone}
    Let $\Lambda_t = \begin{psmallmatrix} X_t \\[2pt] Y_t \end{psmallmatrix}$
 be a smooth one-parameter family of Lagrangian frames in $\op{Lag}(\C^{2d},\omega_d)$.  Suppose $\mc{M}(t)$ is strictly positive (negative) definite. Then the eigenvalues of $ W_{\Lambda_t} $ move strictly (anti-)clockwise on the unit circle as $ t $ increases.
\end{lemma}

Moreover, as Hermitian symplectic group preserves monotone Lagrangian paths, which directly imply the following:  
 \begin{lemma}\label{lem:mono2}
 Let $A \in \mathrm{HSp}(2d)$, and let $\Lambda_t = \begin{psmallmatrix} X_t \\[2pt] Y_t \end{psmallmatrix}$
 be a smooth one-parameter family of Lagrangian frames in $\op{Lag}(\C^{2d},\omega_d)$,
provided \(\mathcal{M}(t)\) is definite.  
  Then the derivative $\partial_t \phi(A \Lambda_t)$ has the same sign as 
  $\partial_t \phi(\Lambda_t)$.
 \end{lemma}
 \begin{proof}
 For any fixed \(B\in\mathrm{HSp}(2d)\)
write \(B\Lambda_t=\begin{psmallmatrix}X_{B,t}\\[2pt]Y_{B,t}\end{psmallmatrix}\)
and set \(\mc{R}_{B,t}:=(X_{B,t}-iY_{B,t})^{-1}\).  Lemma
 \ref{lem:Wdot-phi} yields
 \begin{equation}\label{eq:mono-underA}
 \Omega_{B\Lambda_t} = \mc{R}_{B,t}^*\,\{\Lambda_t^* B^*\mathcal J_{2d} B\,\partial_t\Lambda_t\}\,\mc{R}_{B,t}
 \end{equation}
 Take \(B=I_{2d}\) and \(B=A\). Since \(A\in\mathrm{HSp}(2d)\) satisfies
\(A^*\mathcal J_{2d}A=\mathcal J_{2d}\), hence 
 \[ \partial_t \phi(A \Lambda_t)
= -\frac{1}{\pi}\op{tr}\bigl(\mc{R}_{A,t}^* \mathcal{M}(t)\mc{R}_{A,t}\bigr),\qquad
 \partial_t \phi( \Lambda_t)
 = -\frac{1}{\pi}\op{tr}\bigl(\mc{R}_{I,t}^* \mathcal{M}(t)\mc{R}_{I,t}\bigr).
 \]
As congruence by an invertible matrix preserves definiteness, the result then follows. 
 \end{proof}

 Of greater importance for our purposes, the monotonicity of the cocycle $A^t(\cdot)$ implies the monotonicity of Lagrangian paths $A^t(\cdot)\Lambda$; this in turn implies the monotonicity of the projective action $\phi(A^t(\cdot)\Lambda)$, and ultimately the monotonicity of the fibred rotation number. Indeed, let \(A^t(x)\in C^1(\I\times\T,\mathrm{HSp}(2d))\), for any Lagrangian subspace \(\Lambda= \begin{bsmallmatrix}X\\[2pt]Y\end{bsmallmatrix} \in\Lag(\C^{2d},\omega_d)\),  we define $\Lambda_k(t,x):=A^t_k(x)\begin{psmallmatrix}X\\[2pt]Y\end{psmallmatrix} = \begin{psmallmatrix} X(t,k)\\[2pt] Y(t,k)\end{psmallmatrix}$, and then define 
$\mc R(t,k):=\bigl(X(t,k)-iY(t,k)\bigr)^{-1}$, then we have the following:

\begin{lemma}\label{lem:omega-trace} For every Lagrangian subspace \(\Lambda\in\Lag(\C^{2d},\omega_d)\),  the derivative of the projective action of $A^t_k(x)$ is given by
\begin{equation}\label{der-phi}
    \partial_t\phi[A^t_k(x)](\Lambda)
= \frac{1}{2\pi}\,\op{tr} \bigl(\Omega_{x,\Lambda}(t,k)\bigr),
\end{equation}
 where $\Omega_{x,\Lambda}(t,k): \op{Lag}(\C^{2d},\omega_d) \to \op{Her}_d(\C)$ is given by 
\begin{align}
\Omega_{x,\Lambda}(t,k)
    = -2\,[\mc{R}(t,k)]^*
    \mc{M}_{x,\Lambda}(t,k)
    \mc{R}(t,k) \label{eq:omega-compressed-v1}
\end{align}
where
\[
\mc{M}_{x,\Lambda}(t,k) =\sum_{j=1}^k \Lambda_{j-1}^* [A^t(x_{j-1})^* \mc{J}_{2d}\partial_tA^t(x_{j-1}) ]\Lambda_{j-1} 
\]
with  $x_j:=x+j\alpha$.
\end{lemma}

\begin{proof}
  Lemma
\ref{lem:Wdot-phi} yields 
$\partial_t W_{A^t_k(x)\Lambda}
= i\,W_{A^t_k(x)\Lambda}\,\Omega_{x,\Lambda}(t,k),$
where the Hermitian matrix $\Omega_{x,\Lambda}(t,k)$ is given by 
\begin{align}
\Omega_{x,\Lambda}(t,k)
    = -2\,[\mc{R}(t,k)]^*
    \Biggl\{
    \Lambda_0^*\bigl[A^t_k(x)\bigr]^* \mc{J}_{2d}\,
    {\partial_t}A^t_k(x)\,\Lambda_0\Biggr\}
    \mc{R}(t,k) \label{eq:omega-expanded-v1}.
\end{align}
then \eqref{eq:omega-compressed-v1} follows directly by \eqref{eq:omega-expanded-v1}, the chain rule and Hermitian-symplecticity of $A^t(\cdot)$. 
By construction, $\phi[A^t_k(x)](\Lambda)$ is a continuous real-valued lift of the circle-valued map $\det\big(W_\Lambda^{-1}W_{A^t_k(x)\Lambda}\big)$ via the covering $\pi(s)=e^{2\pi i s}$, then \eqref{der-phi} follows from  Lemma
\ref{lem:Wdot-phi}.
\end{proof}

As a direct consequence, we have have the following:

\begin{corollary}\label{cor:monophi}
Let \(A^t(x)\in C^1(\I\times\T,\mathrm{HSp}(2d))\) be (positively) monotonic. Then for every Lagrangian subspace \(\Lambda\in\Lag(\C^{2d},\omega_d)\) and every integer \(k\ge1\),
\[
\partial_t\phi\bigl[A^t_k(x)\bigr](\Lambda)\;<\;0.
\]
\end{corollary}

\begin{proof}
Let \(u\in\C^d\) be arbitrary, \(u\neq0\), and set 
$v_{j-1}:=\Lambda_{j-1}(t,x)\,u \in \C^{2d}.$
Since \(\Lambda_{j-1}\) is a Lagrangian frame, \(v_{j-1}\) is an isotropic vector.  Hence
\begin{align*}
u^*\bigl(\Lambda_{j-1}^* [A^t(x_{j-1})]^* \mathcal J_{2d}\,\partial_t A^t(x_{j-1})\,\Lambda_{j-1}\bigr)u
&= v_{j-1}^* \bigl([A^t(x_{j-1})]^* \mathcal J_{2d}\,\partial_t A^t(x_{j-1})\bigr) v_{j-1} \\
&= \Psi_{A^t}\bigl(v_{j-1},x_{j-1}\bigr).
\end{align*}
By hypothesis \(\Psi_{A^t}(v,x)>0\) for every isotropic \(v\), therefore the quadratic form on the left is strictly positive for every nonzero \(u\).  Since \(u\) was arbitrary this shows that each compressed matrix
\(\Lambda_{j-1}^* [A^t(x_{j-1})]^* \mathcal J_{2d}\,\partial_t A^t(x_{j-1})\,\Lambda_{j-1}\)
is positive definite.  Summing these positive-definite matrices yields \(\mc{M}_{x,\Lambda}(t,k) \) positive definite, and consequently \(\Omega_{x,\Lambda}=-2\mc{R}^* \mc{M}_{x,\Lambda}\mc{R}\) is negative definite as required.
Since $\Omega_{x,\Lambda}(t,k)$ is negative definite, its trace is negative. Therefore the result follows directly from Lemma \ref{lem:omega-trace}.
\end{proof}

\subsection{Discrete Hellmann-Feynman Identity and premonotonicity of the long-range cocycle}\label{sec:symp-schro}

Consider the eigenvalue equation $
L_{v_d,\alpha,x} u = E u,$
which gives rise to the long-range (transfer) cocycle \((\alpha,L^{E,d})\) where
\[
L^{E,d}(x)=\frac{1}{\hat v_d}
\begin{pmatrix}  
-\hat{v}_{d-1} & \cdots & -\hat{v}_1 & E - 2\cos (2\pi x) - \hat{v}_0 & -\hat{v}_{-1} & \cdots & -\hat{v}_{-d+1} & -\hat{v}_{-d} \\  
\hat{v}_d &  &  &  &  &  &  &  \\  
& & & & & & &  \\  
& & & & & & &  \\  
& & & \ddots & & & &  \\  
& & & & & & &  \\  
& & & & & & &  \\  
& & & & & & \hat{v}_d &  
\end{pmatrix}.
\]
A direct computation \cite{HP} shows that its iterate satisfy 
\[
(\alpha,L^{E,d})^d = (d\alpha,\,L^{E,d}_d) = (d\alpha,\,A^{E,d}),
\]
where \(A^{E,d}(\cdot)\) can be written  as
\[
A^{E,d}(x)
=
\begin{pmatrix}
C_d^{-1}\big(E-V_d(x)\big) & -C_d^{-1}C_d^\ast\\[4pt]
I_d & 0
\end{pmatrix}.
\]

The long-range cocycle plays a distinguished role because it is \emph{symplectic} and it is \emph{premonotonic}.  To describe the symplectic structure, set
$
S_d \;=\;
\begin{pmatrix}
0 & -C_d^\ast\\[4pt]
C_d & 0
\end{pmatrix},$
which is a nondegenerate skew-Hermitian \(2d\times 2d\) matrix defining a Hermitian symplectic form $\psi_d$ on \(\C^{2d}\).  One checks directly \cite{HP} 
\begin{equation}\label{eq:d-symp-cond}
(L^{E,d}(x))^\ast\,S_d\,L^{E,d}(x)=S_d,\qquad (A^{E,d}(x))^\ast\,S_d\,A^{E,d}(x)=S_d.
\end{equation}
 In particular, \((\alpha,L^{E,d})\) and \((d\alpha,A^{E,d})\) are symplectic cocycles.

Recall that a cocycle \((\alpha,A_t)\) is called \emph{premonotonic} if there exists \(n\ge1\) such that its iterate \((\alpha,A_t)^n\) is monotonic. The following lemma records the premonotonicity of the long-range cocycle.

\begin{lemma}[Discrete Hellmann-Feynman Identity]\label{lem:pre-mono}
For any $\alpha \in \mathbb{R}$, the long-range cocycle $(\alpha, L^{E,d})$ is premonotonic, and the $2d$-step cocycle $(2d\alpha, L_{2d}^{E,d})$ is strictly monotonic. More precisely, for any unit vector $v \in \mathbb{C}^{2d}$ and every $x \in \mathbb{T}$, we have
\[
0 \leq \Psi_{L^{E,d}}(v,x) \leq 1 \quad \text{and} \quad 0 \leq \Psi_{A^{E,d}}(v,x) \leq 1,
\]
Furthermore, for the $2d$-step iterate, we have the exact identity:
\begin{equation}\label{eq:monoA^E}
\Psi_{L^{E,d}_{2d}}(v,x) = \|L_d^{E,d}(x) v\|^2 > 0.
\end{equation}
\end{lemma}

\begin{proof}
A direct computation of the exact matrix derivatives shows that for every $d \ge 1$, the following identities hold:
\begin{equation} \label{derivative}
(L^{E,d})^\ast S_d \,\partial_E L^{E,d} = E_{d,d}, 
\quad \text{and} \quad
(A^{E,d})^\ast S_d \,\partial_E A^{E,d} = 
\begin{pmatrix}
I_d & 0 \\
0 & 0_d
\end{pmatrix},
\end{equation}
where $E_{d,d}$ denotes the $2d \times 2d$ matrix with a $1$ at the $(d,d)$ entry and $0$ elsewhere. 
The first statement  follows immediately.

To prove the $2d$-step identity \eqref{eq:monoA^E}, denote $Y_k = L_k^{E,d}(x) v$,  by the chain rule and the symplectic invariance of the transfer matrices, 
\[
\Psi_{L_{2d}^{E,d}}(v, x) = \sum_{k=0}^{2d-1} Y_k^* \left[ (L^{E,d})^* S_d \partial_E L^{E,d} \right] Y_k=    \sum_{k=0}^{2d-1} Y_k^* E_{d,d}Y_k=\sum_{k=0}^{2d-1} |(Y_k)_d|^2 
\]
Because the transfer matrix $L^{E,d}$ acts as a shift operator on the phase space coordinates, the $d$-th component $(Y_k)_d$ sweeps out the sequential physical coordinates of the trajectory, say $y_k$. In particular, 
it is easy to check $Y_d = L_d^{E,d}(x) v$, is composed exactly of these $2d$ contiguous physical components $(y_{2d-1}, \dots, y_0)^\top$. 
Therefore, $\sum_{k=0}^{2d-1} |(Y_k)_d|^2 = \sum_{k=0}^{2d-1} |y_k|^2 = \|L_d^{E,d}(x) v\|^2>0,$ the result follows.
\end{proof}

\begin{remark}
We refer to Lemma \ref{lem:pre-mono} as a Discrete Hellmann-Feynman Identity due to its natural analogue in quantum mechanics. The classical Hellmann-Feynman theorem
$$\frac{\partial \bar{E}}{\partial \lambda} = \int |\psi(x)|^2 \frac{\partial H}{\partial \lambda} dx$$
 establishes that taking the derivative of a system with respect to a parameter yields the spatial integral of the wave function's probability density ($|\psi(x)|^2$). Our algebraic identity perfectly mirrors this physical principle in a discrete setting: the derivative of the transfer matrix's symplectic action with respect to the energy parameter $E$ evaluates exactly to the $\ell^2$-norm (the discrete spatial density, $\sum |y_k|^2$) of the state vector.

\end{remark}

Note by conjugation with  
\[
\mc{P}_d=\begin{pmatrix} C_d & \\ & I_d \end{pmatrix},
\]  
we can transform \(A^{E,d}(\cdot)\)  into the standard Hermitian-symplectic form:  
\begin{equation} \label{eq:hspstruc}
	\mc{P}_d A^{E,d}(\cdot)\mc{P}_d^{-1}=
	\begin{pmatrix}
		(E-V_d(\cdot))C_d^{-1}& -C_d^*\\
		C_d^{-1}& 0
	\end{pmatrix}
	\in \mathrm{HSp}(2d).
\end{equation}

\section{Gap Labeling and gap estimates}\label{sec:periodicapprox}

In this section, we study gap estimates for the finite-range operator with rational frequency
\begin{equation}\label{eq:long-range-op}
(L_{v_d,\alpha,x}u)_n=\sum_{k=-d}^d\hat{v}_k\,u_{n+k}+2\cos\bigl(2\pi(x+n\alpha)\bigr)u_n,
\end{equation}
which will serve as the foundation for the proof of Theorem~\ref{main-thm-analytic} for trigonometric polynomial potentials.

The starting point  is the recently developed Quantitative Avila's Global Theory \cite{GJYZ}, which connects $H_{v_d,\alpha,x}$ with its dual cocycle $(d\alpha, A^{E,d})=(\alpha, L^{E,d})^d$ as follows. Denote $L(E)=L(\alpha,S_E^{v_d})$ the Lyapunov exponent of $(\alpha,S_E^{v_d})$, where we recall $(\alpha,S_E^{v_d}) $ is the Schr\"odinger cocycle related to $H_{V_d,\alpha,x}.$

The following lemma recalls the structural description and hidden subcriticality of the dual center established in \cite{GJ}, building on the quantitative dual framework of \cite{GJYZ}. In the present argument, this structural input is combined with a global symplectic trivialization preserving monotonicity (Proposition~\ref{prop:6.1}) to obtain quantitative estimates for periodic gaps (Proposition~\ref{prop:periodicGAP}).

\begin{lemma}[{\cite{GJ}}]\label{censub}
Let $\alpha\in\mathbb{R}\setminus\mathbb{Q}$ and $E\in\mathbb{R}$. Assume that $\bar{\omega}(\alpha,S_E^{v_d})=1$. 
Then there exists $\delta(E)>0$ such that both $(d\alpha, A^{E,d}(\cdot+i\epsilon))$ and $(\alpha, L^{E,d}(\cdot+i\epsilon))$ are partially hyperbolic with a two-dimensional center bundle for all $|\epsilon|<\delta(E)$. 
In particular, one can take $\delta(E)=\tfrac{L(E)}{2\pi}$ whenever $L(E)>0$.
Moreover, if $E\in\Sigma_{v_d}(\alpha)$ with $L(E)>0$, then for any $|\epsilon|<\tfrac{L(E)}{2\pi}$, we have
\[
L_d\bigl(\alpha,L^{E,d}(\cdot+i \epsilon)\bigr) = L_d(\alpha,L^{E,d}(\cdot)) = 0.
\]
\end{lemma}

Let $v_d$ be a trigonometric polynomial potential of Type I.
Fix $E_{\mathbf{k}}\in\Sigma^+_{v_d,\alpha}$ with
$N_{v_d,\alpha}(E_{\mathbf{k}})\equiv\mathbf{k}\alpha\mod{\mathbb{Z}}$.
Then $\bar{\omega}(\alpha,S_{E_{\mathbf{k}}}^{v_d})=1$ and $L(E_{\mathbf{k}})>0$.
By continuity of the Lyapunov exponent and stability of $T$-acceleration
\cite{GJY}, there exists a compact neighborhood $\I=[i_-,i_+]$ of
$E_{\mathbf{k}}$ such that $L(E)>0$ and
$\bar{\omega}(\alpha,S_E^{v_d})=1$ for all $E\in\I$.
By Lemma~\ref{censub}, the dual cocycle
$(\alpha,L^{E,d}(\cdot+i\epsilon))$ is therefore partially hyperbolic
with a two-dimensional center for all $E\in\I$ and all $|\epsilon|<h_0$,
where $h_0=h_0(\I):=\min_{E\in\I}L(E)/2\pi$.

{By duality, if $\alpha$ is irrational, the integrated density of states (IDS) of the long-range operator \eqref{eq:long-range-op} (denoted by $N^L_{v_d,\alpha}$) coincides with $N_{v_d,\alpha}$. In what follows, we slightly abuse notation by writing $N_{v_d,\alpha,x}$ for the IDS of \eqref{eq:long-range-op}, even when $\alpha$ is irrational. Without loss of generality, we assume
\[
N_{v_d,\alpha}(i_-)<N_{v_d,\alpha}(E_{\mathbf{k}})<N_{v_d,\alpha}(i_+);
\]
if instead $N_{v_d,\alpha}(i_-)\equiv\mathbf{k}\alpha\mod{\mathbb{Z}}$ or $N_{v_d,\alpha}(i_+)\equiv\mathbf{k}\alpha\mod{\mathbb{Z}}$, then the corresponding gap $G_{\mathbf{k}}(\alpha):=\bigcap_{x\in\mathbb{T}}G_{\mathbf{k}}(\alpha,x)$ is already open. Recall that by the gap labeling theorem for the long-range operator \eqref{eq:long-range-op} \cite[Theorem 1.1]{LW}, the $\mathbf{k}$-th gap $G_{\mathbf{k}}(\alpha,x)$ is defined as the level set
\[
G_{\mathbf{k}}(\alpha,x)=\{E\in\mathbb{R}:\;N_{v_{d},\alpha,x}(E)=\mathbf{k}\alpha \mod{\mathbb{Z}}\},
\]
and we say that $G_{\mathbf{k}}(\alpha,x)$ is \emph{collapsed} if it consists of a single point $E_{\mathbf{k}}$.}

By Lemma \ref{lem:pre-mono}, the cocycle $(\alpha,L^{E,d})$ is premonotonic; specifically, its $2d$-th iterate $(2d\alpha, L^{E,d}_{2d})$ is monotonic. 
Although Proposition~\ref{prop:globalmonoframe} is formulated for the standard symplectic structure \(\omega_d\), the result extends naturally to general symplectic structures (in particular, to $\psi_d$ associated with $(\alpha,L^{E,d})$). 
In light of Remark \ref{rem:premono}, we are therefore justified in applying Proposition \ref{prop:globalmonoframe}. 

To state the adapted version of Proposition~\ref{prop:globalmonoframe} for our current setting, we first introduce a generalized Hermitian-symplectic notation.
\begin{definition}
    Let \(\psi\) and \(\omega\) be Hermitian-symplectic forms on \(\C^{2d}\) and \(\C^{2m}\) with associated matrices \(S_\psi\) and \(S_\omega\), respectively. 
    For any matrix \(M\in \C^{2d\times 2m}\), we define
    \[
        \mathrm{HSp}(2d,2m,\psi,\omega) := \bigl\{ M \in \C^{2d\times 2m} \bigm| M^* S_\psi M = S_\omega \bigr\}.
    \]
    We adopt the following notational conventions:
    \begin{itemize}
        \item If \(d=m\), we write \(\mathrm{HSp}(2d,\psi,\omega)\) for short.
        \item If \(d=m\) and \(\psi=\omega\), we simply write \(\mathrm{HSp}(2d,\psi)\).
    \end{itemize}
\end{definition}

With this notation, we obtain the following:
\begin{proposition}\label{prop:6.1}
Let $\alpha\in\mathbb{R}\setminus\mathbb{Q}$, let $v_d$ be a trigonometric
polynomial potential of Type I, and let $E_{\mathbf{k}}\in\Sigma^+_{v_d,\alpha}$.
With $\I$ chosen as above, there exist constants $\delta',h'>0$, a compact neighborhood $\OO\subset\T$ of $\alpha$ and a transformation
\[
B^{E}(\alpha',x)\in C^\omega_{0,h'}(\OO\times\mathbb{T},\mathrm{HSp}(2d,\psi_d,\omega_d))
\]
which is piecewise $C^\omega$ in $E\in\I$, such that
\begin{equation}\label{block-p}
B^{E}(\alpha',x+\alpha')^{-1}L^{E,d}(x)B^{E}(\alpha',x)
= \widehat H^{E}(\alpha',x)\diamond e^{2\pi i\rho^{E}(\alpha',x)}\widehat C^{E}(\alpha',x).
\end{equation}
Here, the components satisfy the following properties:
\begin{itemize}
    \item $\widehat H^{E}(\alpha',x)\in C^\omega_{\delta',0,h'}(\I\times\OO\times\T,\mathrm{HSp}(2d-2))$;
    \item $\rho^{E}(\alpha',x)\in C^\omega_{0,h'}(\OO\times\mathbb{T},\mathbb{R})$ is piecewise $C^\omega$ in $E\in\I$;
    \item $\widehat C^{E}(\alpha',x)\in C^\omega_{0,h'}(\OO\times\mathbb{T},\mathrm{SL}(2,\mathbb{R}))$ is piecewise $C^\omega$ in $E\in\I$ and is premonotonic. More precisely, the iterate
    \begin{equation}\label{iterate-2}
    C^E(\alpha',x):=\widehat{C}^E_{2d}(\alpha',x)
    \end{equation}
    is piecewise monotonic.
\end{itemize}
Moreover, for every $E\in\Sigma(\alpha)\cap\I$, the cocycle $(\alpha,\widehat C^{E,\alpha}(\cdot))$ is subcritical, meaning that for any $|\epsilon|<h'$,
\[
L\bigl(\alpha,\widehat C^{E,\alpha}(\cdot + i\epsilon)\bigr)=L\bigl(\alpha,\widehat C^{E,\alpha}(\cdot)\bigr)=0.
\]
\end{proposition}

\begin{proof}
The proposition follows essentially from \cite[Proposition 9.1]{WAXZ}; the only change is to replace the local symplectification used there (\cite[Theorem 4.2]{WAXZ}) by our global symplectification (Proposition \ref{prop:globalmonoframe}).

We explain the well-definedness of $\rho^E$.
Denote the center block before the scalar decomposition by
$D^E(\alpha',x)$, which is $1$-periodic, and choose $h'$ strictly smaller than the
analytic width supplied by the construction.
At $(E,\alpha')=(E_{\mathbf k},\alpha)$,
both center Lyapunov exponents vanish for $|y|<h'$ by
Lemma~\ref{censub} and Hermitian-symplecticity. Hence
\[
\int_0^1
\log\bigl|\det D^{E_{\mathbf k}}(\alpha,x+iy)\bigr|\,dx
=0,\qquad |y|<h'.
\]
Thus the determinant has zero winding number.
Taking $\OO$ to be an interval, this remains true
throughout $\I\times\OO$ by continuity.
We may therefore choose a $1$-periodic logarithm and set
\[
\rho^E(\alpha',x)
=\frac{1}{4\pi i}\log\det D^E(\alpha',x)
\]
continuously in $(E,\alpha')$ and piecewise analytically in $E$.
This phase is real on the real torus and $1$-periodic.
By the choice of $h'$, it extends continuously to
$\I\times\OO\times\overline{\T_{h'}}$ and is bounded there.
The required real center is then
$\widehat C^E=e^{-2\pi i\rho^E}D^E$.
\end{proof}

From now on we concentrate on a fixed rational frequency $p/q$ (the best approximation of $\alpha$) and a fixed phase $x$; to simplify notation we drop the dependence on $p/q$ and write \eqref{block-p} as
\begin{equation}\label{conjugation-main}
    B^E(x+p/q)^{-1}L^{E,d}(x)B^E(x)
= \widehat{H}^E(x)\diamond e^{2\pi i\rho^E(x)}\widehat{C}^E(x).
\end{equation}
For simplicity, we just assume the transformation $B^E(x)$ is homotopic to a constant. Otherwise,  replace it by $B^E(x)(I_{d-2}\diamond R_{-nx})$, where $n$ is the topological degree of $B^E$. 
Using properties of the periodic cocycle $(p/q,\widehat C^E)$ we obtain the key result below.

\begin{proposition}\label{prop:periodicGAP}
Under the assumption of \ref{prop:6.1}. Let $p/q$ be the best rational approximation to $\alpha$. There exists $q_*>0$, independent of $x$, such that for every $q>q_*$ the following statements hold.
\begin{enumerate}
  \item {For any fixed $x$,  $$1 - N_{v_{d}, \frac{p}{q},x}(E)= \rho_x(\frac{p}{q},\widehat C^{E}).$$}
  As a consequence, for every $E\in G_{\mathbf{k}}(\tfrac{p}{q},x)$, one has
 $ \rho_x(\frac{p}{q},\widehat C^{E})= \frac{\ell}{q}$
where \(\ell\equiv -\mathbf k p\pmod q\). 
  \item There exists constant $C_1=C_1(\I)$ and some endpoint $E_*\in\partial G_{\mathbf{k}}(\tfrac{p}{q},x)$ we have
  \[
  \bigl|G_{\mathbf{k}}(\tfrac{p}{q},x)\bigr|
  \;\ge\; C_1\,
  \min\big\{\big\|\widehat C^{E_*}_{q}(x)-\operatorname{sign}\big(\operatorname{tr}\widehat C^{E_*}_{q}(x)\big)\,I_2\big\|,1\big\}\; g_{\mathbf{k}}^{-1},
  \]
  where
  \[
  g_{\mathbf{k}}
  =\sup_{E\in G_{\mathbf{k}}(\tfrac{p}{q},x)}
  \sum_{k=0}^{q-1}\big\|\lambda_{\max}W_{E,\tfrac{p}{q}}(\cdot)\|_0\,
  \big\|\widehat C^{E}_{k}\big(x+(q-k)\tfrac{p}{q}\big)\big\|^2
  \]
  and $W_{E,p/q}(\cdot)$ is defined in Lemma \ref{lem:LUbd}.
\end{enumerate}
\end{proposition}

{The proof of Proposition \ref{prop:periodicGAP} proceeds in three main steps. First, a consistent normalization of the phase for the transfer cocycle \(L^{E,d}\) is established, connecting it to spectral counting via the integrated density of states; this yields Corollary \ref{cor:rotation-n-cent}, which proves the gap‑labelling identity (1) by relating the rotation number of the reduced cocycle \(\widehat{C}^{E}\) to the gap label \(\mathbf{k}\). Second, the phase of the center cocycle is shown to capture the same gap information as the original cocycle (Corollary \ref{label-pro}).  Third, the size of the gap is estimated. By the strict monotonicity and continuity of the phase, the jump of the phase across the gap boundaries is bounded below by a quantity involving the deviation of \(\widehat{C}^{E}_{q}\) from \(\pm I\). Combining these  yields the lower estimate (2) on the gap length, completing the proof.}

\subsection{Step I. Branch selection for the projective action.} \label{sec:Branch-selection-1}

We begin by choosing a canonical normalization for the phase $\phi$ of the transfer cocycle $L^{E,d}$. The choice is guided by the triviality of the transfer map at large energies. 

Denote by $\Lambda_H$ the horizontal Lagrangian  $\begin{bsmallmatrix}I_d\\[2pt]0\end{bsmallmatrix}$. Let $U_{\e}$ denotes the small neighborhood of $\Lambda_H$ in $\op{Lag}(\C^{2d},\psi_d)$, which can be written as the form $\Lambda=\begin{bsmallmatrix} I_d\\K \end{bsmallmatrix}$ with $\|K\|<\e$.

\begin{lemma}\label{lem:startpoint}
For large $E$, let $ \Lambda\in U_{O(1/E)} $. Then uniformly in $x$, 
\[
 L^{E,d}(x)U_{O(1/E)}=U_{O(1/E)}.
\]
\end{lemma}

\begin{proof}
    Denote $L^{E,d}=\begin{pmatrix}
        A&B\\C&D
    \end{pmatrix}$ and $R=\begin{pmatrix}
        \mathbf{e}_d& \mathbf{e}_1 &\cdots& \mathbf{e}_{d-1}
    \end{pmatrix}$. A direct computation under the frame choice $\Lambda=\begin{pmatrix}
        R\\
        K R
    \end{pmatrix}$   yields
    
    \[
    \begin{aligned}
        AR+BKR&=\begin{pmatrix}
        (\frac{E}{v_d}+*)&*&*&*&\cdots&*\\
        & 1 &0&0&\cdots&0\\
        & & \ddots&0&\cdots&0\\
        & & & 1&\cdots&0\\
        &&&& \ddots &1
    \end{pmatrix}+\begin{pmatrix}
        *&*&*&*&\cdots&*\\
        & 0&0&0&\cdots&0\\
        & & \ddots&0&\cdots&0\\
        & & & 0&\cdots&0\\
        &&&& \ddots &0
    \end{pmatrix}\\
    CR+DKR&=\begin{pmatrix}
        1&0&0&0&\cdots&0\\
        0&0 &0&0&\cdots&0\\
        0&0 &0&0&\cdots&0\\
        \vdots&\vdots &\vdots& \vdots&\ddots&\vdots\\
        0&0 &0&0&\cdots&0\\
    \end{pmatrix}+\begin{pmatrix}
        0&0&0&0&\cdots&0\\
        *&* &*&*&\cdots&*\\
        *&* &*&*&\cdots&*\\
        \vdots&\vdots &\vdots& \vdots&\ddots&\vdots\\
        *&* &*&*&\cdots&*\\
    \end{pmatrix}
    \end{aligned}
    \]
    If $K=O(1/E)$, then $(CR+DKR)(AR+BKR)^{-1}=O(1/E)$. 
\end{proof}

 Lemma \ref{lem:startpoint} permits us to make a specific consistent choice of \[
\phi[L^{E,d}_k(x)] : \R\times \T\times  \operatorname{Lag}(\mathbb{C}^{2d}, \psi_d) \to \mathbb{R}  
\] 
by taking 
\begin{equation}\label{eq:phase-normalization}
   \lim_{E\to +\infty} \phi[L^{E,d}(x)](
\Lambda_H
)=0. 
\end{equation}
and in the following, we just denote  $\phi_{x}^{E}\equiv\phi[L^{E,d}(x)]$, $\phi_{x,k}^{E}\equiv\phi[L^{E,d}_k(x)]$ using this choice. First 
 Lemma \ref{lem:startpoint} implies:

\begin{lemma}
For any $N\geq 1$, then uniformly for $x\in\T$,
$ \lim_{E\to +\infty} \phi_{x,N}^E( \Lambda_H )=0.$
\end{lemma}

\begin{proof}
By the cocycle rule \eqref{eq:cocycle},
$\phi_{x,N}^E( \Lambda_H )=\sum_{m=0}^{N-1}\phi_{x+m\alpha}^E\big(L_m^{E,d}(x)\Lambda_H\big).$
For sufficiently large $E$, Lemma \ref{lem:startpoint} shows that $L^{E,d}(\cdot)$ sends $\Lambda_H$ into $U_{O(1/E)}$.  
By iteration, $L^{E,d}_m(\cdot)$ also sends $\Lambda_H$ into $U_{O(1/E)}$.   
By \eqref{eq:phase-normalization}, the claim follows.
\end{proof}

We now connect the phase with spectral counting.

\begin{lemma}\label{prop:finite-volume-counting}
With the normalization \eqref{eq:phase-normalization}, for every
$N\geq 2$, $E\in\mathbb R$ and $x\in\mathbb T$, we have
\[
\left|
\phi[A_N^{E,d}(x+d\alpha)](\Lambda_H)
-\left(dN-\#\{j:E_j^{dN}\leq E\}\right)
\right|\leq \frac d2,
\]
where the eigenvalues of $\widehat H_x^N$ are counted with
multiplicity. The shift $x+d\alpha$ corresponds to the restriction
$\widehat H_x^N=\widehat H_x^{[1,N]}$.
\end{lemma}

\begin{proof}
By \cite[Section 4]{LW}, one can express the cocycle action on $\Lambda_H$ as
\[
 A_N^{E,d}(x+d\alpha)\begin{pmatrix} I_d \\[2pt] 0 \end{pmatrix}
 =\begin{pmatrix}
   U(E,N+1,x) \\[2pt]
   U(E,N,x)
  \end{pmatrix},
\]
where the block $U(E,N+1,x)$ satisfies
\[
 \det U(E,N+1,x)=0 
 \;\;\Longleftrightarrow\;\; 
 \det(E-\widehat H^N_x)=0.
\]
In particular, Lemma \ref{lemWlambda} implies that  $E$ is an eigenvalue of $\widehat H^N_x$ if and only if $\dim \ker U(E,N+1,x)\ge 1$.  
Equivalently, $-1$ is an eigenvalue of $W_{A^{E,d}_N(x+d\alpha)\Lambda_H}$, with multiplicity equal to $\dim \ker U(E,N+1,x)$.  

Now, by Lemma~\ref{lem:lag-monotone}, for any $N\ge2$, the function
$E \mapsto \phi[A^{E,d}_N(x+d\alpha)](\Lambda_H)$
is strictly decreasing.  
Indeed, as $E$ decreases, all eigenvalues of $W_{A^{E,d}_N(x+d\alpha)\Lambda_H}$ move anti-clockwise along the unit circle, and
\[
 \det W_{A^{E,d}_N(x+d\alpha)\Lambda_H}
   = \exp\!\bigl(2\pi i\,\phi[A^{E,d}_N(x+d\alpha)](\Lambda_H)\bigr).
\]

These crossings of $-1$ by the individual eigenvalues count the
eigenvalues of $\widehat H_x^N$, with multiplicity. With the
normalization at $E=+\infty$, each lifted eigenphase differs
from its crossing count by at most $1/2$. This also holds at
a crossing, whether or not that crossing is included in the count.
Summing over the $d$ eigenphases gives
\[
\left|
\phi[A_N^{E,d}(x+d\alpha)](\Lambda_H)
-\#\{j:E_j^{dN}>E\}
\right|\leq \frac d2.
\]
Since $\widehat H_x^N$ has $dN$ eigenvalues counted with
multiplicity, the assertion follows.
Likewise, crossings of $+1$ correspond to the eigenvalues of
$\widehat H_x^{N-1}$, with multiplicity, since the lower block
is $U(E,N,x)$.
\end{proof}

Therefore, as \(E\) decreases from the largest to the smallest
eigenvalue of \(\widehat H_x^{N-1}\), the image of
\(E \mapsto \phi[A^{E,d}_N(x+d\alpha)](\Lambda_H)\) contains
\([d,d(N-1)]\): at the first $+1$ crossing the total phase is
at most $d$, while at the last it is at least $d(N-1)$.
Similarly, between the largest and smallest eigenvalues of
\(\widehat H_x^N\), the image contains \([d/2,dN-d/2]\).
Dividing the estimate in Lemma~\ref{prop:finite-volume-counting}
by $N$ and letting $N\to\infty$, we obtain
\begin{equation}\label{thm:idseqrot}
    d(1 - N_{v_{d},\alpha,x}(E)) =  \rho_x(E),
\end{equation}
Here the fixed base shift $x\mapsto x+d\alpha$ does not change
the rotation number. Moreover, for each \(x \in \T\), \(\rho_x(E)\equiv \rho_x(d\alpha,A^{E,d})\) is a decreasing function taking values in \([0,d]\).

We now construct a consistent choice of the phase function \(\phi\) for the reduced cocycle \(\widehat{C}^E_q\) that is compatible with the normalization established for \(L^{E,d}\). The construction utilizes conjugation by the symplectic change-of-coordinates matrix \(B^E(x)\) and an adapted family of Lagrangian frames.

 Define the one-parameter families of Lagrangian subspaces:
\[
\Lambda_y := \begin{bmatrix} \cos\pi y \\[2pt] \sin\pi y \end{bmatrix}, \qquad 
\widetilde{\Lambda}_y := \begin{bmatrix} I_{d-1} \\[2pt] \mathbf{0}_{d-1} \end{bmatrix} \diamond \Lambda_y,
\]
where \(\Lambda_y \in \operatorname{Lag}(\mathbb{C}^{2}, \omega_1)\) and \(\widetilde{\Lambda}_y \in \operatorname{Lag}(\mathbb{C}^{2d}, \omega_d)\).

A crucial observation is that equation \eqref{conjugation-main} implies the decomposition:
\[
B^E(x+p/q)^{-1} L^{E,d}(x) B^E(x) \widetilde{\Lambda}_y = \begin{bmatrix} I_{d-1} \\[2pt] \mathbf{0}_{d-1} \end{bmatrix} \diamond \left( \widehat{C}^E(x) \Lambda_y \right).
\]
This allows us to make a consistent choice of the phase function \(\phi\):
\begin{equation}\label{eq:consist-selection}
\phi\left[B^E(x+p/q)^{-1} L^{E,d}(x) B^E(x)\right](\widetilde{\Lambda}_y) = \phi\left[\widehat{C}^E(x)\right](\Lambda_y),
\end{equation}
meaning the branch selection on the left-hand side uniquely determines the corresponding choice on the right. Note that \[
\phi\left[\widehat{C}^E(x+1)\right](\Lambda_y)=\phi\left[\widehat{C}^E(x)\right](\Lambda_y),
\]
Hence $\widehat C^E:\T\to\mathrm{SL}(2,\R)$ is homotopic to a constant, since its projective phase admits a continuous $1$-periodic lift.
As  $B^E(x)$ is homotopic to  constant, the phase function $\phi[B^E(x)](\Lambda)$ is then well-defined as a continuous function on $\mathbb{T}$.

This leads to the following key relations:

\begin{lemma}\label{cen-phi}
One can make a consistent choice of $\phi[\widehat C^E(x)]$ such that 
the following identities hold:
\begin{align}\label{eq:cocycle-rule-N}
\phi[\widehat{C}^E_{Nq}(x)](\Lambda_y) 
= \phi_{x,Nq}^E(B^E(x) \widetilde{\Lambda}_y) + \eta[B^E(x)](\widetilde{\Lambda}_y, \widetilde{\Lambda}_{y_{Nq}}),
\end{align}
where $\widetilde{\Lambda}_{y_k} = [B^E(x_k)^{-1} L^{E,d}_k(x) B^E(x)] \widetilde{\Lambda}_y.$ In particular, we have the following: 
\begin{equation}\label{eq:rule-2}
\phi[\widehat{C}^E_q(x)](\Lambda_y) = \phi_{x,q}^E (B^E(x) \widetilde{\Lambda}_y) + \eta[B^E(x)](\widetilde{\Lambda}_y, B^E(x)^{-1} L^{E,d}_q(x) B^E(x)\widetilde{\Lambda}_y),
\end{equation}
\begin{equation}\label{eq:rule-20}
\phi[C^E_{q}(x)](\Lambda_y) = \phi_{x,2dq}^E(B^E(x) \widetilde{\Lambda}_y) + \eta[B^E(x)](\widetilde{\Lambda}_y, B^E(x)^{-1} A^E_{2q}(x) B^E(x)\widetilde{\Lambda}_y).
\end{equation}
\end{lemma}

\begin{proof}
Define the iterates for $k\geq 0$:
\[
x_k=x+k\tfrac{p}{q}, \quad\Lambda_{y_k} = [\widehat{C}^E_k(x)] \Lambda_y.
\]
We begin by analyzing the product structure:
\begin{equation*}
\begin{aligned}
W_{\widetilde{\Lambda}_{y_{k-1}}}^{-1} W_{\widetilde{\Lambda}_{y_{k}}}^{}
 = W_{\widetilde{\Lambda}_{y_{k-1}}}^{-1} W_{B^E(x_{k-1}) \widetilde{\Lambda}_{y_{k-1}}} 
\cdot W_{B^E(x_{k-1}) \widetilde{\Lambda}_{y_{k-1}}}^{-1} W_{L^{E,d}(x_{k-1}) B^E(x_{k-1}) \widetilde{\Lambda}_{y_{k-1}}}\\ 
\cdot \left[ W_{\widetilde{\Lambda}_{y_{k}}}^{-1} 
W_{B^E(x_k) \widetilde{\Lambda}_{y_{k}}} \right]^{-1}.
\end{aligned}
\end{equation*}
Then follows from this decomposition and \eqref{eq:consist-selection}, we can select $\phi[\widehat{C}^E(x)]$ satisfies
\begin{equation*}
\phi[\widehat{C}^E(x)](\Lambda_y) = \phi_{x}^E (B^E(x) \widetilde{\Lambda}_y) + \phi[B^E(x)](\widetilde{\Lambda}_y)-\phi[B^E(x_1)]( \widetilde{\Lambda}_{y_1}),
\end{equation*}
and in particular, for each \(k \geq 0\):
\begin{equation}\label{eq:rule-10}
\phi[\widehat{C}^E(x_{k-1})](\Lambda_{y_{k-1}}) = \phi_{x_{k-1}}^E (B^E(x_{k-1}) \widetilde{\Lambda}_{y_{k-1}}) + \phi[B^E(x_{k-1})](\widetilde{\Lambda}_{y_{k-1}})-\phi[B^E(x_k)]( \widetilde{\Lambda}_{y_k}),
\end{equation}

To establish \eqref{eq:cocycle-rule-N}, we employ the cocycle rule.
Summing \eqref{eq:rule-10} over \(k = 1\) to \(Nq\) yields:
\begin{align}
\phi[\widehat{C}^E_{Nq}(x)](\Lambda_y) 
&= \sum_{k=1}^{N} \left( \phi_{x,q}^E([L_{q}^{E,d}(x)]^{k-1} B^E(x)\widetilde{\Lambda}_y) + \eta[B^E(x)](\widetilde{\Lambda}_{y_{(k-1)q}}, \widetilde{\Lambda}_{y_{kq}}) \right) \notag\\
&= \phi_{x,Nq}^E(B^E(x) \widetilde{\Lambda}_y) + \eta[B^E(x)](\widetilde{\Lambda}_y, \widetilde{\Lambda}_{y_{Nq}}).
\end{align}
Noting that \(C^E_q(x) = \widehat{C}^E_{2dq}(x)\), the result follows. \end{proof}

We denote this specific choice by \(\widehat{\phi}[\widehat C^E(x)]\) 
and let $\rho_x(p/q,\widehat C^E) $, \(\rho_x(0,\widehat C^E_q)\) be the corresponding fibred rotation number. 
As a direct consequence, Lemma \ref{cen-phi} allows us to give the gap labelling through the center cocycle, which in particular proves $(1)$ of Proposition \ref{prop:periodicGAP}:

\begin{corollary}\label{cor:rotation-n-cent}We have 
$$1 - N_{v_{d}, \frac{p}{q},x}(E)= \rho_x(\frac{p}{q},\widehat C^{E}).$$
As a consequence,  $E\in G_{\mathbf{k}}(p/q,x)\cap \I$ if and only if
$$ \rho_x(\frac{p}{q},\widehat C^{E})= \frac{\ell}{q},  \quad  \quad \rho_x\bigl(0,\,C^E_q\bigr)=2d\ell   $$
where \(\ell\equiv -\mathbf k p\pmod q\). 
\end{corollary}

\begin{proof}
We only need to demonstrate that $$
\rho_x(E) = \frac{1}{2q} \rho_x(0, C^E_q) = d\rho_x\left(\frac{p}{q}, \widehat{C}^E\right).
$$ The second equality follows directly from the identity $C^E_q(x) = \widehat{C}^E_{2dq}(x)$. Therefore, our primary objective is to establish the first equality, beginning with the definition of the rotation number and utilizing the uniformity of the initial Lagrangian. By equation \eqref{eq:cocycle-rule-N}, 
we obtain:
\begin{align*}
\rho_x(E) = \frac{1}{2}\rho_x(2d p/q,A^{E,d}_2) 
&= \lim_{k'\to\infty}\frac{1}{2k'q}\phi[A^{E,d}_{2k'q}(x)](B^E(x)\widetilde\Lambda_y) \\
&= \lim_{k'\to\infty}\frac{1}{2k'q}\phi[L^{E,d}_{2dk'q}(x)](B^E(x)\widetilde\Lambda_y) \\
&= \lim_{k'\to\infty} \big (\frac{1}{2k'q}\phi[{C}^E_{k'q}(x)](\Lambda_y) - \eta[B^E(x)](\widetilde{\Lambda}_y, \widetilde{\Lambda}_{y_{2dk'q}})\big).
\end{align*}
The  result then follows from Lemma~\ref{lem:minmax}. The rest results just follows from \eqref{thm:idseqrot}. 
\end{proof}

\subsection{Step II. Gap labelling through  projective action.}

In the last step, we choice the projective action for the center cocycle \(C^E\) in a manner that maintains compatibility with the normalization established for \(L^{E,d}\). This compatibility is particularly crucial, as it ensures the inheritance of monotonicity properties:

\begin{lemma}\label{mono-center}
For any $y\in[-\tfrac{1}{2},\tfrac{1}{2})$, the function
$\widehat\phi[C^E_q(x)](\Lambda_y)$ is continuous and strictly decreasing on $\I$.
\end{lemma}

\begin{proof}
Note by Proposition \ref{prop:6.1} and  Corollary \ref{cor:monophi}, the result follows. 
 \end{proof}

Moreover, it allows us to obtain gap labelling through  projective action $\widehat\phi[C^E_q(x)](\Lambda_y)$, this finally paves the way to obtain gap estimates in the final step.

Firstly, we make a graph analysis for projective action $\phi_{x,2dq}^E$.

\begin{lemma}\label{lem:rangeofphi-A}
Let \(\mathbf{J}=[a,b]\subseteq\I\). For any $\varepsilon>0$, there exists $\bar q=\bar q(\e)>0$ such that for all $q>\bar q$, one has for every \(\Lambda\in\operatorname{Lag}(\mathbb{C}^{2d},\omega_d)\), the image of the map \(E\mapsto\phi^E_{x,2dq}(\Lambda)\) on \(\mathbf{J}\) contains the interval
\begin{equation}\label{eq:rangeofphiA}
    \big[(1-N_{v_d,\alpha}(b)+\varepsilon)\,2dq,\; (1-N_{v_d,\alpha}(a)-\varepsilon)\,2dq\big].
\end{equation}
\end{lemma}

\begin{proof}
Fix \(\varepsilon>0\). 
By Weyl's inequality for Hermitian matrices, uniformly in $x$
\[
\max_{j}\big|E_j^{2dq}(v_d,p/q)-E_j^{2dq}(v_d,\alpha)\big|
\le \|L^{2dq}_{v_d,p/q,x}-L^{2dq}_{v_d,\alpha,x}\| \le C\,dq|\alpha-p/q|=:\delta(q),
\]
where \(E_j^{2dq}(v_d,\alpha)\) are the eigenvalues of \(L_{v_d,\alpha,x}^{2dq}\). 

Now choose $\delta_*>0$ and then  $q_*$ sufficiently large such that for all $q \ge q_*$, \(\;\delta(q)<\delta_*.\) Moreover, the following hold:
\begin{enumerate}
  \item[(i)] By continuity of \(N_{v,\alpha}\) in \(E\),
  \(\;|N_{v_d,\alpha}(a\pm\delta_*)-N_{v_d,\alpha}(a)|<\varepsilon/6.\)
  \item[(iii)] By the finite-volume IDS approximation (recall \eqref{eq:idsdef}), uniformly in $x$
  \[
  \left|\frac{\#\{j\mid E_j^{2dq}(v_d,\alpha)\le a\pm\delta_*\}}{2dq}-N_{v_d,\alpha}(a\pm\delta_*)\right|<\varepsilon/6.
  \]

\end{enumerate}

Using  the eigenvalue perturbation estimate above we obtain the set inclusion
\[
\{j:\,E_j^{2dq}(v_d,p/q)\le a\}\subseteq\{j:\,E_j^{2dq}(v_d,\alpha)\le a+\delta_*\},
\]
and therefore
\[
\frac{\#\{j \mid E_j^{2dq}(v_d,p/q) \le a\}}{2dq}
\le \frac{\#\{j \mid E_j^{2dq}(v_d,\alpha) \le a+\delta_*\}}{2dq}.
\]
Applying (iii) and (i) in succession yields
\[
\begin{aligned}
\frac{\#\{j \mid E_j^{2dq}(v_d,p/q) \le a\}}{2dq}
< N_{v_d,\alpha}(a+\delta_*) + \frac{\varepsilon}{6} 
< N_{v_d,\alpha}(a) + \frac{\varepsilon}{3} 
< N_{v_d,\alpha}(a) + \frac{\varepsilon}{2}.
\end{aligned}
\]
Analogously, one can find $q_*'>0$ such that for all $q> q_*'$, 
\[
\begin{aligned}
\frac{\#\{j \mid E_j^{2dq}(v_d,p/q) \le b\}}{2dq}
>N_{v_d,\alpha}(b) - \frac{\varepsilon}{2}.
\end{aligned}
\]

Applying Lemma \ref{prop:finite-volume-counting} with $N=2q$,
choosing \(q> \bar q:=\max\{q_*,q_*',2/\varepsilon\}\), and
translating the base point, which is allowed since the preceding
estimates are uniform in $x$, we obtain
\[
\begin{aligned}
    \phi[A^{a,d}_{2q}(x)](\Lambda_H)
    &> (1-N_{v_d,\alpha}(a))\,2dq - \varepsilon\,dq - \frac d2,\\[4pt]
    \phi[A^{b,d}_{2q}(x)](\Lambda_H)
    &< (1-N_{v_d,\alpha}(b))\,2dq + \varepsilon\,dq + \frac d2.
\end{aligned}
\]

Recall that $A^{E,d}_{2q}=L^{E,d}_{2dq}$, so that
$\phi[A^{E,d}_{2q}(x)]=\phi^E_{x,2dq}$ with our normalization.
By Lemma~\ref{lem:minmax}, passing from the horizontal frame $\Lambda_H$ to an arbitrary $\Lambda \in \operatorname{Lag}(\C^{2d},\omega_d)$ changes each endpoint by at most $\pm d$, which implies that 
the image of the map \(E\mapsto\phi^E_{x,2dq}(\Lambda)\) restricted on $\mathbf{J}$ contains the interval
\[
\big[(1-N_{v_d,\alpha}(b))2dq+\varepsilon\,dq+\tfrac{3d}{2},\;
(1-N_{v_d,\alpha}(a))2dq-\varepsilon\,dq-\tfrac{3d}{2}\big].
\]
Since $q>2/\varepsilon$, this interval contains the simpler subinterval
\begin{equation*}
    \big[(1-N_{v_d,\alpha}(b)+\varepsilon)\,2dq,\; (1-N_{v_d,\alpha}(a)-\varepsilon)\,2dq\big].
\end{equation*}
as desired.
\end{proof}

Lemma \ref{lem:rangeofphi-A} allows us to recover the graph of $\widehat \phi[C^E_q(x)]$ from the graph of $\phi_{x,2dq}^E$. Choose \(\varepsilon > 0\) such that
\begin{equation}\label{eq:selectvarepsilon}
    N_{v_d,\alpha}(a) + 2\varepsilon 
    < N_{v_d,\alpha}(b) - 2\varepsilon,
\end{equation}
then for all \(q > \bar{q}(\e)\), we have
\begin{lemma}\label{lem:Zlabel}
 For any $y\in[-\tfrac{1}{2},\tfrac{1}{2})$ and any $E \in \mathbf{J}$, 
\[
\widehat \phi[C^E_{q}(x)](\Lambda_y) \in \mathbb{Z} 
\quad \text{if and only if} \quad 
\phi_{x,2dq}^E\bigl(B^E(x)\widetilde{\Lambda}_y\bigr) \in \mathbb{Z}.
\] 
More precisely, 
\[
\widehat \phi[C^E_{q}(x)](\Lambda_y)
= \phi_{x,2dq}^E\bigl(B^E(x)\widetilde{\Lambda}_y\bigr)
\quad \text{whenever} \quad 
\phi_{x,2dq}^E\bigl(B^E(x)\widetilde{\Lambda}_y\bigr) \in \mathbb{Z}.
\]
\end{lemma}

\begin{proof}

We distinguish the proof into two cases:

\noindent\textbf{Case 1:} Assume \(\widehat{\phi}[C^E_{q}(x)](\Lambda_y) \in \mathbb{Z}\). Then we show that \(\phi_{x,2dq}^E(B^E(x)\widetilde{\Lambda}_y) \in \mathbb{Z}\).

By equation \eqref{conjugation-main}, we have the decomposition:
\[
B^E(x)^{-1} A^{E,d}_{2q}(x) B^E(x) \widetilde{\Lambda}_y = \begin{bmatrix} I_{d-1} \\[2pt] \mathbf{0}_{d-1} \end{bmatrix} \diamond \left( C_q^E(x) \Lambda_y \right) = \widetilde{\Lambda}_y.
\]
This implies that
\[
\eta\left[B^E(x)\right](\widetilde{\Lambda}_y, B^E(x)^{-1} A^{E,d}_{2q}(x) B^E(x)\widetilde{\Lambda}_y) = 0.
\]
Then by Lemma~\ref{cen-phi}, we obtain
\begin{equation}\label{inte}
\phi_{x,2dq}^E(B^E(x)\widetilde{\Lambda}_y) = \widehat{\phi}[C^E_{q}(x)](\Lambda_y) \in \mathbb{Z}.
\end{equation}

\noindent\textbf{Case 2:} Assume \(\widehat{\phi}[C^E_{q}(x)](\Lambda_y) \notin \mathbb{Z}\). Then we show that \(\phi_{x,2dq}^E(B^E(x)\widetilde{\Lambda}_y) \notin \mathbb{Z}\).

Since \(\widehat{\phi}[C^E_{q}(x)](\Lambda_y) \notin \mathbb{Z}\), there exists an integer \(l\) such that \(\widehat{\phi}[C^E_{q}(x)](\Lambda_y) \in (l, l+1)\). By Lemma \ref{lem:minmax}, we have 
\begin{equation} \label{eq:phiB}
    \sup_{(E,x,y)} \left| \eta\left[B^E(x)\right](\widetilde{\Lambda}_y, B^E(x)^{-1} A^E_{2q}(x) B^E(x)\widetilde{\Lambda}_y) \right| < d.
\end{equation} 
Hence, we obtain the estimate
\[
\left| \widehat{\phi}[C^E_{q}(x)](\Lambda_y) - \phi^E_{x,2dq}(B^E(x)\widetilde{\Lambda}_y) \right| < d.
\]
By Lemma \ref{lem:rangeofphi-A}, the image of the mapping \(E \mapsto \widehat{\phi}[C^E_q(x)](\Lambda_y)\) over the compact interval \(\mathbf{J}\) contains an interval of length
\begin{equation}\label{eq:>4d}
    2dq(N_{v_d,\alpha}(b) - N_{v_d,\alpha}(a) - 2\varepsilon) - 2d > 2dq(N_{v_d,\alpha}(b) - N_{v_d,\alpha}(a) - 3\varepsilon) > 4d.
\end{equation}
Therefore, by continuity and monotonicity (as guaranteed by Lemma \ref{mono-center}), there exist at least one of \(E_L, E_R \in \mathbf{J}\) with \(E_L < E < E_R\) such that
\[
\widehat{\phi}[C^{E_L}_q(x)](\Lambda_y) = l + 1, \qquad
\widehat{\phi}[C^{E_R}_q(x)](\Lambda_y) = l.
\]
Without loss of generality, we may assume that \(E_R\) exists in \(\mathbf{J}\).

We divide the proof into two subcases.

\smallskip
\noindent\textbf{Subcase 2.1:} \(E_L \in [a, b]\).
Fix a pair \(E_L < E < E_R\). For any fixed \(x\), consider the two-parameter family
\[
U(s,t) := W_{B^s(x)\widetilde{\Lambda}_y}^{-1} W_{B^s(x) (I_{2d-2} \diamond C^t_q(x)) \widetilde{\Lambda}_y}, \quad (s,t) \in [E_L, E_R] \times [E_L, E_R],
\]
which is continuous in \((s,t)\). Let \(\arg\det U(s,t)\) denote a continuous lift of the circle-valued function \(\det U(s,t)\). By construction and the identity
\[
W_{B^E(x)\widetilde{\Lambda}_y}^{-1} W_{A^E_{2q}(x) B^E(x)\widetilde{\Lambda}_y} = W_{B^E(x)\widetilde{\Lambda}_y}^{-1} W_{B^E(x) (I_{2d-2} \diamond C^E_{q}(x)) \widetilde{\Lambda}_y},
\]
we may choose
\[
\arg\det U(t,t) = \phi_{x,2dq}^t(B^t(x)\widetilde{\Lambda}_y).
\]
In particular, by equation \eqref{inte} from Case 1, we obtain
\[
\arg\det U(E_R, E_R) = \widehat{\phi}[C^{E_R}_q(x)](\Lambda_y) = l,
\]
and similarly, \(\arg\det U(E_L, E_L) = l+1\).

Since \(\widehat{\phi}[C^{E_R}_q(x)](\Lambda_y) = l\), the cocycle \(C^{E_R}_q(x)\) fixes \(\Lambda_y\), which implies that \(U(s, E_R) = I_{2d}\) for all \(s\). By continuity in \(s\), it follows that
\[
\arg\det U(s, E_R) \equiv l \quad \text{for all } s \in [E_L, E_R].
\]
From the monotonicity in the reduced parameter (i.e., for fixed \(s\), the map \(t \mapsto \arg\det U(s,t)\) is strictly decreasing; see Lemma~\ref{lem:mono2} and Lemma~\ref{mono-center}), we deduce that for every \(s \in [E_L, E_R]\) and every \(t \in (E_L, E_R)\),
\[
l = \arg\det U(s, E_R) < \arg\det U(s, t) < \arg\det U(s, E_L) = l+1.
\]
Taking \(s = t\) yields
\[
l < \arg\det U(t,t) = \phi_{x,2dq}^t(B^t(x)\widetilde{\Lambda}_y) < l+1
\]
for every \(t \in (E_L, E_R)\). In particular, for \(t = E\), we obtain
$\phi_{x,2dq}^E(B^E(x)\widetilde{\Lambda}_y) \in (l, l+1).$

\smallskip
\noindent\textbf{Subcase 2.2:} \(E_L < a\).
In this case, we extend the interval \([a, b]\) to \([a-\delta, b] \subset \I\) such that for all \(E \in [a-\delta, E_R]\), we still have \(\widehat{\phi}[C^E_{q}(x)](\Lambda_y) \in [l, l+1)\). Set \(a_* = \widehat{\phi}[C^{a-\delta}_{q}(x)](\Lambda_y)\).

Define \(\widetilde{C}^E_{q}(x) = R_{f_E} C^E_{q}(x)\), where  \(f_E\) is a smooth decreasing function satisfying:
\begin{enumerate}
    \item For any \(E \geq a\), \(f_E = 0\).
    \item For any \(E < a - \frac{\delta}{2}\), \(f_E = l+1 - a_* + \delta\).
\end{enumerate}
One can verify that \(\widehat{\phi}[\widetilde{C}^E_{q}(x)](\Lambda_y)\) is still strictly decreasing, and there exists \(\widetilde{E}_L \in (a-\delta, a)\) such that \(\widehat{\phi}[\widetilde{C}^{\widetilde{E}_L}_{q}(x)](\Lambda_y) = l+1\), since
$
\lim_{E \to a-\delta} \widehat{\phi}[\widetilde{C}^E_{q}(x)](\Lambda_y) > l+1.$

Now, we replace \(C^t_{q}\) by \(\widetilde{C}^t_{q}\) and let \(\widetilde{C}^t_{q}\) play the role of \(C^t_{q}\) in Subcase 2.1. For any \(s \in [\widetilde{E}_L, E_R]\), a similar continuity argument shows that
\[
\arg \det U(s, E_R) \equiv l, \quad \arg \det U(s, \widetilde{E}_L) \equiv l+1.
\]
Then for any \(\widetilde{E}_L \leq s \leq E_R\) and \(\widetilde{E}_L < t < E_R\), we have
\[
l < \arg \det U(s, t) < l+1.
\]
Since \(\widetilde{C}^E_{q}(x) = C^E_{q}(x)\) for \(a \leq s < E_R\), it follows that
$\phi^E_{x,2dq}(B^E(x)\widetilde{\Lambda}_y) \in (l, l+1),$ thus we finish the whole proof. \qedhere
	\end{proof}
\begin{figure}[htbp]
    \centering
\begin{tikzpicture}[scale=0.9]
    \def\m{1}
    \def\N{6}
    \pgfmathsetmacro{\ymax}{\m*\N}
    \def\a{1.5}
    \def\delta{0.5}
    
    \pgfmathsetmacro{\xmin}{-\a-\delta}
    \pgfmathsetmacro{\xmax}{\a+\delta}
    \def\Amin{0.05}

    \pgfmathsetmacro{\sraw}{ln(2*\ymax - 1)} 
    \pgfmathsetmacro{\s}{(\sraw + 0.8)/\a}

    \pgfmathsetmacro{\Bleft}{\ymax/(1+exp(-\s*\a))}
    \pgfmathsetmacro{\Bright}{\ymax/(1+exp(\s*\a))}

    \pgfmathsetmacro{\marginLeft}{\Bleft - (\ymax - 0.5)}
    \pgfmathsetmacro{\marginRight}{(0.5*\m) - \Bright}
    \pgfmathsetmacro{\marginMin}{min(\marginLeft,\marginRight)}
    \pgfmathsetmacro{\A}{min(\Amin,0.8*max(\marginMin,0.01))}

    \pgfmathdeclarefunction{B}{1}{\pgfmathparse{\ymax/(1+exp(\s*#1))}}
    \pgfmathdeclarefunction{R}{1}{\pgfmathparse{B(#1) + \A * sin(360 * B(#1))}}

    \begin{axis}[
        scale only axis,
        width=15cm,
        height=7cm,
        axis lines=middle,
        xlabel=$E$,
        ylabel={\qquad\qquad{\color{ blue}$\phi_{x,2q}^E(B^E(x)\widetilde \Lambda_y) $ } \quad {\color{red} $\widehat{\phi} [C^E_{q}(x)](\Lambda_y)$}},
        xmin=\xmin,
        xmax=\xmax,
        ymin=-0.5,
        ymax=\ymax+0.6,
        clip=false,
        xtick={-\a, \a},
        xticklabels={$a$, $b$},
        ytick=\empty,
        extra y ticks={\ymax, \ymax-0.5, 0.5*\m},
        extra y tick labels={, $(1-N_{v_d,\alpha}(a)-\varepsilon)2dq$, $(1-N_{v_d,\alpha}(b)+\varepsilon)2dq$},
    ]

    \addplot[thick, blue, domain=\xmin:\xmax, samples=200, smooth] {B(x)};
    \addplot[thick, red, domain=\xmin+0.4:\xmax-0.4, samples=400, smooth] {R(x)};

    \draw[dashed] (axis cs: \xmin, \ymax) -- (axis cs: \xmax, \ymax);
    \draw[dashed] (axis cs: \xmin, \ymax-0.5) -- (axis cs: -\a, \ymax-0.5);
    \draw[dashed] (axis cs: \a, 0.5*\m) -- (axis cs: \xmax, 0.5*\m);
    \draw[dashed] (axis cs: -\a, 0) -- (axis cs: -\a, \ymax-0.5);
    \draw[dashed] (axis cs: \a, 0) -- (axis cs: \a, 0.5*\m);

    \end{axis}

\end{tikzpicture}
    \caption{Relationship between $\phi_{x,2q}^E(B^E(x)\widetilde \Lambda_y)$ and  $\widehat{\phi} [C^E_{q}(x)](\Lambda_y)$.}
    \label{fig:my_curves}
\end{figure}
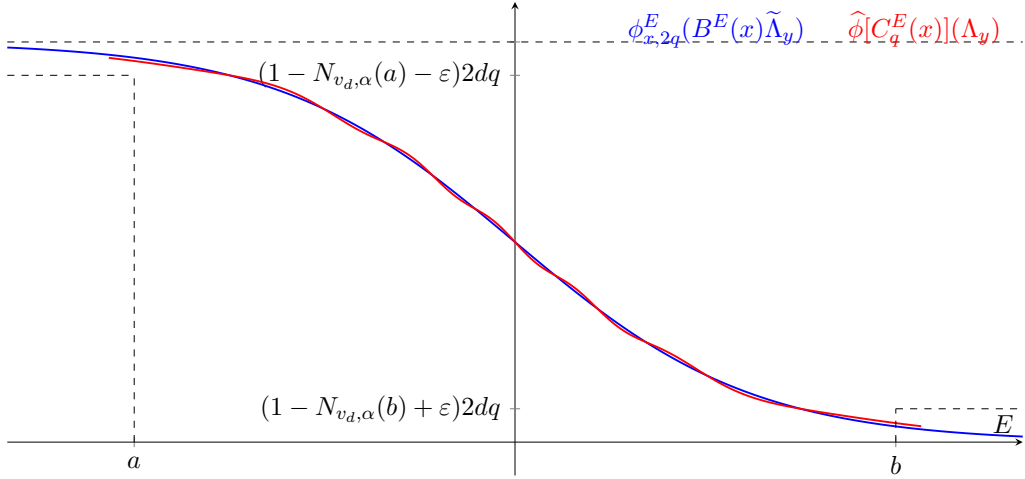

This lemma is crucial, as it allows to obtain gap labelling  through the projecitve action of $\widehat\phi[C^E_q(x)](\Lambda_y)$:

\begin{corollary}\label{label-pro}
For any $x,$ $E\in G_{\mathbf{k}}(p/q,x)\cap \mathbf{J}$ if and only if
there exists $y\in[-1/2,1/2)$, 
$$\widehat\phi\bigl[C^E_q(x)\bigr](\Lambda_y)=\phi^E_{x,2dq}\bigl(B^E(x)\widetilde\Lambda_y\bigr)=2d\ell,$$
where \(\ell\equiv -\mathbf k p\pmod q\).
\end{corollary}
\begin{proof}
  Direct consequence of Corollary  \ref{cor:rotation-n-cent}  and Lemma \ref{lem:Zlabel}.
\end{proof}

Take an interval \(\mathbf J=[a,b]\subset\I\) and choose \(\varepsilon>0\) so small such that
\begin{equation}\label{eq:definedelta-polished}
    N_{v_d,\alpha}(a)+2\varepsilon < N_{v_d,\alpha}(E_{\mathbf k}) < N_{v_d,\alpha}(b)-2\varepsilon.
\end{equation}
By Lemma \ref{lem:rangeofphi-A} and Lemma \ref{lem:Zlabel}, for every \(y\in[-\tfrac12,\tfrac12)\) the image of the map
\(E\mapsto \widehat\phi\bigl[C^E_q(x)\bigr](\Lambda_y)\) on \(\mathbf J\) contains the interval
\[
\Big[(1-N_{v_d,\alpha}(b)+\varepsilon)\,2dq,\; (1-N_{v_d,\alpha}(a)-\varepsilon)\,2dq\Big],
\]
provided $q>\bar q(\e)$. 
In particular, every integer multiple \(2d\ell\) lying in this interval is attained as a value of \(\widehat\phi[C^E_q(x)](\Lambda_y)\) for some energy \(E\in\mathbf J\).
Hence by \eqref{eq:>4d} there exist integers \(\ell\) with \(1\le \ell\le q-1\) and pairwise disjoint subintervals (which might be collapsed) \(G^\ell(p/q,x)\subset\mathbf J\) (ordered from right to left in \(\mathbf J\)) such that for every \(E\in G^\ell(p/q,x)\)
\begin{equation}\label{eq:varepsilonlabel}
    \widehat\phi\bigl[C^E_q(x)\bigr](\Lambda_y) \;=\; 2d\ell 
 \in \Big[(1-N_{v_d,\alpha}(b)+\varepsilon)\,2dq,\; (1-N_{v_d,\alpha}(a)-\varepsilon)\,2dq\Big].
\end{equation}
for some $y\in[-1/2,1/2).$
Write \(\ell_m\) and \(\ell_M\) for the minimal and maximal integers \(\ell\) for which \eqref{eq:varepsilonlabel} holds. By Corollary \ref{label-pro}, the intervals 
$\op{int}G^\ell(p/q,x)$ with \(\ell\equiv -\mathbf k p\pmod q\) are just  spectral gaps $G_{\mathbf k}(\frac{p}{q},x)$, predicted by the Gap Labelling Theorem. The complement of $\bigcup_{\ell=\ell_m}^{\ell_M}
G^\ell(p/q,x)$ with respect to $ \mathbf J $ contains  a union of $\ell_M-\ell_m+2$ open intervals
$\Delta^\ell$ such that the image of $\widehat\phi[C^E_{q}(x)]$ is contained in $(2d\ell,2d(\ell+1))$, $\ell=\ell_m-1,\dots,\ell_M$. We call the closure $\overline{\Delta^l}$ the $\ell$-th band of the spectrum.

\subsection{Step III. Gap estimates}

Thus we only need to study $ \op{int}G^\ell(p/q,x) $. Indeed, we can make a more geometric argument, which further explain why $ \op{int}G^\ell(p/q,x) $ are spectral gaps. 
let $a(E)=\op{tr} {C^{E}_{q}(x)}\in\R  $, and 
$\Delta^\ell_s $ be the open intervals such that the image of $\widehat\phi[C^E_{q}(x)]$ is contained in $(2d\ell+s,2d\ell+s+1)$, $s=0,\dots,2d-1$. For $E \in \Delta^\ell_s$, 
$\Lambda_y \mapsto C^{E}_{q}(x) \Lambda_y$   has no fixed point so $|a(E)|<2$, and
$\widehat\phi[C^E_q(x)](\Lambda_y)$ strictly decreasing then implies that
$a(E)$ monotonic.  Since at the boundary of
$\Delta^\ell_s$, $\Lambda_y \mapsto C^{E}_{q}(x) \Lambda_y$ has fixed
points (which implies $a(E)$ reach $\pm 2$), we conclude that $a|_{\Delta^\ell_s}$ is a diffeomorphism 
to $(-2,2)$. 
The set $\bigcup_{\ell=\ell_m}^{\ell_M} \op{int} G^\ell(p/q,x)$ consists of the set of $E\in \mathbf J$ such that
$\Lambda_y \mapsto C^{E}_{q}(x) \Lambda_y$ has exactly two fixed points, so that
$C_{q}^{E}(x)$ is hyperbolic.

For $\ell_m \leq \ell \leq \ell_M$ and $E \in \partial G^\ell(p/q,x) $,
$\op{tr} C^{E}_{q}(x)=\pm 2$ on $\partial
G^\ell(p/q,x) $, and $G^\ell(p/q,x)$ is collapsed if
$C^{E}_{q}(x)=\pm I_2 $.  In the non-collapsed case, we will give 
 quantitative
estimates of the gaps, before giving its proof, we state a technical lemma, which will also be used in the case $\beta(\alpha)=0$. 

\begin{lemma}\label{lem:esti-deri-C_q} 
For $C^t\in C^1(\I\times\T,\mathrm{HSP}(2))$, the following inequalities hold: 
\begin{align} 
    \av{\partial_t\phi[C^t_k(x)](\Lambda)}&\leq \frac{1}{\pi}\sup\limits_{\substack{\|v\|=1,\\ x\in\T}} \av{\Psi_{C^t}(v,x)}  \cdot\sum_{m=1}^{k}\bigl\|C^t_{m}(x+(k-m)\alpha)^{-1}\bigr\|^2,\label{eq:kstep-derivative}
\end{align}
Moreover, if $ \inf_{\|v\|=1,\,x\in\T} {\Psi_{C^t}(v,x)}>0 $, then
\begin{align}
     -\partial_t\phi[C^t_k(x)](\Lambda)
    &\leq \frac{1}{\pi}\sup_{\substack{\|v\|=1,\\ x\in\T}} {\Psi_{C^t}(v,x)}\cdot\sum_{m=1}^{k}\bigl\|C^t_{m}(x+(k-m)\alpha)^{-1}\bigr\|^{2},\label{eq:onestep-derivative-1}\\ 
    -\partial_t\phi[C^t_k(x)](\Lambda) &\geq \frac{1}{\pi}\inf_{\substack{\|v\|=1,\\ x\in\T}} {\Psi_{C^t}(v,x)} \cdot\sum_{m=1}^{k}\bigl\|C^t_{m}(x+(k-m)\alpha)\bigr\|^{-2}.\label{eq:onestep-derivative}
\end{align}
\end{lemma}
\begin{proof}
Fix an isotropic frame vector $\Lambda\in \C^2$, and write
\(\Lambda_\ell:=C^t_\ell(x)\Lambda\). 
By Lemma \ref{lem:omega-trace}, one obtains the exact identity
\begin{equation}\label{deri-phi-center}
    -\partial_t\phi\bigl[C^t_k(x)\bigr](\Lambda)
= \frac{1}{\pi}\sum_{\ell=0}^{k-1}
\frac{\Lambda_\ell^*\,\bigl[C^t(x_\ell)\bigr]^*J\,\partial_t C^t(x_\ell)\,\Lambda_\ell}
{\|\Lambda_k\|^2},
\qquad x_\ell:=x+\ell\alpha.
\end{equation}
Taking absolute values and using the definition of \(\Psi_{C^t}\) gives
\[
\bigl|\partial_t\phi[C^t_k(x)](\Lambda)\bigr|
\leq \frac{1}{\pi}\sup_{\|v\|=1,x\in\T}\bigl|\Psi_{C^t}(v,x)\bigr|
\sum_{\ell=0}^{k-1}\frac{\|\Lambda_\ell\|^2}{\|\Lambda_k\|^2}.
\]

Now by the cocycle factorization $
\Lambda_k = C^t_{k-\ell}(x+\ell\alpha)\,\Lambda_\ell,$
\[
\|C^t_{k-\ell}(x+\ell\alpha)\|^{-2}
\le
\frac{\|\Lambda_\ell\|^2}{\|\Lambda_k\|^2}
\le
\|C^t_{k-\ell}(x+\ell\alpha)^{-1}\|^2.
\]
Substituting this bound into the sum and reindexing with \(m=k-\ell\) yields the desired bound of \eqref{eq:kstep-derivative}.  Inequalities \eqref{eq:onestep-derivative-1} and \eqref{eq:onestep-derivative} are similar. \qedhere
\end{proof}

Once we have this, now we finish the proof of gap estimates. 

\begin{lemma}\label {lem:gapsestimate}

Let $\ell_m \leq \ell \leq \ell_M$ and let $E_\ell \in \partial G^\ell(p/q,x)$.  Then 
\begin{equation}
|G^\ell(p/q,x)| \geq C_1 \min \{\|{\widehat C^{E_\ell}_{q}(x)}-\op{sign}(\op{tr} {\widehat C^{E_\ell}_{q}(x)}) I_2 \|,1\} g^{-1}_\ell
\end{equation}
for some constant $C_1=C_1(\I,\OO) $ and with
\begin{equation}
g_\ell= \sup_{E \in G^l}
\sum_{k=0}^{q-1} \|\lambda_{\max}W_{E,\tfrac{p}{q}}(\cdot)\|_0\|\widehat C^{E}_k(x+(q-k)\tfrac{p}{q})\|^2.
\end{equation}

\end{lemma}

\begin{proof}
If $G^\ell(p/q,x)$ is collapsed, it means $C_q^E(x)=\pm I_2$, then  Corollary \ref{cor:rotation-n-cent} implies $\widehat C_q^E(x)=\pm I_2$. For non-collapsed case,
assume $E_\ell$ is the left boundary of $G^\ell(p/q,x)$. Then the image
$ \widehat\phi[C^{E_\ell}_q(x)](\Lambda_y)$ is an interval $[2d\ell,2d\ell+\epsilon_\ell]$, then there exists $R(x)\in \rmm{SO}(2,\R)$ such that \[
    R(x)^{-1}\widehat{C}^{E_\ell}_{2dq}(x)R(x)=\pm \begin{pmatrix}
        1& e\\ 0&1
    \end{pmatrix}
    \]
    for some $e<0$. 
    Since $
    {C}^{E_\ell}_{q}(x)=\widehat{C}^{E_\ell}_{2dq}(x)=[\widehat{C}^{E_\ell}_{q}(x)]^{2d}$,  and these matrices commute, it follows that 
    \[
    R(x)^{-1}\widehat{C}^{E_\ell}_{q}(x)R(x)=\pm \begin{pmatrix}
        1& e/2d\\ 0&1
    \end{pmatrix}
    \]
    Then the image $ \widehat\phi[\widehat C^{E_\ell}_q(x)](\Lambda_y)$ is an interval $[\ell,\ell+\widehat\epsilon_\ell]$ where
\begin{equation}
C_2^{-1} \leq
\frac {\widehat \epsilon_\ell} {\min \{\|{\widehat C^{E_\ell}_{q}(x)}-\op{sign}(\op{tr} {\widehat C^{E_\ell}_{q}(x)} ) I_2\|,1\}} \leq C_2
\end{equation}
for some absolute constant $C_2$.

Let $E_\ell'$ denote the other endpoint of $G^\ell(p/q,x)$. By construction of the labelling, at $E_\ell'$ the corresponding projective image attains the value $\ell$ from above, so
$
\sup_{y}\widehat\phi\big[\widehat C^{E_\ell'}_q(x)\big](\Lambda_y)=\ell.
$
Consequently the difference of the projective images at the two endpoints is at least $\widehat\varepsilon_\ell$:
\[
\widehat\phi\big[\widehat C^{E_\ell}_q(x)\big](\Lambda)-\widehat\phi\big[\widehat C^{E_\ell'}_q(x)\big](\Lambda)
\;\ge\; \widehat\varepsilon_\ell .
\]

By  Lemma \ref{mono-center}, the map $E\mapsto \widehat\phi\big[\widehat C^{E}_{q}(x)\big](\Lambda)$ is continuous on $\I$ and piecewise analytic; let $\{i_s\}_{s=1}^m$ be the (finitely many) points of non-differentiability in the relevant interval.
Extremal $s_1$, $s_2$ with $\bigcup\limits_{s=s_1}^{s_2}[i_s,i_{s+1}]\subset G^\ell(p/q,x)  $. 
Then we have \begin{align*}
    \widehat \phi[\widehat C^{E_\ell}_{q}(x)]&(\Lambda)-\widehat \phi[\widehat C^{E_\ell'}_{q}(x)](\Lambda)\\
    &= \widehat \phi[\widehat C^{E_\ell}_{q}(x)](\Lambda)-\widehat \phi[\widehat C^{i_{s_1}}_{q}(x)](\Lambda)+\widehat \phi[\widehat C^{i_{s_2+1}}_{q}(x)](\Lambda)-\widehat \phi[\widehat C^{E_\ell'}_{q}(x)](\Lambda)\\
    &\qquad+\sum_{s=s_1+1}^{s_2} \widehat \phi[\widehat C^{i_s}_{q}(x)](\Lambda)-\widehat \phi[\widehat C^{i_{s+1}}_{q}(x)](\Lambda)  \\
    &= 
    \int_{E_\ell}^{i_{s_1}}\partial_E\widehat \phi[\widehat C^{E}_{q}(x)](\Lambda)dE+
    \int_{i_{s_2+1}}^{E_\ell'}\partial_E\widehat \phi[\widehat C^{E}_{q}(x)](\Lambda)dE\\
    &\qquad+  \sum_{s=s_1+1}^{s_2}
    \int_{i_{s}}^{i_{s+1}}\partial_E\widehat \phi[\widehat C^{E}_{q}(x)](\Lambda)dE\\
    &\leq \sup_{E\in G^\ell(p/q,x)\setminus \{i_s\}_{s=1}^m }\av{\partial_E\widehat \phi[\widehat C^{E}_{q}(x)](\Lambda)} (E_\ell'-E_\ell)
\end{align*}

By Lemma \ref{lem:esti-deri-C_q}, we have \[
\begin{aligned}
     \sup_{y }&\sup_{E\in G^\ell(p/q,x)\setminus \{i_s\}_{s=1}^m}  \av{\partial_E\widehat \phi[\widehat C^{E}_{q}(x)](\Lambda_y)}\\
     &\qquad\qquad \leq \frac{1}{\pi}\sup_{E\in G^\ell(p/q,x)\setminus \{i_s\}_{s=1}^m}  \sup_{\|v\|=1} \av{\Psi_{\widehat C^E}(v,x)}  \sum_{k=0}^{q-1}\|\widehat C^{E}_k(x+(q-k)\tfrac{p}{q})\|^2
\end{aligned}
\]
Although we lost the monotonicity here, but by Lemma \ref{lem:LUbd} and Lemma \ref{lem:pre-mono}, we still have estimate
\begin{equation}\label{eq:LipschitzControl}
    \begin{aligned}
        |\Psi_{\widehat C^E}(u,x)| 
    &< C(\I,\OO)\|\lambda_{\max}W_{E,\tfrac{p}{q}}(\cdot)\|_0 \sup_{\|v\|=1}\av{\Psi_{L^E}(v,x)}
      + C(\I,\OO)\|\widehat C^E(\cdot)\|_0^2 \epsilon
    \end{aligned}
\end{equation}
where $\epsilon$ can be taken arbitrarily small. Therefore
\[
\begin{aligned}
     \sup_{y\in \R/\Z }\sup_{E\in G^\ell(p/q,x)\setminus \{i_s\}_{s=1}^m}  \av{\partial_E\widehat \phi[C^{E}_{q}(x)](\Lambda_y)} & \leq C_4(\I,\OO) \widehat g_{l}
\end{aligned}
\]
for some $C_4(\I,\OO) >0$. \qedhere
\end{proof}

\section{All frequency and Dimension-Free Aubry Duality}\label{sec:aubry}
\label{sec:Dim-free-Aubry}

In this section, we develop a quantitative, all-frequency, dimension-free Aubry duality for Schr\"odinger operators with trigonometric potential, which generalizes the previous duality argument for almost Mathieu operators \cite[Lemma 6.4]{avila2016dry}. 

Given $v\in C^\omega_{\epsilon}(\T,\R)$, so $|\hat{v}_k|\leq e^{-2\pi\epsilon |k|}$. For a vector-valued function \(f: \T \to \C^{2d}\) with spatial components indexed by \(k \in \{-d, \dots, d-1\}\) (ordered as \(f(x) = \big(f(d-1, x), f(d-2, x), \dots, f(0, x), f(-1, x), \dots, f(-d, x)\big)^\top\)), define the weighted analytic norm
\[
\tnorm{f}_\epsilon = \sum_{k=-d}^{d-1} e^{-\xi |k|} \|f(k, \cdot)\|_\epsilon
\]
for some fixed \(\xi \ll 2\pi\epsilon\).  Then we have the following:

\begin{proposition}\label{dim-duality}
Let \(\alpha \in \R \setminus \Q\). For every \(\epsilon > 0\), \(c > 0\) and \(K>0\), there exist \(\delta_* > 0\), \(q_* > 0\), and a subsequence \(p_n/q_n\) of the best approximations of \(\alpha\) with the following property. Let \(p_n/q_n\) satisfy \(q_n > q_*\) and let \(E \in [-3 - \inf |v_d|, 3 + \sup |v_d|]\). Then there exist no \(\rho \in C^\omega_\epsilon(\T, \R)\) with $\|\rho\|_\epsilon\leq K$ and \(U_d \in C^\omega_\epsilon(\T, \mathrm{HSp}(2d, 2, \psi_d, \omega_1))\) such that
\begin{equation} \label{ests-new}
\tnorm{L^{E,d}(\cdot) U_d(\cdot) - U_d(\cdot + \alpha) \left(e^{2\pi i \rho(\cdot)} I_2\right)}_{\epsilon} \leq e^{(-c + \delta_*) q_n}, \quad \tnorm{U_d}_\epsilon \leq e^{\delta_* q_n}.
\end{equation}
Moreover, if \(\beta(\alpha) > 0\) and the subsequence \(p_n/q_n\) satisfies \(q_{n+1} > e^{(\beta - o(1))q_n}\), then the same statement holds upon replacing \(\alpha\) by \(p_n/q_n\) in \eqref{ests-new}.
\end{proposition}

\begin{proof}

In all subsequent arguments, \(o(q_n)\) terms are uniform in the truncation dimension \(d\) and in phases satisfying $\|\rho\|_\epsilon\leq K$, with $K$ fixed. The small exponential losses below can be made arbitrarily small by first choosing $\delta_*$ sufficiently small and then taking $n$ sufficiently large.
We proceed by contradiction. Let \(p/q = p_n/q_n\) and \(E\) satisfy the assumptions, and suppose there exist \(\rho\) and \(U_d\) satisfying \eqref{ests-new}. Let \(u_+^d(\cdot, x)\) and \(u_-^d(\cdot, x)\) denote the first and second columns of \(U_d(x)\), respectively. Then
\[
\tnorm{L^{E,d} u_{\pm}^d - e^{2\pi i\rho} u_{\pm}^d(\cdot + \alpha)}_\epsilon \leq e^{(-c + o(1)) q}.
\]

The first step is an averaging argument via cohomological equation approximation, split into two cases:

\medskip
\noindent\textbf{Case 1:} $\beta(\alpha)=0$. For this case, one can always find $\psi\in C^\omega_{\epsilon/2}(\R/\Z,\R)$ solve the cohomological equation
$\psi(x+ \alpha)-\psi(x) = \rho(x)-\hat{\rho}_0 $ with estimates
\[
\|{\psi}(x) \|_{\epsilon/2} \leq O(1) \|\rho\|_\epsilon.
\]
\noindent\textbf{Case 2:} $\beta(\alpha)>0$. For this case, the cohomological equation cannot be solved completely. 
However, its \emph{dominant part} can still be solved, 
providing a good approximation for further analysis.

The following lemma can be viewed as a refinement of \cite[Lemma~6.1]{GJY}. The improved control is needed to preserve the exponential accuracy of the input errors in Proposition~\ref{dim-duality}.

\begin{lemma}\label{lem:cohomologicaleq}
    For $\alpha\in\R\backslash\Q$ with $\beta(\alpha) > 0$ and $\rho \in C^\omega_\epsilon(\mathbb{T},\mathbb{R})$, there exist sequences $p_n/q_n$ with $q_{n+1}>e^{(\beta-o(1))q_n}$ and $\psi_n \in C^{\omega}_{\epsilon/2}(\mathbb{T},\mathbb{R})$ such that
    \begin{equation}\label{eq:psix1}
        \|{\psi_n}(x) \|_{\epsilon/2} \leq o(q_n) \|\rho\|_\epsilon
    \end{equation}
    \begin{equation}\label{eq:psix2-irra}
        \|\psi_n(x + \alpha) - \psi_n(x) - \rho(x) + \hat{\rho}_0 \|_{\epsilon/2} < e^{-c^* q_n} \|\rho\|_\epsilon
    \end{equation}
    \begin{equation}\label{eq:psix2}
        \|\psi_n(x + p_n/q_n) - \psi_n(x) - \rho(x) + \hat{\rho}_0 \|_{\epsilon/4} < e^{-c_* q_n} \|\rho\|_\epsilon
    \end{equation}
    where $c^*=c^*(\epsilon)>0$, $c_* = c_*(\epsilon, \|\rho\|_\epsilon, \beta) > 0$.
\end{lemma}

\begin{proof}
    We first recall a small divisor estimate:
    \begin{lemma}[\cite{HY}]\label{hy}
        For any $0 < |k| < \frac{q_n}{6}$ with $k \notin \{ l q_{n-1} : l \in \mathbb{Z} \}$, $
        \|k\alpha\|_{\mathbb{T}} \geq \frac{1}{7q_{n-1}}.$    \end{lemma}

 Since $\beta(\alpha) > 0$, there exists a subsequence $\{n_k\}$ such that
    \[
    \left| \alpha - \frac{p_{n_k}}{q_{n_k}} \right| < \exp\left( -(\beta - o(1)) q_{n_k} \right), \quad o(1) \to 0 \text{ as } k \to \infty.
    \]
    For each $k \geq 2$, define $\psi_k \in C^\omega (\mathbb{R}/\mathbb{Z}, \mathbb{R})$ as the solution to 
    \begin{equation}\label{eq:holo}
        \psi_k(x + \alpha) - \psi_k(x) = \sum_{0 < |k| \leq q_{n_k}/6} \hat{\rho}_k e^{2\pi i k x}.
    \end{equation}
    Let $q' = q_{n_{k-1}}$ and $q'' = q_{n_k}$. We consider two cases based on the gap between indices:

    \textbf{Case 1: $n_{k-1} + 1 < n_k$.} 
    Split the summation in \eqref{eq:holo} at $N = q_{n_{k-1}+1}/6$. The coefficient estimates are:
    \[
    |(\widehat{\psi_k})_k| \leq \begin{cases} 
        7q' |\hat{\rho}_k| & 0 < |k| \leq N, \, k \notin \mathbb{Z}q' \\
        2q_{n_{k-1}+1} |\hat{\rho}_k| & 0 < |k| \leq N, \, k \in \mathbb{Z}q' \\
        7q_{n_k-1} |\hat{\rho}_k| & N < |k| \leq q''/6, \, k \notin \mathbb{Z}q_{n_k-1} \\
        2q'' |\hat{\rho}_k| & N < |k| \leq q''/6, \, k \in \mathbb{Z}q_{n_k-1}
    \end{cases}
    \]
    This yields:
    \begin{equation}\label{eq:ps}
        \|\psi_k\|_{\epsilon/2} \leq C(\epsilon) \left(7q' + 2q_{n_{k-1}+1} e^{-\pi\epsilon q'} + 7q_{n_k-1} e^{-\pi\epsilon N} + 2q'' e^{-\pi\epsilon q_{n_k-1}} \right) \|\rho\|_\epsilon \leq o(q'')\|\rho\|_\epsilon
    \end{equation}
    since $q' \ll q_{n_{k-1}+1} \leq q''$ and the exponentials dominate.

    \textbf{Case 2: $n_{k-1} + 1 = n_k$.} 
    Then $q' = q_{n_k-1}$ and:
    \begin{equation}\label{eq:ps-1}
        \|\psi_k\|_{\epsilon/2} \leq C(\epsilon) \left(7q' + 2q'' e^{-\pi\epsilon q'} \right) \|\rho\|_\epsilon \leq o(q'')\|\rho\|_\epsilon
    \end{equation}

    By \eqref{eq:holo}  and Fourier decay:
   \begin{equation}\label{eq:trunc}
        \|\psi_k(\cdot + \alpha) - \psi_k(\cdot) - \rho + \hat{\rho}_0 \|_{\epsilon/2} < e^{-c^* q''} \|\rho\|_\epsilon, \quad c^* = c^*(\epsilon) > 0.
    \end{equation}
    Combining \eqref{eq:ps}, \eqref{eq:ps-1}, \eqref{eq:trunc}, and the mean value theorem:
    \begin{align*}
        &\|\psi_{k}(x + p_{n_k}/q_{n_k}) - \psi_{k}(x) - \rho(x) + \hat{\rho}_0 \|_{\epsilon/4} \\
        &\qquad \leq \|\psi_k(\cdot + \alpha) - \psi_k(\cdot) - \rho + \hat{\rho}_0 \|_{\epsilon/2} + \|\psi_k'\|_{\epsilon/4} \cdot |\alpha - p_{n_k}/q_{n_k}| \\
        &\qquad < \left( e^{-c^*  q_{n_k}} + q_{n_k} e^{-(\beta - o(1)) q_{n_k}} \right) \|\rho\|_\epsilon \\
        &\qquad < e^{-c_* q_{n_k}} \|\rho\|_\epsilon
    \end{align*}
    for sufficiently large $k$, where $c_* = \min(c^*(\epsilon), \beta/2) > 0$. Relabeling the subsequence as $p_n/q_n$ completes the proof.
\end{proof}
\begin{remark}
For a bounded family of phases $\|\rho\|_{\epsilon}\leq K$ with $K$ fixed, the sublinear bound on $\psi_n$ gives $\|e^{\pm2\pi i\psi_n}\|_{\epsilon/2}\leq e^{o(q_n)}$. Thus the averaging transformation preserves exponentially small input errors at the same scale $q_n$, for both the irrational shift and its rational approximant.
\end{remark}

The argument for the irrational shift \(\alpha\) is analogous to the rational case below. We focus on the case \(\beta(\alpha) > 0\) and prove the statement for the rational shift \(p/q\) (the extension to \(\alpha\) follows identically).

Let \(\psi \in C^\omega_{\epsilon/2}(\T, \R)\) be the approximating solution from Lemma \ref{lem:cohomologicaleq} for \(p/q\). Define the gauge-transformed center bundle basis 
$\bar{U}_d(x) = e^{2\pi i \psi(x)} U_d(x),$ 
with columns \(\bar{u}_{\pm}^d(k, x) = e^{2\pi i \psi(x)} u_\pm^d(k, x)\). Lemma \ref{lem:cohomologicaleq} imply that 
\[
\tnorm{\bar{u}_\pm^d}_{\epsilon/2} \leq e^{o(q) \|\rho\|_\epsilon} \tnorm{u_\pm^d}_{\epsilon/2} = e^{o(q)},
\]
\[
\tnorm{L^{E,d} \bar{u}_\pm^d - e^{2\pi i \hat{\rho}_0} \bar{u}_\pm^d(\cdot + p/q)}_{\epsilon/4} \leq e^{(-c + o(1)) q} + e^{-c_* q} e^{o(q)} \leq e^{(-c' + o(1)) q},
\]
where \(c' = \min\{c, c_*\}\).

By  the structure of $L^{E,d}$ and  the definition of the weighted norm, this gives component-wise error bounds:
    \[
    \|\bar{u}_\pm^d(k+1, x) - e^{2\pi i \hat{\rho}_0} \bar{u}_\pm^d(k, x + p/q)\|_{\epsilon/4} \leq e^{\xi |k|} e^{(-c' + o(1)) q}, \qquad -d\leq k\leq d-2,
    \]
 while for the top component \(k = d-1\), we have
 \begin{eqnarray}
 \label{top}
  &&  \left\|(E - 2\cos x)\bar{u}_\pm^d(0, x) - \sum_{j=-d}^{d-1} \hat{v}_j \bar{u}_\pm^d(j, x) - \hat{v}_d e^{2\pi i \hat{\rho}_0} \bar{u}_\pm^d(d-1, x + p/q)\right\|_{\epsilon/4}\\ 
  \nonumber & \leq& 
    | \hat{v}_d | e^{\xi (d-1)} e^{(-c' + o(1)) q}.
    \end{eqnarray}

Using a telescoping sum on the shift error bounds, we relate all components of \(\bar{u}_\pm^d\) to the 0-th component: for all \(m \in \{-d, \dots, d-1\}\),
\begin{equation}\label{eq:component-rel}
\|\bar{u}_\pm^d(m, x) - e^{2\pi i \hat{\rho}_0 m} \bar{u}_\pm^d(0, x + m p/q)\|_{\epsilon/4} \leq \left( \sum_{j=0}^{|m|-1} e^{\xi |j|} \right) e^{(-c' + o(1)) q}.
\end{equation}
Substituting the component relation \eqref{eq:component-rel} into \eqref{top}, we derive the reduced equation for the 0-th component:
\begin{align*}
 \Big\| (E - 2\cos 2\pi x) \bar u_{\pm}^0(x) &- \sum_{k=-d}^{d} \hat{v}_k e^{2\pi i \hat \rho_0 k} \bar u_{\pm}^0(x+k p/q) \Big\|_{\epsilon/4} \\
&\le \Bigg( \sum_{k=-d}^{d} |\hat{v}_k| \Big( \sum_{i=0}^{|k|-1} e^{\xi i} \Big) \Bigg) e^{(-c'+o(1))q} < M e^{(-c' + o(1)) q} 
\end{align*}
where \(M\) is a constant independent of \(d\) (using the analytic decay of \(v\) and \(\xi \ll 2\pi\epsilon\)).

We now prove a uniform lower bound on the \(L^2\) norm of the 0-th component:

\begin{lemma}\label{lem:lowerb-u0-new}
$\|u_+^0\|_{L^2(\T)} \geq e^{-o(q)}$, uniformly in the truncation $d$.
\end{lemma}

\begin{proof}
By the symplectic normalization $U_d^*S_dU_d=\mathcal J_2$,
\[
1 = -\psi_d(u_+, u_-) = u_+^{\text{top},*} C_d^* u_-^{\text{bot}} - u_+^{\text{bot},*} C_d u_-^{\text{top}},
\]
where $u^{\text{top}} = (u^{d-1}, u^{d-2}, \dots, u^0)^\top \in \C^d$ and $u^{\text{bot}} = (u^{-1}, u^{-2}, \dots, u^{-d})^\top \in \C^d$. Using the Toeplitz structure of $C_d$ and reindexing sums with $m = d-(j-i)$ for $1 \leq m \leq d$, this simplifies exactly to the truncated Wronskian:
\[
1 = \sum_{m=1}^d \sum_{l=0}^{m-1} \left(\overline{\hat{v}_m }  \overline{u_+^l} u_-^{l-m} - \hat{v}_m u_-^l\overline{u_+^{l-m}}  \right),
\]
which further implies that
\[
1 \leq \sum_{m=1}^d |\hat{v}_m| \sum_{l=0}^{m-1} \left( \left| \int_\T \overline{u_+^l} u_-^{l-m} dx \right| + \left| \int_\T \overline{u_-^l} u_+^{l-m} dx \right| \right).
\]
Apply the Cauchy-Schwarz inequality to each integral, using the Fourier decay  of $\hat{v}_m$ and the weighted norm bound $\|u_-^k\|_{L^2} \leq \|u_-^k\|_{L^\infty} \leq e^{\xi |k| + o(q)}$,  we have:
\[
1 \leq e^{o(q)} \sum_{m=1}^d e^{-2\pi\epsilon m} \sum_{l=0}^{m-1} \left( \|u_+^l\|_{L^2} e^{\xi |l-m|} + e^{\xi l} \|u_+^{l-m}\|_{L^2} \right).
\]
Define $\delta = 2\pi\epsilon - \xi > 0$,  reindexing the double sum and using the geometric series bound $\sum_{m > |k|} e^{-\delta m} \lesssim e^{-\delta |k|}$, we obtain 
\begin{equation}\label{eq:weighted-lb}
\sum_{k=-d}^{d-1} e^{-\delta |k|} \|u_+^k\|_{L^2} \geq e^{-o(q)}.
\end{equation}

On the other hand, \eqref{eq:component-rel} implies that \[
\sum_{i=0}^{|k|} c_i^{-1} = \sum_{i=0}^{|k|} e^{\xi i} \lesssim e^{\xi |k|},
\]
since $\xi > 0$ is fixed. Substituting back gives the sharp, uniform error bound:
\begin{equation}\label{eq:component-error}
\left| \|u_+^k\|_{L^2} - \|u_+^0\|_{L^2} \right| \lesssim e^{\xi |k|} e^{(-c' + o(1))q} \quad \forall k \in \{-d, \dots, d-1\}.
\end{equation}

Consequently, we have estimate \[
\left| \sum_{k=-d}^{d-1} e^{-\delta |k|} \|u_+^k\|_{L^2} - \|u_+^0\|_{L^2} \sum_{k=-d}^{d-1} e^{-\delta |k|} \right| \leq  \sum_{k=-d}^{d} e^{-\delta |k|} \left| \|u_+^k\|_{L^2} - \|u_+^0\|_{L^2} \right| \leq e^{(-c'' + o(1))q},
\]
for fixed $c'' = \min(c', \delta - \xi) > 0$.
Substitute the lower bound \eqref{eq:weighted-lb}:
\[
e^{-o(q)} \leq C_\delta \|u_+^0\|_{L^2} + e^{(-c'' + o(1))q},
\]
the desired result follows.
\end{proof}

Consider the Fourier series of the 0-th component: \(\bar{u}_\pm^d(0, x) = \sum_{k \in \Z} \widehat{(\bar{u}_\pm^d(0))}_k e^{2\pi i k x}\). From the reduced spectral equation, we get the Fourier space recurrence:
\[
\left|(E - v_d(kp/q + \hat{\rho}_0)) \widehat{(\bar{u}_\pm^d(0))}_k - \left( \widehat{(\bar{u}_\pm^d(0))}_{k+1} + \widehat{(\bar{u}_\pm^d(0))}_{k-1} \right) \right| \leq e^{(-c' + o(1)) q} e^{-\pi \epsilon |k|/2},
\]
where \(v_d(x) = \sum_{m=-d}^d \hat{v}_m e^{2\pi i m x}\). Define the Fourier coefficient vectors:
\[
\U_k^\pm = \begin{pmatrix} \widehat{(\bar{u}_\pm^d(0))}_k \\ \widehat{(\bar{u}_\pm^d(0))}_{k-1} \end{pmatrix}, \quad \U_k = \begin{pmatrix} \U_k^+ & \U_k^- \end{pmatrix} \in \C^{2 \times 2}.
\]
Then the recurrence becomes:
\begin{align}
\|\U_k^\pm\| &\leq e^{o(q)} e^{-\pi \epsilon |k|/2}, \label{dec} \\
\left\|\A(kp/q + \hat{\rho}_0) \U_k^\pm - \U_{k+1}^\pm\right\| &\leq e^{(-c' + o(1)) q} e^{-\pi \epsilon |k|/2}, \label{jump}
\end{align}
where the transfer matrix is \(\A(x) = \begin{pmatrix} E - v_d(x) & -1 \\ 1 & 0 \end{pmatrix}\).

By Parseval's identity and Lemma \ref{lem:lowerb-u0-new}:
\[
\sum_{k \in \Z} \|\U_k^\pm\|^2 = 2 \int_{\T} |\bar{u}_\pm^d(0, x)|^2 dx = 2 \int_{\T} |u_\pm^d(0, x)|^2 dx \geq e^{-o(q)}.
\]
Thus there exist \(k_\pm\) with \(|k_\pm| \leq o(q)\) such that \(\|\U_{k_\pm}^\pm\| \geq e^{-o(q)}\). Let \(C_0\) be a constant satisfying \(\ln \|\A\|_0 \leq C_0\). Using the recurrence \eqref{jump}, we propagate the lower bound: for any small \(\delta_0 > 0\),
\begin{equation} \label{U}
\|\U_k^\pm\| \geq e^{(-C_0 \delta_0 - o(1)) q}, \quad |k| \leq (\delta_0 - o(1)) q.
\end{equation}

From \eqref{dec} and \eqref{jump}, we estimate the difference of determinants:
\begin{equation} \label{jump1}
|\det \U_k - \det \U_{k+1}| \leq e^{(-c' + o(1)) q} e^{-\pi \epsilon |k|}.
\end{equation}
Summing \eqref{jump1} over all \(k \in \Z\) gives the uniform bound:
\begin{equation} \label{U1}
|\det \U_k| \leq e^{(-c' + o(1)) q} \quad \forall k.
\end{equation}

For \(|k| \leq (\delta_0 - o(1)) q\), decompose \(\U_k^+ = \gamma_k \U_k^- + \V_k\) with \(\langle \V_k, \U_k^- \rangle = 0\). From the lower bound \eqref{U} and determinant bound \eqref{U1}:
\begin{align*}
e^{(-C_0 \delta_0 - o(1)) q} &\leq |\gamma_k| \leq e^{(C_0 \delta_0 + o(1)) q}, \\
\|\V_k\| &\leq e^{(-c' + C_0 \delta_0 + o(1)) q}.
\end{align*}
Using the recurrence \eqref{jump}, we find \(|\gamma_{k+1} - \gamma_k| \leq e^{(-c' + C_0 \delta_0 + o(1)) q}\) for \(|k| \leq (\delta_0 - o(1)) q\). Summing this difference gives:
\[
\|\U_k^+ - \gamma_0 \U_k^-\| \leq e^{(-c' + C_0 \delta_0 + o(1)) q}, \quad |k| \leq (\delta_0 - o(1)) q.
\]
By Fourier inversion, this implies the spatial bound:
\[
\|\bar{u}_+^d(0, \cdot) - \gamma_0 \bar{u}_-^d(0, \cdot)\|_0 = \|u_+^d(0, \cdot) - \gamma_0 u_-^d(0, \cdot)\|_0 \leq e^{(-\delta_1 + o(1)) q}
\]
for some \(\delta_1 > 0\), choosing \(\delta_0\) sufficiently small.

Finally, using the symplectic form invariance:
\[
-1 = \psi_d(u_+^d(\cdot, x), u_-^d(\cdot, x)) = \psi_d(u_+^d(\cdot, x) - \gamma_0 u_-^d(\cdot, x), u_-^d(\cdot, x)).
\]
Repeating the argument from Lemma \ref{lem:lowerb-u0-new} with \(u_+^d\) replaced by \(u_+^d - \gamma_0 u_-^d\), we conclude that:
\[
\|u_+^d(0, \cdot) - \gamma_0 u_-^d(0, \cdot)\|_0 \geq e^{-o(q)},
\]
which contradicts the previous upper bound. This completes the proof.
\end{proof}

\section{Proof of Theorem \ref{main-thm-analytic} (trigonometric polynomial potentials)}\label{sec:endpf}

As a preliminary attempt,  we finish the proof of Theorem \ref{main-thm-analytic} for trigonometric polynomial potentials, and distinguish the proof into two cases: 

\subsection{Case 1:  $\beta(\alpha)>0$.}

In this case, the main idea is to use periodic approximation. Take $|\alpha - p/q| = o(1)$. Recall \[
G_{\mathbf{k}}(\alpha):=\bigcap_{x\in\mathbb{T}}G_{\mathbf{k}}(\alpha,x)=\bigcap_{x\in\mathbb{T}}\{E\in\mathbb{R}:\;N_{v_d,\alpha,x}(E)=\mathbf{k}\alpha \mod{\mathbb{Z}}\},
\]  The following $1/2$-H\"older continuity of joint-gaps for the long-range operator is a key:
\begin{lemma}\label{cor:rational-ap}
If $|G_{\mathbf{k}}(\alpha)| > 0$, then for $|\alpha - \alpha'|$ sufficiently small (depending only on $G_{\mathbf{k}}(\alpha)$ and $v$), we have $|G_{\mathbf{k}}(\alpha')| > 0$ and
\[
|G_{\mathbf{k}}(\alpha) - G_{\mathbf{k}}(\alpha')| \leq C(v_d) |\alpha - \alpha'|^{1/2}.
\]
\end{lemma}
\begin{remark}
    In the context of the Schr\"odinger case, this result was established by \cite{AvS1991measure}. We assume that Lemma \ref{cor:rational-ap} is known to the community; however, due to the absence of an exact reference, we provide a proof in Appendix \ref{12holder} for the sake of completeness.
\end{remark}

Now, by Lemma \ref{cor:rational-ap}, it suffices to demonstrate that $ G_{\mathbf{k}}(p/q) $ admits a subexponential estimate for sufficiently large $ q $. 
Recall from \eqref{eq:definedelta-polished} that we choose $ \varepsilon > 0 $ such that
$$
1 - N_{v_d,\alpha}(E_{\mathbf{k}}) \in (1 - N_{v_d,\alpha}(b) + 2\varepsilon,\, 1 - N_{v_d,\alpha}(a) - 2\varepsilon).
$$
Then, there exists $ q > q_*(\varepsilon) > 0 $ such that $ |\alpha - p/q| = o(1) $ and
$$
-\mathbf{k}p/q \pmod{\mathbb{Z}} = \mathbf{l}/q \in (1 - N_{v_d,\alpha}(b) + \varepsilon,\, 1 - N_{v_d,\alpha}(a) - \varepsilon).
$$
Consequently, the $ \mathbf{l} $-th gap $ G^\mathbf{l}(p/q,x) = G_{\mathbf{k}}(p/q,x) $ satisfies $ |\mathbf{k}| = o(q) $. Such a gap (possibly collapsed) falls within the scope of Section \ref{sec:periodicapprox}.

Let $\I$ be chosen as in the beginning of Section \ref{sec:periodicapprox}.
By Proposition \ref{prop:6.1}, we have
$B^{E}(\alpha', x) \in C^\omega_{0,h'}(\OO \times \mathbb{T}, \mathrm{HSp}(2d))$, which gives for any $|\alpha - p/q| = o(1)$,
\begin{equation}\label{conj2}
B^E(p/q, x + p/q)^{-1} L^{E,d}(x) B^E(p/q, x) = \widehat{H}^E(p/q, x) \diamond e^{2\pi i \rho^E(p/q, x)} \widehat{C}^E(p/q, x),
\end{equation}
with the estimate for any $0<\epsilon_0\leq h'$
\begin{equation}\label{esti-b}
\|B^E(p/q, x)\|_{\epsilon_0} \leq 2 \|B^E(\alpha, x)\|_{\epsilon_0},
\end{equation}
Moreover, by Proposition \ref{prop:6.1}, for any $E \in \Sigma(\alpha) \cap \I$, $(\alpha, \widehat{C}^E(\alpha, x))$ is subcritical, then we  have 
\begin{equation}
\left\| \widehat{C}^E_q(\alpha, \cdot) \right\|_{\epsilon_0} \leq e^{o(q)}
\end{equation}
uniformly over $E$ in $\Sigma(\alpha) \cap \I$.

By the joint continuity and the compactness of $\I\times \overline{\T_{\epsilon_0}}$, one can approximate $\widehat{C}^E(\alpha, x)$ via $\widehat{C}^E(p/q, x)$ uniformly in $(E,x)$. Then a subadditive argument (see, for instance, Lemma 3.1 of \cite{AYZ}) shows that
\begin{equation}
\left\| \widehat{C}^E(p/q, \cdot + (q-1)\alpha) \cdots \widehat{C}^E(p/q, \cdot
+ \alpha) \widehat{C}^E(p/q, \cdot) \right\|_{\epsilon_0} \leq e^{o(q)}
\end{equation}
uniformly over $E$ in $\Sigma(\alpha) \cap \I$, if $|\alpha - p/q| = o(1)$.
Continuity of the spectrum and general upper semicontinuity (again, Proposition 3.1 of \cite{AYZ}) then gives
\begin{equation}\label{subex1}
\max_{0 \leq k \leq q} \left\| \widehat{C}_k^E(p/q, \cdot) \right\|_{\epsilon_0} \leq e^{o(q)}
\end{equation}
if $|\alpha - p/q| = o(1)$ and $E$ is $o(1)$-close to $\Sigma(p/q)\cap \I$. 

We first show
\begin{lemma}\label{lem:gap-esti}
For  every  $x\in\T$,   $|G_{\mathbf k}(p/q,x)| \geq e^{-o(q)}$.
\end{lemma}
\begin{proof}
Fix $x\in\T$. By Proposition  \ref{prop:periodicGAP},
a key goal is to show that if  $E_* \in \partial G^{\mathbf{l}}({p/q},x)$, then
\begin{equation}\label{upp}
\left\| \widehat{C}_q^{E_*}(p/q, x) - \pm I_2 \right\| \geq e^{-o(q)},
\end{equation}
where $\widehat{C}_k^E(p/q, x) = \widehat{C}^E(p/q, x + (k-1)p/q) \cdots \widehat{C}^E(p/q, x)$ denotes the iteration under frequency $p/q$.

We  now apply the following quantitative version of the solution to the almost reducibility conjecture in the case of Liouvillean frequencies \cite{AviAlmost, avila2016dry}:

\begin{theorem}[\cite{avila2016dry}] \label {exp}

For every $0<\epsilon<\epsilon_0$, $c_0>0$ and $\delta>0$,
there exist $\delta_*>0$ and $q_*>0$ with the following property.
Let $p/q \in \Q$ and
$A \in C^\omega_{\epsilon_0}(\R/\Z,\mathrm{SL}(2,\R))$ be such that $q>q_*$ and
\begin{equation} \label {sizek}
\sup_{0 \leq k \leq q} \frac {1} {q} \ln \|A_k\|_{\epsilon_0} \leq \delta_*,
\end{equation}
where $A_k(x)=A(x+(k-1) p/q) \cdots A(x)$.  If $\|A_q(0)-I\| \leq e^{-c_0
q}$ with $I \in \{I_2,-I_2\}$ then there exists $\bar B \in
C^\omega_\epsilon(\R/\Z,\rmm{SL}(2,\R))$ and $R \in \rmm{SO}(2,\R)$
such that $\|\bar B\|_\epsilon \leq e^{\delta q}$, $R^q=I$, and $$\tilde
A(x)=\bar B(x+p/q)A(x)\bar B(x)^{-1}$$ satisfies
$\|\tilde A-R\|_0<e^{(-c_0+\delta) q}$ and $\|\tilde
A-R\|_\epsilon<e^{(-\gamma c_0+\delta) q}$ with $\gamma=1-\frac {\epsilon}
{\epsilon_0}$.  Moreover, $\bar B$ may be chosen to be homotopic to a constant.

\end{theorem}

Fix $0<\epsilon<\epsilon_0$ and $c_0 > 0$ small. By \eqref{subex1} and Theorem \ref{exp}, if 
$
\left\| \widehat{C}_q^{E_*}(p/q, x) - \pm I_2 \right\| \leq e^{-c_0 q},$
then there exists $\bar{B} \in C^\omega_\epsilon(\mathbb{R}/\mathbb{Z}, \mathrm{SL}(2, \mathbb{R}))$ homotopic to a constant such that $\|\bar{B}\|_\epsilon \leq e^{o(q)}$ and
\[
\tilde{A}(x) = \bar{B}(x + p/q) \widehat{C}^{E_*}(p/q, x) \bar{B}(x)^{-1}
\]
satisfies $\|\tilde{A} - R\|_\epsilon \leq e^{(-\gamma c_0 + o(1)) q}$ for some $\gamma > 0$ and $R^q = \pm I_2$. By calculating the fibred rotation number, the matrix $R$ must in fact be $R_{{\mathbf l}/2q}$ or $R_{({\mathbf l}-q)/2q}$.
Now define on $2\T$
\[
\bar{B}^{E_*}(p/q, x) = B^{E_*}(p/q, x) \left( I_{2d-2} \diamond \bar{B}(x)^{-1} R_{-\mathbf k x/2} \right).
\]
Since $\mathbf l\equiv-\mathbf k p\pmod q$, \eqref{conj2} implies that, for one fixed choice of sign,
\[
\left\|  L^{E_*,d}(\cdot) \left.\bar{B}^{E_*}(p/q, \cdot)\right|_{E^c_{E_*}} - \left.\bar{B}^{E_*}(p/q, \cdot + p/q)\right|_{E^c_{E_*}}
\pa{\pm e^{2\pi i \rho^{E_*}(p/q, \cdot)} I_2 }\right\|_{\epsilon} \leq e^{(-\gamma c_0 + o(1)) q},
\]
where the restriction denotes the two center columns of the frame.
Moreover, $|\mathbf{k}| = o(q)$ and \eqref{esti-b} imply that $\|\bar{B}^{E_*}(p/q, \cdot)\|_{\epsilon} \leq e^{o(q)}$. Absorbing the constant sign into the scalar phase, we apply Proposition \ref{dim-duality} on $2\T$ to obtain a contradiction, proving \eqref{upp}.

Now if $|G^{\mathbf l}(p/q)| = o(1)$, then for $E \in G_{\mathbf k}(p/q,x)$, we have \eqref{subex1}, then
\eqref{esti-b} implies that 
\[
g_{\mathbf k} \leq \sup_{E \in G_{\mathbf k}(p/q,x)} \|B^E(p/q, \cdot)|_{E^c_E}\|_0^2 \cdot q \cdot e^{o(q)} \leq C(\I, \OO) \cdot q \cdot e^{o(q)} \leq e^{o(q)}.
\]
By Proposition  \ref{prop:periodicGAP} and \eqref{upp}, this implies that $|G_{\mathbf k}(p/q,x)| \geq e^{-o(q)}$, as desired.
\end{proof}

Once we have this, then the desired result immediately follows from the following result: 

\begin{lemma}\label{lem:local-stability}
If $|G_{\mathbf k}(p/q)| = o(1)$, then there exists $c>0$ such that for any $x,x'\in\T$,
\[
|G_{\mathbf k}(p/q,x) \Delta G_{\mathbf k}(p/q,x')| < e^{-cq}.
\]
As a consequence, we have for all $x\in\T$,
\[
|G_{\mathbf k}(p/q)| > |G_{\mathbf k}(p/q,x)| - e^{-cq}.
\]
\end{lemma}

\begin{proof}
 Let $E_+^{\mathbf{k}}(x)$ be the right boundary of $ G_{\mathbf k}(p/q,x) $, that is \[
    E_+^{\mathbf{k}}(x):=\{E: \sup_{y}\widehat \phi[\widehat C_{2dq}^{E_+^{\mathbf{k}}(x)}(x)](\Lambda_y)=2d\mathbf{l} \}. 
    \] It is clear $E_+^{\mathbf{k}}(x)$ is continuous in $x\in\T$. Let  $E_0=\min_{x\in\T}E_+^{\mathbf{k}}(x)$, $E_1=\max_{x\in\T}E_+^{\mathbf{k}}(x)$.
    For any $E$ $o(1)$-close to $[E_0,E_1]$, by \eqref{subex1} and  $ \op{tr} \widehat{C}_{2dq}^E(p/q,x)$ is $1/q$-periodic, one has
    \[
    t(E,x):=\op{tr} \widehat{C}_{2dq}^E(p/q,x)=\sum_{k\in\Z}a_{kq}(E)e^{2\pi i kq x}
    \]
    with $|a_{kq}(E)|\leq e^{o(q)}e^{-2\pi kq\epsilon}$. Therefore, for any $0<\epsilon_1<\epsilon_0$, there exists $c_1>0$ such that \begin{equation}\label{eq:trace-ecq}
        \|t(E,\cdot)-\widehat t_0(E)\|_{\epsilon_1}< e^{-c_1q}.
    \end{equation}
    Let $x_1$ be such that $E_1=E_+^{\mathbf{k}}(x_1)$. By the definition of $E_1$, it follows that \(\op{tr}\widehat C^{E_1}_{2dq}(p/q,x_1)=\pm 2\). Without loss of generality, we assume it equals to $2$.

    By \eqref{eq:trace-ecq} and \eqref{upp}, we have \[
    \|\widehat C_{2dq}^{E_1}- I_2\|_{\epsilon_1}>e^{-o(q)}, \qquad \|\det (\widehat C_{2dq}^{E_1}-I_2) \|_{\epsilon_1}< e^{-cq}.
    \]
    We need apply the following partial solution of Avila's almost reducible conjecture: \begin{lemma}[Lemma 4.6 of \cite{AviAlmost}] \label{lem:arcparabolic}For every $0<\epsilon<\epsilon_1<\epsilon_0$, $c_0>0$ and $\delta>0$,
there exist $\delta_*>0$ and $q_*>0$ with the following property.
Let $p/q \in \Q$ and
$A \in C^\omega_{\epsilon_0}(\R/\Z,\mathrm{SL}(2,\R))$ be such that $q>q_*$ and \eqref{sizek} holds.  If $\|A_q(\cdot)\pm I_2\|_{\epsilon_1} > e^{-\delta
q}$ and $\|\det (A_q\pm I_2)\|_{\epsilon_1} < e^{-c_0q}$, then there exists $\bar B \in
C^\omega_\epsilon(\R/\Z,\rmm{SL}(2,\R))$ with $\|\bar B\|_\epsilon \leq e^{C\delta q}$,  such that $$\|\bar B(x+p/q)A(x)\bar B(x)^{-1} \pm \begin{psmallmatrix}
    1& b \\ 0& 1
\end{psmallmatrix} \|_\epsilon<e^{-\gamma c_0 q}$$ with $|b|<e^{-\delta q}$,  $C>1,\gamma<1$. 
     \end{lemma} 
\begin{remark}
    The proof is an adaptation of Lemma 4.6 in \cite{AviAlmost}. It proceeds by setting $W = A_q \pm I_2$ and making minor adjustments to the parameters to fit our context. Additionally, we eliminate the non-zero Fourier modes up to order $q$.
\end{remark}
By Lemma \ref{lem:arcparabolic}, there exists $\bar B\in C^\omega_\epsilon(\T,\rmm{PSL}(2,\R))$ such that \[
    \bar B(x+p/q)\widehat C^{E_1}(x) \bar B(x)^{-1} = e^{F(x)}\begin{pmatrix}
        1 & b_0\\
        0& 1
    \end{pmatrix}
    \]
    with estimates $|b_0|<e^{-o(q)}$ and $\|F\|_\epsilon<e^{-c_2q}$ for some $c_2>0$. A direct computation (c.f. \cite{AFK2011Kam}) shows that 
 \begin{equation}\label{eq:conjugacy}
\widehat C^{E_1}_{2dq}(x) = e^{\widetilde{F}(x)} M(x), \quad \text{where } M(x) := B(x)^{-1} \begin{pmatrix} 1 & b_1 \\ 0 & 1 \end{pmatrix} B(x).
\end{equation}
    with  $\|\tilde{F}\|_\epsilon<e^{-c_3q}$ for some $c_3>0$.

Note that $P \Lambda_0 = \Lambda_0$ for the parabolic component $P = \begin{pmatrix} 1 & b_1 \\ 0 & 1 \end{pmatrix}$. Let $y_*(x)$ be defined such that $B(x)\Lambda_{y_*(x)}=\Lambda_0$. Consequently, we have the invariance $M(x)\Lambda_{y_*(x)}=\Lambda_{y_*(x)}$.

In light of \eqref{eq:conjugacy}, we fix a branch of the projective action $\widehat{\phi}$ satisfying the following properties:
\[
\begin{aligned}
    \widehat{\phi}[M(x)](\Lambda_{y_*(x)}) = 2d\mathbf{l}, \quad \text{and} \quad
    \widehat{\phi}[e^{\widetilde{F}(x)}](\Lambda) \in (-1/2, 1/2) \quad \text{for all } \Lambda.
\end{aligned}
\]
This choice of branch guarantees the additive decomposition:
\[
\widehat\phi[\widehat C^{E_1}_{2dq}(x)](\Lambda_{y}) = \widehat{\phi}[M(x)](\Lambda_{y}) + \widehat{\phi}[e^{\widetilde{F}(x)}](M(x)\Lambda_{y}).
\]
In particular, for the direction $\Lambda_{y_*(x)}$, this yields:
\[
\widehat\phi[\widehat C^{E_1}_{2dq}(x)](\Lambda_{y_*(x)}) = 2d\mathbf{l} + \widehat\phi[e^{\widetilde{F}(x)}](\Lambda_{y_*(x)}).
\]
The phase map $A \mapsto \widehat\phi[A](\Lambda)$ is Lipschitz continuous with respect to the matrix norm for $A$. Therefore, we obtain the lower bound:
\[
 \widehat\phi[\widehat C^{E_1}_{2dq}(x)](\Lambda_{y_*(x)}) \ge 2d\mathbf{l} - C e^{-c_3 q}.
\]

By Corollary~\ref{cor:monophi}, for any fixed $x$ and $\Lambda_y$, the phase $\widehat \phi$ is strictly monotone in $E$. Let us define the minimal derivative magnitude
\[
D_{\min} := \inf_{E \in [E_1-\delta_0, E_1]\setminus\{i_s\} } -\partial_E \widehat{\phi}[\widehat{C}_{2dq}^E(x)](\Lambda_{y_*(x)}),
\]
where $\delta_0$ is a suitable small constant. By Lemmas \ref{lem:LUbd} and \ref{lem:esti-deri-C_q}, we have $D_{\min} > e^{-o(q)}$. Set $\delta := \frac{Ce^{-c_3 q}}{D_{\min}}$.
Thus, for any $x\in\mathbb{T}$, the monotonicity implies:
\[
\sup_{y} \widehat{\phi}[\widehat{C}_{2dq}^{E_1-\delta}(x)](\Lambda_y)
\geq \widehat\phi[\widehat C^{E_1-\delta}_{2dq}(x)](\Lambda_{y_*(x)})
> \widehat\phi[\widehat C^{E_1}_{2dq}(x)](\Lambda_{y_*(x)}) + D_{\min} \delta
> 2d\mathbf{l}.
\]
Since $E_+^{\mathbf k}(x)$ is the energy where the phase supremum equals $2d\mathbf{l}$, it must lie within the interval $(E_1-\delta, E_1]$. Consequently, the distance is bounded by:
\[
|E_1 - E_+^{\mathbf k}(x)| \le \delta \le \frac{C e^{-c_3 q}}{e^{-o(q)}} \le e^{-c q}.
\]
This proves the exponential clustering of the right gap endpoints. An analogous analysis applied to the left boundary of $G_{\mathbf{k}}(p/q, x)$ yields the same estimate for the left endpoints, which completes the proof.
\end{proof}

\subsection{Case 2:  $\beta(\alpha)=0$.}

Again by Proposition \ref{prop:6.1}, we obtain
\begin{equation}\label{conj1}
    B^E( x + \alpha)^{-1} L^{E,d}(x) B^E(x) = \widehat{H}^E(x) \diamond e^{2\pi i \rho^E(x)} \widehat{C}^E(x),
\end{equation}
where $\rho^E( x) \in C^\omega(\mathbb{T}, \mathbb{R})$, 
and $\widehat{C}^E( x) \in C^\omega(\mathbb{T}, \mathrm{SL}(2, \mathbb{R}))$
are both piecewise $C^\omega$ with respect to $E \in \I$. Moreover, $\widehat{C}^E$ is premonotonic (which means a finite iterates of $\widehat C^E$ is monotonic) and for any $E\in\Sigma(\alpha)\cap\I$, $C^E$ is subcritical.

For the case $\beta(\alpha)=0$, a qualitative result is enough:
\begin{theorem}[\cite{avila2008absolutely,Avi2023KAM}]\label{prop-main-dio}
Let $\alpha\in\R\backslash \Q$ with $\beta(\alpha)=0$, 
$k_0\in\Z$ and  suppose that  $(\alpha, A)$ is subcritical and
$\mathbf{l}=2\rho_f(\alpha,A)-k_0\alpha \in \Z$. Then $(\alpha,A)$ is reducible, i.e.
there exists $B(x)\in C^\omega(\T,\mathrm{PSL}(2,\R))$, $c \in \R$ such that 
\begin{eqnarray}
(-1)^{\mathbf{l}} B(x+\alpha)^{-1}A(x)B(x)=I_2+c\mathfrak{L} .
\end{eqnarray}
\end{theorem}

Now consider $E_*\in\partial G_{\mathbf{k}}(\alpha)$, by Theorem \ref{prop-main-dio} and Proposition \ref{dim-duality}, 
without loss of generality, there exist $\bar B\in C^\omega(2\T,\rmm{SL}(2,\R)$,  $c>0$ such that 
\[
\bar B(x+2d\alpha)^{-1} \widehat C^{E_*}_{2d}(x) \bar B(x)= \begin{pmatrix}
    1& c\\
    0&1
\end{pmatrix}
\]
Meanwhile, set
\[
M^E(x):=\bar B(x+2d\alpha)^{-1}\,\widehat C^{E}_{2d}(x)\,\bar B(x),
\qquad
\Lambda_y:=\begin{pmatrix}\cos(\pi y)\\\sin(\pi y)\end{pmatrix}.
\]
Then we have the following: 
\begin{lemma}\label{lem:conefield}
Let \(\alpha\in\R\setminus\Q\) with \(\beta(\alpha)=0\). Fix a small
\(\varepsilon_0\in(0,\tfrac14)\) 
and set
\[
D_{\max}:=\sup_{(E,x,y)\in\I \backslash \{i_s\} \times\T\times[0,\varepsilon_0]}
\bigl(-\partial_E\phi[M^E(x)](\Lambda_y)\bigr).
\]
Then we have estimate
$$|G_{\mathbf{k}}(\alpha)|>\frac{\min\{c,1\}\tan^2(\pi\varepsilon_0)}{20D_{\max}}.$$
\end{lemma}

\begin{proof}
Consider the cone
$\mathcal C:=\{\Lambda_y:\; y\in[0,\varepsilon_0]\}.$ We make the consistent choice of  \(\phi\) so that \(\phi[M^{E_*}(x)](\Lambda_0)=0\)
for all \(x\). By Corollary~\ref{cor:monophi} the function
\(E\mapsto\phi[M^E(x)](\Lambda_y)\) is strictly decreasing, hence for
every \(x\)
\begin{equation}\label{eq:lower-end}
\phi[M^{E}(x)](\Lambda_{0})>0\qquad\text{for }E<E_*.
\end{equation}

For the upper endpoint we compute the projective angle at \(E_*\). For
any \(y\in[0,\varepsilon_0]\),
\[
\phi[M^{E_*}(x)](\Lambda_{y})
=\frac{1}{\pi}\arctan\!\Big(\frac{\sin(\pi y)}{\cos(\pi y)+c\sin(\pi y)}\Big)-y.
\]
Elementary estimates 
give the uniform bound
\begin{equation*}
-y< \phi[M^{E_*}(x)](\Lambda_{y})< -\frac{\min\{c,1\}}{20}\,\tan^2(\pi y)\qquad(y\in(0,\varepsilon_0]).
\end{equation*}
By the definition of \(D_{\max}\) we have for \(E< E_*\)
\[
\phi[M^{E}(x)](\Lambda_{y})
\le \phi[M^{E_*}(x)](\Lambda_{y}) + D_{\max}|E-E_*|.
\]
Now if we define $
\delta:=\frac{\min\{c,1\}\tan^2(\pi\varepsilon_0)}{20D_{\max}},$ this implies that for
every \(E\in(E_*-\delta,E_*)\) and every \(x\),
\begin{equation}\label{eq:upper-end}
    \phi[M^{E}(x)](\Lambda_{\varepsilon_0})
\le \phi[M^{E_*}(x)](\Lambda_{\varepsilon_0}) + D_{\max}\delta
<0.
\end{equation}

As $M^E$ induces orientation-preserving homeomorphism on $S^1$, by \eqref{eq:lower-end} and \eqref{eq:upper-end},  the projective image of the entire cone \(\mathcal C\) under \(M^E(x)\)
is strictly contained in the interior of \(\mathcal C\), uniformly in
\(x\). By the conefield criterion, this implies that \((2d\alpha,M^E(x))\)  is uniformly hyperbolic.
By \eqref{conj1}, this implies that $(\alpha,L^{E,d})$ is uniformly hyperbolic, the desired result follows. 
\end{proof}

\begin{remark}\label{mp-co}
   This proof was proposed in \cite{avila2016dry},  it recovers  \cite[Theorem 6.2]{LYZZ2017Asymptotics}, where the proof is based on the Moser-P\"oschel argument.
\end{remark}

\section{Proof of Theorem \ref{main-thm-analytic} (analytic potential)}\label{analytic-proof}

\subsection{Transition to the analytic case}

In the preceding attempts, we established that the DTMP holds for trigonometric polynomial potentials. For a general real-analytic potential $v$, however, the dual operator has infinite range, and a corresponding finite-dimensional dual cocycle cannot be directly defined. A natural strategy is therefore to approximate $v$ by its Fourier truncations $v_d(x)=\sum_{j=-d}^d\hat v_j e^{2\pi ijx}$. Since each endpoint of the gap with label $\mathbf{k}$ changes by at most $\|v-v_d\|_0$, the corresponding gap for $v$ is open provided
\[
|G_{\mathbf{k}}(v_d,\alpha)|>2\|v-v_d\|_0.
\]

The gap-opening mechanism established for $v_d$ in the previous section gives a positive lower bound for each fixed $d$, but further control of its dependence on $d$ is needed for this comparison. The relevant bounds involve the center frame and the energy derivatives of the projective action. Using the analytic regularity of $v$ and the decay of its Fourier coefficients, we track this dependence through the global symplectic reduction, together with uniform control of the scalar phases. Combined with dimension-independent Aubry duality and the cone argument below, these estimates yield gap lower bounds of the form $e^{-o(q_n)}$ at suitable approximation scales $q_n$, which dominate the exponentially small truncation errors.

 To explicitly track the dependence on the potential, let $L^v(E) = L(\alpha, S_E^v)$ and $L^v_\varepsilon(E)$ denote the real and complexified Lyapunov exponents of the Schr\"odinger cocycle $(\alpha, S_E^v)$, respectively. Let $v$ be an analytic potential of Type I, and let $E_{\mathbf{k}} \in \Sigma^{+}_{v,\alpha}$ be the maximal energy satisfying $N_{v,\alpha}(E_{\mathbf{k}}) \equiv \mathbf{k}\alpha \pmod{\mathbb{Z}}$. At this boundary, we have $\bar{\omega}(\alpha,S_{E_{\mathbf{k}}}^{v})=1$ and $L^v(E_{\mathbf{k}})>0$.   By the continuity of the Lyapunov exponent \cite{bj} and the stability of $T$-acceleration \cite{GJY}, these properties are robust under both energy perturbation and potential approximation. First, there exists a compact connected neighborhood $\I=[i_-,i_+]$ of $E_{\mathbf{k}}$ such that $L^{v}(E)>0$ and $\bar{\omega}(\alpha,S_{E}^{v})=1$ uniformly for all $E \in \I$. We may assume the strict inequalities $N_{v,\alpha}(i_-) < N_{v,\alpha}(E_{\mathbf{k}}) < N_{v,\alpha}(i_+)$ hold, as otherwise the gap $G_{\mathbf{k}}(v,\alpha)$ is trivially open. Furthermore, because $v$ is of Type I, we can ensure $\I$ contains at least one spectral gap.\footnote{Otherwise, enlarge $\I$ to cover the maximal connected component of $\Sigma_{v,\alpha}$ containing $E_{\mathbf{k}}$, and then extend it slightly beyond that component.}  Second, applying the same continuity and stability results to the approximations, there exists $d_0>0$ such that for all $d>d_0$, the truncated cocycle satisfies $L^{v_d}(E)>0$ and $\bar{\omega}(\alpha,S_{E}^{v_d})=1$ throughout $\I$. Consequently, if we let $E_\mathbf{k}^d$ denote the maximal energy satisfying $N_{v_d,\alpha}(E_\mathbf{k}^d) \equiv \mathbf{k}\alpha \pmod{\mathbb{Z}}$, this boundary is strictly confined within $\I$, yielding $N_{v_d,\alpha}(i_-) < N_{v_d,\alpha}(E_\mathbf{k}^d) < N_{v_d,\alpha}(i_+)$. Furthermore, we have the following observations:

 \begin{lemma}\label{lem:E_kd}
    We have $|E_{\mathbf{k}}^d - E_{\mathbf{k}}| \leq  \|v_d - v\|_0$.
\end{lemma}
\begin{proof}
Let $\delta = \|v_d - v\|_0$. Then,  $v(x) - \delta \le v_d(x) \le v(x) + \delta$ for each $x \in \mathbb{T}$ and therefore
\begin{equation}\label{eq:operator_order}H_{v,\alpha,x} - \delta I \le H_{v_d,\alpha,x} \le H_{v,\alpha,x} + \delta I \quad \text{for each } x \in \mathbb{T}.\end{equation}
This leads to the inequalities $N_{v,\alpha}(E - \delta) \le N_{v_d,\alpha}(E) \le N_{v,\alpha}(E + \delta)$ for all $E \in \mathbb{R}$.
 Indeed, from \(H_{v,\alpha,x} - \delta I \le H_{v_d,\alpha,x}\) we obtain \((H_{v,\alpha,x}-\delta I)^{\Lambda} \le H_{v_d,\alpha,x}^{\Lambda}\) for any $\Lambda\subset \Z$ (see Section \ref{sec:ids-revisited} for the definition of $H_{v_d,\alpha,x}^{\Lambda}$). By the Courant–Fischer min-max principle, their eigenvalues (arranged non-decreasingly) satisfy
\[
\lambda_k((H_{v,\alpha,x}-\delta I)^{\Lambda}) \le \lambda_k(H_{v_d,\alpha,x}^{\Lambda}) \quad \text{for each } k=1,...,|\Lambda| .
\]
Hence, it follows that \(\#\{k \mid \lambda_k(H_{v,\alpha,x}^{N}) \le E+\delta\} \ge \#\{k \mid \lambda_k(H_{v_d,\alpha,x}^{N}) \le E\}\) (See Section \ref{sec:ids-revisited} for the definition of $H_{v,\alpha,x}^{N}$). Dividing by \(N\) and taking the limit \(N\to\infty\) preserves the inequality, giving \(N_{v,\alpha}(E+\delta) \ge N_{v_d,\alpha}(E)\). The inequality \(N_{v_d,\alpha}(E) \ge N_{v,\alpha}(E-\delta)\) follows analogously from \(H_{v_d,\alpha,x} \le H_{v,\alpha,x} + \delta I\).

Now, set $c = N_{v,\alpha}(E_\mathbf{k})$. To establish the upper bound on $E_{\mathbf{k}}^d - E_{\mathbf{k}}$, let $\epsilon > 0$ be arbitrary and evaluate the lower bound $N_{v_d,\alpha}(E) \ge N_{v,\alpha}(E - \delta)$ at the shifted energy $E = E_{\mathbf{k}} + \delta + \epsilon$. This yields$$N_{v_d,\alpha}(E_{\mathbf{k}} + \delta + \epsilon) \ge N_{v,\alpha}(E_{\mathbf{k}} + \epsilon).$$Since $E_{\mathbf{k}}$ is defined as the maximal energy for which the integrated density of states equals $c$, the function $N_{v,\alpha}$ must strictly exceed $c$ at any larger energy argument. Thus, we have $N_{v,\alpha}(E_{\mathbf{k}} + \epsilon) > c$, which implies $N_{v_d,\alpha}(E_{\mathbf{k}} + \delta + \epsilon) > c$. By the definition of $E_{\mathbf{k}}^d$ as the maximal energy satisfying $N_{v_d,\alpha}(E_{\mathbf{k}}^d) = c$, any energy at which the IDS strictly exceeds $c$ must lie strictly to the right of $E_{\mathbf{k}}^d$. This forces the inequality $E_{\mathbf{k}}^d < E_{\mathbf{k}} + \delta + \epsilon$. Sending $\epsilon \to 0^+$ then dictates\begin{equation}\label{eq:dist_upper}E_{\mathbf{k}}^d - E_{\mathbf{k}} \le \delta.\end{equation}

Conversely, to obtain the matching lower bound, we evaluate the upper bound $N_{v,\alpha}(E + \delta) \ge N_{v_d,\alpha}(E)$ at $E = E_{\mathbf{k}}^d + \epsilon$ for an arbitrary $\epsilon > 0$, obtaining$$N_{v,\alpha}(E_{\mathbf{k}}^d + \delta + \epsilon) \ge N_{v_d,\alpha}(E_{\mathbf{k}}^d + \epsilon).$$By the identical maximality argument applied to the truncated system, we have $N_{v_d,\alpha}(E_{\mathbf{k}}^d + \epsilon) > c$, which immediately ensures $N_{v,\alpha}(E_{\mathbf{k}}^d + \delta + \epsilon) > c$. The maximality of $E_{\mathbf{k}}$ for the limit operator then requires $E_{\mathbf{k}} < E_{\mathbf{k}}^d + \delta + \epsilon$. Letting $\epsilon \to 0^+$ yields $E_{\mathbf{k}} - E_{\mathbf{k}}^d \le \delta.$ Combining this with the symmetric bound \eqref{eq:dist_upper}, we  complete the proof.\end{proof}

\subsection{Estimates on the central frame}

To establish quantitative estimates for the spectral gaps, the primary step is to bound the central frame of the conjugacy $B^{E,d}$. The construction of $B^{E,d}$ relies on isolating the center bundle and obtaining explicit quantitative bounds for its corresponding frame that precisely track the dependence on the truncation scale $d$. As a necessary preliminary, we state an adapted version of Theorem 6.3 from \cite{GJ}, which provides this required quantitative control over the center bundle. Let ${h_0=h_0(\I):=\min_{E\in\I} L^{v}(E)/(4\pi)>0}$.

\begin{theorem}\label{theorem-main1general}
{
Let $\alpha\in\R\backslash\Q$, $v\in C_{h_1}^\omega(\T,\R)$, and $\I$ a compact interval such that
$\bar{\omega}(\alpha,S_E^v)=1$ and $L^v(E)>0$ for all $E\in\I$.
Then there exist $\delta>0$ and $d_0$ such that, for every $d\geq d_0$
with $\widehat v_d\ne0$, there are
$\mathcal{O}^{E,d}\in C_{\delta,h_0}^\omega(\I\times\T,\rmm{HSp}(2d\times 2,\psi_d,\omega_1))$, and
$\mathcal M^{E,d}\in C_{\delta,h_0}^\omega(\I\times\T,\rmm{HSp}(2))$ such that
\[
L^{E,d}(x)\mathcal{O}^{E,d}(x)
=\mathcal{O}^{E,d}(x+\alpha)\mathcal M^{E,d}(x),
\]
where $\mathcal M^{E,d}\rightarrow\mathcal M^E$ for some
$\mathcal M^E\in C_{\delta,h_0}^\omega(\I\times\T,\rmm{HSp}(2))$
in the $\|\cdot\|_{\I,h_0}$-norm topology.
Moreover, for sufficiently large $d$ and every $E\in\I$, we have
\[
L(\alpha,\mathcal M^{E,d}(\cdot+i\e))
=L(\alpha,\mathcal M^{E,d}(\cdot)),
\qquad \forall\,|\e|<h_0.
\]
If, in addition, $E\in\Sigma_{v_d,\alpha}\cap\I$, then
\[
L(\alpha,\mathcal M^{E,d}(\cdot+i\e))=0,
\qquad \forall\,|\e|<h_0.
\]
Furthermore, if we write
\[
\mathcal{O}^{E,d}(x)=\begin{pmatrix}
\mathcal{O}^{E,d}(x,d-1) &
\mathcal{O}^{E,d}(x,d-2) &
\cdots &
\mathcal{O}^{E,d}(x,-d)
\end{pmatrix}^\top,
\]
then there exist $0<\gamma<h_1$ and $K\geq1$, depending only on
$\alpha,v,h_1,\I$ and independent of $d,d'$, such that,
after increasing $d_0$, the following estimates hold:
\begin{enumerate}[label=(H\arabic*), leftmargin=*, itemsep=3pt]
\item For every $\eta>0$, there exist $r_\eta,C_\eta>0$ and
$d_\eta\geq d_0$, depending only on $\alpha,v,h_1,\I,\eta$, such that
\[
\|\mathcal{O}^{z,d}(\cdot,k)\|_{h_0}\leq C_\eta e^{\eta|k|},
\qquad -d\leq k\leq d-1,
\]
whenever $z\in\I_\delta$, $\operatorname{dist}(z,\Sigma_{v,\alpha})<r_\eta$
and $d\geq d_\eta$.
\item For $z\in\I_\delta$ with
$\operatorname{dist}(z,\Sigma_{v,\alpha})=o(1)$ and $d'>d\geq d_0$,\footnote{More precisely, there exists $r=r(\alpha,v,h_1,\I)>0$, independent of $d,d'$, such that the estimate holds whenever $\operatorname{dist}(z,\Sigma_{v,\alpha})<r$.}
\[
\|\mathcal{O}^{z,d}(\cdot,k)-\mathcal{O}^{z,d'}(\cdot,k)\|_{h_0}
\leq K\|v_d-v_{d'}\|_{h_1-\gamma}
\leq e^{-\gamma d},
\qquad -d\leq k\leq d-1.
\]
\item $\|\mathcal M^{E,d}-\mathcal M^E\|_{\I,h_0}
\leq K\|v_d-v\|_{h_1-\gamma}$ for $d\geq d_0$.
\end{enumerate}
}
\end{theorem}

\begin{proof}
{
The proof is adapted from Theorem 6.3 in \cite{GJ}, with some necessary adjustments.
For the reader's convenience, we provide a sketch, following the
notation and labels there. We explain the energy dependence and
the convergence estimates. Constants may depend on
$\alpha,v,h_1,\I$, but not on $d,d'$ or the indices $i,j,k,\ell$.
We use $C$ for constants that may change from line to line;
additional dependence is indicated by subscripts, as in $C_m$
and $C_\eta$. All auxiliary choices below are fixed in terms of
these data. The common constant $K$ in (H2)--(H3) is chosen at
the end. Truncation degrees are taken through indices with
$\widehat v_d\ne0$, as appropriate to the infinite-range case.
}

Note that, for real $E\in\I$ and all sufficiently large $d$,
the proof of Lemma 6.3 of \cite{GJ} gives a two-dimensional
center, holomorphic on $\T_{3h_0/2}$.
For a normalized center frame
$H_d\in\rmm{HSp}(2d\times2,\psi_d,\omega_1)$, (6.29) gives
\[
G_d=H_d(H_d^\sharp S_dH_d)^{-1}H_d^\sharp
=-H_d\mathcal J_2H_d^\sharp.
\]
Here $F^\sharp(z,\theta)=F(\bar z,\bar\theta)^*$.
This expression is independent of the choice of frame and satisfies
\[
G_dS_dG_d=G_d,\qquad
\operatorname{tr}(G_dS_d)=2,\qquad G_d^\sharp=-G_d.
\]
On the other hand, the resolvent construction at the begining can extends $G_d$ holomorphically
to a common domain $\I_{2\delta}\times\T_{3h_0/2}$, with
$\delta,d_0$ independent of $d$.
By the identity theorem, the above identities hold throughout
this domain. Since $S_d$ is invertible and $(G_dS_d)^2=G_dS_d$,
\[
\operatorname{rank}G_d=\operatorname{rank}(G_dS_d)
=\operatorname{tr}(G_dS_d)=2.
\]

Indeed, by (6.9)--(6.14), stability and compactness,
choose finitely many real-centered energy disks $U_\nu$ and
shifts $\varepsilon_\nu<h_1$ such that
\[
\overline{\I_{2\delta}}\subset\bigcup_{\nu=1}^{N}U_\nu,
\qquad
\varepsilon_1(E)<\varepsilon_\nu<\varepsilon_2(E)
\quad (E\in U_\nu\cap\R).
\]
We take $\delta>0$ small enough that the real energies remain
Type I and $L^v(E)/(2\pi)>3h_0/2$ on these neighborhoods.
The corresponding dual resolvents are holomorphic and uniformly bounded on a
neighborhood of $\overline{U_\nu}\times\overline{\T_{3h_0/2}}$
for all $d\geq d_0$.
Their weighted differences (6.16) give local holomorphic extensions of
$G_d$. These agree on overlaps by the identity theorem, since
they agree for real energies, and hence define the common extension.
The resolvent identity and (6.16), (6.28) give,
for $z\in U_\nu$,
\begin{equation}\label{eq:gj-selected-column-bounds}
\begin{aligned}
\|u_d^\ell(z,\cdot,j)\|_{3h_0/2}
&\leq Ce^{2\pi\varepsilon_\nu|j-\ell|},\\
\|u_d^\ell(z,\cdot,j)-u_{d'}^\ell(z,\cdot,j)\|_{3h_0/2}
&\leq Ce^{2\pi\varepsilon_\nu|j-\ell|}
\|v_d-v_{d'}\|_{\varepsilon_\nu}.
\end{aligned}
\end{equation}
Here the difference of the shifted dual operators is Fourier
equivalent to multiplication by
$(v_d-v_{d'})(\cdot\pm i\varepsilon_\nu)$, so its norm is at most
$\|v_d-v_{d'}\|_{\varepsilon_\nu}$.
For (6.23)--(6.24), matching the half-line square-summable solutions
at the source gives an $\ell^2$ solution with right-hand side
$\delta_\ell$. The shifted operators are already invertible, so
uniqueness identifies this solution with the resolvent column.
With $u_d^\ell$ defined as in (6.16), direct substitution of
(6.26)--(6.27) into (6.23)--(6.24) gives, with the convention (6.29),
\[
(G_d)_{j,\ell}(z,\theta)=-u_d^\ell(z,\theta,j).
\]
This corrects the overall sign in (6.28).
The expression is independent of the admissible shift, and the
same estimates apply to the intrinsic $G_d$ defined above.

{
We use the same balance between spatial growth and Fourier decay
as in the weighted estimates of Section~\ref{sec:Dim-free-Aubry},
in particular the truncated Wronskian in
Lemma~\ref{lem:lowerb-u0-new}.
Put $\varepsilon_{\max}=\max_\nu\varepsilon_\nu<h_1$ and fix
$0<\gamma<(h_1-\varepsilon_{\max})/4$ for the Fourier tail estimates.
For each fixed $m$, \eqref{eq:gj-selected-column-bounds} gives
a holomorphic limit $G^m$ of the central blocks $G_d^m$, with
\begin{equation}\label{eq:weighted-center-convergence}
\|G_d^m(z,\cdot)-G^m(z,\cdot)\|_{3h_0/2}
\leq C_m\|v_d-v\|_{\varepsilon_{\max}},
\qquad z\in\I_{2\delta}.
\end{equation}
The compatible limits $G^m$ define an infinite matrix
$G=(G_{j,\ell})_{j,\ell\in\Z}$ whose central blocks are
\[
G^m=(G_{j,\ell})_{-m\leq j,\ell\leq m-1}.
\]
Since $\operatorname{rank}G_d^m\leq2$, all $3\times3$ minors of
$G^m$ vanish by convergence. Hence every finite block of $G$
has rank at most two.
}

The symplectic form has entries
\[
(S_d)_{j,\ell}
=(\mathbf1_{j<0}-\mathbf1_{\ell<0})
\widehat v_{\ell-j}\mathbf1_{|\ell-j|\leq d}.
\]
Let $S_{j,\ell}$ denote the same expression without the Fourier cutoff.
Only pairs crossing the origin contribute, and for these
$|j-\ell|=|j|+|\ell|$. Thus the same summation as in
Lemma~\ref{lem:lowerb-u0-new} gives, for fixed $i,k$,
\[
\begin{split}
&\sum_{\max\{|j|,|\ell|\}>N}
|(G_d)_{i,j}(S_d)_{j,\ell}(G_d)_{\ell,k}|\\
&\qquad\leq Ce^{2\pi\varepsilon_{\max}(|i|+|k|)}
\sum_{m>N}m e^{-2\pi(h_1-\varepsilon_{\max})m}.
\end{split}
\]
Here entries are set to zero whenever a row or column index lies
outside $\{-d,\ldots,d-1\}$. The bound is uniform on the common
complex domain.
The tail of $\sum_{j,\ell}(G_d)_{j,\ell}(S_d)_{\ell,j}$ satisfies
the same bound without the factor depending on $i,k$.
Consequently, the identities $G_dS_dG_d=G_d$ and
$\operatorname{tr}(G_dS_d)=2$ from (6.29) pass to absolutely
convergent sums:
\begin{equation}\label{eq:weighted-limit-algebra}
\begin{aligned}
\sum_{j,\ell\in\Z}G_{i,j}S_{j,\ell}G_{\ell,k}&=G_{i,k},\\
\sum_{j,\ell\in\Z}G_{j,\ell}S_{\ell,j}&=2.
\end{aligned}
\end{equation}
These identities force $G$ to have rank two: rank zero is excluded
by the second identity, while in rank one the first would make
the second sum equal to one.
By compactness, a fixed central window contains a nonzero
$2\times2$ minor at every point of the closed common parameter
domain. We may therefore choose $n_0$, independently of $d$, such
that $G^{n_0}$ has rank two throughout that domain.
For real parameters, $iG^{n_0}$ has signature $(1,1)$, inherited
from the finite-dimensional center. This gives the fixed block
and convergence needed.

We now continue as in Lemma 6.6 and (6.33)--(6.36).
The kernel of $G^{n_0}$ has a holomorphic complement on the common
Stein domain; Proposition~\ref{lem:ajslem} gives a periodic
holomorphic frame $e=(e_1,e_2)$ for this complement.
With $L_d^{n_0}$ denoting the corresponding columns, set
\[
O'_d=L_d^{n_0}e,\qquad
\Omega'=-e^\sharp G^{n_0}e.
\]
By the same argument as in Lemma~\ref{lem:Key-holom-extend},
the identity $(G^{n_0})^\sharp=-G^{n_0}$ and the choice of $e$
imply that $\Omega'$ is nondegenerate.
Lemma 6.1 of \cite{GJ}, as used in (6.36), supplies a periodic
analytic initial normalization $P_0$ at a fixed $E_{\mathrm{ref}}\in\I$
(with a constant change of basis to match our sign convention).
Solve
\[
\partial_zP=-\tfrac12(\Omega')^{-1}(\partial_z\Omega')P,
\qquad P(E_{\mathrm{ref}},\theta)=P_0(\theta).
\]
On the symmetric simply connected energy neighborhood this gives
a jointly holomorphic invertible $P$, and differentiation shows
$P^\sharp\Omega'P=\mathcal J_2$.
We now restrict to $\I_\delta\times\T_{h_0}$, where the fixed
frames and their required inverses are uniformly bounded.
All subsequent estimates are on this domain unless otherwise stated.
By (6.35), set
\[
K_d=-\mathcal J_2^{-1}P^\sharp e^\sharp G_d^{n_0}eP,
\qquad \mathcal O^{z,d}=O'_dP K_d^{-1/2}.
\]
Here $K_d\to I_2$ uniformly, and the power-series branch gives
an analytic, uniformly Lipschitz correction. The identity
$K_d^\sharp\mathcal J_2=\mathcal J_2K_d$ implies
$(\mathcal O^{z,d})^\sharp S_d\mathcal O^{z,d}=\mathcal J_2$.

The fixed-block estimate \eqref{eq:weighted-center-convergence}
also holds between degrees $d,d'$, with
$\|v_d-v_{d'}\|_{\varepsilon_{\max}}$ on the right.
It passes to $K_d^{-1/2}$ and to the fixed blocks of
$\mathcal O^{z,d}$. Recover $\mathcal M^{z,d}$ from one such block
and its spatial shift as in (6.40): after multiplication by
$e^\sharp(z,\theta+\alpha)$, the coefficient matrix tends to
$-\Omega'P$, which is uniformly invertible. Consequently,
\[
\|\mathcal M^{z,d}-\mathcal M^{z,d'}\|_{h_0}
\leq C\|v_d-v_{d'}\|_{\varepsilon_{\max}}.
\]
Since $\varepsilon_{\max}<h_1-\gamma$, letting $d'\to\infty$ proves (H3).
For real energies, the invariance and the stated Lyapunov exponent
properties follow from the identification with the dual center,
as in \cite{GJ}. The conjugacy identity extends to complex energies
by the identity theorem.

For (H1), fix $\eta>0$. At every
$E_0\in\Sigma_{v,\alpha}\cap\overline{\I_{2\delta}}$,
the first turning point is zero, so choose
\[
0<\varepsilon<\min\{\varepsilon_2(E_0),\eta/(4\pi),(h_1-\gamma)/4\}.
\]
By stability and compactness, there are finitely many neighborhoods
$V_\nu$, $0<r_\eta<\delta$, and $d_\eta$ such that
\[
\{z\in\I_\delta:\operatorname{dist}(z,\Sigma_{v,\alpha})<r_\eta\}
\subset\bigcup_{\nu=1}^{N_\eta}V_\nu.
\]
Each chosen shift remains admissible on its $V_\nu$ for
$d\geq d_\eta$, and \eqref{eq:gj-selected-column-bounds} gives
\[
\|u_d^\ell(z,\cdot,k)\|_{h_0}
\leq C_\eta e^{\eta|k-\ell|/2},
\qquad z\in\bigcup_{\nu=1}^{N_\eta}V_\nu,\quad d\geq d_\eta.
\]
By (6.33), the rows of the chosen frame are
\[
\mathcal O^{z,d}(\cdot,k)
=-\bigl(u_d^{n_0-1}(z,\cdot,k),\ldots,u_d^{-n_0}(z,\cdot,k)\bigr)
eP K_d^{-1/2}.
\]
The number of selected columns is fixed, and $eP K_d^{-1/2}$ is
uniformly bounded independently of $d$. Hence
\begin{equation}\label{eq:complex-energy-row-bound}
\|\mathcal O^{z,d}(\cdot,k)\|_{h_0}
\leq C_\eta e^{\eta|k|/2},
\end{equation}
which proves (H1). If the spectral set is empty, the near-spectrum
assertions are vacuous for a sufficiently small neighborhood.

For (H2), fix $0<\eta<\pi(h_1-\gamma)$ in terms of $h_1,\gamma$
and use the preceding finite cover. The resulting $C_\eta$ is now
absorbed into $C$. On each member, with its corresponding shift
$\varepsilon$, the resolvent identity and (6.33) give
\[
\begin{aligned}
\|O'_d(\cdot,k)-O'_{d'}(\cdot,k)\|_{h_0}
&\leq C e^{\eta|k|/2}\|v_d-v_{d'}\|_\varepsilon,\\
\|K_d-K_{d'}\|_{h_0}
&\leq C\|v_d-v_{d'}\|_\varepsilon.
\end{aligned}
\]
Subtracting $\mathcal O^{z,d}=O'_dP K_d^{-1/2}$ and using the
row bound from (H1) for the normalization error therefore gives
\begin{equation}\label{eq:center-row-difference-growth}
\|\mathcal O^{z,d}(\cdot,k)-\mathcal O^{z,d'}(\cdot,k)\|_{h_0}
\leq C e^{\eta|k|/2}\|v_d-v_{d'}\|_\varepsilon.
\end{equation}
Since $w=v_d-v_{d'}$ has Fourier support in $|k|>d$,
\[
\|w\|_\varepsilon
\leq C e^{-2\pi(h_1-\gamma-\varepsilon)d}\|w\|_{h_1-\gamma}.
\]
The constant is uniform because $\varepsilon<(h_1-\gamma)/4$.
For $|k|\leq d$, our choice of $\eta$ gives
$2\pi(h_1-\gamma-\varepsilon)-\eta/2>\pi(h_1-\gamma)$.
Thus the decay absorbs the row growth, proving the first bound
in (H2) on the fixed neighborhood
$\operatorname{dist}(z,\Sigma_{v,\alpha})<r_\eta$.
Finally, the Fourier decay of $v$ gives
\[
\|v_d-v_{d'}\|_{h_1-\gamma}\leq C e^{-2\pi\gamma d}.
\]
Choose $K\geq1$ to dominate the constants in (H2)--(H3), and then
increase $d_0$ to include the threshold just chosen and to ensure
\[
K\|v_d-v_{d'}\|_{h_1-\gamma}\leq e^{-\gamma d},
\qquad d'>d\geq d_0.
\]
The thresholds $d_\eta$ in (H1) may be increased accordingly.
All these estimates are uniform for $d'>d$ in the stated domains.\qedhere

\end{proof}

Having established the uniform estimates for the center frame $\mathcal{O}^{E,d}$, we are now in a position to assemble the full conjugacy. The following proposition embeds this center frame into a global symplectification, providing the complete conjugacy $B^{E,d}$ and the quantitative bounds necessary for our gap estimates.

\begin{proposition}\label{prop:10.3}
Let $\alpha\in \R\backslash\Q$, $v\in C_{h_1}^\omega(\T,\R)$, and $\I$ as above. Choose compact intervals
\[
\I_{d_*}=[E_{\mathbf k}-r_{d_*},E_{\mathbf k}+r_{d_*}]
\subset\I,\qquad r_{d_*}\to0,
\]
We choose $r_{d_*}\to0$ sufficiently slowly that $r_{d_*}\geq e^{-\sqrt{d_*}}$ and
$|E_{\mathbf k}^{d}-E_{\mathbf k}|<r_{d_*}/2$ for all $d\geq d_*$,
for sufficiently large $d_*$. For every such $d_*$, there exists $\epsilon_*(\I,d_*)$, such that for any $0<\epsilon<\epsilon_*$, there exist $
B^{E,d}_{\epsilon,d_*}(x)\in C^\omega_{h_0}(\mathbb{T},\mathrm{HSp}(2d,\psi_d,\omega_d))$ which is
piecewise $C^\omega$ in $E\in\I$, such that
\begin{equation}\label{block-p-analytic}
B^{E,d}_{\epsilon,d_*}(x+\alpha)^{-1}L^{E,d}(x)B^{E,d}_{\epsilon,d_*}(x)
= \widehat H^{E,d}(x)\diamond e^{2\pi i\rho^{E,d}_{\epsilon,d_*}(x)}\widehat C^{E,d}_{\epsilon,d_*}(x).
\end{equation}
Here the components satisfy:
\begin{itemize}
    \item $\widehat H^{E,d}(x)\in C^\omega_{\delta',h_0}(\I\times\T,\mathrm{HSp}(2d-2))$ for some $\delta'>0$;
    \item $\rho^{E,d}_{\epsilon,d_*}(x)\in C^\omega_{h_0}(\mathbb{T},\mathbb{R})$, $\widehat C^{E,d}_{\epsilon,d_*}(x)\in C^\omega_{h_0}(\mathbb{T},\mathrm{SL}(2,\mathbb{R}))$ are piecewise $C^\omega$ in $E\in\I$;
    \item $\widehat C^{E,d_*}_{\epsilon,d_*}(x)\in C^\omega_{h_0}(\mathbb{T},\mathrm{SL}(2,\mathbb{R}))$ is premonotonic. More precisely, the iterate
    \[
    C^{E,d_*}_{\epsilon,d_*}(x):=(\widehat{C}^{E,d_*}_{\epsilon,d_*})_{2d_*}(x)
    \]
    is piecewise monotonic.
\end{itemize}
For sufficiently large $d$ and any $E\in\I$, we have
\[
L\bigl(\alpha,\widehat C^{E,d}_{\epsilon,d_*}(\cdot + i\varepsilon)\bigr)=L\bigl(\alpha,\widehat C^{E,d}_{\epsilon,d_*}(\cdot)\bigr),\qquad \forall\,|\varepsilon|<h_0,
\]
and, if $E\in\Sigma(v_d,\alpha)\cap\I$, then the cocycle $(\alpha,\widehat C^{E,d}_{\epsilon,d_*})$ is subcritical, i.e.,
\[
L\bigl(\alpha,\widehat C^{E,d}_{\epsilon,d_*}(\cdot + i\varepsilon)\bigr)=0,\qquad \forall\,|\varepsilon|<h_0.
\]
Moreover, for fixed $d_*$ and $\epsilon$, the phases can be chosen compatibly so that
\[
\|\rho^{E,d}_{\epsilon,d_*}-\rho^E_{\epsilon,d_*}\|_{\I,h_0}
+\|\widehat C^{E,d}_{\epsilon,d_*}-\widehat C^E_{\epsilon,d_*}\|_{\I,h_0}
\longrightarrow0,
\]
where $\rho^E_{\epsilon,d_*}\in C^\omega_{h_0}(\mathbb T,\R)$ and
$\widehat C^E_{\epsilon,d_*}\in C^\omega_{h_0}(\mathbb T,\rmm{SL}(2,\R))$
are piecewise $C^\omega$ in $E\in\I$.

If we denote the matrix formed by the $d$-th and $2d$-th columns
of $B^{E,d}_{\epsilon,d_*}(x)$ by
\[
U^{E,d}_{\epsilon,d_*}(x)
=
\begin{pmatrix}
U^{E,d}_{\epsilon,d_*}(x,d-1) &
U^{E,d}_{\epsilon,d_*}(x,d-2) &
\cdots &
U^{E,d}_{\epsilon,d_*}(x,-d)
\end{pmatrix}^\top,
\]
then there exist $\gamma>0$ and $K>0$, independent of $d$, $d_*$
and $\epsilon$, such that the following estimates hold on
$\I_{d_*}$ for $d\geq d_*$ and sufficiently large $d_*$.
\begin{enumerate}[label=(H\arabic*), leftmargin=*, itemsep=3pt]
\item For each $-d\leq k\leq d-1$, 
\[
\|U^{E,d}_{\epsilon,d_*}(x,k)\|_{\I_{d_*},h_0}
\leq K e^{o(|k|)}.
\]

\item Let $d'>d$. For each $-d\leq k\leq d-1$,
\[
\|U^{E,d}_{\epsilon,d_*}(x,k)
-U^{E,d'}_{\epsilon,d_*}(x,k)\|_{\I_{d_*},h_0}
\leq K\|v_d-v_{d'}\|_{h_1-\gamma}
\leq e^{-\frac{\gamma}{2}d}.
\]

\item $\|\widehat C^{E,d}_{\epsilon,d_*}
-\widehat C^E_{\epsilon,d_*}\|_{\I_{d_*},h_0}
\leq K\|v_d-v\|_{h_1-\gamma}$.

\item For $E,E'\in\I_{d_*}$,
\[
\|\widehat C^E_{\epsilon,d_*}
-\widehat C^{E'}_{\epsilon,d_*}\|_{h_0}
\leq K|E-E'|.
\]

\item For $d=d_*$, uniformly in $\Lambda$, we have
\[
\inf_{E\in\I_{d_*}\backslash\{i_s\}}
-\partial_E\phi[\widehat C^{E,d_*}_{\epsilon,d_*}(x)](\Lambda)
>-K\epsilon.
\]
\end{enumerate}
For $d\geq d_*$ and sufficiently large $d_*$, the phases also satisfy
\[
\|\rho^{E,d}_{\epsilon,d_*}\|_{\I_{d_*},h_0}
+\|\rho^E_{\epsilon,d_*}\|_{\I_{d_*},h_0}\leq K,
\]
with $K$ independent of $d,d_*$ and $\epsilon$.
\end{proposition}

\begin{proof}

The frames in the hyperbolic part do not affect the quantitative estimates. Thus, one can select any frame given in Proposition \ref{prop:globalsympframe}. For the center part, we select the frames $\mathcal{O}^{E,d}$  explicitly constructed in Theorem \ref{theorem-main1general}. Denote the extended frames by $\mathcal{B}^{E,d}(x)$.  
Then, the proposition follows essentially from \cite[Proposition 9.1]{WAXZ}; the only change is to replace the local symplectification used there (\cite[Theorem 4.2]{WAXZ}) by our global symplectification (Proposition \ref{prop:globalmonoframe}). By Theorem \ref{theorem-main1general}, the symplectic structure used here is $\mathcal{B}^{E,d}(\bar{x})^* S_d\mathcal{B}^{E,d}(x) = \mathcal{J}_{2d} $. 

For each fixed reference dimension $d_*$, we apply
Proposition~\ref{prop:globalmonoframe} on the interval $\I$.
We take $E_{\mathbf k}$ as the initial point of the parallel transport,
with the original center frame there, and propagate in both energy
directions. We include the endpoints of $\I_{d_*}$ among the
partition points. Thus the transport to a point of $\I_{d_*}$
only uses local segments contained in $\I_{d_*}$.
The mesh size is smaller than $\epsilon$, and its admissible upper
bound may depend on $d_*$. The resulting two-dimensional gauge is
then used for every $d$; it is held fixed when $d\to\infty$.
The construction is global on $\I$, whereas the constants in
(H1)--(H5) are estimated only on $\I_{d_*}$.
We suppress $\alpha$ in the notation, as the irrational $\alpha$ is
fixed throughout. The following simple observation gives the
dimension-independent pairing estimates: 
\begin{lemma}
\label{lem:psi-bound}
Let $\{\hat v_k\}_{k\in \Z}$ satisfy $|\hat v_k| \le C_0 e^{-\gamma |k|}$ for some $C_0,\gamma>0$,
and define $S_d$ as before.  
Suppose for each $d$, vectors $u^{d}, v^{d} \in \mathbb{C}^{2d}$ satisfy, for some $\varepsilon < \gamma$ and $C>0$,
\[
|u^{d}(k)| \le C e^{\varepsilon |k|}, \qquad |v^{d}(k)| \le C e^{\varepsilon |k|}, \qquad
k = d-1,\dots,0,-1,\dots,-d .
\]
Then there exists a constant $K = K(C_0,C,\gamma,\varepsilon)$ independent of $d$ such that
\[
\bigl|\psi_d(u^{d}, v^{d})\bigr| \le K .
\]
\end{lemma}
\begin{proof}
    A straightforward computation.
\end{proof}
In the following quantitative estimates, $s,t\in\I_{d_*}$
and $d\geq d_*$. Since $E_{\mathbf k}\in\Sigma_{v,\alpha}$ and
$r_{d_*}\to0$, these energies approach the fixed spectrum
$\Sigma_{v,\alpha}$ as $d_*\to\infty$. Hence (H1) and (H2) of
Theorem~\ref{theorem-main1general} apply uniformly for $d\geq d_*$,
once $d_*$ is sufficiently large.
For every $\eta>0$, the row bounds are $C_\eta e^{\eta|k|}$,
uniformly for $d\geq d_*$ once $d_*$ is sufficiently large.
Set $\mc{O}^{t,d}(x)=\begin{pmatrix}
    u_{d}^{t,d}(x)& u_{-d}^{t,d}(x)
\end{pmatrix}$. Combine Theorem \ref{theorem-main1general}(H1) with the replaced symplectic relation \[
 \mathcal{B}^{t,d}(x)^{-1}=-\mathcal{J}_{2d}\mathcal{B}^{t,d}(\bar{x})^* S_d,
 \] to \eqref{eq:representQ}, Lemma \ref{lem:psi-bound} shows that  \[Q^{t,d}_s(x)=\begin{pmatrix}
\psi_d\bigl(u_{-d}^{s,d}(\bar{x}),\, u_d^{t,d}(x)\bigr) & \psi_d\bigl(u_{-d}^{s,d}(\bar{x}),\, u_{-d}^{t,d}(x)\bigr) \\[4pt]
-\psi_d\bigl(u_d^{s,d}(\bar{x}),\, u_d^{t,d}(x)\bigr) & -\psi_d\bigl(u_d^{s,d}(\bar{x}),\, u_{-d}^{t,d}(x)\bigr)
\end{pmatrix}.\] Thus the entries of $Q_s^{t,d}$ are uniformly bounded by
Lemma~\ref{lem:psi-bound}.

Fix $b>0$ smaller than the Fourier decay rate in
Lemma~\ref{lem:psi-bound}. By the uniform interpretation of
\eqref{eq:complex-energy-row-bound}, choose a fixed $\rho>0$ so that
the row bound $C_b e^{b|k|}$ holds for $|z-E_{\mathbf k}|<2\rho$
and all sufficiently large $d$.
For sufficiently large $d_*$, we have
$\I_{d_*}\subset(E_{\mathbf k}-\rho,E_{\mathbf k}+\rho)$.
Applying Lemma~\ref{lem:psi-bound} to the same inner products for
complex $t$ gives
\[
\sup_{s\in\I_{d_*},\ |z-E_{\mathbf k}|<2\rho}
\|Q_s^{z,d}\|_{h_0}\leq C,\qquad d\geq d_*.
\]
Since $z\mapsto Q_s^{z,d}$ is holomorphic, Cauchy's estimate yields
\[
\sup_{s,t\in\I_{d_*}}
\|\partial_t^j Q_s^{t,d}\|_{h_0}
\leq j!C\rho^{-j},\qquad j=1,2.
\]
These constants are independent of $d$ and $d_*$.
Since $Q_s^{s,d}=I_2$, the local radius in
Lemma~\ref{lem:detQ} can consequently be chosen independently
of $d$ and $d_*$, also for complex $t$ near $s$.
We restrict the real local segments to $\I_{d_*}$.

Recall \eqref{eq:represent-of-Phi-2}, for any $w'\in E^c_s(d)$,
\begin{equation}
    \Phi_{t,s}^{d}(w') = (I_{2d} - \mathbb{P}_{E_{s}^c(d)}) \, \mathcal{O}^{t,d}(x)  \, Q_s^{t,d}(x)^{-1} z .
\end{equation}
where $z=(\psi_d(u_{-d}^{s,d}(\bar{x}), w'), -\psi_d(u_d^{s,d}(\bar{x}), w'))^\top$. Then by \eqref{eq:ProjeE_c}, 
\[
\begin{aligned}
    \widehat{u}_{\pm d}^{t,d}(x) = \widehat{\nabla}_{t,s}^{d} u_{\pm d}^{s,d}(x)&=u_{\pm d}^{s,d}(x)+ \mc V - u_d^{s,d}(x)\psi_d (u_{-d}^{s,d}(\bar{x}),\mc V) + u_{-d}^{s,d}(x) \psi_d (u_d^{s,d}(\bar{x}),\mc V),
\end{aligned}
\]
where $\mc V=\mathcal{O}^{t,d}(x)  \, Q_s^{t,d}(x)^{-1} \begin{pmatrix}
    \psi_d(u_{-d}^{s,d}(\bar{x}), u_{\pm d}^{s,d}(x))\\
    -\psi_d(u_{d}^{s,d}(\bar{x}), u_{\pm d}^{s,d}(x))
\end{pmatrix}$. The preceding representation and Lemma~\ref{lem:psi-bound}
preserve the fixed exponential row bounds on a fixed complex
neighborhood of $s$, where $Q_s^{t,d}$ is uniformly invertible.
Applying Lemma~\ref{lem:psi-bound} to the entries of $J_s$ on
this neighborhood, and then Cauchy's estimate on a smaller
neighborhood, yields the uniform constants needed in
Lemma~\ref{lem:N_s(t)}. Here the adjoint factors are extended
by reflection in the energy variable, so the pairings are holomorphic. Moreover, since $N_s(t,x)$ is close to $I_2$,
Theorem~\ref{theorem-main1general}(H1) and the graph-transform
formula imply that there exists $K_1$ independent of $d$ and $d_*$
such that
\[
\|\bigl(\widetilde{\nabla}_{t,s}^d u_{\pm d}^{s,d}\bigr)(k,\cdot)\|_{h_0} \leq K_1 e^{o(|k|)} \text{ for } -d\leq k \leq d-1.
\]
Here the estimate is uniform over the admissible local pairs
$s,t\in\I_{d_*}$. By \eqref{eq:R_tsalpha},
Lemma~\ref{lem:psi-bound} bounds the holomorphic inner products
on the same fixed complex neighborhoods. Cauchy's estimate
then shows that $\eta$ in Lemma~\ref{lem:localR_t} is independent
of $d$ and $d_*$ on these local segments.
Other constants \(\bar C'\) in \eqref{eq:partial-J_s}, \(\delta_2\) in Lemma~\ref{lem:N_s(t)}, \(\eta'\) in Lemma~\ref{lem:quanti-tildeB}, \(C'''\) in Lemma~\ref{lem:LUbd}
are independent of $d$ and $d_*$ on these windows in the same manner. 

For given $d_*$, take any $0<\epsilon<\epsilon_*$ for some sufficiently small $\epsilon_*(\I,d_*)<\delta_2$, by \eqref{eq:4.28} and \eqref{deri-phi-center}, (H5) holds for
$d=d_*$ and $E\in\I_{d_*}$ after excluding the finite set
$\{i_s(\epsilon,d_*)\}$ where differentiability fails.
No monotonicity assertion is made for $d\neq d_*$.
The global premonotonicity at $d=d_*$ follows from the fixed-dimensional
construction on $\I$, with the admissible mesh bound depending on $d_*$.
By the total holonomy \eqref{eq:total-holonomy}, set
\[
U^{E,d}_{\epsilon,d_*}(x)=\mathcal{O}^{E,d}(x)\widetilde{R}^{\epsilon,d_*}_E(x).
\]
Set
\[
D^{E,d}_{\epsilon,d_*}(x)
:=\widetilde{R}^{\epsilon,d_*}_E(x+\alpha)^{-1}
\mc M^{E,d}(x)\widetilde{R}^{\epsilon,d_*}_E(x),
\]
and define $D^E_{\epsilon,d_*}$ analogously with $\mc M^E$.
Both are $1$-periodic, since the center cocycles and the holonomy gauge are $1$-periodic.
For fixed $\epsilon,d_*$, the gauge and its inverse are bounded
in the $\|\cdot\|_{\I,h_0}$-norm, so
$D^{E,d}_{\epsilon,d_*}\to D^E_{\epsilon,d_*}$ in this norm by
Theorem~\ref{theorem-main1general}(H3).
By the same argument as in Proposition~\ref{prop:6.1}, using
Theorem~\ref{theorem-main1general} at $E=E_{\mathbf k}^d$,
the determinants have zero winding number throughout $\I$,
and the same holds for the limit by convergence.
Their determinants therefore admit $1$-periodic logarithms, and the
square roots admit real-analytic phases on $\T$.
Dividing by compatible square roots gives the real centers
\[
\widehat C^{E,d}_{\epsilon,d_*}
=\frac{D^{E,d}_{\epsilon,d_*}}{\sqrt{\det D^{E,d}_{\epsilon,d_*}}},
\qquad
\widehat C^E_{\epsilon,d_*}
=\frac{D^E_{\epsilon,d_*}}{\sqrt{\det D^E_{\epsilon,d_*}}},
\]
and hence
\begin{equation}\label{eq:limitCd*}
\begin{aligned}
D^{E,d}_{\epsilon,d_*}(x)
&=e^{2\pi i\rho^{E,d}_{\epsilon,d_*}(x)}
\widehat C^{E,d}_{\epsilon,d_*}(x),\\
D^E_{\epsilon,d_*}(x)
&=e^{2\pi i\rho^E_{\epsilon,d_*}(x)}
\widehat C^E_{\epsilon,d_*}(x).
\end{aligned}
\end{equation}
The Hermitian-symplectic identity bounds the inverses of the centers
in the same strip. With compatible branches, division by the square
root of the determinant is locally analytic in the $h_0$-norm.
Thus $\widehat C^{E,d}_{\epsilon,d_*}\to\widehat C^E_{\epsilon,d_*}$
in the $\|\cdot\|_{\I,h_0}$-norm. All these choices are piecewise
analytic in energy.

On $\I_{d_*}$, the holonomy starts at $E_{\mathbf k}$,
and its norm and inverse norm are bounded by $e^{\eta|\I|}$.
The estimates (H1-H2) follow from \eqref{eq:total-holonomy} and
Theorem~\ref{theorem-main1general}(H1-H2).
The same gauge bounds and Theorem~\ref{theorem-main1general}(H3)
give the estimate in (H3) for $D^{E,d}_{\epsilon,d_*}-D^E_{\epsilon,d_*}$.
The energy derivative bound for the gauge in Lemma~\ref{lem:localR_t}
and the analyticity of $\mc M^E$ near $E_{\mathbf k}$ likewise give
\[
\|D^E_{\epsilon,d_*}-D^{E'}_{\epsilon,d_*}\|_{h_0}
\leq K|E-E'|,\qquad E,E'\in\I_{d_*},
\]
by integration across the finitely many mesh points.
Since the matrices and their inverses are uniformly bounded on
$\I_{d_*}$, division by the compatible square roots transfers these
two estimates directly to the real centers, proving (H3-H4)
with constants independent of $d,d_*$ and $\epsilon$.

Finally, we record the corresponding phase estimates.
At $E_{\mathbf k}$, choose the same limiting phase for all
$\epsilon,d_*$, since the holonomy is the identity there,
and extend this choice continuously in energy.
The extra phase-strip margin in the proof of
Theorem~\ref{theorem-main1general} bounds this initial phase in
the $h_0$-norm. Continuing the logarithm through finitely many
energy neighborhoods gives a bounded phase on $\I$ for fixed
$\epsilon,d_*$.
The common gauge cancels in the determinant ratio:
\[
\frac{\det D^{E,d}_{\epsilon,d_*}}{\det D^E_{\epsilon,d_*}}
=\frac{\det\mc M^{E,d}}{\det\mc M^E}.
\]
Thus, for compatible finite-dimensional phases,
\[
\rho^{E,d}_{\epsilon,d_*}-\rho^E_{\epsilon,d_*}
=\frac{1}{4\pi i}\operatorname{Log}
\frac{\det\mc M^{E,d}}{\det\mc M^E},
\]
where $\operatorname{Log}$ is the branch near $1$ with
$\operatorname{Log}1=0$. The bounds on $\mc M^E$ and its inverse,
together with Theorem~\ref{theorem-main1general}(H3), give
\[
\|\rho^{E,d}_{\epsilon,d_*}-\rho^E_{\epsilon,d_*}\|_{\I,h_0}
\leq K\|\mc M^{E,d}-\mc M^E\|_{\I,h_0}
\longrightarrow0.
\]
Here $K$ is independent of $d,d_*$ and $\epsilon$.
On $\I_{d_*}$, the preceding Lipschitz bound for $D^E_{\epsilon,d_*}$
also makes $\det D^E_{\epsilon,d_*}/\det D^{E_{\mathbf k}}_{\epsilon,d_*}$
uniformly close to $1$. With the same branch of $\operatorname{Log}$,
\[
\begin{aligned}
\rho^E_{\epsilon,d_*}-\rho^{E_{\mathbf k}}_{\epsilon,d_*}
&=\frac{1}{4\pi i}\operatorname{Log}
\frac{\det D^E_{\epsilon,d_*}}{\det D^{E_{\mathbf k}}_{\epsilon,d_*}},\\
\|\rho^E_{\epsilon,d_*}-\rho^{E_{\mathbf k}}_{\epsilon,d_*}\|_{h_0}
&\leq K|E-E_{\mathbf k}|\leq Kr_{d_*}.
\end{aligned}
\]
The common initial phase satisfies
$\|\rho^{E_{\mathbf k}}_{\epsilon,d_*}\|_{h_0}\leq K_0$,
with $K_0$ independent of $\epsilon,d_*$.
The triangle inequality therefore gives
\[
\|\rho^E_{\epsilon,d_*}\|_{\I_{d_*},h_0}\leq K_0+Kr_{d_*}.
\]
Combining this with the phase convergence estimate yields
\[
\begin{aligned}
&\|\rho^{E,d}_{\epsilon,d_*}\|_{\I_{d_*},h_0}
+\|\rho^E_{\epsilon,d_*}\|_{\I_{d_*},h_0}\\
&\quad\leq
\|\rho^{E,d}_{\epsilon,d_*}-\rho^E_{\epsilon,d_*}\|_{\I_{d_*},h_0}
+2\|\rho^E_{\epsilon,d_*}\|_{\I_{d_*},h_0}\\
&\quad\leq K\|\mc M^{E,d}-\mc M^E\|_{\I,h_0}
+2K_0+2Kr_{d_*}\leq K,
\end{aligned}
\]
after enlarging $K$, for $d\geq d_*$ and sufficiently large $d_*$.
Here we use $r_{d_*}\to0$ and convergence of $\mc M^{E,d}$ uniformly
over the tail $d\geq d_*$. The constant $K$ is independent of
$d,d_*$ and $\epsilon$, without further reducing the strip $\T_{h_0}$.
\end{proof}

\smallskip

 Fix $d_*$ arbitrarily and take $0<\epsilon<\epsilon_*(\I,d_*)$ small enough so that Proposition \ref{prop:10.3} applies.
By Aubry duality, the IDS of the finite-range dual operator $L_{v_d,\alpha,x}$ coincides with $N_{v_d,\alpha}$. Consequently, we will also view $N_{v,\alpha}$ as the IDS of the infinite-range dual operator $L_{v,\alpha,x}$. 
According to \cite[Appendix A]{LW}, the rotation numbers of the dual and conjugated cocycles are related by the topological degree of the conjugacy. Because the hyperbolic part $\widehat{H}^{E,d}(x)$ and the phase term $e^{2\pi i\rho_{\epsilon,d_*}^{E,d}(x)}$ do not contribute, evaluating at $E_{\mathbf{k}}^d$ yields: $$\rho(\alpha,S_{E_{\mathbf{k}}^d}^{v_d}) = 1-N_{v_d}(E_{\mathbf{k}}^d) = \rho(\alpha,L^{E_{\mathbf{k}}^d,d}) = \rho(\alpha,\widehat C^{E_{\mathbf{k}}^d,d}_{\epsilon,d_*})+\frac{\deg B^{E_{\mathbf{k}}^d,d}_{\epsilon,d_*}\cdot \alpha}{2} \pmod{\mathbb{Z}}.$$ By Lemma \ref{lem:E_kd}, $E_{\mathbf{k}}^d \to E_{\mathbf{k}}$ as $d \to \infty$. Combining this with the uniform convergence of the cocycles established in Proposition \ref{prop:10.3}(H3), we obtain $\lim\limits_{d\to\infty} \rho(\alpha,\widehat C^{E_{\mathbf{k}}^d,d}_{\epsilon,d_*}) = \rho(\alpha,\widehat C^{E_{\mathbf{k}}}_{\epsilon,d_*}).$ Because the mapping $E \mapsto B^{E,d}_{\epsilon,d_*}$ is continuous on $\I$, its topological degree is locally constant with respect to $E$ and can be denoted simply by $n_B^{\epsilon,d_*}(d)$. Taking the limit as $d \to \infty$ in the relations above, the convergence of the limit boundaries forces the sequence $n_B^{\epsilon,d_*}(d) \cdot \alpha \pmod{\mathbb{Z}}$ to converge to a constant $c_* \in [0,1)$. This directly implies:
$$\rho(\alpha,S_E^v) \equiv \rho(\alpha,\widehat C^{E}_{\epsilon,d_*})+c_*(\epsilon,d_*) \mod \Z.$$
The following lemma establishes the rigidity of this limiting shift $c_*(\epsilon,d_*)$: the presence of a spectral gap forces it to be an exact integer multiple of the frequency $\alpha$.

\begin{lemma}\label{lem:rotn}
    Suppose $\I$ contains at least one spectral gap, then there exists an integer $n_*(\epsilon,d_*)$ such that for any $E\in \I$, $$\rho(\alpha,\widehat C^{E}_{\epsilon,d_*})\equiv\rho(\alpha,S_E^v)+n_*(\epsilon,d_*)\alpha \mod \Z.$$
For all sufficiently large $d$, the same integer satisfies
\begin{equation}\label{eq:rotn-finite}
\rho(\alpha,\widehat C^{E,d}_{\epsilon,d_*})
\equiv\rho(\alpha,S_E^{v_d})+n_*(\epsilon,d_*)\alpha
\pmod{\mathbb Z},\qquad E\in\I.
\end{equation}
The truncation threshold can be chosen independently of the admissible
reference dimension and mesh.
\end{lemma}
\begin{proof}
By hypothesis, the interval $\I$ contains a spectral gap. Let us fix an energy $E_0$ strictly within this gap. The Schrödinger cocycle $(\alpha,S^v_{E_0})$ is therefore uniformly hyperbolic. By the robustness of uniform hyperbolicity, the approximating cocycle $(\alpha,S^{v_d}_{E_0})$ is also uniformly hyperbolic for all sufficiently large $d$, which implies $E_0 \notin \Sigma_{v_d,\alpha}$.
Consequently, Proposition \ref{prop:10.3} guarantees that for any $|\varepsilon| < h_0$, the complexified Lyapunov exponent of the transformed cocycle satisfies$$L\big(\alpha,\widehat C^{E_0}_{\epsilon,d_*}(\cdot + i\varepsilon)\big) = \lim_{d\to\infty} L\big(\alpha,\widehat C^{E_0,d}_{\epsilon,d_*}(\cdot + i\varepsilon)\big) = \lim_{d\to\infty} L\big(\alpha,\widehat C^{E_0,d}_{\epsilon,d_*}\big).$$ Furthermore, applying the quantitative global theory \cite[Theorem 1.2]{GJYZ}, we obtain the strict lower bound $\lim_{d\to\infty} L\big(\alpha,\widehat C^{E_0,d}_{\epsilon,d_*}\big) 
> 0.$ Since its complexified Lyapunov exponent is constant and strictly positive on the strip $|\varepsilon| < h_0$, Avila's global theory \cite{avila0} ensures that the limiting cocycle $(\alpha,\widehat C^{E_0}_{\epsilon,d_*})$ is regular and uniformly hyperbolic.

Because both cocycles evaluated at $E_0$ are uniformly hyperbolic, their fibered rotation numbers are locked to the frequency module. Therefore, there exist integers $m$ and $m_*$ such that $\rho(\alpha,S_{E_0}^v) \equiv m\alpha \pmod{\mathbb{Z}}$ and $\rho(\alpha,\widehat C^{E_0}_{\epsilon,d_*}) \equiv m_*(\epsilon,d_*)\alpha \pmod{\mathbb{Z}}$. Substituting these identities into our general relation yields $m\alpha \equiv m_*(\epsilon,d_*)\alpha + c_*(\epsilon,d_*) \pmod{\mathbb{Z}}.$ Setting $n_*(\epsilon,d_*) = m_*(\epsilon,d_*) - m$, we conclude that $c_*(\epsilon,d_*) \equiv -n_*(\epsilon,d_*)\alpha \pmod{\mathbb{Z}}$. Substituting this back into the relation for an arbitrary $E \in \I$ yields $\rho(\alpha,\widehat C^{E}_{\epsilon,d_*}) \equiv \rho(\alpha,S_E^v) + n_*(\epsilon,d_*)\alpha \pmod{\mathbb{Z}}$.

For the finite-dimensional assertion, first fix the reference dimension
and mesh. At $E_0$, uniform convergence and openness of uniform
hyperbolicity place $\widehat C^{E_0,d}_{\epsilon,d_*}$ and
$\widehat C^{E_0}_{\epsilon,d_*}$ in the same uniformly hyperbolic
neighborhood for all sufficiently large $d$. Their rotation numbers
therefore coincide. Likewise, the potential segment joining $v_d$ to
$v$ remains gapped at $E_0$ for sufficiently large $d$, so
$\rho(\alpha,S_{E_0}^{v_d})=\rho(\alpha,S_{E_0}^{v})$ modulo $\mathbb Z$.
The finite-dimensional rotation-number relation preceding the lemma
has an energy-independent shift. Evaluating it at $E_0$ identifies
this shift as $n_*(\epsilon,d_*)\alpha$, proving \eqref{eq:rotn-finite}
on all of $\I$.
Finally, this stabilization can be checked in the common center
coordinates before applying the holonomy gauge in \eqref{eq:limitCd*}.
For each reference dimension and mesh, the same gauge conjugates
both the finite and limiting centers; hence its degree contributes
identically to their rotation numbers. The equality at $E_0$, and
thus the truncation threshold, is independent of these later choices.
\end{proof}

As a consequence of Lemma \ref{lem:rotn},  ARC implies that 
$(\alpha,\widehat C^{E_\mathbf{k}}_{\epsilon,d_*})$ (resp. $(\alpha,\widehat C^{E_\mathbf{k}^d,d}_{\epsilon,d_*})$) is almost reducible to a parabolic cocycle:
To establish a unified almost reducibility estimate that bridges the Diophantine ($\beta(\alpha)=0$) and Liouvillean ($\beta(\alpha)>0$) regimes, we naturally parameterize both cases via the sequence of best rational approximations, letting $(q_n)_{n \in \mathbb{N}}$ be the sequence of continued fraction denominators of the frequency $\alpha$. For any given scale $q_n$, we define the truncation degree $d_n$ as the maximal integer $d_n \le q_n$ such that the Fourier coefficient $\hat{v}_{d_n} \neq 0$. As all Fourier coefficients in the intermediate range $(d_n, q_n]$ are exactly zero,  one conclude $\|v_{d_n} - v\|_\gamma \le e^{-\kappa q_n}$ for some uniform $\kappa > 0$. This scale-matching allows us to synthesize the two cases into a single, unified lemma.

\begin{lemma}\label{lem:unified-ar}Let $\alpha \in \mathbb{R} \setminus \mathbb{Q}$. There exist $0 < \gamma < h_0$, a sequence of conjugacies $\bar B_{q_n}^{\epsilon,d_*}(x) \in C^\omega_\gamma(\mathbb{T},\mathrm{PSL}(2,\mathbb{R}))$ with $\deg \bar B_{q_n}^{\epsilon,d_*} = 2(n_*(\epsilon,d_*)-\mathbf{k})$ and constants $|c_{q_n}| > e^{-o(q_n)}$ such that for  sufficiently large $n$, we have for all $d\geq\max\{d_n,d_*\}$,
\begin{equation}\label{eq:unified-ar}(-1)^{\mathbf{l}} \bar B_{q_n}^{\epsilon,d_*}(x+\alpha)^{-1}\widehat{C}^{E_\mathbf{k}^{d},d}_{\epsilon,d_*}(x)\bar B_{q_n}^{\epsilon,d_*}(x)=I_2+c_{q_n}\mathfrak{L} + F_{q_n}^{\epsilon,d_*}(x) + Q^{d}_{\epsilon,d_*}(x)\end{equation}
with the following quantitative estimates hold for some uniform $\kappa > 0$:
\begin{equation}\label{estis}
    \|\bar B_{q_n}^{\epsilon,d_*}\|_\gamma \le e^{o(q_n)}, \qquad \|F_{q_n}^{\epsilon,d_*}\|_\gamma \le e^{-\kappa q_n}, \qquad \|Q^{d}_{\epsilon,d_*}\|_\gamma \le \min \{ e^{-\kappa q_n}, e^{-\kappa  d}\}
\end{equation}\end{lemma}
\begin{remark}
    If \(\beta(\alpha) > 0\), we will select the subsequence \(p_n/q_n\) satisfies \(q_{n+1} > e^{(\beta - o(1))q_n}\).
\end{remark}
\begin{proof}
By Proposition \ref{prop:10.3}, the limiting cocycle $(\alpha,\widehat C^{E_{\mathbf{k}}}_{\epsilon,d_*})$ is subcritical, and Lemma~\ref{lem:rotn} implies that $\rho(\alpha,\widehat C^{E_{\mathbf{k}}}_{\epsilon,d_*})=(n_*(\epsilon,d_*)-\mathbf{k})\alpha \mod \Z$.
At $E_{\mathbf k}$ the holonomy is the identity. Hence \eqref{eq:limitCd*} gives
\[
\widehat C^{E_{\mathbf{k}}}_{\epsilon,d_*}(x)
=e^{-2\pi i\rho^{E_{\mathbf k}}_{\epsilon,d_*}(x)}\mc M^{E_{\mathbf k}}(x).
\]
By the common phase choice in Proposition~\ref{prop:10.3}, this real center
cocycle is independent of $\epsilon,d_*$ and of the truncation dimension $d$.
Thus if $\beta(\alpha) = 0$, 
by Theorem \ref{prop-main-dio} applied to this fixed real center cocycle, there exists 
 $\bar B^{\epsilon,d_*}\in C^\omega_\gamma(2\T,\mathrm{SL}(2,\R))$ with $\deg \bar B^{\epsilon,d_*}=2(n_*(\epsilon,d_*)-\mathbf{k})$ and a constant $c\in \R$ such that 
\[
(-1)^{\mathbf{l}}\bar B^{\epsilon,d_*}(x+\alpha)^{-1} \widehat C^{E_\mathbf{k}}_{\epsilon,d_*}(x) \bar B^{\epsilon,d_*}(x)= I_2+c\mathfrak{L}.
\]
Take $q_n=q_n(\mathbf{k},n_*)$ sufficiently large such that $\|\bar B^{\epsilon,d_*}\|_\gamma \le e^{o(q_n)}$, then we set  
 $F_{q_n}^{\epsilon,d_*}(x) \equiv 0$, $\bar B_{q_n}^{\epsilon,d_*} \equiv \bar B^{\epsilon,d_*}$ and $c_{q_n} \equiv c$. The conjugacy and $c$ are chosen once for the fixed real center, so the estimates are independent of $\epsilon,d_*$.

If $\beta(\alpha) > 0$, we replace Theorem~\ref{prop-main-dio} by the following quantitative version of ARC:  \begin{theorem}[\cite{avila2016dry}]\label{prop-main-liouville}
Let $\alpha\in\R\backslash \Q$ with $\beta(\alpha)>0$ and $k_0\in\Z$. Suppose that $(\alpha, A)$ is subcritical and $\mathbf{l}=\rho(\alpha,A)-k_0\alpha \in \Z$. Then  there exists $B_{q_n}(x)\in C^\omega_\gamma(\T,\mathrm{PSL}(2,\R))$ with $\deg B_{q_n}=2k_0$ and $c_{q_n}\in\R$ such that 
\begin{eqnarray}
(-1)^{\mathbf{l}} B_{q_n}(x+\alpha)^{-1}A(x)B_{q_n}(x)=I_2+c_{q_n}\mathfrak{L} + F_{q_n}(x),
\end{eqnarray}
with the following quantitative estimates hold for some uniform $\kappa>0$:
\[
    \|B_{q_n}\|_\gamma \leq e^{o(q_n)}, \qquad 
    \|F_{q_n}\|_\gamma \leq e^{-\kappa q_n}.
\]
\end{theorem}
Applying this theorem to the same fixed real center, we likewise choose
the conjugacies, parabolic coefficients and remainders independently of
$\epsilon,d_*$. In both cases, we retain the superscripts $\epsilon,d_*$
in the conjugacies and remainders for consistency with the subsequent notation.
We then focus on the approximating cocycle $(\alpha,\widehat{C}^{E_{\mathbf{k}}^d,d}_{\epsilon,d_*})$, which yields \eqref{eq:unified-ar}. 
By Proposition~\ref{prop:10.3}, $\widehat C^E_{\epsilon,d_*}$ is piecewise $C^\omega$ in $E\in\I$, and estimate (H3) holds. Combining this with Lemma~\ref{lem:E_kd} and Proposition~\ref{prop:10.3}(H4), we obtain there is a uniform constant $C>0$  such that
\[
\|Q^d_{\epsilon,d_*}\|_\gamma < \|\bar B_{q_n}^{\epsilon,d_*}\|_\gamma^2\bigl(\|\widehat C^{E_{\mathbf{k}}^d,d}_{\epsilon,d_*}- \widehat C^{E_{\mathbf{k}}^d}_{\epsilon,d_*}\|_\gamma + \|\widehat C^{E_{\mathbf{k}}^d}_{\epsilon,d_*}-\widehat C^{E_{\mathbf{k}}}_{\epsilon,d_*}\|_\gamma\bigr) < C\|\bar B_{q_n}^{\epsilon,d_*}\|_\gamma^2 \|v_d-v\|_{h_1-\gamma}.
\]
By our choice of $d_n$, we have for all $d\geq\max\{d_n,d_*\}$, $C\|v_d-v\|_{h_1-\gamma} <\min\{ e^{-\kappa' q_n}, e^{-\kappa 'd}\}$. The estimates follows by picking a smaller common $\kappa$.

We now arrive to show $|c_{q_n}| > e^{-o(q_n)}$ uniformly, we assume for contradiction that $|c_{q_n}| \le e^{-\delta q_n}$ for some $\delta > 0$.  Note that the center bundle frame $U^{E,d}_{\epsilon,d_*}(x)$ obeys
$$L^{E,d}(x) U^{E,d}_{\epsilon,d_*}(x) = U^{E,d}_{\epsilon,d_*}(x+\alpha) e^{2\pi i\rho_{\epsilon,d_*}^{E,d}(x)}\widehat{C}^{E,d}_{\epsilon,d_*}(x).$$
Taking $E = E_\mathbf{k}^{d}$ and $d \geq\max\{d_n,d_*\}$, we define the conjugated extended frame:$$\widetilde U_{d}(x) := U^{E_\mathbf{k}^{d},d}_{\epsilon,d_*}(x) \bar B_{q_n}^{\epsilon,d_*}(x).$$Substituting the normal form \eqref{eq:unified-ar} into the block-diagonalization yields:$$L^{E_\mathbf{k}^{d},d}(x) \widetilde U_{d}^{\epsilon,d_*}(x) = \widetilde U_{d}^{\epsilon,d_*}(x+\alpha) e^{2\pi i \rho_{\epsilon,d_*}^{E_\mathbf{k}^{d},d}(\cdot)} \left( \pm  I_2 \pm \big(c_{q_n}\mathfrak{L} + F_{q_n}^{\epsilon,d_*}(x) + Q^{d}_{\epsilon,d_*}(x)\big) \right).$$  The discrepancy is entirely governed by
\begin{equation}\label{eq:err-matrix}
\pm  \widetilde U_{d}^{\epsilon,d_*}(x+\alpha) e^{2\pi i \rho_{\epsilon,d_*}^{E_\mathbf{k}^{d},d}(\cdot)} \big(c_{q_n}\mathfrak{L} + F_{q_n}^{\epsilon,d_*}(x) + Q^{d}_{\epsilon,d_*}(x)\big).\end{equation}
To evaluate \eqref{eq:err-matrix} against Proposition \ref{dim-duality}, we calculate the weighted spatial-analytic norm $\tnorm{\cdot}_\gamma$.
Since $E_{\mathbf k}^d\in\I_{d_*}$, the phase bound in
Proposition~\ref{prop:10.3} gives
\[
\|\rho_{\epsilon,d_*}^{E_{\mathbf k}^d,d}\|_\gamma\leq K,
\]
uniformly in $d,d_*$ and $\epsilon$.
In particular, the scalar factor in \eqref{eq:err-matrix} is uniformly bounded.
Moreover,
Proposition \ref{prop:10.3} (H1) ensures that the spatial components of the original frame grow sub-exponentially: $\|U_{\epsilon,d_*}^{E_\mathbf{k}^{d},d}(\cdot,k)\|_\gamma \le e^{o(|k|)}$. For any fixed weight parameter $\xi > 0$, the exponential decay  dominates the sub-exponential growth. Therefore, we obtain a finite bound:$$\tnorm{U_{\epsilon,d_*}^{E_\mathbf{k}^{d},d}}_\gamma = \sum_{k=-d}^{d-1} e^{-\xi|k|} \|U_{\epsilon,d_*}^{E_\mathbf{k}^{d},d}(k, \cdot)\|_\gamma \le \sum_{k \in \mathbb{Z}} e^{-\xi|k| + o(|k|)} =: C_\xi < \infty.$$ Here, $C_\xi$ is a global constant  independent of  $d$. As a consequence: $$\tnorm{\widetilde U^{\epsilon,d_*}_{d}}_\gamma \le \tnorm{U_{\epsilon,d_*}^{E_\mathbf{k}^{d},d}}_\gamma \|\bar B_{q_n}^{\epsilon,d_*}\|_\gamma \le C_\xi \, e^{o(q_n)} \le e^{o(q_n)}.$$
Combining the hypothesis $|c_{q_n}| \le e^{-\delta q_n}$ alongside 
the quantitative bounds \eqref{estis} it yields, we have
$$\tnorm{L^{E_\mathbf{k}^{d},d}(\cdot) \widetilde U_{d}^{\epsilon,d_*}(\cdot) - \widetilde U^{\epsilon,d_*}_{d}(\cdot+\alpha) \big(\pm e^{2\pi i \rho_{\epsilon,d_*}^{E_\mathbf{k}^{d},d}(\cdot)}  I_2\big)}_\gamma  \le e^{o(q_n)} e^{-\delta' q_n} \le e^{-\frac{\delta'}{2} q_n},$$
for some $\delta'$, this directly contradicts Proposition \ref{dim-duality}. The result follows. \end{proof}

\begin{remark}

After a constant diagonal conjugacy, we may assume throughout that
$e^{-o(q_n)}<|c_{q_n}|\leq 1$, retaining the notation above.
Indeed, the normal form and \eqref{estis} give $|c_{q_n}|\leq e^{o(q_n)}$.
Set $a_{q_n}=\max\{1,\sqrt{|c_{q_n}|}\}$ and
$D_{q_n}=\operatorname{diag}(a_{q_n},a_{q_n}^{-1})$.
Replacing $\bar B_{q_n}$ by $\bar B_{q_n}D_{q_n}$ replaces $c_{q_n}$ by
$c_{q_n}/a_{q_n}^2$ and conjugates both remainders by $D_{q_n}$.
Since $\|D_{q_n}^{\pm1}\|\leq e^{o(q_n)}$, all estimates in \eqref{estis}
remain valid after decreasing $\kappa$ if necessary. This conjugacy is
independent of $E$ and has degree zero.

\end{remark}

\subsection{Derivative estimates for projective action}

Now we give finer estimates concerning derivative estimates for projective action, and distinguish the proof into two steps.\\

\smallskip
\noindent \textbf{Step 1: Finer monotonicity of $L^{E,d}$.}
Let $u = (u_1, u_2)^\top \in \mathbb{R}^2$ be a unit vector, and define its natural embedding into $\mathbb{R}^{2d}$ by $$\bar{u} := (0, \dots, 0, u_1, 0, \dots, 0, u_2)^\top \in \mathbb{R}^{2d},$$ where $u_1$ and $u_2$ occupy the $d$-th and $2d$-th entries, respectively. Next, we define the following transformed vectors:$$\bar{v} := B^{E,d}(x)\bar{u}, \quad \bar{w} := L^{E,d}(x)B^{E,d}(x)\bar{u}, \quad \text{and} \quad \bar{x} := B^{E,d}(x+\alpha)^{-1}L^{E,d}(x)B^{E,d}(x)\bar{u}.$$By construction, the vectors $\bar{v}$ and $\bar{w}$ are isotropic with respect to $(\mathbb{C}^{2d},\psi_d)$, while $\bar{x}$ is isotropic with respect to $(\mathbb{C}^{2d},\omega_d)$. Recall that $$\Psi_{L^{E,d}}(\bar v, x) = \psi_d\bigl(L^{E,d}(x)\bar{v}, (\partial_E L^{E,d})(x)\bar{v}\bigr),$$ where $\psi_d$ is the symplectic form introduced in Section \ref{sec:symp-schro}, and $\omega_d$ denotes the standard Hermitian symplectic form on $\mathbb{C}^{2d}$.

\begin{lemma}\label{lem:subexplowb}
  There exists \(D(\I)>0\), independent of \(d\), such that for every $E\in\I_{d_*}$ and $d\geq d_*$,
  \begin{equation}\label{eq:upper bound}
    \sup_{\substack{\|u\|=1\\ x\in\mathbb{T}}}\Psi_{L^{E,d}}(\bar v, x)=\sup_{\substack{\|u\|=1\\ x\in\mathbb{T}}}\psi_d\bigl(L^{E,d}(x)\bar{v},\, (\partial_E L^{E,d})(x)\bar{v}\bigr)\le D .
  \end{equation}
  Furthermore, for $E\in\I_{d_*}$ and $d\geq d_*$, with $d_*$ sufficiently large, then
  \begin{equation}\label{eq:lower bound}
    \inf_{\substack{\|u\|=1\\ x\in\mathbb{T}}}\Psi_{L^{E,d}_{2d}}(\bar v, x)=\inf_{\substack{\|u\|=1\\ x\in\mathbb{T}}}
    \psi_d\bigl(L^{E,d}_{2d}(x)\bar{v},\, (\partial_E L^{E,d}_{2d})(x)\bar{v}\bigr)\ge e^{-o(d)} .
  \end{equation}
\end{lemma}

\begin{proof}
  From the Discrete Hellmann-Feynman Identity (Lemma~\ref{lem:pre-mono}), set \(Y_k = L_k^{E,d}(x)\bar v\). Then
 $ \Psi_{L_k^{E,d}}(\bar v, x) = \sum_{j=0}^{k-1} |(Y_j)_d|^2,
  $  and in particular \(\sum_{j=0}^{2d-1} |(Y_j)_d|^2 = \|L_d^{E,d}(x)\bar v\|^2\).
  For \(k=1\) and \(k=2d\),
  \[
  \Psi_{L^{E,d}}(\bar v, x) = |B^{E,d}(x,0)\bar u|^2, \qquad
  \Psi_{L_{2d}^{E,d}}(\bar v, x) = \|L_d^{E,d}(x)\bar v\|^2.
  \]
  The upper bound~\eqref{eq:upper bound} then follows from (H1) in Proposition \ref{prop:10.3}.

  For the lower bound \eqref{eq:lower bound}, using the block-diagonal conjugation~\eqref{block-p-analytic}, we can write
  \[
  L_d^{E,d}(x)\bar v = \exp\left(2\pi i \sum_{j=0}^{d-1} \rho^{E,d}(x+j\alpha)\right)\cdot U^{E,d}(x+d\alpha) \, \widehat C_d^{E,d}(x) \, u.
  \]
  Consequently, the symplectic form evaluates to
  \[
  \Psi_{L_{2d}^{E,d}}(\bar v, x)
  = \bigl\|U^{E,d}(x+d\alpha) \widehat C_d^{E,d}(x) u\bigr\|^2
  = u^* M(x) u,
  \]
  where \(M(x) = \widehat C_d^{E,d}(x)^* G(x+d\alpha) \widehat C_d^{E,d}(x)\) and \(G(x) = U^{E,d}(x)^* U^{E,d}(x)\).

To bound the smallest eigenvalue of \(M(x)\) from below, we denote the columns of \(U^{E,d}(x)\) by \(u_+,u_-\).
Applying an orthogonal decomposition to $u_-$, we write \(u_- = \lambda u_+ + w\) with \(w \perp u_+\), which implies \(\det G(x) = \|u_+\|^2 \|w\|^2\).
  Since the frame is symplectic, we have \(u_+^* S_d u_- = -1\). Substituting the orthogonal decomposition into this relation yields
  \(1 = |u_+^* S_d w| \le \|u_+\| \|S_d\|_{\text{op}} \|w\|\).
  Exponential decay of the hopping terms yields \(\|S_d\|_{\text{op}} \le C\). This immediately implies that 
  \(\det M(x) = \det G(x) \ge C^{-2} > 0\).

For $E\in\I_{d_*}$, (H1) of Proposition~\ref{prop:10.3} gives
$\|G(x)\|_{\text{op}}\leq e^{o(d)}$.
The limiting center $(\alpha,\widehat C^{E_{\mathbf k}}_{\epsilon,d_*})$
is subcritical and independent of $\epsilon,d_*$.
By (H3-H4) of Proposition~\ref{prop:10.3},
\[
\|\widehat C^{E,d}_{\epsilon,d_*}
-\widehat C^{E_{\mathbf k}}_{\epsilon,d_*}\|_{h_0}
\leq K\|v_d-v\|_{h_1-\gamma}+Kr_{d_*}\longrightarrow0
\]
uniformly for $E\in\I_{d_*}$ and $d\geq d_*$ as $d_*\to\infty$.
The uniform subadditive estimate and continuity of a fixed finite
iterate therefore give $\|\widehat C_d^{E,d}\|_0\leq e^{o(d)}$.
Here the exponent loss is made arbitrarily small by first increasing
$d_*$, uniformly for $d\geq d_*$ and all admissible meshes.
  Thus the largest eigenvalue of $M$ is bounded by  \(\lambda_{\max}(M) \le \|\widehat C_d^{E,d}\|^2 \|G\|_{\text{op}} \le e^{o(d)}\).
  Because the determinant is uniformly bounded away from zero, we obtain that for the smallest eigenvalue of $M$
  \[
  \lambda_{\min}(M) = \frac{\det M(x)}{\lambda_{\max}(M)} \ge \frac{C^{-2}}{e^{o(d)}} \ge e^{-o(d)}.
  \]
Since \(u^* M u \ge \lambda_{\min}(M)\)  for any unit vector \(u\), the infimum in~\eqref{eq:lower bound} follows immediately.
  \qedhere
\end{proof}

\smallskip
\noindent \textbf{Step 2: Finer derivative estimates of the projection.}  Based on Lemma \ref{lem:unified-ar}, we define the following conjugated cocycles for any $E \in \I$:
\[
\begin{aligned}
    M^{E,q}_{\epsilon,d_*}(x)&=\bar B_q^{\epsilon,d_*}(x+\alpha)^{-1} \widehat C^{E}_{\epsilon,d_*}(x) \bar B_q^{\epsilon,d_*}(x)-F_q^{\epsilon,d_*}(x),\\
    M^{E,d,q}_{\epsilon,d_*}(x)&=\bar B_q^{\epsilon,d_*}(x+\alpha)^{-1} \widehat C^{E,d}_{\epsilon,d_*}(x) \bar B_q^{\epsilon,d_*}(x).
\end{aligned}
\]
As the conjugacy $\widehat C^{E,d}_{\epsilon,d_*}(x)$ is continuous and piecewise premonotonic with respect to $E$, this property is inherited by the projective phase action  $\phi[M^{E,d,q}_{\epsilon,d_*}(x)](\Lambda)$. Let $\{i_s(\epsilon,d_*)\}$ denote the finite set of non-differentiability points. The following lemma establishes uniform bounds on these phase derivatives on $\I_{d_*}$, for $d=d_*$, away from these points.

\begin{lemma}\label{lem:derivate-bound}
Let $d=d_*$ be sufficiently large and $0<\epsilon<\epsilon_*(\I,d)$ be sufficiently small. There exists $D=D(\I)>0$, independent of $d$ and $\epsilon$, such that uniformly in $\Lambda$:
\begin{align}
    -\mc D_q\epsilon \le \inf_{E\in \I_{d_*}\setminus \{i_s\} } -\partial_E\phi[M^{E,d,q}_{\epsilon,d_*}(x)](\Lambda) \le \sup_{E\in \I_{d_*}\setminus \{i_s\} } -\partial_E\phi[M^{E,d,q}_{\epsilon,d_*}(x)](\Lambda) \leq \mc D_q,  \label{eq:deriv-bounds-2} 
\end{align}
where 
$\mc D_q := D \|\bar B_q^{\epsilon,d_*}\|_0^2 \le e^{o(q)}$. Furthermore, at points of differentiability, for $E\in\I_{d_*}$, then
\begin{align}
    -\partial_E\phi[(M^{E,d,q}_{\epsilon,d_*})_{2d}(x)](\Lambda) \ge e^{-o(d)}\mc D_q^{-1}.
\end{align}
\end{lemma}
\begin{remark}
    The key observation here is that although the one-step projective action of $M_{\epsilon,d_*}^{E,d,q}(x)$ is not necessarily strictly monotonic, its failure to be monotonic is strictly controlled and could be arbitrarily small. 
\end{remark}

\begin{proof}
At points of differentiability, since $\bar{B}_q^{\epsilon,d_*} \in \mathrm{Sp}(2, \mathbb{R})$ preserves the Hermitian symplectic inner product $\psi=\omega_2$, and $\bar{B}_q^{\epsilon,d_*}$ is locally constant with respect to $E$, 
by Lemma \ref{lem:omega-trace} (more precisely, \eqref{deri-phi-center}), we have uniformly in $\Lambda$:
\begin{align}
 \label{eq:upperb-M}    
    \sup_{E\in \I_{d_*}\setminus \{i_s\} } -\partial_E\phi[M^{E,d,q}_{\epsilon,d_*}(x)](\Lambda) &\leq& \frac{1}{\pi}\max_{s\in\{1,-1\}}\left\{\|\bar{B}_q^{\epsilon,d_*}\|_0^{2s}\|\widehat C^{E,d}_{\epsilon,d_*}\|_0^{2s} \sup_{\substack{\|v\|=1,\\ x\in\T}} \Psi_{\widehat C^{E,d}_{\epsilon,d_*}}(v,x) \right\}, \\
\label{eq:lowerb-M}
    \inf_{E\in \I_{d_*}\setminus \{i_s\} } -\partial_E\phi[M^{E,d,q}_{\epsilon,d_*}(x)](\Lambda)&\geq& \frac{1}{\pi}\min_{s\in\{1,-1\}}\left\{\|\bar{B}_q^{\epsilon,d_*}\|_0^{2s}\|\widehat C^{E,d}_{\epsilon,d_*}\|_0^{2s} \inf_{\substack{\|v\|=1,\\ x\in\T}} \Psi_{\widehat C^{E,d}_{\epsilon,d_*}}(v,x) \right\}.
\end{align}
Combining Lemma \ref{lem:subexplowb}, Lemma \ref{lem:LUbd}
(more precisely, \eqref{eq:4.28}), and Proposition \ref{prop:10.3}(H5),
with \eqref{eq:upperb-M} and \eqref{eq:lowerb-M}, gives
\eqref{eq:deriv-bounds-2} on $\I_{d_*}$.
The constant $D$ is independent of $d$ and $\epsilon$, by the uniform
bounds in Proposition \ref{prop:10.3} on this interval.

For the $2d$-step lower bound, apply the chain identity
\eqref{eq:Psi-C-chain} to the $2d$-step intertwining relation.
The two endpoint frame terms are bounded in absolute value by
\[
K\epsilon\left(1+\|(\widehat C^{E,d}_{\epsilon,d_*})_{2d}(x)u\|^2\right),
\qquad \|u\|=1,
\]
Here we use the local frame estimate underlying \eqref{eq:4.28}
and its dimension-independent constants on $\I_{d_*}$.
By Lemma \ref{lem:subexplowb}, the main term is bounded below by
$\mu_d=e^{-o(d)}$, uniformly on this interval. For fixed $d=d_*$,
choose $\epsilon$ sufficiently small that
\[
K\epsilon\left(1+\sup_{E\in\I_{d_*},\,x\in\T}
\|(\widehat C^{E,d}_{\epsilon,d_*})_{2d}(x)\|^2\right)<\mu_d/2.
\]
The supremum is bounded uniformly over the admissible meshes at fixed $d$.
Thus the remaining quadratic form is at least $\mu_d/2$.
On the near-spectrum interval, the $2d$-step center norm is $e^{o(d)}$;
applying the projective derivative formula gives the lower bound
$e^{-o(d)}$ before conjugation by $\bar B_q$.
Since $\bar B_q$ is independent of $E$, its projective derivative costs
at most the factor $\|\bar B_q\|_0^2$, giving the asserted estimate.
Continuity and piecewise analyticity allow integration across the finitely
many mesh points.
\end{proof}

\subsection{Unified Cone Field Argument via Pre-monotonicity}
We now develop a quantitative cone argument to prove uniform lower bounds for the spectral gaps of $L^{E,d}$, completing the proof of Theorem~\ref{main-thm-analytic}.   Unlike the classical finite-range argument (Lemma \ref{lem:conefield}), where one applies the monotone iterate $\widehat C^E_{2d}$ directly to obtain an invariant cone, quantitative analytic approximation requires controlling the truncation error at every intermediate step. Since the one-step cocycle is only \emph{pre}-monotonic, the projective phase may temporarily drift in the wrong direction. The key idea is to quantitatively bound this transient excursion, ensuring that the orbit remains inside the cone until the monotonicity of the $2d$-th iterate restores cone invariance.

\begin{proposition}\label{prop:conefield-analytic}
Let \(\alpha\in\R\setminus\Q\). 
Fix a sufficiently small \(\varepsilon_0\in(0,\tfrac14)\). Then for sufficiently large $n$,  the spectral gap $G_{\mathbf{k}}(v_{d_n},\alpha)$ satisfies the estimate
\begin{equation}\label{poly-es}
    |G_{\mathbf{k}}(v_{d_n},\alpha)| > \min\left\{\frac{|c_{q_n}| \pi\varepsilon_0^2}{4\mathcal{D}_{q_n}},\frac{r_{d_n}}{2}\right\}.
\end{equation}
As a consequence,  \begin{equation}\label{ana-es}
    |G_{\mathbf{k}}(v,\alpha)| > \frac12\min\left\{\frac{|c_{q_n}| \pi\varepsilon_0^2}{4\mathcal{D}_{q_n}},\frac{r_{d_n}}{2}\right\}.
\end{equation}
\end{proposition}

{ 
\begin{proof}
As \eqref{ana-es} follows directly from \eqref{poly-es} and Lipschitz  continuity with respect to the potential, we only need to prove \eqref{poly-es}.    Fix $q_n$ sufficiently large, determined as in Lemma \ref{lem:unified-ar}. Fix $d = d_n = d_*$ and take $0 < \epsilon < \epsilon_*(\I, d_n)$ sufficiently small. In what follows, we omit $d_*$ and $\epsilon$ from the notation.   For definiteness, assume that $c_{q_n}>e^{-o(q_n)}$. Consider the cone field $\mathcal{C}=\{\Lambda_y : y\in[0,\e_0]\}$ in the projective line, where $\Lambda_y=(\cos(\pi y), \sin(\pi y))^\top$.

We pick the lifts so that $\phi[M^{E_\mathbf{k},q_n}(x)](\Lambda_0)=0$. By the cocycle rule, this implies $\phi[M^{E_\mathbf{k},q_n}_{k}(x)](\Lambda_0)=0$ for all $k\geq 1$, and ensures $\phi[M^{E_\mathbf{k}^{d_n},d_n,q_n}(x)](\Lambda_0)\in (-1/2,1/2)$.

Define the one-step phase map for the limit and the truncated systems, respectively: \[
\begin{aligned}
    f(E,x,\e) &:= \phi[M^{E,q_n}(x)](\Lambda_{\e})+\e,\\
    f^{d_n}(E,x,\e) &:= \phi[M^{E,d_n,q_n}(x)](\Lambda_{\e})+\e.
\end{aligned}
\]
At the precise limit boundary point $E_\mathbf{k}$, an exact computation yields for any $y \in [-\e_0, \e_0]$:
\[
f(E_\mathbf{k}, x, y) = \frac{1}{\pi}\arctan\!\left(\frac{\sin(\pi y)}{\cos(\pi y)+c_{q_n}\sin(\pi y)}\right).
\]
Since $0<c_{q_n}\leq1$, reducing $\e_0$ independently of $n$ if necessary, a Taylor expansion provides uniform quadratic bounds on the parabolic drift. For the positive phase, we obtain
\begin{equation}\label{eq:approxf-1}
    f(E_\mathbf{k},x,\e) \leq \e - \frac{c_{q_n}}{2}\pi \e^2 \qquad \text{for } \e \in [0, \e_0],
\end{equation}
while for the negative phase, we have the corresponding bound:
\begin{equation}\label{eq:approxf-2}
    f(E_\mathbf{k},x,\e) \geq \e - 2c_{q_n}\pi \e^2 \qquad \text{for } \e \in [-\e_0, 0).
\end{equation}
Next, we track the evolution of boundary by considering the iterated phases
\[
\begin{aligned}
\e_0^{d_n}(E,x) &:= \e_0, & \e_{k+1}^{d_n}(E,x) &:= f^{d_n}(E, x+k\alpha, \e_k^{d_n}),\\
\hat\e_0^{d_n}(E,x) &:= 0, & \hat\e_{k+1}^{d_n}(E,x) &:= f^{d_n}(E, x+k\alpha, \hat\e_k^{d_n}).
\end{aligned}
\]
Recall that $M^{E_\mathbf{k}^{d_n},d_n,q_n}(x) =M^{E_\mathbf{k},q_n}(x)+ Q^{d_n}(x)+F_{q_n}(x)$, standard perturbation of the projective action yields:
\[
    f^{d_n}(E_\mathbf{k}^{d_n}, x, \e) = f(E_\mathbf{k},x,\e) 
    + \mathcal{O}(\|Q^{d_n}\|_0+\|F_{q_n}\|_0) = f(E_\mathbf{k},x,\e)+ 
    \mathcal{O}(e^{-\kappa {q_n}}).
\]

With $\delta_0$ defined below, the choice of $r_{d_n}$ ensures that
$(E_\mathbf{k}^{d_n}-\delta_0,E_\mathbf{k}^{d_n})\subset\I_{d_n}$.
Moreover, $r_{d_n}\geq e^{-\sqrt{d_n}}\geq e^{-\sqrt{q_n}}$, so
$\delta_0\geq e^{-o(q_n)}$ and it dominates the exponential errors.
Fix an energy $E \in (E_\mathbf{k}^{d_n}-\delta_0, E_\mathbf{k}^{d_n})$ and denote $\delta := E_\mathbf{k}^{d_n} - E \in (0, \delta_0)$. Integrating the uniform derivative bounds \eqref{eq:deriv-bounds-2} over the energy segment in $\I_{d_n}$ yields the perturbed map estimates:
\begin{equation}\label{eq:perturbed-bounds}
    f^{d_n}(E_\mathbf{k}^{d_n}, x, \e) -\mathcal{D}_{q_n} \epsilon \delta \le f^{d_n}(E,x,\e) \le f^{d_n}(E_\mathbf{k}^{d_n}, x, \e) + \mathcal{D}_{q_n} \delta
\end{equation}
where $\mc D_{q_n}:=D \|\bar B_{q_n}\|^2$. Combining this with the truncation error yields
\begin{equation}\label{eq:perturbed-bounds-2}
    f(E_\mathbf{k},x,\e) - \mathcal{D}_{q_n} \epsilon \delta - C_0 e^{-\kappa q_n} \le f^{d_n}(E,x,\e) \le f(E_\mathbf{k},x,\e) + \mathcal{D}_{q_n} \delta + C_0 e^{-\kappa q_n}.
\end{equation}
where $C_0 > 0$ is a constant. 

Crucially, the uniform control parameter  $\epsilon<\epsilon_*(\I,d_n)$ from Lemma \ref{lem:derivate-bound} can be chosen arbitrarily small. For our purpose, we explicitly fix the parameter
$  \epsilon \leq \frac{2|c_{q_n}|\pi\e_0^2}{2\mathcal{D}_{q_n} \delta_0 (1+2c_{q_n}\pi\e_0)^{2{d_n}+1}}$. With this choice, we establish the following key geometric control: 
\begin{lemma}\label{lem:orbit}
Setting $\delta_0 := \min\left\{\frac{|c_{q_n}| \pi\varepsilon_0^2}{4\mc D_{q_n}},\frac{r_{d_n}}2\right\}$, we have  $-\e_{0} < \e_k^{d_n}(E,x) < \e_{0}$ for all $1 \leq k \leq 2d_n$.   
\end{lemma}
\begin{proof}
We proceed by induction, splitting the finite orbit into two possible regimes. Let $1 \le k_* \le 2{d_n}+1$ be the first index for which $\e_{k_*}^{d_n} < 0$.

\noindent\textbf{Case 1 ($1 \le k < k_*$):}
In this regime, $\e_k^{d_n} \ge 0$. Because $M^{E_\mathbf{k}^{d_n},{d_n},q_n}(x) \in \mathrm{SL}(2, \mathbb{R})$, the  projective map $f^{d_n}(E_\mathbf{k}^{d_n}, x, \cdot)$ is strictly increasing on $[0, \e_0]$. Applying \eqref{eq:perturbed-bounds}, \eqref{eq:perturbed-bounds-2} inductively alongside \eqref{eq:approxf-1} yields
\[
    \e_{k+1}^{d_n} \le f^{d_n}(E_\mathbf{k}^{d_n},x+k\alpha,\e_k) + \mc D_{q_n}\delta \leq f^{d_n}(E_\mathbf{k}^{d_n},x+k\alpha,\e_0) + \mc D_{q_n}\delta \le \e_0 - \frac{c_{q_n}}{2}\pi\e_0^2 + \mc D_{q_n}\delta +C_0e^{-\kappa q_n}.
\]
Because we chose $\delta < \delta_0$, we have $\mc D_{q_n}\delta <\mc D_{q_n}\delta_0 \leq \frac{c_{q_n}}{4}\pi\e_0^2<\frac{c_{q_n}}{2}\pi\e_0^2-C_0 e^{-\kappa q_n}$ (for $n$ sufficiently large). Thus $\e_{k+1}^{d_n} < \e_0$.

\noindent\textbf{Case 2 ($k_* \le k \le 2d_n$):}
For the initial step into the negative region $(-\e_0, 0)$, since $\e_{k_*-1}^{d_n} \ge 0$, the lower bound in \eqref{eq:perturbed-bounds-2} ensures that the excursion is bounded by $\e_{k_*}^{d_n} \ge -\mc{D}_{q_n}\epsilon\delta-C_0e^{-\kappa q_n}$. 
While the orbit remains negative, we couple  \eqref{eq:approxf-2} with \eqref{eq:perturbed-bounds-2} to find
\[
    \e_{k+1}^{d_n} \ge \e_k^{d_n} - 2c_{q_n}\pi(\e_k^{d_n})^2 - \mc D_{q_n}\epsilon\delta-C_0e^{-\kappa q_n}.
\]
We define the non-negative sequence $y_m := -\e_{k_*+m}^{d_n} $, which has the initial condition $y_0 \le \mc D_{q_n}\epsilon\delta+C_0e^{-\kappa q_n} $ and satisfies the recurrence
\[
    y_{m+1} \le y_m + 2c_{q_n}\pi y_m^2 +\mc D_{q_n}\epsilon\delta+C_0e^{-\kappa q_n}.
\]
We msut verify that  $y_m < \e_0$ for all $0 \le m \le 2{d_n} - k_*$. Assuming inductively that $y_m < \e_0$, we bound the quadratic term via $2c_{q_n}\pi y_m^2 < 2c_{q_n}\pi\e_0 y_m$, which yields a linear upper bound:
\[
    y_{m+1} < (1 + 2c_{q_n}\pi\e_0) y_m + \mc D_{q_n}\epsilon\delta+C_0e^{-\kappa q_n}.
\]
Iterating this linear bound from $y_0$ gives
\[
    y_m < (\mc D_{q_n}\epsilon\delta+C_0e^{-\kappa q_n}) \sum_{j=0}^{m} (1+2c_{q_n}\pi\e_0)^j < (\mc D_{q_n}\epsilon\delta_0+C_0e^{-\kappa q_n}) \frac{(1+2c_{q_n}\pi\e_0)^{2{d_n}+1}}{2c_{q_n}\pi\e_0}<\e_0,
\]
where, taking $\e_0$ sufficiently small in terms of $\kappa$, the contribution of $e^{-\kappa q_n}$ after the geometric summation is still exponentially small, since $d_n\leq q_n$ and $e^{-o(q_n)}<c_{q_n}\leq1$.
This proves that the orbit never escapes through the lower boundary ($-\e_0 < \e_k^{d_n}$). If the perturbation eventually pushes the orbit back into the non-negative region $[0, \e_0]$, the uniform bounds of Case 1 simply resume. \end{proof}

By the preceding induction, the finite orbit is strictly confined to $(-\e_0, \e_0)$. This local control allows us to evaluate the $2{d_n}$-step phase map on the boundaries of the cone $\mathcal{C} = \{\Lambda_y : y \in [0, \e_0]\}$. 

At the upper boundary, Lemma \ref{lem:orbit} guarantees that the phase maps strictly downward: $\e_{2{d_n}}^{d_n}(E,x)
< \e_0$. At the lower boundary, the unperturbed map $f(E_\mathbf{k},x,\cdot)$ fixes $0$. Since its phase derivative is $1+O(\e_0)$ near $0$, the same finite-step perturbation estimate gives $\sup_x|\hat\e_{2d_n}^{d_n}(E_\mathbf{k}^{d_n},x)|\leq e^{-\kappa_2q_n}$ for some $0<\kappa_2<\kappa$. By Lemma \ref{lem:derivate-bound}, the $2{d_n}$-step iterate $M^{E,{d_n},q_n}_{2{d_n}}$ is strictly monotonic with respect to $E$. More precisely, combining \eqref{eq:lower bound}, \eqref{eq:LUbd} and \eqref{eq:lowerb-M}, and fixing $0<\kappa_1<\kappa_2$, for sufficiently large $n$ and $E < E_\mathbf{k}^{d_n}- e^{-\kappa_1 q_n}$ the boundary is pushed strictly upward:
\begin{equation}\label{eq:newfordio}
    \hat \e_{2{d_n}}^{d_n}(E,x)> \hat \e_{2{d_n}}^{d_n}(E_\mathbf{k}^{d_n},x)+ e^{-o({d_n})}\|\bar B_{q_n}\|^{-2} e^{-\kappa_1 {q_n}}  > -e^{-\kappa_2 q_n}+ e^{-o({d_n})}e^{-o(q_n)} e^{-\kappa_1 {q_n}}>0.
\end{equation}
Because $M^{E,d_n,q_n}_{2d_n}(x)$ induces orientation-preserving homeomorphism on the projective line, the behavior at the boundaries implies that
$
    M^{E,d_n,q_n}_{2d_n}(x)\mathcal{C} \subset \mathrm{int}(\mathcal{C}).
$
This strict forward invariance of the cone field proves that the perturbed cocycle $(\alpha, M^{E,d_n,q_n})$ is uniformly hyperbolic and its fibered rotation number $\rho(\alpha,M^{E,d_n,q_n})=0$ for all $E \in (E_\mathbf{k}^{d_n} - \delta_0, E_\mathbf{k}^{d_n}-e^{-\kappa_1 q_n})$. Recall $\deg \bar B_{q_n} = 2(n_*-\mathbf{k})$, the cocycle $(\alpha,\widehat C^{E,d_n})$ is also uniformly hyperbolic, with $\rho(\alpha,\widehat C^{E,d_n})=(n_*-\mathbf{k})\alpha\mod\Z$.
By \eqref{eq:rotn-finite}, for this fixed sufficiently large $n$,
\[
\rho(\alpha,S_E^{v_{d_n}})=-\mathbf{k}\alpha\pmod{\mathbb Z},
\qquad E\in(E_\mathbf{k}^{d_n}-\delta_0,
E_\mathbf{k}^{d_n}-e^{-\kappa_1q_n}).
\]
Thus $N_{v_{d_n},\alpha}$ equals its value at $E_\mathbf{k}^{d_n}$
on this interval. By monotonicity of the IDS, it has the same value
throughout $(E_\mathbf{k}^{d_n}-\delta_0,E_\mathbf{k}^{d_n})$.
This interval therefore lies in the gap with label $\mathbf{k}$
\cite{avila2016dry}. The cone inequalities remain strict at
$E=E_\mathbf{k}^{d_n}-\delta_0$, so uniform hyperbolicity extends
slightly to its left, proving \eqref{poly-es}.

Apply the order argument in Lemma~\ref{lem:E_kd} to both the minimal
and maximal energies at this IDS level. Each gap endpoint changes
by at most $\|v-v_{d_n}\|_0$, and consequently
\[
|G_{\mathbf{k}}(v,\alpha)|
\geq |G_{\mathbf{k}}(v_{d_n},\alpha)|-2\|v-v_{d_n}\|_0
>\delta_0-2\|v-v_{d_n}\|_0>\delta_0/2>0.
\]
The last inequality holds for a fixed sufficiently large $n$, since
$\delta_0\geq e^{-o(q_n)}$ whereas
$\|v-v_{d_n}\|_0\leq e^{-\kappa q_n}$, after decreasing $\kappa$
if necessary. This proves \eqref{ana-es}.
\end{proof}

\begin{remark}
    Establishing the intermediate lower bound $-\e_{0} < \e_k^{d_n}(E,x)$ is  essential for this argument. Because the one-step cocycle lacks strict monotonicity in $E$, the perturbation can initially push the phase into the negative region ($\e < 0$). Here, the unperturbed parabolic fixed point acts as a repeller, actively accelerating the negative drift. Explicitly bounding this excursion ensures the orbit does not escape the domain of validity of our local quadratic estimates,  before the strictly monotonic $2d_n$-th iterate can rescue the phase and push it back into the positive cone, and ensure the survival of the gap cone. 
\end{remark}

\begin{remark}
The assumption $c_{q_n}>0$ is consistent with the choice of the right endpoint
$E_{\mathbf{k}}^{d_n}$. Indeed, if $c_{q_n}<0$, reflecting both the phase
and the energy about $E_{\mathbf{k}}^{d_n}$ gives the same cone argument on
$(E_{\mathbf{k}}^{d_n}+e^{-\kappa_1q_n},E_{\mathbf{k}}^{d_n}+\delta_0)$.
By \eqref{eq:rotn-finite}, the IDS would equal its value at
$E_{\mathbf{k}}^{d_n}$ on this interval, contradicting the maximality of
$E_{\mathbf{k}}^{d_n}$. Hence $c_{q_n}>0$ for sufficiently large $n$.
\end{remark}

}

	\appendix

\section{Proof of  Lemma \ref{prop:regu-bundle}}\label{regu-bund}

    First by compactness of $\{\alpha_0\}\times\I\times \overline{\T_{h}}$, to show Lemma \ref{prop:regu-bundle} it suffices to show the following claim: for any $t_0\in \I$ any $x_0+iy_0\in \T_{h_0}$, there exists a neighborhood $\mathcal V\subset \CC\times \R\times \CC$ of $(t_0,\alpha_0,x_0+iy_0)$ such that for any $(t,\alpha, x+iy)\in \V$ the cocycle $(\alpha, A_t(x+iy))$ is partially hyperbolic and $E^\ast_{t,\alpha}(x+iy)$ holomorphically depends on $t $ and $x+iy$, continuously depends on $\alpha$, $\ast\in \{u,s,c, cs,cu\}$\footnote{$E^{cs}:=E^c\oplus E^s, E^{cu}:=E^c\oplus E^u$.}.

We first show the claim for $\ast=u$, the proof of $\ast=cu$ is the same. Then we consider the inverse cocycle we get the proof of the claim for $s$ and $cs$ bundle. The last step is to use the fact that $E^c$ is the intersection of $E^{cs}$ and $E^{cu}$ with holomorphicity on $t,x+iy$ and continuity on $\alpha$.

We denote by $u^{t_0,\alpha_0}(x+iy_0)$ the unstable bundle at $x+iy_0$ of the cocycle $(\alpha_0, A_{t_0}(\cdot+i y_0))$,  as an invariant section of the bundle $\T\times G(d-m, 2d)$. Consider a cone field $U=U_\epsilon$ which is defined by $\{x, B(u^{t_0,\alpha_0}(x+iy_0),\epsilon)\}_{x\in \T} $, $ B(u(x),\epsilon)\subset G(d-m,2d)$ for $\epsilon$ sufficiently small. 

Partially hyperbolicity of $(\alpha_0, A_{t_0}(\cdot +iy_0))$ implies the existence of $n$ such that the cocycle $(n\alpha_0, A_{t_0}^n(\cdot +iy_0))$ maps $\overline{U}$ into $U$. In particular by continuity for $(t,\alpha, z)$ in a small neighborhood $\V\subset \C\times \R\times \C $ of $( t_0, \alpha_0, iy_0)$, $A_{t}^n(x-n\alpha +z)$ also maps $\overline{U}(x-n\alpha)$ into $U(x)$ for any $x$, hence by classical cone field criterion, $(n\alpha, A_{t}^n(\cdot +z))$ is $(d-m)$-dominated \cite{AJS2014}. Therefore as \cite[Lemma 6.2]{AJS2014}, $E^u(x)$ for $(\alpha, A_{t}(\cdot +z))$ at $x$ is just  $$\lim_{k\to \infty}A^{kn}_{t}(x-kn\alpha+z)\cdot u(x-kn\alpha). 
$$
(The key is that here we allow $\alpha$ to vary! This is actually quite classical in the theory of smooth dynamical systems, cf. the proof of  \cite{CP}, Corollary 2.8.) For each $k$, $(t,\alpha, z)\mapsto A^{kn}_{t}(x-kn\alpha+z)\cdot u(x-kn\alpha)$ is joint continuous in $(t,\alpha,z)$ and holomorphic in $t,z$, taking values in a small ball $U(x)$, then by Montel theorem the limit function $E^u_{t,\alpha}(z)$ is also holomorphic depends on $t, z$. 
Joint continuity of $E^u$ in $(t,\alpha,z)$ follows from uniform convergence on compact sets.

\section{$1/2$-Hölder continuity of spectra}\label{12holder}

Let $L_{\alpha,\theta}\in \mathcal{L}(\ell^2(\mathbb{Z}))$ be the family of finite-range operators
\begin{equation}\label{defLgeneral}
	(L_{\alpha,\theta} u)_n \;:=\; \sum_{k=-d}^{d} \hat{v}_k u_{n-k} + w(\theta + n\alpha)\, u_n,
\end{equation}
where the hopping coefficients $\hat v_k\in\C$ are fixed and satisfy $\hat v_{-k}=\overline{\hat v_k}$, and $w\in C^1(\T,\R)$.  For $\alpha\in\R$ set
	$\sigma_\alpha \;:=\; \bigcup_{\theta\in\T} \sigma\big(L_{\alpha,\theta}\big).$
    
 \begin{lemma}\label{Lcontinuity}
There exist constants $\tilde{C}(\hat{v}, w), C(\hat{v}, w) > 0$ such that if $\lambda \in \sigma_\alpha$ and $|\alpha - \alpha'| < \tilde{C}(\hat{v}, w)$, then there exists $\lambda' \in \sigma_{\alpha'}$ such that
\begin{align*}
|\lambda - \lambda'| \le C(\hat{v}, w) |\alpha - \alpha'|^{1/2}.
\end{align*}
\end{lemma}

\begin{proof}
Let $L \ge 1$ be given. There exists $\phi_L \in \ell^2(\mathbb{Z})$ and $\theta$ such that
\begin{align}
\|(L_{\alpha,\theta} - \lambda)\phi_L\| \le \frac{1}{L} \|\phi_L\|.
\end{align}
Let $\eta_{j,L}$ be the test function centered at $j$,
\begin{align*}
\eta_{j,L}(n) =
\begin{cases}
1 - {|n-j|}/{L}, & |n-j| \le L, \\
0, & |n-j| > L.
\end{cases}
\end{align*}
Then the overlap identity holds:
\begin{align}\label{sumeta}
\sum_j (\eta_{j,L}(n))^2 = a_L := 1 + \frac{(L-1)(2L - 1)}{3L}
\end{align}
and independent of $n$. Using this,
\begin{align}\label{C1}
\sum_j \|\eta_{j,L}(L_{\alpha,\theta} - \lambda) \phi_L\|^2 = a_L \|(L_{\alpha,\theta} - \lambda) \phi_L\|^2 \le \frac{a_L}{L^2} \|\phi_L\|^2 = \frac{1}{L^2} \sum_j \|\eta_{j,L} \phi_L\|^2.
\end{align}
Since $\|u+v\|^2\leq 2\|v\|^2+2 \|u\|^2$, by \eqref{C1}, we get
\begin{align}\label{C2}
\sum_j \| (L_{\alpha,\theta} - \lambda) \eta_{j,L} \phi_L \|^2 &\le 2 \sum_j \| \eta_{j,L} (L_{\alpha,\theta} - \lambda) \phi_L \|^2 + 2 \sum_j \| [\eta_{j,L}, L_{\alpha,\theta}] \phi_L \|^2.
\end{align}
For the commutator term, observe:
\begin{align*}
[\eta_{j,L}, L_{\alpha,\theta}] \phi(n) &= \eta_{j,L}(n)(L_{\alpha,\theta} \phi)(n) - (L_{\alpha,\theta}(\eta_{j,L} \phi))(n) \\
&= \sum_{k=-d}^d \hat{v}_k \left[ \eta_{j,L}(n) \phi(n-k) - \eta_{j,L}(n-k) \phi(n-k) \right] \\
&= \sum_{k=-d}^d \hat{v}_k (\eta_{j,L}(n) - \eta_{j,L}(n - k)) \phi(n - k),
\end{align*}
which implies
\begin{align*}
\|[\eta_{j,L}, L_{\alpha,\theta}] \phi\|^2 \le \sum_n \left| \sum_{k=-d}^d \hat{v}_k (\eta_{j,L}(n) - \eta_{j,L}(n - k)) \phi(n - k) \right|^2.
\end{align*}
Using the inequality $|\eta(n) - \eta(n-k)| \le \frac{|k|}{L}$ and Cauchy-Schwarz, we estimate:
\begin{align*}
\|[\eta_{j,L}, L_{\alpha,\theta}] \phi\|^2 \le \frac{C_v}{L^2} \|\phi\|^2,
\end{align*}
where $C_v = (2d+1) \sum_{k=-d}^d |\hat{v}_k|^2 k^2$.

For  $L>L_0$, summing over $j$, 
\begin{align*}
\sum_j \| [\eta_{j,L}, L_{\alpha,\theta}] \phi_L \|^2 \le \frac{3C_v}{L} \|\phi_L\|^2 \le \frac{3C_v}{L a_L} \sum_j \|\eta_{j,L} \phi_L\|^2.
\end{align*}
Substituting into \eqref{C2}, and noting $a_L \sim \tfrac{2}{3}L$, we obtain:
\begin{align*}
\sum_j \| (L_{\alpha,\theta} - \lambda) \eta_{j,L} \phi_L \|^2 \le \frac{C_1(v)}{L^2} \sum_j \|\eta_{j,L} \phi_L\|^2.
\end{align*}
Hence for certain $j$, $\psi := \eta_{j,L} \phi_L \neq 0$ and satisfies
\begin{align}\label{approxeig}
\| (L_{\alpha,\theta} - \lambda) \psi \| \le \frac{\pa{C_1(v)}^{1/2}}{L} \|\psi\|.
\end{align}

Given $\alpha^\prime$ near $\alpha$, take $\theta'$ such that $\theta + j\alpha = \theta' + j\alpha'$. On $\mathrm{supp}(\psi)$, we have:
\begin{align*}
|w(\theta + n\alpha) - w(\theta' + n\alpha')| \le \|w'\|_{\infty} L |\alpha - \alpha'|.
\end{align*}
So $\| (L_{\alpha',\theta'} - L_{\alpha,\theta}) \psi \| \le C_2(w) L |\alpha - \alpha'| \|\psi\|$.
Combining with \eqref{approxeig}:
\begin{align*}
\| (L_{\alpha',\theta'} - \lambda) \psi \| \le \left(\frac{C_1(v)^{1/2}}{L} + C_2(w) L |\alpha - \alpha'| \right) \|\psi\|.
\end{align*}
Optimizing by setting $L \sim |\alpha - \alpha'|^{-1/2} > L_0$ yields:
\begin{align*}
\| (L_{\alpha',\theta'} - \lambda) \psi \| \le C(v,w) |\alpha - \alpha'|^{1/2} \|\psi\|.
\end{align*}
Thus, $\lambda$ lies within $C(v,w) |\alpha - \alpha'|^{1/2}$ of $\sigma(L_{\alpha',\theta'}) \subset \sigma_{\alpha'}$. \end{proof}

\smallskip 

\subsection*{Proof of Lemma \ref{cor:rational-ap}}
    Since $L_{\alpha,\theta}$ has no spectrum in $G^m(\alpha)=(E^m_-(\alpha),E^m_+(\alpha))$, then by Lemma \ref{Lcontinuity}, $L_{\alpha',\theta}$ has no spectrum in \[
    (E^m_-(\alpha)+C(v,w)|\alpha-\alpha'|^{1/2},E^m_+(\alpha)-C(v,w)|\alpha-\alpha'|^{1/2})
    \]
    which is non-empty (for $|\alpha-\alpha'|$ small). A simple continuity argument shows that the fibred rotation number unchanged. Symmetry implies the result. \qed

\section*{Acknowledgements} 
The authors would like to thank Svetlana Jitomirskaya for useful discussions. The third author would like to thank Qiongling Li for useful discussions on Proposition \ref{lem:ajslem}. A major part of this work was completed while the first author was a Ph.D. student at Nankai University, during which he was supported by NSFC grant (123B2005). Xianzhe Li was supported by NSF DMS-2052899, DMS-2155211, Simons 
896624 and an AMS-Simons Travel Grant. Disheng Xu was supported by National Key R$\&$D Program of China No. 2024YFA1015100, NSFC grant (12526201).  Qi Zhou was supported by NSFC grant (12531006,12526201) and  Nankai Zhide Foundation.

\section*{AI disclosure}
The authors used AI tools to assist with literature searches, proof verification and checking, figure preparation, and manuscript editing. All the main ideas in this paper came from the human authors. Perhaps the clearest sign of this is that only a human mathematician would be stubborn enough to pursue such a convoluted route.

\color{black}


\begin{thebibliography}{99}

\bibitem{Anosov1967} 
Anosov, D. V., Tangential fields of transversal foliations in $Y$-systems, \textit{Math. Notes} \textbf{2} (1967), 818--823.

\bibitem{AA} 
Argentieri, F. and Avila, A., In preparation.



\bibitem{avila2008absolutely} 
Avila, A., The absolutely continuous spectrum of the almost Mathieu operator, preprint (2008), arXiv:0810.2965.

\bibitem{AviAlmost} 
Avila, A., Almost reducibility and absolute continuity I, preprint (2010), arXiv:1006.0704.

\bibitem{avila0} 
Avila, A., Global theory of one-frequency Schr\"odinger operators, \textit{Acta Math.} \textbf{215} (2015), 1--54.

\bibitem{Avi2023KAM} 
Avila, A., KAM, Lyapunov exponents, and the spectral dichotomy for typical one-frequency Schr\"odinger operators, preprint (2023), arXiv:2307.11071.

\bibitem{ABD} 
Avila, A., Bochi, J., and Damanik, D., Opening gaps in the spectrum of strictly ergodic Schr\"odinger operators, \textit{J. Eur. Math. Soc.} \textbf{14} (2012), 61--106.

\bibitem{AFK2011Kam}
A.~Avila, B.~Fayad, and R.~Krikorian, A {KAM} scheme for {SL}(2, $\mathbb{R} $) cocycles with {L}iouvillean frequencies, \textit{Geom. Funct. Anal.} \textbf{21} (2011), no. 5, 1001--1019.

\bibitem{AJ05} 
Avila, A. and Jitomirskaya, S., The Ten Martini Problem, \textit{Ann. of Math.} \textbf{170} (2009), 303--342.

\bibitem{AJ08} 
Avila, A. and Jitomirskaya, S., Almost localization and almost reducibility, \textit{J. Eur. Math. Soc.} \textbf{12} (2010), 93--131.

\bibitem{AJS2014} 
Avila, A., Jitomirskaya, S., and Sadel, C., Complex one-frequency cocycles, \textit{J. Eur. Math. Soc.} \textbf{16} (2014), 1915--1935.

\bibitem{AK06} 
Avila, A. and Krikorian, R., Reducibility or non-uniform hyperbolicity for quasiperiodic Schr\"odinger cocycles, \textit{Ann. of Math.} \textbf{164} (2006), 911--940.

\bibitem{AK} 
Avila, A. and Krikorian, R., Monotonic cocycles, \textit{Invent. Math.} \textbf{202} (2015), 271--331.

\bibitem{ALSZ2024abominable} 
Avila, A., Last, Y., Shamis, M., and Zhou, Q., On the abominable properties of the almost Mathieu operator with well-approximated frequencies, \textit{Duke Math. J.} \textbf{173} (2024).

\bibitem{AYZ} 
Avila, A., You, J., and Zhou, Q., Sharp phase transitions for the almost Mathieu operator, \textit{Duke Math. J.} \textbf{166} (2017), 2697--2718.

\bibitem{avila2016dry} 
Avila, A., You, J., and Zhou, Q., Dry Ten Martini problem in the non-critical case, preprint (2024), arXiv:2306.16254v2.

\bibitem{AOS} 
Avron, J. E., Osadchy, D., and Seiler, R., A topological look at the quantum Hall effect, \textit{Phys. Today} \textbf{56} (2003), 38--42.

\bibitem{AS1983Almost} 
Avron, J. and Simon, B., Almost periodic Schr\"odinger operators II. The integrated density of states, \textit{Duke Math. J.} \textbf{50} (1983), 369--391.

\bibitem{AvS1991measure} 
Avron, J. E., van Mouche, P., and Simon, B., On the measure of the spectrum for the almost Mathieu operator, \textit{Comm. Math. Phys.} \textbf{139} (1991), 215--215.

\bibitem{BBL} 
Band, R., Beckus, S., and Loewy, R., The Dry Ten Martini Problem for Sturmian Hamiltonians, preprint (2024), arXiv:2402.16703.

\bibitem{BLT} 
Bellissard, J., Lima, R., and Testard, D., Almost periodic Schr\"odinger operators, in \textit{Mathematics + Physics: Lectures on Recent Results}, World Scientific, Singapore, 1985, 1--64.

\bibitem{BS} 
Bellissard, J. and Simon, B., Cantor spectrum for the almost Mathieu equation, \textit{J. Funct. Anal.} \textbf{48} (1982), 408--419.

\bibitem{BB} 
Berti, M. and Biasco, L., Forced vibrations of wave equations with non-monotone nonlinearities, \textit{Ann. I. H. Poincar\'e-AN} \textbf{23} (2006), 439--474.

\bibitem{BouGreens} 
Bourgain, J., \textit{Green's Function Estimates for Lattice Schr\"odinger Operators and Applications}, Princeton University Press, Princeton, 2004.

\bibitem{BG} 
Bourgain, J. and Goldstein, M., On nonperturbative localization with quasi-periodic potential, \textit{Ann. of Math.} \textbf{152} (2000), 835--879.

\bibitem{BJ02} 
Bourgain, J. and Jitomirskaya, S., Absolutely continuous spectrum for 1D quasiperiodic operators, \textit{Invent. Math.} \textbf{148} (2002), 453--463.

\bibitem{bj} 
Bourgain, J. and Jitomirskaya, S., Continuity of the Lyapunov exponent for quasiperiodic operators with analytic potential, \textit{J. Stat. Phys.} \textbf{108} (2002), 1203--1218.

\bibitem{BrinStuck} 
Brin, M. and Stuck, G., \textit{Introduction to Dynamical Systems}, Cambridge University Press, Cambridge, 2002.

\bibitem{CEY} 
Choi, M. D., Elliott, G. A., and Yui, N., Gauss polynomials and the rotation algebra, \textit{Invent. Math.} \textbf{99} (1990), 225--246.

\bibitem{CP} 
Crovisier, S. and Potrie, R., Introduction to partially hyperbolic dynamics, Lecture Notes for a minicourse at ICTP (2015).

\bibitem{DGY} 
Damanik, D., Gorodetski, A., and Yessen, W., The Fibonacci Hamiltonian, \textit{Invent. Math.} \textbf{206} (2016), 629--692.


\bibitem{De} 
Deimling, K., \textit{Nonlinear Functional Analysis}, Springer-Verlag, Berlin, 1985.

\bibitem{E92} 
Eliasson, L. H., Floquet solutions for the one-dimensional quasiperiodic Schr\"odinger equation, \textit{Comm. Math. Phys.} \textbf{146} (1992), 447--482.

\bibitem{Ekeland} 
Ekeland, I., \textit{Convexity Methods in Hamiltonian Mechanics}, Springer, Berlin, 1990.

\bibitem{Forstneric2017} 
Forstnerič, F., \textit{Stein Manifolds and Holomorphic Mappings: The Homotopy Principle in Complex Analysis}, 2nd edn, Springer, Berlin, 2017.


\bibitem{GeARC}
Ge, L., On the almost reducibility conjecture, \textit{Geom. Funct. Anal.} \textbf{34} (2024), 32--59.

\bibitem{GJ} 
Ge, L. and Jitomirskaya, S. Y., Intrinsic symplectic structure and sharp arithmetic universality, preprint (2026), arXiv:2407.08866v2.

\bibitem{GJY} 
Ge, L., Jitomirskaya, S. Y., and You, J., Kotani theory, Puig's argument, and stability of The Ten Martini Problem, preprint (2023), arXiv:2308.09321.

\bibitem{Martini}
Ge, L., Jitomirskaya, S. Y., and You, J., The Robust Ten Martini Problem, arxiv:2609.15812v1

\bibitem{GJY2}
Ge, L., Jitomirskaya, S. Y., and You, J., In preparation.  

\bibitem{GJYZ} 
Ge, L., Jitomirskaya, S. Y., You, J., and Zhou, Q., Multiplicative Jensen's formula and quantitative global theory of one-frequency Schr\"odinger operators, \textit{Forum Math. Pi}, \textbf{14} (2026), e1.

\bibitem{GeWXu} 
Ge, L., Wang, Y., and Xu, J., The Dry Ten Martini Problem for $C^2$ cosine-type quasiperiodic Schr\"odinger operators, preprint (2025), arXiv:2503.06918.

\bibitem{GY2020Arithmetic} 
Ge, L. and You, J., Arithmetic version of Anderson localization via reducibility, \textit{Geom. Funct. Anal.} \textbf{30} (2020), 1370--1401.

\bibitem{GYZ2024Structured} 
Ge, L., You, J., and Zhou, Q., Sharp Spectral Gaps, Arithmetic Localization, and Reducibility via Resonance Analysis, arXiv:2407.05490v2.

\bibitem{GJLS1997Duality} 
Gordon, A. Y., Jitomirskaya, S., Last, Y., and Simon, B., Duality and singular continuous spectrum in the almost Mathieu equation, \textit{Acta Math.} \textbf{178} (1997), 169--183.

\bibitem{GBSS} 
Gottlob, E., Borgnia, D. S., Slager, R. J., and Schneider, U., Quasiperiodicity protects quantized transport in disordered systems without gaps, \textit{PRX Quantum} \textbf{6} (2025), 020359.

\bibitem{GS1} 
Goldstein, M. and Schlag, W., H\"older continuity of the integrated density of states for quasi-periodic Schr\"odinger equations and averages of shifts of subharmonic functions, \textit{Ann. of Math.} \textbf{154} (2001), 155--203.

\bibitem{GS11} 
Goldstein, M. and Schlag, W., On resonances and the formation of gaps in the spectrum of quasi-periodic Schr\"odinger equations, \textit{Ann. of Math.} \textbf{173} (2011), 337--475.

\bibitem{HS1}
Han, R. and Schlag, W., Nonperturbative localization on the strip and Avila's almost reducibility conjecture, \textit{Forum Math. Pi} \textbf{14} (2026), e20.

\bibitem{HP} 
Haro, A. and Puig, J., A Thouless formula and Aubry duality for long-range Schr\"odinger skew-products, \textit{Nonlinearity} \textbf{26} (2013), 1163--1189.

\bibitem{Ha} 
Harper, P. G., Single band motion of conduction electrons in a uniform magnetic field, \textit{Proc. Phys. Soc. London A} \textbf{68} (1955), 874--892.

\bibitem{HatcherVB} 
Hatcher, A., \textit{Vector Bundles and K-Theory}, Version 2.2 (2017).

\bibitem{HS} 
Helffer, B. and Sj\"ostrand, J., Semiclassical analysis for Harper's equation. III. Cantor structure of the spectrum, \textit{M\'em. Soc. Math. France} \textbf{39} (1989), 1--124.

\bibitem{HY} 
Hou, X. and You, J., Almost reducibility and non-perturbative reducibility of quasi-periodic linear systems, \textit{Invent. Math.} \textbf{190} (2012), 209--260.

\bibitem{howardmaslov} 
Howard, P., Latushkin, Y., and Sukhtayev, A., The Maslov index for Lagrangian pairs on $\mathbb{R}^{2n}$, \textit{J. Math. Anal. Appl.} \textbf{451} (2017), 794--821.

\bibitem{HOWARD2016} 
Howard, P. and Sukhtayev, A., The Maslov and Morse indices for Schr\"odinger operators on [0,1], \textit{J. Differ. Equ.} \textbf{260} (2016), 4499--4549.


\bibitem{JM2012}
Jitomirskaya, S. Y. and Marx, C.~A.,
 Analytic quasi-periodic {Schr\"odinger} operators and rational frequency approximants,
 {\it Geom. Funct. Anal.} {\bf 22} (2012), no.~5, 1407--1443.

\bibitem{gaplabel} 
Johnson, R. and Moser, J., The rotation number for almost periodic potentials, \textit{Comm. Math. Phys.} \textbf{84} (1982), 403--438.

\bibitem{Kac} 
Kac, M., Public communication at 1981 AMS Annual Meeting.

\bibitem{Klitzing1980} 
Klitzing, K. v., Dorda, G., and Pepper, N., New method for high-accuracy determination of the fine-structure constant based on quantized Hall resistance, \textit{Phys. Rev. Lett.} \textbf{45} (1980), 494--497.

\bibitem{L} 
Last, Y., Zero measure spectrum for the almost Mathieu operator, \textit{Comm. Math. Phys.} \textbf{164} (1994), 421--432.

\bibitem{L1} 
Last, Y., Spectral theory of Sturm-Liouville operators on infinite intervals: a review of recent developments, in \textit{Sturm-Liouville Theory}, Birkh\"auser, Basel, 2005, 99--120.

\bibitem{LYZZ2017Asymptotics} 
Leguil, M., You, J., Zhao, Z., and Zhou, Q., Asymptotics of spectral gaps of quasi-periodic Schr\"odinger operators, \textit{Camb. J. Math.} \textbf{12} (2024), 753--830.

\bibitem{L}
Li, X. Dry Ten Martini Problem in the Subcritical Type I Regime. 
preprint (2026), arXiv:2609.23467.


\bibitem{LW} 
Li, X. and Wu, L., The fibred rotation number for ergodic symplectic cocycles and its applications: I. Gap Labelling Theorem, \textit{Math. Z.} \textbf{311} (2025), Art. 53.

\bibitem{LYZ2024Exact} 
Li, X., You, J., and Zhou, Q., Exact local distribution of the absolutely continuous spectral measure, \textit{Proc. Lond. Math. Soc.}, \textbf{132} (2026), e70150.

\bibitem{LZ} 
Li, X. and Zhou, Q., In preparation.

\bibitem{LY2015} Liu W. and Yuan, X., Spectral gaps of almost Mathieu operators in the exponential regime, \textit{J. Fractal Geom.} \textbf{2} (2015), no.~1, 1--51.


\bibitem{Long} 
Long, Y., \textit{Index Theory for Symplectic Paths with Applications}, vol. 207, Birkh\"auser, Basel, 2012.

\bibitem{MS1974} 
Milnor, J. W. and Stasheff, J. D., \textit{Characteristic Classes}, Annals of Mathematics Studies, No. 76, Princeton University Press, Princeton, 1974.

\bibitem{MP84} 
Moser, J. and P\"oschel, J., An extension of a result by Dinaburg and Sinai on quasiperiodic potentials, \textit{Comment. Math. Helv.} \textbf{59} (1984), 39--85.

\bibitem{OA} 
Osadchy, D. and Avron, J. E., Hofstadter butterfly as quantum phase diagram, \textit{J. Math. Phys.} \textbf{42} (2001), 5665--5671.

\bibitem{Pe} 
Peierls, R., Zur Theorie des Diamagnetismus von Leitungselektronen, \textit{Z. Phys.} \textbf{80} (1933), 763--791.

\bibitem{P} 
Puig, J., Cantor spectrum for the almost Mathieu operator, \textit{Comm. Math. Phys.} \textbf{244} (2004), 297--309.

\bibitem{P06} 
Puig, J., A nonperturbative Eliasson's reducibility theorem, \textit{Nonlinearity} \textbf{19} (2006), 355--376.

\bibitem{R} 
Rauh, A., Degeneracy of Landau levels in crystals, \textit{Phys. Status Solidi B} \textbf{65} (1974), 131--135.

\bibitem{Shabat92} 
Shabat, B. V., \textit{Introduction to Complex Analysis, Part II: Functions of Several Variables}, Translations of Mathematical Monographs, vol. 110, American Mathematical Society, Providence, 1992.

\bibitem{Sin} 
Sinai, Ya. G., Anderson localization for one-dimensional difference Schr\"odinger operator with quasiperiodic potential, \textit{J. Stat. Phys.} \textbf{46} (1987), 861--909.

\bibitem{SS} 
Sontz, S. B., \textit{Principal Bundles: The Classical Case}, Springer, Berlin, 2015.

\bibitem{TKNN1982} 
Thouless, D. J., Kohmoto, M., Nightingale, M. P., and den Nijs, M., Quantised Hall conductance in a two-dimensional periodic potential, \textit{Phys. Rev. Lett.} \textbf{49} (1982), 405--408.

\bibitem{WXZ} 
Wang, D., Xu, D., and Zhou, Q., Hyperbolicity for one-frequency analytic quasi-periodic (Hermitian)-symplectic cocycles, preprint (2025), arXiv:2503.15281.

\bibitem{WAXZ} 
Wang, Z., Xu, D., and Zhou, Q., Subordinacy theory for long-range operators: Hyperbolic geodesic flow insights and monotonicity theory, preprint (2025), arXiv:2506.23098.

\bibitem{Xu2019Density} 
Xu, D., Density of positive Lyapunov exponents for symplectic cocycles, \textit{J. Eur. Math. Soc.} \textbf{21} (2019), 3143--3190.

\bibitem{You}You, J., Some problems in quasi-periodic Schr\"odinger operators, J. Math. Phys. {\bf 67} (2026), 062704.

\end{thebibliography}
\end{document}